\documentclass{amsbook}

\usepackage{amssymb}
\usepackage{amsfonts}
\usepackage{amsmath}
\usepackage{lastpage}
\usepackage[legalpaper,bookmarks=true,colorlinks=true,linkcolor=blue,citecolor=blue]{hyperref}
\usepackage{graphicx}%
\usepackage{fancyhdr}
\usepackage{color}
\usepackage[mathlines]{lineno}
\usepackage{lscape}
\usepackage{epsfig}
\usepackage{natbib}
\usepackage{geometry}
\usepackage{tgbonum}
\usepackage{imakeidx}
\usepackage{imakeidx}
\makeindex[columns=3, title=Alphabetical Index]

\fontfamily{qcr}\selectfont

\newtheorem{theorem}{Theorem}
\theoremstyle{plain}

\newtheorem{lemma}{Lemma}

\newtheorem{proposition}{Proposition}

\numberwithin{equation}{section}

\newcommand{\Bin}{\bigskip \noindent}

\newcommand{\Ni}{\noindent}

\begin{document}
\Large
\pagenumbering{roman}
\begin{center}

\huge \textbf{Gane Samb Lo \\ Aladji Babacar Niang \\ Lois Chinwendu Okereke}\\
\vskip 6cm
\Huge \textbf{A Course on Elementary Probability Theory} \\

\vskip 6cm

\huge \textit{\textbf{Statistics and Probability African Society (SPAS) Books Series}.\\
\textbf{Saint-Louis, Calgary, Alberta. 2020}}.\\

\bigskip \Large  \textbf{DOI} : http://dx.doi.org/10.16929/sts/2020.001\\
\bigskip \textbf{ISBN:} 9798582099772 
\end{center}

\newpage
\begin{center}
\huge \textbf{SPAS TEXTBOOKS SERIES}
\end{center}

\bigskip \bigskip

\Large

 \begin{center}
 \textbf{GENERAL EDITOR of SPAS EDITIONS}
 \end{center}

\vskip 0.7cm
\noindent \textbf{Prof Gane Samb LO}\\
gane-samb.lo@ugb.edu.sn, gslo@ugb.edu.ng\\
Gaston Berger University (UGB), Saint-Louis, SENEGAL.\\
African University of Science and Technology, AUST, Abuja, Nigeria.\\

\begin{center}
\Large \textbf{HONORARY EDITOR}
\end{center}

\noindent \textbf{Hamet Seydi}.\\
Universit\'es of SENEGAL.\\
hseydi@gmail.com\\

\begin{center}
\Large \textbf{ASSOCIATED EDITORS}
\end{center}

\noindent \textbf{Brahim MEZERDI}.\\
mezerdi@univ-biskra.dz, bmezerdi@yahoo.fr\\
Universit\'e Biskra, Algeria\\

\noindent \textbf{Blaise SOME}\\
some@univ-ouaga.bf\\
Chairman of LANIBIO, UFR/SEA\\
Ouaga I Pr Joseph Ki-Zerbo University.\\

\begin{center}
\Large \textbf{ADVISORS}
\end{center}

\bigskip
\noindent \textbf{Mohamed Ahsanullah}\\
ahsan@rider.edu\\
Rider University, Lawrence, USA.\\

\noindent \textbf{Dr Diam Ba}\\
dba@statpas.org\\
Gaston Berger University, Senegal.\\

\noindent \textbf{Tchilabalo Abozou KPANZOU}\\
kpanzout@yahoo.fr, kpanzouta@statpas.org\\
Kara University, Togo.\\

\newpage
\noindent \Large \Large \textbf{List of published books}\\

\Ni Collection : SPAS BOOKS SERIES.\\

\Ni 1. A Collection of Papers in Mathematics and Related Sciences. Hamet SEYDI, Gane Samb LO, and Aboubakary DIAKHABY (eds.) (2018)\\
Doi : 10.16929/sbs/2018.100. ISBN:  978-2-9559183-0-2. Available at:\\
Euclid.org :https://projecteuclid.org/euclid.spaseds/1569509457
Amazon.com : Search COLLECTION-PAPERS-MATHEMATICS-RELATED-SCIENCES\\

\bigskip
\noindent \textbf{For the latest updates, visit our official website :}\\

\begin{center}
www.statpas.org/spaseds/
\end{center}

\newpage
\noindent \textbf{Library of Congress Cataloging-in-Publication Data}\\

\noindent Gane Samb LO, 1958-;\\
Aladji Babacar Niang 1994-;
Lois Chinwendu Okereke 1989-\\

\noindent A Course on Elementary Probability Theory.\\

\noindent SPAS textbooks Series, 2020.\\

%\noindent Copyright \copyright Statistics and Probability African Society (SPAS).\\

\noindent \textit{DOI} : 10.16929/sts/2020.001\\

\noindent \textit{ISBN} 9798576736386.

\newpage

\noindent \textbf{Author : Gane Samb LO}\\
\bigskip

\bigskip
\noindent \textbf{Emails}:\\
\noindent gane-samb.lo@ugb.edu.sn, ganesamblo@ganesamblo.net, gslo@aust.edu.ng\\

\bigskip
\noindent \textbf{Url's}:\\
\noindent www.ganesamblo@ganesamblo.net\\
\noindent www.statpas.net/@ganesamblo\\

\bigskip \noindent \textbf{Affiliations}.\\
Main affiliation : University Gaston Berger, UGB, SENEGAL.\\
African University of Sciences and Technology, AUST, Abuja, Nigeria.\\
Affiliated as a researcher to : LSTA, Pierre et Marie Curie University, Paris VI, France.\\

\noindent \textbf{Teaches or has taught} at the graduate level in the following universities:\\
Saint-Louis, Senegal (UGB)\\
Abuja - Nigeria (AUST)\\
Banjul - Gambia (TUG)\\
Bamako - Mali (USTTB)\\
Ouagadougou - Burkina Faso (UJK)\\
African Institute of Mathematical Sciences, Mbour, SENEGAL, AIMS.\\
Franceville - Gabon\\
Kara, (Togo), undergraduate level\\

\newpage
\noindent \textbf{General acknowledgment}.\\

\Ni As main author, I wish to acknowledge the help of many people.  I have been working for almost thirty years, mainly at the Saint-Louis Gaston University, since 1991. I also taught in many African Universities, especially at the master degree level. I worked in the administration of the university at all levels (head of department, dean faculty, vice-president). I also supervised more than twenty Ph.D. theses and several master dissertations. As well, I animated the LERSTAD, a research group I created around 1992 and ran it for years with a regular weekly seminar.\\

\Ni The books I am writing are the outcomes of all these activities. I worked with amazing people, younger colleagues and Ph.D. students who became later high profiled researchers in Africa, Europa, Canada and the United States, and other parts of the world. For the last few years, I have been teaching in Nigeria, in the African University of Sciences and Technology of the Nelson Mandela Institute, Abuja, Nigeria.\\

\Ni The books of this series in Mathematics in general,  and in Random Analysis (Probability Theory and Statistics and their applications) are written in English since we want to reach a bigger public. But French versions will be published for books regularly used in undergraduate education. A broader presentation of our series of book can be found in the general preface, in page \pageref{generalPreface}.\\

\Ni I want to thank many people involved in the publication process of our books, particularly this one:\\

\Ni (1) Students who followed this course at Universit\'e Gaston Berger and at Universit\'e Dakar-Bourguiba (SENEGAL) during years. The course I taught them has become the book you have in your hand.\\

\Ni (2) Members of my research teams (LERSTAD, IMHOTEP) who are asked to read all our books and who regularly take part in the editing (Drs. Tchilabola A. Kpanzou [Togo]; Harouna Sangar\'e, Soumaila Demb\'el'\'e, Mouminou Diallo [Mali]; Modou Ngom, Diam B\^a, Amadou Dadhi\'e Ba; Mrs Gorgui Gning, Cherif Mamadou Moctar Traor\'e, etc.)\\

\Ni (3) In 2018, my former students in Measure Theory and Integration in Abuja : Tagbo Innocent Aroh, Lois Chinwendu Okereke, Abubakar Adamu, Aicha Adam Aminu, Chidiebere Eze Leonard Eze, helped in polishing the English text and checked the mathematical formulas.\\

\Ni (4) My collaborators Aladji Babacar Niang and Lois Chinwendu Okereke, who conducted a thorough editing of the book. I am convinced that their work tremendously improved this book. As a result, I consider them co-authors with at least twenty-five percent of the co-authorship for each of them.\\

\bigskip \noindent \textbf{Acknowledgment of Funding}.\\

\Ni The author acknowledges continuous and various support from the authorities of Gaston Berger University. It is true that in Senegalese public universities, the weekly teaching charge for full professors is five hours per week. The remaining time is devoted for supervision and research activities. Accordingly, scholars have all the means to realize their research activities : personal office, equipment with computers, printers, fax machines, ink, paper, internet connection, etc., and funding for participation in conferences.\\

\Ni In that sense, I acknowledge that writing the books of this series is implicitly funded by the university and the state of Senegal. I express my most sincere appreciations to their authorities.

%\maketitle
\newpage
\Ni \textbf{Abstract of the book}.\\
\renewcommand\footrulewidth{0.5pt}

\Ni  This book introduces to the theory of probabilities from the beginning. Assuming that the reader possesses the normal mathematical level acquired at the end of the secondary school, we aim to equip him with a solid basis in probability theory. The theory is preceded by a general chapter on counting methods. Then, the theory of probabilities is presented in a discrete framework. Two objectives are sought. The first is to give the reader the ability to solve a large number of problems related to probability theory, including application problems in a variety of disciplines. The second is to prepare the reader before he takes course on the mathematical foundations of probability theory. In this later book, the reader will concentrate more on mathematical concepts, while in the present text, experimental frameworks are mostly found. If both objectives are met, the reader will have already acquired a definitive experience in problem-solving ability with the tools of probability theory and at the same time he is ready to move on to a theoretical course on probability theory based on the theory of Measure and Integration. The book ends with a chapter that allows the reader to begin an intermediate course in mathematical statistics.\\


\noindent \textbf{AMS 2010 Classification Subjects :} 60GXX; 62GXX.

\newpage
\Ni \textbf{R\'esum\'e de l'ouvrage}. \\

\Ni Cet ouvrage introduit \`a la th\'eorie des probabilit\'es depuis le d\'ebut. En supposant que le lecteur poss\`ede le niveau math\'ematique normal acquis \`a la fin du lyc\'ee, nous ambitionnons de le doter d'une base solide en th\'eorie des probabilit\'es. L'expos\'e de la th\'eorie est pr\'ec\'ed\'e d'un chapitre g\'en\'eral sur les m\'ethodes de comptage. Ensuite, la th\'eorie des probabilit\'es est pr\'esent\'ee dans un cadre discret. Deux objectifs sont recherch\'es. Le premier est de donner au lecteur la capacit\'e \`a r\'esoudre un grand nombre de probl\`emes li\'es \`a la th\'eorie des probabilit\'es, y compris les probl\`emes d'application dans une vari\'et\'e de disciplines. Le second \'etait de pr\'eparer le lecteur avant qu'il n'aborde l'ouvrage sur les fondements math\'ematiques de la th\'eorie des probabilit\'es. Dans ce dernier ouvrage, le lecteur se concentrera davantage sur des concepts math\'ematiques tandis que dans le pr\'esent texte, il se trouvent surtout des cadres exp\'erimentaux. Si les deux objectifs sont atteints, le lecteur aura d\'ej\`a acquis une exp\'erience d\'efinitive en capacit\'e de r\'esolution de probl\`emes de la vie r\'eelle avec les outils de la th\'eorie des probabilit\'es et en même temps, il est prêt \`a passer \`a un cours th\'eorique sur les probabilit\'es bas\'ees sur la th\'eorie de la mesure et l'int\'egration. Le livre se termine par par un chapitre qui permet au lecteur de commencer un cours interm\'ediaire en statistiques math\'ematiques.\\

%\maketitle

\newpage
\bigskip \noindent \textbf{Dedication}.\\

\begin{center}
To our beloved and late sister Khady Kane LO\\
27/07/1953 - 7/11/1988
\end{center}

\begin{figure}[h]
	\centering
		\includegraphics[width=0.80\textwidth]{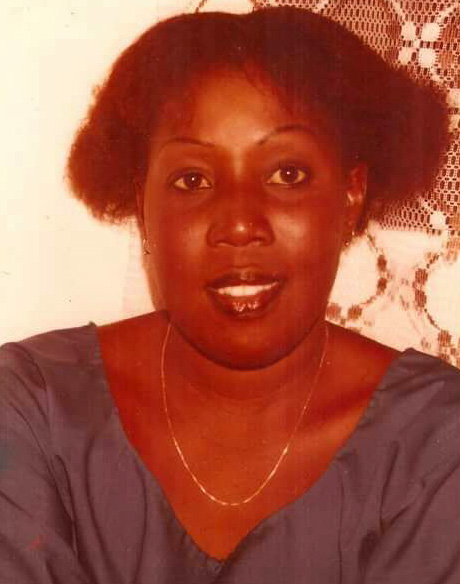}
	\caption{The ever smiling young lady}
	\label{proba01.crv.fig1}
\end{figure}

\tableofcontents
%\listoftables
\frontmatter
\mainmatter
\Large

\chapter*{General Preface} \label{generalPreface}

\noindent \textbf{This textbook} is part of a series whose ambition is to cover broad part of Probability Theory and Statistics.  These textbooks are intended to help learners and readers, both of of all levels, to train themselves.\\

\noindent As well, they may constitute helpful documents for professors and teachers for both courses and exercises.  For more ambitious  people, they are only starting points towards more advanced and personalized books. So, these texts are kindly put at the disposal of professors and learners.

\bigskip \noindent \textbf{Our textbooks are classified into categories}.\\

\noindent \textbf{A series of introductory  books for beginners}. Books of this series are usually accessible to student of first year in 
universities. They do not require advanced mathematics.  Books on elementary probability theory and descriptive statistics are to be put in that category. Books of that kind are usually introductions to more advanced and mathematical versions of the same theory. The first prepare the applications of the second.\\

\noindent \textbf{A series of books oriented to applications}. Students or researchers in very related disciplines  such as Health studies, Hydrology, Finance, Economics, etc.  may be in need of Probability Theory or Statistics. They are not interested in these disciplines  by themselves, rather in the need to apply the findings of these disciplines as tools to solve their specific problems. So adapted books on Probability Theory and Statistics may be composed to focus on the applications of such fields. A perfect example concerns the need of mathematical statistics for economists who do not necessarily have a good background in Measure Theory.\\

\noindent \textbf{A series of specialized books on Probability theory and Statistics of high level}. This series begin with a book on Measure Theory, its counterpart of probability theory, and an introductory book on topology. On that basis, we will have, as much as possible,  a coherent presentation of branches of Probability theory and Statistics. We will try  to make it self-contained, as much as possible, so that anything we need will be in the series.\\

\noindent Finally, \textbf{research monographs} close this architecture. The architecture should be so large and deep that the readers of monographs booklets will find all needed theories and inputs in it.\\

\bigskip \noindent We conclude by saying that, with  only an undergraduate level, the reader will  open the door of anything in Probability theory and statistics with \textbf{Measure Theory and integration}. Once this course validated, eventually combined with two solid courses on topology and functional analysis, he will have all the means to get specialized in any branch in these disciplines.\\

\bigskip \noindent Our collaborators and former students are invited to make live this trend and to develop it so  that  the center of Saint-Louis becomes or continues to be a renowned mathematical school, especially in Probability Theory and Statistics.

\chapter*{Preface of the first edition 2018}

The current series of Probability Theory  \index{probability theory} and Statistics  \index{statistics} are based on two {introductory  books for beginners} : A Course of Elementary probability Theory and A course on Descriptive Statistics.\\

\noindent All the more or less advanced probability courses are preceded by this one. We strongly recommend you do not skip it. It has the tremendous advantage of making the feel reader the essence of probability theory  \index{probability theory} by using extensively random experiences. The mathematical concepts come only after a complete description of a random experience.\\

\noindent  This book introduces the theory of probabilities from the beginning. Assuming that the reader possesses the normal mathematical level acquired at the end of the secondary school, we aim to equip him with a solid basis in probability theory.  \index{probability theory} The theory is preceded by a general chapter on counting methods. Then, the theory of probabilities is presented in a discrete framework.\\

\noindent Two objectives are sought. The first is to give the reader the ability to solve a large number of problems related to probability theory,  \index{probability theory} including application problems in a variety of disciplines. The second is to prepare the reader before he approached the textbook on the mathematical foundations of probability theory. In this book, the reader will concentrate more on mathematical concepts, while in the present text, experimental frameworks are mostly found. If both objectives are met, the reader will have already acquired a definitive experience in problem-solving ability with the tools of probability theory and at the same time he is ready to move on to a theoretical course on probability theory based on the theory of measurement and integration.\\

\noindent The book ends with a chapter that allows the reader to begin an intermediate course in mathematical statistics \index{statistics} .\\

\noindent I wish you a pleasant reading and hope receiving your feedback.\\

\noindent To my late and beloved sister Khady Kane Lo(1953-  ).\\

\noindent Saint-Louis, Calgary, Abuja, Bamako, Ouagadougou, 2017.

\chapter*{Introduction} \label{proba01.intro}

\noindent There exists a tremendous number of random phenomena in nature, real life and experimental sciences.\\

\noindent Almost everything is random in nature : whether, occurrences of rain and their durations, number of double stars in a region of the sky, lifetimes of plants, of humans, and of animals, life span of a radioactive atom, phenotypes of offspring of plants or any biological beings, etc.\\

\noindent The general theory states that each phenomena has a structural part (that is deterministic) and a random part (called the error or the deviation).\\

\noindent Randomness also appears as conceptual experiments : tossing  \index{tossing} a coin once or 100 times, throwing three dice, arranging a deck of cards, matching two decks, playing roulette, etc.\\

\noindent Every day human life is subject to randomness : waiting times for buses, traffic, number of calls on a telephone, number of busy lines in a communication network, sex of a newborn, etc.\\

\noindent The reader is referred to \cite{feller1} \index{Feller} for a more diverse and rich set of examples.\\

\noindent The quantitative study of random phenomena is the objective of Probability Theory  \index{probability theory} and Statistics  \index{statistics} Theory. Let us give
two simple examples to briefly describe each of these two disciplines.

\noindent In Probability Theory,  \index{probability theory} one assigns a chance of realization to a random event before its realization,
taking into account the available information.\\

\noindent {\bf Example }: A \textit{good coin}, that is, a homogenous  \index{homogenous} and well balanced, is tossed. Knowing that the coin cannot stand on its rim, there is a 50\%  chances of having a Tail \index{tail} .\\

\noindent We base our conclusion on the lack of any reason to favor one of the possible outcomes : head  \index{head} or tail.  \index{tail} So we convene that these outcomes
are equally probable and then get 50\% chances for the occurring for each of them.\\

\noindent Let us start with an example. Suppose that we have a coin and we do not know any thing of the material structure of the coin. In particular, we doubt that the coin is homogenous \index{homogenous} .\\

\noindent We decide to toss it repeatedly and to progressively monitor the occurring frequency  \index{frequency} of the head.  \index{head} We denote
by $N_{n}$ the number of heads  \index{head} obtained after $n$ tossing  \index{tossing} and define the frequency  \index{frequency} of the heads by

\begin{equation*}
F_{n}=\frac{N_{n}}{n}.
\end{equation*}

\Bin It is conceivable to say that the stability of $F_n$ in the neighborhood of some value $p \in ]0,1[$ is an
important information about the structure of the coin.\\

\noindent In particular, if we observe that the frequency  \index{frequency} $F_n$ does not deviate from $p=50\%$ more than $\varepsilon=0.001$ whenever
$n$ is greater than $100$, that is

\begin{equation*}
\left| F_{n}-1/2\right| \leq 10^{-3}
\end{equation*}

\Bin  for $n\geq 100$, we will be keen to say that the probability of occurrence of the head  \index{head} is $p=50\%$ and, by this,
we accept that that the coin is fair.\\

\noindent Based on the data (also called the statistics),  \index{statistics} we \textbf{have estimated} the probability of having a head  \index{head} at $p=50\%$ and
accepted the pre-conceived idea (hypothesis) that the coin is \textit{good} in the sense of homogeneity.\\

\noindent The reasoning we made and the method we applied are perfect illustrations of Statistical Methodology : estimation and
model or hypothesis validation from data.\\

\noindent \textbf{In conclusion},  Statistics  \index{statistics} Theory enables the use of the data (also called statistics or observations), to estimate the law of
a random phenomenon and to use that law to predict the future of the same phenomenon (\textit{inference}) or to predict any other phenomenon that seems identical to it (\textit{extrapolation}).\\

\noindent \textbf{(1)} The discipline of Statistics  \index{statistics} Theory and that of Probability Theory  \index{probability theory} are two ways of treating
the same random problems.\\

\noindent \textbf{(2)} The first is based primarily on the data to draw conclusions .\\

\noindent \textbf{(3)} The second is based on theoretical, and mathematical considerations to establish theoretical formulas.\\

\noindent\textbf{(4)}  Nevertheless, the Statistics  \index{statistics} discipline may be seen as the culmination of Probability Theory \index{probability theory} .\\

\newpage

\bigskip \noindent This book is an introduction to Elementary Probability Theory.  \index{probability theory} It is the result of many years of teaching
the discipline in Universities and High schools, mainly in Gaston Berger Campus of Saint-Louis, SENEGAL.\\

\noindent It is intended to those who have never done it before. It focuses on the essential points of the theory.
It is particularly adapted for undergraduate students in the first year of high schools.\\

\noindent The basic idea of this book consists of introducing Probability Theory,  \index{probability theory} and the notions of events and random variables
in discrete probability spaces.  \index{probability space} In such spaces, we discover as much as possible at this level, the fundamental
properties and tools of the theory.\\

\noindent So, we do not need at this stage the elaborated notion of $\sigma$-algebras or fields. This useful and brilliant method has already been used, for instance, in Billingsley \cite{billingsleyMT} \index{Billingsley} for deep and advanced probability problems.\\

\noindent The book will finish by a bridge chapter towards a medium course of Mathematical Statistics.  \index{statistics} In this chapter, we use an
analogy method to express the former results in a general shape, that will be the first chapter of the a fore mentioned course.\\

\noindent The reader will have the opportunity to master the tools he will be using in this course,  with the coming course on the mathematical foundation of Probability Theory,  \index{probability theory} which will be an element of our Probability and Statistics  \index{statistics} series. For such an advanced course, still to come, the reader will have to get prepared by a course of Measure Theory of Integration.\\

\noindent In this computer dominated world, we are lucky to have very powerful and free software like \textit{R} and \textit{Scilab}, to cite only the most celebrated. We seized this tremendous opportunity to use numerical examples using R software throughout the text.\\

\noindent The remainder of the book is organized as following.\\

\bigskip \noindent Chapter \ref{proba01.combin} is a quick but sufficient introduction to combinatorial analysis.
The student interested in furthering knowledge in this subject  can refer to the course on general algebra.\\

\bigskip \noindent Chapter \ref{proba01.pm} is devoted to an introduction to probability measures  \index{probability measure} in discrete spaces.
The notion of equi-probability  \index{equi-probability} is also dealt with, there.\\

\bigskip \noindent Conditional  \index{conditional} probability and independence  \index{independence} of random events are addressed in Chapter \ref{proba01.cpi}.\\

\bigskip \noindent Chapter \ref{proba01.rv} introduces random variables, and presents a review of the usual probability law \index{probability law} s.\\

\bigskip \noindent Chapter \ref{proba01.param} is devoted to the computations of the parameters of the usual laws.
Mastering the results of this chapter and those in Chapter \ref{proba01.rv} is in fact of great importance for higher level courses.\\

\bigskip \noindent A window on dependence methods will be presented in Chapter  \ref{proba01.randomcouples}.\\

\bigskip \noindent Distribution functions of random variables are studied in Chapter \ref{proba01.crv}, which introduces continuous random variable \index{continuous Random Variable} s.\\

%\bigskip \noindent Chapter \ref{proba01.summary} is an extrapolation of the whole book and constitutes an opening into a course of Mathematical Statistics \index{statistics} .\\

%\bigskip \noindent A chapter \ref{proba01.npm} devoted ot numerical probability theory  \index{probability theory} is expected soon to.\\

\chapter{Elements of Combinatorics}  \index{combinatoric} \label{proba01.combin}

Here, we are concerned with counting cardinalities of subsets of a reference set $\Omega$,  by following specific rules. We begin by a general counting principle.\\

\noindent The cardinality of a set $E$ is the number of its elements, denoted by $Card(E)$ or $\#(E)$. For an infinite set $E$, we have

$$
Card(E)=\#(E)=+\infty.
$$

\Ni Let us consider the following example : A student has two (2) skirts and four (4) pants in his closet. He decides to pick at random
a skirt and pants to dress. In how many ways can he dress by choosing a skirt and pants? Surely, he has 2 choices for a skirt and for each of these choices, he has four possibilities to pick pants. In total, he has

$$
2 \times 4
$$

\bigskip \noindent ways to dress.\bigskip

\noindent We applied the following general counting principle.\\

\noindent \textbf{Proposition 2.1}. Suppose that the set $\Omega$ of size $n$ can be partitioned into $n_{1}$ subsets
$\Omega _{i}$, $i = 1, . . .,n_{1}$ of same size; and that each of these subsets $\Omega _{i}$ can be split into $n_{2}$ subsets $\Omega _{ij}$, j = 1, . . ., $n_{2}$, $i = 1, . . ., n_{1}$, of same size, and that each of the $\Omega _{ij}$ can be divided into $n_{3}$ subsets
$\Omega _{ijh}$, $h = 1, . . ., n_{3}$, $j = 1, . . ., n_{2}$, $i = 1, . . ., n_{1}$ of same size also. Suppose that we may proceed like that
up an order $k$ with $n_k$ subsets with common size $a$.\\

\bigskip \noindent Then the cardinality of $\Omega$ is given by\\

\begin{equation*}
n=n_{1}  \times n_{2} \times \times ... \times n_{k} \times a.
\end{equation*}

\bigskip

\noindent \textbf{Proof}. Denote by $B_{h}$ the cardinality of a subset generated at step $h$, for $h=0,...,k$. A step $h$,
we have $n_h$ subsets partitioned into $n_{h+1}$ subsets of same size. Then we have

\begin{equation*}
B_{h}=n_{h+1} \times B_{h+1} \text{ for all } 0 \leq h\leq k-1,
\end{equation*}

\Bin with $B_{0} = n$, $B_{k+1} = a$. Now, the proof is achieved by induction in the following way

\begin{eqnarray*}
n&=&B_{0}=n_{1} \times B_{1}\\
&=&n_{1} \times n_{2}  B_{2}\\
&=&n_{1} \times n_{2} \times n_{3} B_{3}\\
&...& \\
&=&n_{1} \times n_{2} \times n_{3} \times .... \times n_k \times B_k,
\end{eqnarray*}

\bigskip \noindent with $B_k=a$.\\

\noindent Although this principle is simple, it is a fundamental tool in combinatorics  \index{combinatoric} : divide and count.\\

\noindent But it is not always applied in a so simple form. Indeed, partitioning is the most important skill to develop in order to apply the principle successfully.\\

\noindent  This is what we will be doing in all this chapter.

\noindent Let $E$ be a set of $n$ elements with  $n\geq 1$.\\

\noindent \textbf{Definition 2.1}. A $p$-tuple of $E$ is an ordered subset of $E$ with $p$ distinct elements of $E$.
A $p$-tuple is called a $p$-permutation  \index{permutation} or $p$-arrangement  \index{arrangement} of elements of $E$ and the number of $p$-permutations of elements of $E$ is denoted
as $A_{n}^{p}$ (we read $p$ before $n$). It is also denoted by

$$
(n)_p.
$$

\bigskip \noindent We have the result.\\

\noindent \textbf{Theorem 2.1}. For all  $1\leq  p \leq n$, we have

\begin{equation*}
A_{n}^{p}=n(n-1)(n-2)......(n-p+1).
\end{equation*}

\bigskip \noindent \textbf{Proof}. Set $E = \left\{ x_{1},...,x_{n}\right\} $. Let $\Omega$ be the class of all ordered subsets of $p$
elements of $E$. We are going to apply the counting principle to $\Omega$.\\

\noindent It is clear that $\Omega$ may be divided into $n = n_{1}$ subsets $F_{i}$, where each $F_{i}$ is the class of ordered subsets of $\Omega$
with first elements $x_{i}$. Since the first element is fixed to $x_i$, the cardinality of $F_i$ is the number of ordered subsets of $E \setminus \{x_i\}$, so that the classes $F_{i}$ have a common cardinality which is the number of $(p-1)$-permutations  \index{permutation} from a set of $(n-1)$ elements. We have proved that

\begin{equation}
A_{n}^{p}=n \times A_{n-1}^{p-1}, \label{proba01.f01}
\end{equation}

\bigskip \noindent for any $1\leq p \leq n$. We get by induction

$$
A_{n}^{p}=n \times (n-1) \times A_{n-2}^{p-2}
$$

\Ni and after $h$ repetitions of  (\ref{proba01.f01}), we arrive at

$$
A_{n}^{p}=n \times (n-1) \times (n-2) \times ... \times (n-h) \times A_{n-h-1}^{p-h-1}.
$$

\bigskip \noindent For $h=p-2$, we have

$$
A_{n}^{p}=n \times (n-1) \times (n-2) \times ... \times (n-p+2) \times A_{n-p+1}^{1}.
$$

\bigskip \noindent And, clearly, $A_{n-p+1}^{1}=n-p+1$ since $A_{n-p+1}^{1}$ is the number of singletons form a set of $(n-p+1)$ elements.

\bigskip \noindent \textbf{Remark}. Needless to say, we have $A_{n}^{p}=0$ for $p\geq n+1$.\\

\noindent Here are some remarkable values of $A_{n}^{p}$. For any positive integer $n\geq 1$, we have\\

\noindent \textbf{(i)}  $A_{n}^{0} = 1$.\\

\noindent \textbf{(ii)} $A_{n}^{1} = n$.\\

\bigskip \noindent From an algebraic point of view, the numbers $A_{n}^{p}$ also count the number of injections.

\bigskip \noindent We begin by recalling the following algebraic definitions :\\

\noindent \textit{A function} from a set $E$ to a set $F$ is a correspondence from $E$ to $F$ such that each element of $E$ has at most one image in $F$.\\

\noindent \textit{A mapping} from a set $E$ to a set $F$ is a correspondence from $E$ to $F$ such that each element of $E$ has exactly one image in $F$.\\

\noindent \textit{An injection} from a set $E$ to a set $F$  is a mapping from $E$ to $F$ such that each any two distinct elements of $E$ have distinct images in $F$.\\

\noindent \textit{A surjection} from a set $E$ to a set $F$  is a mapping from $E$ to $F$ such that each element of $F$ is the image of at least one element of $E$.\\

\noindent \textit{A bijection} from a set $E$ onto a set $F$  is a mapping from $E$ to $F$ such that each element of $F$ is the image of one and only element of $E$.\\

\noindent  If there is an injection (respectively a bijection) from $E$ to $F$, then we have the inequality: $Card(E) \leq Card(F)$ (respectively, the equality $Card(E) = Card(F)$).\\

\noindent \textbf{The number of $p$-arrangements  \index{arrangement} from $n$ elements ($p \leq n$) is the number of injections from a set of $p$ elements to a set of $n$ elements}.\\

\bigskip \noindent The reason is the following. Let $E=\{y_{1},...,y_{p}\}$ and $F=\{x_{1},x_{2,}...,x_{n}\}$ be sets with $Card(E) = p$ and $Card(F) = n$ with $p \leq n$. Forming an injection $f$ from $E$ on $F$ is equivalent to choosing a $p$-tuple $(x_{i_{1}}...x_{i_{p}})$ in $F$
and to set the following correspondence $f(y_{h}) = x_{i_{h}}$, $h=1,...,p$. So, we may find as many injections from
$E$ to $F$ as $p$-permutations  \index{permutation} of elements of from $F$.\\

\noindent Thus, $A_{n}^{p}$ is also the number of injections from a set of $p$ elements on a set of $n$ elements.\\

\bigskip \noindent \textbf{The number of mappings  \index{number of mapping} from a set of $p$ elements to a set of $n$ elements is $n^p$}.\\

\noindent Indeed, let $E=\{y_{1},...,y_{p}\}$ and $F=\{x_{1},x_{2,}...,x_{n}\}$ be sets with $Card(E) = p$ and $Card(F) = n$ with no relation between $p$ and $n$. Forming a mapping $f$ from $E$ on $F$ is equivalent to choosing, for any $x\in E$, one arbitrary element of $F$ and to assign it to $x$ as its image. For the first element $x_{1}$ of $E$, we have $n$ choices, $n$ choices also  for the second $x_{2}$, $n$ choices also for the third $x_{3}$, and so forth. In total, we have  choices to form a mapping from $E$ to $F$.\\

\noindent \textbf{Example on vote casting in Elections}. The $p$ members of some population are called to mandatory cast a vote for one of $n$ candidates at random. The number of possible outcomes is the number of injections from a set of $p$ elements to a set of $n$ elements : $n^p$.

\noindent \textbf{Definition 2.2.}
A permutation  \index{permutation} or an ordering of the elements of $E$ is any ordering of all its elements. If $n$ is the cardinality of $E$, the number of permutations of the element of $E$ is called : \textit{n factorial}, and denoted $n!$ (that is $n$ followed with an exclamation point).\\

\bigskip \noindent Before we give the properties of \textit{n factorial}, we point out that a permutation  \index{permutation} of $n$ objects is a $n$-permutation of $n$ elements.\\

\noindent \textbf{Theorem 2.2.} For any $n \geq 1$, we have\\

\noindent \textbf{(i)}  $n! = n (n-1) (n-2) ...... 2 \times  1$.\\

\noindent \textbf{(ii)} $n! = n (n-1)!$.\\

\bigskip

\noindent \textbf{Proof.} By remembering that a permutation  \index{permutation} of $n$ objects is an $n$-permutation of $n$ elements, we may see that Point (i)  is obtained for $p=n$ in the Formula of Theorem 2.1. Point (ii) obviously derives from (i) by induction.\\

\noindent \textbf{Exercise 2.1.} What is the number of bijections between two sets of common cardinalily $n\geq 1$?\\

\noindent \textbf{Exercise 2.2.} Check that for any $1 \leq  p \leq n$,\\
$$
A_{n}^{p} = n! / (n-p)!.
$$

Let $\mathcal{E}$ be collection of $n$ distinct objects. Suppose that $\mathcal{E}$ is partitioned into $k$ sub-collections
$\mathcal{E}_{j}$ of respective sizes $n_{1}$, . . .$n_{k}$ the following properties apply :\\

\noindent \textbf{(i)} Two elements of same sub-collection are indistinguishable between them.\\

\noindent \textbf{(ii)} Two elements from two different sub-collections are distinguishable one from the other.\\

\bigskip \noindent What is the number of permutations  \index{permutation} of $\mathcal{E}$? Let us give an example.\\

\noindent Suppose we have $n=20$ balls of the same form (not distinguishable by the form, not at sight nor by touch), with $n_1=6$ of them of red color, $n_2=9$ of blue color and $n_3=5$ of green color. And we can distinguish them only by their colors. In this context, two balls of the same color are the same for us, from our sight. An ordering of these $n$ balls in which we cannot distinguish the balls of same color is called a permutation  \index{permutation} with repetition or a visible \textit{permutation}.\\

\noindent Suppose that we have realized a permutation  \index{permutation} of these $n$ balls. Permuting for example only red balls between them does not change the appearance of the global permutation. Physically, it has changed but not visibly meaning  \index{mean} from our sight. Such a permutation may be described as \textit{visible}, or \textit{permutation with repetition}.\\

\noindent In fact, any of the $n!$ real permutations  \index{permutation} represents all the \textit{visible} repetitions where the $6$ red balls are permuted between them, the $9$ blue balls are permuted between them and the $5$ green balls are permuted between them. By the counting principle, a real permutation represents exactly

$$
9! \times 6! \times 5!
$$

\bigskip \noindent  permutations  \index{permutation} with repetition. Hence, the number of permutations with repetition of these $20$ balls is

$$
\frac{20!}{9! \times 6! \times 5!}.
$$

\bigskip \noindent Now, we are going to do the same reasoning in the general case.\\

\noindent As we said before, we have two types of permutations  \index{permutation} here :\\

\noindent \textbf{(a)} The real or physical permutations,  \index{permutation} where the $n$ elements are supposed to be distinguishable.\\

\noindent \textbf{(b)}  The permutations  \index{permutation} with repetition in which we cannot distinguish between the elements of a same sub-collection.\\

\bigskip

\noindent \textbf{Theorem 2.3}. The number of permutations  \index{permutation} with repetition of a collection $\mathcal{E}$ of $n$ objects partitioned into $k$
sub-collections $\mathcal{E}_{i}$ of size $n_{i}$, $i=1,...,k$, such that only elements from different sub-collections are distinguishable between them, is given by
\begin{equation*}
B(n_{1},n_{2},...,n_{k})=\frac{n!}{n_{1}!n_{2}!......n_{k}!}.
\end{equation*}

\bigskip \noindent \textit{Terminology}. The numbers of permutations  \index{permutation} with repetition are also called \textit{multinomial coefficients},  \index{multinomial coefficient} in reference Formula (\ref{proba01.ForMMonial01}).\\

\noindent \textbf{Proof.} Let $B(n_{1},n_{2},...,n_{k})$ be the number of permutations  \index{permutation} with repetition. Consider a fixed real permutation. This real permutation corresponds exactly to all permutations with repetition obtained by permuting the $n_{1}$ objects of $\mathcal{E}_{1}$ between them,
the $n_{2}$ objects of $\mathcal{E}_{2}$ between them, $\cdots$, and the $n_{k}$ objects of $\mathcal{E}_{k}$ between them. And we obtain $n_{1}!n_{2}!......n_{k}!$ permutations  \index{permutation} with repetition corresponding to the same real permutation. Since this is true for any real permutation which generates
$$
n_{1}! \times n_{2}!\times ... \times n_{k}!
$$

\noindent visible permutations,  \index{permutation} we get that

\begin{equation*}
(n_{1}!n_{2}!.....n_{k}!)\times B(n_{1},n_{2},...,n_{k}) =n!.
\end{equation*}

\bigskip \noindent This gives the result in the theorem.\\

\noindent Let us define the numbers of combinations  \index{combination} as follows.\\

\noindent \textbf{Definition 2.3.} Let $E$ be a set of $n$ elements. A combination  \index{combination} of $p$ elements of $E$ is a subset of $E$ of size $p$.\\

\noindent In other words, any subset $F$ of $E$ is a $p$-combination  \index{combination} of elements of $E$ if and only if $Card(F)=p$.\\

\noindent It is important to remark that we can not have combinations  \index{combination} of more that $n$ elements in a set of $n$ elements.\\

\noindent We have :\\

\noindent \textbf{Theorem 2.4.} The number of combinations  \index{combination} of $p$ elements from $n$ elements is given by \\
$$
\left(\begin{tabular}{c} $n$ \\ $p$ \end{tabular} \right)=\frac{n!}{p!(n-p)!}.
$$

\bigskip

\noindent \textbf{Proof}. Denote by
$$
\left(\begin{tabular}{c} $n$ \\ $p$ \end{tabular} \right)
$$

\Bin the number of combinations  \index{combination} of $p$ elements from $n$ elements.\\

\noindent The collection of $p$-permutations  \index{permutation} is exactly obtained by taking all the orderings of the elements of the combinations  \index{combination} of  $p$ elements from $E$. Each combination gives $p!$ $p$-permutations of $E$. Hence, the cardinality of the collection of $p$-permutations is exactly $p!$ times that of the class of $p$-combinations of $E$, that is

\begin{equation*}
p!\text{ }\left(\begin{tabular}{c} $n$ \\ $p$ \end{tabular} \right)=A_{n}^{p}.
\end{equation*}

\noindent By using Exercise 2.1 above, we have

\begin{equation*}
\left(\begin{tabular}{c} $n$ \\ $p$ \end{tabular} \right)=\frac{A_{n}^{p}\text{ }}{p!}=\frac{n!\text{ }}{\text{ }%
p!(n-p)!.}.
\end{equation*}

\bigskip \noindent We also have this definition :\\

\noindent \textbf{Definition 2.4.} The numbers $\left(\begin{tabular}{c} $n$ \\ $p$ \end{tabular} \right)$, are also called binomial coefficients  \index{binomial coefficient} because of Formula (\ref{proba01.ForBin}) below.

\bigskip

\noindent The urn model plays a very important role in discrete Probability Theory  \index{probability theory} and Statistics,  \index{statistics} especially in sampling theory. The simplest example of urn model is the one where we have a number of balls, distinguishable or not, of different colors.\\

\bigskip \noindent We are going to apply the concepts seen above in the context of urns.\\

\bigskip \noindent Suppose that we want to draw $r$ balls at random from an urn containing $n$ balls that are distinguishable by touch (where touch means  \index{mean} hand touch).\\

\bigskip \noindent We have two ways of drawing.\\

\noindent \textbf{(i)} \textit{Drawing without replacement}. This means  \index{mean} that we draw a first ball and we keep it out of the urn. Now, there
are $(n-1)$ balls in the urn. We draw a second and we have $(n-2)$ balls left in the urn. We repeat this procedure until we have the $p$ balls. Of course, $p$ should be less or equal to $n$.\\

\noindent \textbf{(ii)}  \textit{Drawing with replacement}. This means  \index{mean} that we draw a ball and take note of its identity or its characteristics (that are studied) and put it back in the urn. Before each drawing, we have exactly $n$ balls in the urn. A ball can be drawn several times.\\

\noindent It is clear that the drawing model \textbf{(i)} is exactly equivalent to the following one :\\

\noindent \textbf{(i-bis)} We draw $p$ balls at the same time, simultaneously, at once.\\

\bigskip \noindent Now, we are going to see how the $p$-permutations  \index{permutation} and the $p$-combinations  \index{combination} occur here, by a series of questions and answers.\\

\noindent \textbf{Questions.} Suppose that an urn contains $n$ distinguishable balls. We draw $p$ balls. In how many ways can the drawing occur? Or what is the number of possible outcomes in term of subsets formed by the $p$ drawn balls?\\

\noindent \textbf{Solution 1.} If we draw the $p$ balls \textbf{without replacement} and \textbf{we take into account the order}, the number of possible outcomes is the number of $p$-permutation \index{permutation} s.\\

\noindent \textbf{Solution 2.} If we draw the $p$ balls \textbf{without replacement} and \textbf{we do not take the order into account or there is no possible ordering}, the number outcomes is the number of $p$-combinations  \index{combination} from $n$.\\

\noindent \textbf{Solution 3.} If we draw the $p$ balls \textbf{with replacement}, the number of outcomes is $n^{p}$, the number of mapping \index{number of mapping} s.

\bigskip \noindent \textbf{Needless to say, the ordering is always assumed if we proceed by a drawing with replacement}.\\

\noindent \textbf{Please, keep in mind these three situations that are the basic keys in Combinatoric \index{combinatoric} s}.\\

\noindent Now, let us explain the solutions before we continue.\\

\noindent \textbf{Proof of Solution 1}. Here, we draw the $p$ balls one by one. We have $n$ choices for the first ball. Once this ball is out, we have $(n-1)$ remaining balls in the urn and we have $(n-1)$ choices for the second ball. Thus, we have

$$
n \times (n-1)
$$

\bigskip \noindent possible outcomes to draw two ordered balls. For three balls, $p=3$, we have the number

$$
n \times (n-1) \times (n-2).
$$

\bigskip \noindent \textbf{Remark} that for $p=1,2,3$, the number of possible outcomes is

$$
n \times (n-1) \times (n-2) \times .... \times (n-p+1) = A_n^p.
$$

\bigskip \noindent You will not have any difficulty to get this by induction.\\

\noindent \textbf{Proof of Solution 2}. Here, there is no ordering. So we have to divide the number of ordered outcomes $A_n^p$ by $p!$ to
get the result.\\

\noindent \textbf{Proof of Solution 3}. At each step of the $p$ drawing, we have $n$ choices. At the end, we get $n^p$ ways to draw $p$ elements.\\

\noindent We are going to devote a special subsection to the numbers of combinations  \index{combination} or binomial coefficient \index{binomial coefficient} s.\\

\noindent Here, we come back to the number of combinations  \index{combination} that we call Binomial Coefficients  \index{binomial coefficient} here. First, let us state their main properties.\\

\bigskip

\noindent \textbf{Main Properties}.\\

\noindent \textbf{Proposition 2.1.} We have\\
\noindent \textbf{(1)} $\left(\begin{tabular}{c} $n$ \\ $0$ \end{tabular} \right) = 1$ for all $n \geq  0$.\\
\noindent \textbf{(2)} $\left(\begin{tabular}{c} $n$ \\ $1$ \end{tabular} \right) = n$ for all $n \geq  1$.\\
\noindent \textbf{(3}  $\left(\begin{tabular}{c} $n$ \\ $p$ \end{tabular} \right)=\left(\begin{tabular}{c} $n$ \\ $n-p$ \end{tabular} \right)$, for $0 \leq  p\leq  n$.\\
\noindent \textbf{(4)}
$\left(\begin{tabular}{c} $n-1$ \\ $p-1$ \end{tabular} \right)+\left(\begin{tabular}{c} $n-1$ \\ $p$ \end{tabular} \right)=\left(\begin{tabular}{c} $n$ \\ $p$ \end{tabular} \right)$, for all $1 \leq  p \leq n$.\\

\bigskip

\noindent \textbf{Proof.} Here, we only prove Point (4). The other points are left to the reader as exercises.\\

\begin{eqnarray*}
\left(\begin{tabular}{c} $n-1$ \\ $p-1$ \end{tabular} \right)+\left(\begin{tabular}{c} $n-1$ \\ $p$ \end{tabular} \right)&=&\left\{ \frac{(n-1)!}{(p-1)!(n-p)!}\right\}+\left\{ \frac{(n-1)!}{p!(n-p-1)!}\right\}\\
&=&p\left\{ \frac{(n-1)!}{p!(n-p)!} \right\} +(n-p)\left\{ \frac{(n-1)!}{p!(n-p)!}\right\}\\
&=&\left\{ \frac{(n-1)!}{p!(n-p)!}\right\} \left\{p+n-p\right\}\\
&=&\frac{n(n-1)!}{p!(n-p)!}\\
&=&\frac{n!}{p!(n-p)!}=\left(\begin{tabular}{c} $n$ \\ $p$ \end{tabular} \right). \ \ QED
\end{eqnarray*}

\bigskip

\noindent \textbf{Pascal's Triangle}.  \index{Pascal's triangle} We are going to reverse the previous way by giving some of these properties as characteristics of the binomial coefficients.  \index{binomial coefficient} We have\\

\noindent \textbf{Proposition 2.3.} The formulas (i) and (ii):

\noindent \textbf{(i)} $\left(\begin{tabular}{c} $n$ \\ $0$ \end{tabular} \right) = \left(\begin{tabular}{c} $n$ \\ $n$ \end{tabular} \right) =1$ for all $n \geq  0$;\\
\noindent \textbf{(ii)}
$\left(\begin{tabular}{c} $n-1$ \\ $p-1$ \end{tabular} \right)+\left(\begin{tabular}{c} $n-1$ \\ $p$ \end{tabular} \right)=\left(\begin{tabular}{c} $n$ \\ $p$ \end{tabular} \right)$, for all $1 \leq  p \leq n$\\

\Bin entirely characterize the binomial coefficient \index{binomial coefficient} s
$$
\left(\begin{tabular}{c} $n$ \\ $p$ \end{tabular} \right), \ \ 1 \leq p \leq  n.
$$

\bigskip \noindent \textbf{Proof.} We give the proof by using Pascal's triangle.  \index{Pascal's triangle} Point (ii) gives the following clog rule (\textit{r\`egle du sabot in French})

\begin{equation*}
\begin{tabular}{|l|l|l|}
\hline
n / p & p-1 & p \\ \hline
n-1 & $\left(\begin{tabular}{c} $n-1$ \\ $p-1$ \end{tabular} \right)$ & $\left(\begin{tabular}{c} $n-1$ \\ $p$ \end{tabular} \right)$ \\ \hline
n &  & $\left(\begin{tabular}{c} $n$ \\ $p$ \end{tabular} \right)=\left(\begin{tabular}{c} $n-1$ \\ $p-1$ \end{tabular} \right)+\left(\begin{tabular}{c} $n-1$ \\ $p$ \end{tabular} \right)$ \\ \hline
\end{tabular}
\end{equation*}

\bigskip \noindent or more simply

\begin{equation*}
\begin{tabular}{|l|l|}
\hline
$u$ & $v$ \\ \hline
& $u+v$ \\ \hline
\end{tabular}
\end{equation*}

\bigskip

\noindent With that rule, we may construct the Pascal's triangle  \index{Pascal's triangle} of the numbers $\left(\begin{tabular}{c} $n$ \\ $p$ \end{tabular} \right)$.

\bigskip

\begin{equation*}
\begin{tabular}{|c|c|c|c|c|c|c|c|c|c|}
\hline
$n / p$ & 0 & 1 & 2 & 3 & 4 & 5 & 6 & 7 & 8 \\
\hline
0 & \textbf{1} &  &  &  &  &  &  &  &  \\
\hline
1 & \textbf{1}=\textit{u} & \textbf{1}=\textit{v} &  &  &  &  &  &  &  \\
\hline
2 & \textbf{1} & $\overset{u+v}{2}$ & \textbf{1} &  &  &  &  &  &  \\
\hline
3 & \textbf{1} & \textit{3} & \textit{3} & \textbf{1} &  &  &  &  &  \\
\hline
4 & \textbf{1} & \textit{4=u} & \textit{6=v} & \textit{4} & \textbf{1} &  &
&  &  \\
\hline
5 & \textbf{1} & \textit{5} & $\overset{u+v}{10}$ & \textit{10} & \textit{5}
& \textbf{1} &  &  &  \\
\hline
6 & \textbf{1} &  &  &  &  &  & \textbf{1} &  &  \\
\hline
7 & \textbf{1} &  &  &  &  &  &  & \textbf{1} &  \\
\hline
8 & \textbf{1} &  &  &  &  &  &  &  & \textbf{1}\\
\hline
\end{tabular}
\end{equation*}

\bigskip

\noindent The reader is asked to continue to fill this triangle himself. Remark that filling the triangle only requires the first column ($n=0$), the diagonal ($p=n$) and the clog rule. So points (i) and (ii) are enough to determine all the binomial coefficients.  \index{binomial coefficient} This leads to the following conclusion.\\

\noindent \textbf{Proposition 2.4}. Any array of integers $\beta (p,n)$,  $0 \leq  p \leq n$ such that\\

\noindent \textbf{(i)}  $\beta (n,n)$=$\beta (0,n)$= 1, for all $n \geq 0$,\\

\noindent \textbf{(ii)} $\beta(p - 1, n - 1) + \beta(p, n - 1) = \beta(p, n),$ for 1 $\leq  p \leq n$.\\

\noindent is exactly the array of binomial coefficients,  \index{binomial coefficient} that is
$$
\beta(p,n)=\left(\begin{tabular}{c} $n$ \\ $p$ \end{tabular} \right),  \ \ 0 \leq  p \leq n.
$$

\bigskip \noindent We are going to visit the  Newton's formula \index{Newton's formula} .\\

\noindent \textbf{The Newton's Formula \index{Newton's formula} }.\\

\noindent Let us apply the result just above to the power $n$ of a sum of two scalars in a commutative ring ($\mathbb{R}$ for example).

\bigskip

\begin{theorem} \label{theo2.4} For any $(a, b)\in \mathbb{R}^2$, for any $n\geq 1$, we have
\begin{equation}
(a+b)^n=\sum_{p=0}^{n}\ \left(\begin{tabular}{c} $n$ \\ $p$ \end{tabular} \right) a^{p\ }b^{n-p}. \label{proba01.ForBin}
\end{equation}
\end{theorem}

\bigskip

\noindent \textbf{Proof}. Since $\mathbb{R}$ is a commutative ring, we know that $(a+b)^{n}$ is a polynomial in $a$ and $b$ and it is written as a linear combination  \index{combination} of terms  $a^{p} \times b^{n-p}$, $p = 0$,...,$n$. We have the formula

\begin{equation}
(a+b)^{n}=\sum_{p=0}^{n}\ \beta (p,\ n)\ a^{p\ }b^{n-p}. \label{proba01.ForBin}
\end{equation}

\Bin It will be enough to show that the array $\beta(p,\ n)$ is actually that of the binomial coefficients.  \index{binomial coefficient} To begin, we write

\begin{equation*}
(a+b)^{n}=(a+b)\times (a+b)\times ...\times (a+b).
\end{equation*}

\bigskip \noindent  From there, we see that $\beta(p, n)$ is the number of choices of $a$ or $b$ in each factor $(a + b)$ such that $a$ is chosen $p$ times and $b$ is chosen $(n - p)$ times. Thus $\beta (n,\ n) = 1$ since $a^{n}$ is obtained in the unique case where $a$ is chosen in each factor $(a + b)$. Likely $\beta (0,\ n) = 1$, since this corresponds to the monomial $b^{n}$, that is the unique case where $b$ is chosen in each case. So, Point (i) is proved for the array $\beta(\cdot, \cdot)$. Next, we have

\begin{eqnarray*}
(a+b)^{n} &=&(a+b)(a+b)^{n-1}\\
&=&(a+b)\times \sum_{p=0}^{n-1}\ \beta(p,\ n-1) a^{p}\ b^{n-p-1}\\
&=&(a+b) \left(\cdots +\beta (p-1,\ n-1)\ a^{p-1}b^{n-p} \right.\\
&+& \left. \cdots + \beta (p,\ n-1) a^{p}\ b^{n-p-1}+ \cdots \right)
\end{eqnarray*}

\bigskip

\noindent This means  \index{mean} that, when developing $(a + b)(a + b)^{n-1}$, the term $a^{p} \times b^{n-p}$ can only come out \\

\noindent (1) either from the product of $b$ by $a^{p} \times b^{n-p-1}$ of the binomial $(a + b)^{n-1}$,\\

\noindent (2) or from the product of $a$ by $a^{p-1} \times b^{n-p}$ of the binomial $(a + b)^{n-1}$.\\

\noindent Then we see that for $ \leq  p \leq n$, we get

\begin{equation*}
\beta (p-1,\ n-1)+\beta (p,\ n-1)=\beta (p,\ n).
\end{equation*}

\bigskip

\noindent We conclude that the array $\beta (p, n)$ fulfills Points (i) and (ii) above. Then, this array is that of the binomial coefficients.  \index{binomial coefficient} QED.\\

\bigskip \noindent \textbf{Remark}. The name of binomial coefficients  \index{binomial coefficient} comes from this Newton's formula \index{Newton's formula} .\\

\bigskip \noindent \textbf{Multiple Newton's Formula \index{Newton's formula} }.\\

\noindent We are going to generalize the Newton's Formula  \index{Newton's formula} from dimension $k=2$ to an arbitrary dimension $k\geq 2$. Then binomial coefficients  \index{binomial coefficient} will be replaced by the numbers of permutations  \index{permutation} with repetition.\\

\noindent Let $k\geq 2$ and let us be given real numbers $a_1$,$a_2$, ..., and  $a_k$ and let $n\geq 1$ be a positive number. Consider

$$
\Gamma_n =\{(n_1,...,n_k), n_1 \geq 0,...,n_k \geq 0, n_1 + ... + n_k=n \}.
$$

\Bin We have

\begin{equation}
(a_1+...+a_k)^n=\sum_{(n_1,...,n_k) \in \Gamma_n} \frac{n!}{(n_1! \times n_2 ! \times ... \times n_k !)} a_{1}^{n_1} \times a_{2}^{n_2} \times ... \times a_{k}^{n_k}, \label{proba01.ForMMonial01}
\end{equation}

\bigskip \noindent that we may write in a more compact manner in

\begin{equation}
\left(\sum_{i=1}^{k} a_i\right)^n=\sum_{(n_1,...,n_k) \in \Gamma_n} \frac{n!}{\prod_{i=1}^{k} n_i !} \prod_{i=1}^{k} a_{i}^{n_i}. \label{proba01.ForMMonial02}
\end{equation}

\bigskip \noindent We will see how this formula is important for the multinomial law, which in turn is so important in Statistics \index{statistics} .\\

\bigskip \noindent \textbf{Proof}. Let us give a simple proof of it.\\

\noindent To develop $(\sum_{i=1}^{k} a_i)^n$, we have to multiply $(\sum_{i=1}^{k} a_i)$ by itself $n$ times.  By the distributivity of the sum with respect to the product, the result will be a sum of products

$$
z_1 \times z_2 \times ... \times z_n,
$$

\Bin where each $z_i$ is one of the $a_1$, $a_2$, ..., $a_k$. By commutativity, each of these products is of the form

\begin{equation}
a_{1}^{n_1} \times a_{2}^{n_1} \times ... \times a_{k}^{n_k}, \label{proba01.multi01}
\end{equation}

\noindent where $(n_1,...,n_k) \in \Gamma_n$. And for a fixed $(n_1,...,n_k) \in \Gamma_n$,  the product (\ref{proba01.multi01}) is the same as all products

$$
z_1 \times z_2 \times ... \times z_n,
$$

\Bin in which we have $n_1$ of the $z_i$ identical to $a_1$, $n_2$ identical to $a_2$, ..., and $n_k$ identical to  $a_k$. These products correspond to  the permutations  \index{permutation} with repetition of $n$ elements such that $n_1$ are identical, $n_2$ are identical, ..., and $n_k$ are identical. Then, each product (\ref{proba01.multi01}) occurs

$$
\frac{n!}{(n_1! \times n_2 ! \times ... \times n_k !)}
$$

\Bin times in the expansion. This puts an end to the proof.\\

\noindent The number $n!$ grows and becomes huge very quickly. In many situations, it may be handy to have an asymptotic equivalent formula.\\

\Ni This formula is the Sterling's one and it is given as follows : \\

\begin{equation*}
n!=(2 \pi  n)^{\frac{1}{2}}\ \left(\frac{n}{e}\right)^{n}\ \exp (\theta _{n}),
\end{equation*}

\bigskip \noindent with, for any $\eta>0$, for $n$ large enough,

\begin{equation*}
\ \left| \theta _{n}\right| \leq \frac{1+\eta}{12n}.
\end{equation*}

\bigskip \noindent This implies, in particular, that

\begin{equation*}
n!\sim (2 \pi n)^{\frac{1}{2}}\ \left(\frac{n}{e}\right)^{n},
\end{equation*}

\noindent as $n\ \rightarrow \ \infty $.\\

\noindent One can find several proofs (See \cite{feller1} \index{Feller}, page 52, for example). In this textbook, we provide a proof in the lines of the one in \cite{valiron} \index{Valiron}, pp. 167, that is based on Wallis integrals.  \index{Wallis integral} This proof is exposed in Chapter \ref{proba01_appendix} (which is an appendix), Section
\ref{proba01_appendix_stirling}.\\

\noindent We think that a student in first year of University will be interested by an application of the course on Riemann integration.\\

 %proba01.combin
\chapter{Introduction to Probability Measures}  \index{probability measure} \label{proba01.pm}

Assume that we have a perfect die whose six faces are numbered from $1$ to $6$. We want to toss it twice. Before we toss it, we know that the outcome will be a couple $(i,j)$, where $i$ is the number that will appear first and $j$ the second.\\

\noindent We always keep in mind that, in probability theory,  \index{probability theory} we will be trying to give answers about events that have not occurred yet. In the present example, the possible outcomes form the set

\begin{equation*}
\Omega =\left\{ 1,2,...,6\right\} \times \left\{
1,2,...,6\right\} =\left\{ (1,1),(1,2),...,(6,6)\right\}.
\end{equation*}

\bigskip \noindent $\Omega $ is called the \textit{sample} space of the experiment or the \textit{probability space}.  \index{probability space} Here, the size of the set $\Omega$ is finite and is exactly $36$. Parts or subsets of $\Omega$ are called events. For example,\\

\noindent (1) \{(3, 4)\} is the event : \textit{Face 3 comes out in the first tossing  \index{tossing} and Face 4 in the second},\\

\noindent (2)  A = $\left\{ (1,1),(1,2),(1,3),(1,4),(1,5),(1,6)\right\} $ is the event  : \textit{1 comes out in the first tossing \index{tossing} }.\\

\noindent Any element of $\Omega$, as a singleton, is an elementary event.  \index{elementary event} For instance $\{(1,1)\}$ is the elementary event :
\textit{Face 1 appears in both tossing \index{tossing} }.\\

\noindent In this example, we are going to use the perfectness of the die, and the regularity of the geometry of the die, to the conviction that \\

\noindent (1) All the elementary events  \index{elementary event} have equal chances of occurring, that is one chance out of 36.\\

\noindent (2)  Each event $A$ of  $\Omega$ has a number of chances of occurring, which equal to its cardinality.\\

\noindent We recall that an event $A$ occurs if and only if the occurring elementary event  \index{elementary event} is in $A$.\\

\bigskip

\noindent Denote $\mathbb{P}(A)$ to be the fraction of the number of chances of occurrence of $A$ over the total number of chances $n=36$, i.e.,

\begin{equation*}
\mathbb{P}(A)=\frac{Card (A)}{36}\text{.}
\end{equation*}

\bigskip \noindent Here, we say that $\mathbb{P}(A)$ is the probability that the event $A$ occur after the tossing \index{tossing} s.\\

\noindent We may easily check the following facts.\\

\noindent \textbf{(1)}   $\mathbb{P}(\Omega)=1$ and $0\leq \mathbb{P}(A) \leq 1$ for all $A\subseteq \Omega $.\\

\noindent \textbf{(2)}  For all $A$, $B$, parts of $\Omega$ such that $A \cap B = \emptyset $, we have

\begin{equation*}
\mathbb{P}(A \cup B) = \mathbb{P}(A)+\mathbb{P}(B).
\end{equation*}

\bigskip \noindent \textbf{Notation} : if $A$ and $B$ are disjoint,  \index{joint} we adopt the following convention and write :

$$
A \cup B=A+B.
$$

\Bin As well, if $(A_n)_{n\geq 0}$ is a sequence of pairwise disjoint  \index{joint} events, we write

$$
\bigcup_{n\geq 0} A_n = \sum_{n\geq 0} A_n.
$$

\bigskip \noindent We summarize this by saying : \noindent \textbf{We may use the symbol + (plus) in place of the symbol $\cup$ (union), when the sets are mutually  disjoint,  \index{joint} that is pairwise disjoint.}\\

\noindent The so-defined mapping $\mathbb{P}$ is called a probability measure  \index{probability measure} because of (1) and (2) above.\\

\noindent  If the space $\Omega$ is infinite, (2) is written as follows.\\

\noindent \textbf{(2)} For any sequence of events $(A_n)_{n\geq 0}$ pairwise disjoint,  \index{joint} we have

$$
\mathbb{P}\left(\sum_{n\geq 0} A_n\right) = \sum_{n\geq 0} \mathbb{P}(A_n)
$$

\bigskip \noindent \textbf{Terminology}.\\

\noindent \textbf{(1)} The events $A$ and $B$ are said to be mutually exclusive if $A \cap B=\emptyset$. In other words, the events $A$ and $B$ cannot occur simultaneously.\\

\noindent \textbf{(2)} If we have that $\mathbb{P}(A)=0$, we say that the event $A$ is impossible, or that $A$ is a null-set with respect to
$\mathbb{P}$.\\

\noindent \textbf{(3)} If we have $\mathbb{P}(A)=1$, we say that the event $A$ a sure event with respect to $\mathbb{P}$, or that the event $A$ holds the probability measure  \index{probability measure} $\mathbb{P}$.\\

\bigskip \noindent \textbf{Nota-Bene}. For any event $A$, the number $\mathbb{P}(A)$ is a probability (that $A$ occur). But the application, that is the mapping, $\mathbb{P}$ is called a probability measure \index{probability measure} .\\

\noindent Now we are ready to present the notion of probability measures.  \index{probability measure} But we begin with discrete ones.

\bigskip

\bigskip \noindent Let  $\Omega$ be a set with finite cardinality or with infinite countable cardinality. Let $\mathcal{P}$($\Omega)$ be the class of parts of $\Omega $.\\

\noindent \textbf{Definition 1.} An mapping $\mathbb{P}$, defined from $\mathcal{P}$($\Omega $) to $\left[ 0,1\right]$ is called a \textit{probability measure}  \index{probability measure}  on $\Omega$ if and only if :\\

\noindent \textbf{(1)} $\mathbb{P}(\Omega) = 1.$\\

\noindent \textbf{(2)} For all sequences of pairwise disjoint  \index{joint} events $(A_{i},i\in \mathbb{N})$ of $\Omega$, we have

\begin{equation*}
\mathbb{P}\left( \sum_{i}A_{i}\right) =\sum_{i}\mathbb{P}(A_{i}).
\end{equation*}

\bigskip

\noindent We say that the triplet ($\Omega $, $\mathcal{P}$($\Omega $), $\mathbb{P}$) is a probability space \index{probability space} .\\

\noindent \textbf{Terminology}. Point (2) means  \index{mean} that $\mathbb{P}$ is additive on $\mathcal{P}$($\Omega $).  We refer to it under the name of \textit{additivity}.\\

\noindent Let $(\Omega, \mathcal{P}(\Omega), \mathbb{P})$ a probability space.  \index{probability space} We have the following properties. Each of them will be proved just after its statement.\\

\noindent \textbf{(A)}  $\mathbb{P}(\emptyset ) = 0$.\\

\noindent \textbf{Proof : }. By additivity, we have

$$
\mathbb{P}(\emptyset) =\mathbb{P}(\emptyset +\emptyset)  =\mathbb{P}(\emptyset) + \mathbb{P}(\emptyset) = 2 \mathbb{P}(\emptyset).
$$

\bigskip \noindent  Then, we get $\mathbb{P}(\emptyset)=0.$\\

\noindent \textbf{(B)}  If $(A,B)\in \mathcal{P}(\Omega)^2$ and if $A \subseteq  B$,
then $\mathbb{P}(A) \leq \mathbb{P}(B)$.\\

\noindent \textbf{Proof.} Recall the definition of the difference of subsets :

$$
B \setminus A=\{x \in \Omega, x \in B \text{ and } x \notin A\}=B \cap A^c.
$$

\bigskip \noindent Since $A \subset B$, we have

\begin{equation*}
B=(B \setminus A)+A.
\end{equation*}

\bigskip \noindent Thus, by additivity,
\begin{equation*}
\mathbb{P}(B)=\mathbb{P}((B \setminus A)+A)=\mathbb{P}(B \setminus A)+\mathbb{P}(A).
\end{equation*}

\bigskip \noindent It follows that
\begin{equation*}
\mathbb{P}(B)-\mathbb{P}(A)=\mathbb{P}(B\setminus A)\geq 0.
\end{equation*}

\bigskip \noindent \textbf{(C)}  If $(A,B)\in \mathcal{P}(\Omega)^2$ and if $A \subseteq $ B,
through $\mathbb{P}(B\setminus A)=\mathbb{P}(B)-\mathbb{P}(A)$.\\

\noindent \textbf{Proof.} This is already proved through (B).\\

\noindent \textbf{(D)} (\textbf{Continuity Property of a Probability measure})  \index{probability measure} Let $(A_{n})_{n\geq 0}$ be a non decreasing sequence of subsets of $\Omega$ with limit $A$, that is
:\\

\noindent (1) For all $n \geq $ 0, A$_{n}$ $\subseteq  A_{n+1}$\\

\noindent and\\

\noindent (2) $\cup _{n\geq 0\text{ }}$A$_{n} = A$.\\

\noindent Then  $\mathbb{P}$(A$_{n}$) $\uparrow \mathbb{P}(A)$ as $n \uparrow \infty$.\\

\bigskip \noindent \textbf{Proof : } Since the sequence $(A_{j})_{j\geq 0}$ is non-decreasing, we have

\begin{equation*}
A_{k}=A_{0}+(A_{1}\setminus A_{0})+(A_{2}\setminus
A_{1})+...........+(A_{k}\setminus A_{k-1}),
\end{equation*}

\bigskip \noindent for all $k\geq 1$. Finally, we have

\begin{equation*}
A =A_{0}+(A_{1}\setminus A_{0})+(A_{2}\setminus
A_{1})+...........+(A_{k}\setminus A_{k-1})+.......
\end{equation*}

\bigskip \noindent Denote $B_{0}=A_{0}$ and  $B_{k}=A_{k}\setminus A_{k-1}$, for $k \geq  1$. By using the additivity of $\mathbb{P}$,
we get

$$
\mathbb{P}(A)=\sum_{j\geq 0}\mathbb{P}(B_{j})=\lim_{k\rightarrow \infty }\sum_{0\leq j\leq k}\mathbb{P}(B_{j}).
$$

\bigskip \noindent But

\begin{eqnarray*}
\sum_{0\leq j\leq k}\mathbb{P}(B_{j})&=&\mathbb{P}(A_{0})+\sum_{1\leq j\leq k}\mathbb{P}\left( A_{j}\setminus A_{j-1}\right)\\
&=&\mathbb{P}(A_{0})+\sum_{1\leq j\leq k}\mathbb{P}(A_{j})-\mathbb{P}(A_{k-1})\\
&=&\mathbb{P}(A_{0})+(\mathbb{P}(A_{1})-\mathbb{P}(A_{0}))+(\mathbb{P}(A_{2})\\
&-&\mathbb{P}(A_{1}))+...+(\mathbb{P}(A_{k})-\mathbb{P}(A_{k-1}))\\
&=&\mathbb{P}(A_{k}).
\end{eqnarray*}

\bigskip \noindent We arrive at

\begin{equation*}
\mathbb{P}(A)=\sum_{j\geq 0}\mathbb{P}(B_{k})=\lim_{k\rightarrow
\infty }\sum_{0\leq j\leq k}\mathbb{P}(B_{j})\text{=}\lim_{k\rightarrow
\infty }\mathbb{P}(A_{k}).
\end{equation*}

\bigskip \noindent Hence, we have

\begin{equation*}
\lim_{n\rightarrow \infty }\mathbb{P}(A_{n})=\mathbb{P}(A).
\end{equation*}

\bigskip \noindent Taking the complements of the sets of Point (D), leads to the following point.\\

\noindent \textbf{(E)} (\textbf{Continuity Property of a Probability measure})  \index{probability measure} Let $(A_{n})_{n\geq 0}$ be a sequence of non-increasing subsets of $\Omega $
to $A$, that is,\\

\noindent (1)   For all $n \geq $ 0, $A_{n+1}$ $\subseteq  A_{n}$,\\

\noindent and\\

\noindent (2)  $\cap _{n\geq 0}A_{n} = A$.\\

\noindent Then $\mathbb{P}(A_{n}) \downarrow \mathbb{P}(A)$ when  $n \uparrow \infty $.

\bigskip

\noindent At this introductory level, we usually work with discrete probabilities defined on an enumerable space $\Omega$. A
probability measure  \index{probability measure} on such a space is said to be discrete.

\noindent Let $Card (\Omega) \leq  Card(\mathbb{N})$, meaning  \index{mean} that $\Omega$ is enumerable, meaning also that we may
write $\Omega$ in the form : $\Omega= \left\{ \omega _{1},\omega _{2},.....\right\}$.\\

\noindent The following theorem allows to build a probability measure  \index{probability measure} on discrete spaces.\\

\noindent \textbf{Theorem}. Defining discrete  probability measure  \index{probability measure} $\mathbb{P}$ on $\Omega$, is equivalent to providing numbers
$p_{i}$, $1\leq i$, such that $0\leq p_{i} \leq 1$ and $p_{1}+p_{2}+...=1$ so that, for any subset of $A$ of $\Omega$,

\begin{equation*}
\mathbb{P}(A)=\sum_{\omega _{i}\in A}p_{i}.
\end{equation*}

\bigskip

\noindent \textbf{Proof.} Let $\mathbb{P}$ be a probability measure  \index{probability measure} on
$\Omega=\left\{ \omega _{i},i\in I\right\} $, I $\subseteq $ $\mathbb{N}$. Denote

$$
\mathbb{P}(\left\{ \omega _{i}\right\} ) = p_{i}, i \in  I.
$$

\bigskip \noindent We have

\begin{equation}
\forall \ (i\in I),\text{ }0\leq p_{i}\leq 1  \label{prop01}
\end{equation}

\Bin and

\begin{equation}
\mathbb{P}(\Omega)=\sum_{i\in I} \mathbb{P}(\{\omega_i\})=\sum_{i\in I}p_{i}=1  \label{prop02}.
\end{equation}

\noindent Moreover, if $A = \left\{ \omega _{i_{1}},\omega _{i_{2}}, ...,\omega _{i_{j}},..., j\in J\right\}\subseteq \Omega $, with
$i_{j} \in  I$, $J \subseteq I$, we get by additivity of  $\mathbb{P}$,

\begin{equation*}
\mathbb{P}(A)=\sum_{j \in J} \mathbb{P}(\left\{ \omega _{ij}\right\})=\sum_{j \in J}p_{i_{j}} =
\sum_{\omega _{i}\in A}p_{i}.
\end{equation*}

\bigskip \noindent It is clear that the knowledge of the numbers $(p_{i})_{i\in I}$ allows to compute the probabilities
$\mathbb{P}(A)$ for all subsets $A$ of $\Omega $.\\

\noindent Conversely, suppose that we are given numbers  $(p_{i})_{i\in I\text{ }}$ such that Equations (\ref
{prop01}) and (\ref{prop02}) hold. Then the mapping $\mathbb{P}$ defined on $\mathcal{P}(\Omega)$ by
\begin{equation*}
\mathbb{P}(\left\{ \omega _{i_{1}},...,\omega _{i_{k}}\right\}
)=\sum_{j=1}^{k}p_{i_{j}}
\end{equation*}

\bigskip \noindent is a probability measure  \index{probability measure} on $\mathcal{P}(\Omega)$.

\newpage

\bigskip \noindent The notion of equi-probability  \index{equi-probability} is very popular in Probability Theory  \index{probability theory} on \textbf{finite} sample spaces. It means  \index{mean} that on the finite sample space with size $n$, all the individual events have equal probability of $1/n$ to occur.\\

\noindent This happens in a fair lottery : if the lottery is based on picking $k=7$ numbers out of $n=40$ fixed numbers, all choices have the same probability of winning the max lotto.\\

\noindent We have the following rule.\\

\noindent \textbf{Theorem}. If a probability measure  \index{probability measure} $\mathbb{P}$, that is defined on a finite sample set $\Omega$ with cardinality $n\geq 1$ with
$\Omega$= $\left\{ \omega _{1},...,\omega _{n}\right\}$, assigns the same probability to all the elementary events,  \index{elementary event} then
for all $1 \leq i \leq n$, we have,
\begin{equation*}
\mathbb{P}(\left\{ \omega _{i}\right\} )=\frac{1}{n}
\end{equation*}

\noindent and for $A \in \mathcal{P}(\Omega)$,

\begin{equation*}
\mathbb{P}(A)=\frac{Card (A)}{n},
\end{equation*}

\Bin \textit{i.e.}

\begin{equation*}
\mathbb{P}(A)=\frac{Card (A)}{Card (\Omega)}
\end{equation*}

\Bin and finally

\begin{equation}
\mathbb{P}(A)=\frac{\text{number of favorable cases}}{\text{number  \index{number of favorable cases} of possible cases}}. \label{proba03}
\end{equation}

\bigskip \bigskip

\noindent \textbf{Proof}. Suppose that for any $i\in \left[1,\text{ }n \right]$, $\mathbb{P}(\left\{\omega_{i}\right\} )=p_{i}=p$, where $p$ is a
constant number between $0$ and $1$. This leads to $1=p_{1}+...+p_{n}=np$. Then

$$
p=\frac{1}{n}.
$$

\bigskip \noindent Let $A=\left\{ \omega_{i_{1}},...,\omega _{i_{k}}\right\}$. We have,

\begin{equation*}
\mathbb{P}(A)=p_{i_{1}}+...+p_{i_{k}}=p+....+p=\frac{k}{n}=\frac{Card (A)}{Card (\Omega)}.
\end{equation*}

\bigskip \noindent \textbf{Remark}. In real situations, such as the lottery and the dice tossing,  \index{tossing} the equi-probability  \index{equi-probability} hypothesis is intuitively deduced, based on symmetry, geometry and logic properties. For example, in a new wedding couple, we use logic to say that : there  is no reason that having a girl as a first child is more likely than having a boy as a first child, and vice-verse. So we conclude that the probability of having a first child girl is one half.\\

\noindent In the situation of equi-probability,  \index{equi-probability} computing probabilities becomes simpler. It is reduced to counting problems, based on the results of chapter 2.\\

\noindent In this situation, everything is based on Formula (\ref{proba03}).\\

\noindent This explains the importance of Combinatorics  \index{combinatoric} in discrete Probability Theory \index{probability theory} .\\

\noindent \textbf{Be careful}. Even if equi-probability  \index{equi-probability} is popular, the contrary is also very common.\\

\noindent \textbf{Example}. A couple wants to have three children. The space is

$$
\Omega=\left\{GGG, GGB, GBG, GBB, BGG, BGB, BBG, BBB \right\}.
$$

\Bin In the notation above, GGG is the elementary event  \index{elementary event} that the couple has three girls, GBG is the event that the couple has first a girl, next a boy and finally a girl, etc.\\

\noindent We suppose the  eight individual events, that are characterized by the gender of the first, and the second and the third child, have equal probabilities of occurring.\\

\noindent Find the probability that each of the following events happen :\\

\noindent (1) the couple has at least one boy.\\
\noindent (2) there is no girl older than a boy.\\
\noindent (3) the couple has exactly one girl.\\

\bigskip

\noindent \textbf{Solution}. Because of equi-probability,  \index{equi-probability} we only have to compute the cardinality of each event and, next
use Formula (\ref{proba03}).\\

\noindent (1) The event A=(the couple has exactly on boy) is :

$$
A=\left\{GGB, GBG, GBB, BGG, BGB, BBG, BBB \right\}.
$$

\bigskip \noindent Then,

$$
\mathbb{P}(A)=\frac{\text{number of favorable cases}}{\text{number  \index{number of favorable cases} of possible cases}}=\frac{7}{8}
$$

\bigskip \noindent (2) The event B=(there is no girl older than a boy) is :

$$
B=\left\{GGG, BGG,  BBG, BBB \right\}.
$$

\noindent Then

$$
\mathbb{P}(A)=\frac{\text{number of favorable cases}}{\text{number  \index{number of favorable cases} of possible cases}}=\frac{4}{8}=\frac{1}{2}
$$

\bigskip \noindent (3) The event C=(The couple has exactly one girl) is :

$$
C=\left\{GBB,BGB, BBG \right\}.
$$

\noindent Then

$$
\mathbb{P}(A)=\frac{\text{number of favorable cases}}{\text{number  \index{number of favorable cases} of possible cases}}=\frac{3}{8}.
$$
 
 %proba01.pm
\chapter{Conditional  \index{conditional} Probability and Independence}  \index{independence} \label{proba01.cpi}

\bigskip

Suppose that we are tossing  \index{tossing} a die three times and considering the outcomes in the order of occurring. The set of all individual events $\Omega $ is the set of triplets $(i, j, k)$, where $i$, $j$ and $k$ are, respectively the face that comes out in the first, in the second and in the third tossing, that is

\begin{equation*}
\Omega =\left\{ 1,2,...,6\right\} ^{3}=\left\{ (i, \ j, \ k), \ 1\leq i, \j, \ k \leq 6 \right\} .
\end{equation*}

\bigskip

\noindent Denote by $A$ the event : {\bf the sum of the three numbers $i$, $j$ and $k$ is six (6)} and by $B$ the event : {\bf the number $1$ appears in the first tossing}.  \index{tossing} We have

\begin{eqnarray}
A &=&\left\{ (i,j,k) \in \{1,2,...,6\}^3, \text{ } i+j+k=6 \right\} \notag\\
&=&\left\{ (1,1,4), (1,2,3), (1,3,2), (1,4,1),(2,1,3),(2,2,2),(2,3,1),(3,1,2),(3,2,1),(4,1,4)\right\} \notag
\end{eqnarray}

\noindent and

\begin{eqnarray}
B &=&\left\{ (i,j,k) \in \{1,2,...,6\}^3, \text{ } i=1 \right\} \notag \\
&=& \{1 \} \times \{1,2,...,6\} \times \{1,2,...,6\} \notag
\end{eqnarray}

\bigskip \noindent Remark that $Card(A)=10$ and $Card(B) = 1 \times 6 \times 6=36$.\\

\noindent Suppose that we have two observers named \textbf{Observer 1} and \textbf{Observer 2}.\\

\noindent \textbf{Observer 1} tosses the die three times and gets the outcome. Suppose the event $A$ occurred. \textbf{Observer 1} knows that $A$ is realized.\\

\noindent \textbf{Observer 2}, who is somewhat far from \textbf{Observer 1}, does not know. But \textbf{Observer 1} let him know that the event $A$ occurred.\\

\noindent Now, given this information, \textbf{Observer 2} is willing to know the probability that $B$ has occurred.\\

\noindent In this context, the event $B$ can not occur out of $A$. The event $A$ becomes the set of individual events, the sample space, with respect to \textbf{Observer 2}.\\

\noindent Then, the event $B$ occurs if and only if $A\cap B$ occurs. Since we are in an equiprobability experience, from the point of view of \textbf{Observer 2}, the probability that $B$ occurs \textbf{given} $A$ already occurred, is
\begin{equation*}
\frac{Card(B\cap A)}{CardA}.
\end{equation*}

\bigskip

\noindent This probability is the conditional  \index{conditional} probability of $B$ \textbf{given} $A$, denoted by $\mathbb{P(}B/A)$ or $\mathbb{P}_{A}(B)$ :

\begin{equation*}
\mathbb{P}_{A}(B)=\frac{Card(A\cap B)}{Card(A)}=\frac{\mathbb{P}(A\cap B)}{\mathbb{P}(A)}.
\end{equation*}

\noindent In the current case,

$$
A \cap B =\left\{ (1,1,4), (1,2,3), (1,3,2), (1,4,1)\right\}
$$

\noindent and then

\begin{equation*}
\mathbb{P}_{A}(B)=\frac{4}{10}=\frac{2}{5},
\end{equation*}

\bigskip \noindent which is different of the unconditional  \index{conditional} probability of $B$ :

\begin{equation*}
\mathbb{P}(B)=\frac{6}{6\times 6 \times 6 \times}=\frac{1}{36}.
\end{equation*}

\bigskip \noindent Based on that example, we may give the general definitions pertaining of the conditional  \index{conditional} probability concept.\\

\noindent \textbf{Theorem 1.} (Definition). Let $(\Omega, \mathcal{P}(\Omega), \mathbb{P})$ be a probability space.  \index{probability space} For any event
$A$ such that $\mathbb{P}(A)>0$, the application
\begin{equation*}
\begin{array}{cccc}
\mathbb{P}_{A}: & \mathcal{P}(\Omega ) & \longmapsto & \left[ 0,1\right] \\
& B & \hookrightarrow & \mathbb{P}_{A}(B)=\mathbb{P(B}/A)=\mathbb{P}(A\cap
B)/\mathbb{P}(A)
\end{array}
\end{equation*}

\noindent is a probability measure.  \index{probability measure} It is supported by $A$, meaning  \index{mean} that we have $\mathbb{P}_{A}(A)=1$. The application $\mathbb{P}_{A}$  is called the conditional  \index{conditional} probability \textbf{given} $A$.\\

\noindent \textbf{Proof}. Let $A$ satisfy $\mathbb{P}(A)>0$. We have for all $B \in \Omega$,

\begin{equation*}
A\cap B\subseteq A.
\end{equation*}

\Bin Thus

\begin{equation*}
\mathbb{P}(A\cap B)\leq \mathbb{P}(A)
\end{equation*}

\Bin and next,

\begin{equation*}
\mathbb{P}_{A}(B)\leq 1.
\end{equation*}

\bigskip \noindent It is also clear that $A$ is a support of the probability measure  \index{probability measure} $\mathbb{P}_{A}$, meaning  \index{mean} that

$$
\mathbb{P}_{A}(A)=1,
$$

\Bin \noindent since

$$
\mathbb{P}_{A}(A)=\mathbb{P}(A\cap A)/\mathbb{P}(A)=\mathbb{P}(A)/\mathbb{P}(A)=1.
$$

\bigskip \noindent It is also obvious that

\begin{equation*}
\mathbb{P}_{A}\geq 0.
\end{equation*}

\bigskip \noindent Further if $B=\sum_{j\geq 1} B_{j}$, we have

\begin{equation*}
A\cap \left(\sum_{j\geq 1} B_{j}\right)=\sum_{j\geq 1} A\cap B_{j}.
\end{equation*}

\Bin By applying the additivity of $\mathbb{P}$, we get

\begin{equation*}
\mathbb{P} \left(A\cap \left(\sum_{j\geq 1} B_{j}\right) \right)=\sum_{j\geq 1} \mathbb{P}(A\cap B_{j}).
\end{equation*}

\Bin By dividing by $\mathbb{P}$(A), we arrive at

\begin{equation*}
\mathbb{P}_{A}\left( \sum_{j\geq 1} B_{j} \right)=\sum_{j\geq 1} \mathbb{P}_{A}(B_{j}).
\end{equation*}

\Bin Hence  $\mathbb{P}_{A}$ is a probability measure \index{probability measure} .\\

\bigskip

\noindent \textbf{Theorem 2}. For all events $A$ and $B$ such that $\mathbb{P}(A)>\ 0$,
we have

\begin{equation}
\mathbb{P}(A\cap B)=\mathbb{P}(A)\times \mathbb{P}\left( B/A\right). \label{PC01}
\end{equation}

\Bin
\noindent Moreover, for any family of events $A_{1}$, $A_{2}$, . . ., $A_{n}$, we have

\begin{eqnarray}
&&\mathbb{P}(A_{1}\cap A_{2}\cap ........A_{n}) \label{PC02}\\
&=&\mathbb{P}(A_{1})\times
\mathbb{P}\left( A_{2}/A_{1}\right) \times \mathbb{P}\left(
A_{3}/A_{1}\cap A_{2}\right) \times ...\times \mathbb{P}\left(
A_{n}/A_{1}\cap A_{2}\cap ...\cap A_{n-1}\right), \notag
\end{eqnarray}

\Bin
\noindent with the convention that $\mathbb{P}_{A}(B)=0$ for any event $B$, whenever we have $\mathbb{P}(A) = 0$.\\

\noindent Formula (\ref{PC02}) is the \textbf{progressive conditioning formula} which is very useful when dealing with Markov chains.\\

\bigskip \noindent \textbf{Proof}. The first formula (\ref{PC01}) is a rephrase of the conditional  \index{conditional} probability for $\mathbb{P}\left(A\right)>0$.\\

\noindent Next, we do understand from the example that, \textbf{given} an impossible event $A$, no event $B$ can occur since $A\cap B$ is still impossible and $\mathbb{P}\left( B/A\right)=0$ for any event $B$ when $\mathbb{P}\left(A\right)=0$. So, Formula (\ref{PC01}) still holds when $\mathbb{P}\left(A\right)=0$.\\

\noindent As to the second formula (\ref{PC02}), we get it by iterating Formula (\ref{PC01}) $n$ times :

\begin{equation*}
\mathbb{P}(A_{1}\cap A_{2}\cap ...\cap A_{n})=\mathbb{P}(A_{1}\cap A_{2}\cap
...\cap A_{n-1})\times \mathbb{P}\left( A_{n}/A_{1}\cap
A_{2}\cap ...\cap A_{n-1}\right).
\end{equation*}

\bigskip

\noindent We are going to introduce a very important formula. This formula is very important in the applications of Probability Theory \index{probability theory} .

\bigskip

\bigskip Consider a partition $E_{1}, E_{2} ,. . ., E_{k}$ of $\Omega $, that is, the $E_i$ are disjoint  \index{joint} and satisfies the relation
\begin{equation*}
\sum_{1\leq i\leq k}E_{i}=\Omega .
\end{equation*}

\noindent \textbf{The causes \index{causes} }.\\

\noindent If we have the partition of $\Omega$ in the form $\sum_{1\leq i\leq k}E_{i}=\Omega$, we call the events $E_i$'s the causes  \index{causes} and the numbers
$\mathbb{P}(E_i)$ the \textit{prior probabilities \index{prior probabilities} }.\\

\noindent Consider an arbitrary event $B$, we have by the distributivity of the intersection over the union, that

$$
B=\Omega \cap B=(E_{1}+E_{2}+...+E_{k})\cap B=E_{1}\cap B+....+E_{k}\cap B.
$$

\Bin
\noindent From the formula

\begin{equation}
B=E_{1}\cap B+....+E_{k}\cap B, \label{proba01.cpi.DecompCausis}
\end{equation}

\Bin \noindent we say that : for $B$ to occur, each cause $E_i$ contributes by the part $E_{i} \cap B$. The denomination of the $E_i$ as causes  \index{causes} follows from this fact.\\

\noindent The first important formula is the following.\\

\noindent \textbf{Total Probabilities Formula \index{Total Probabilities Formula} }.\\

\noindent Suppose that the sample space $\Omega$ is portioned into causes  \index{causes} $E_1$, ..., $E_k$, then for any event $B$,

\begin{equation}
\mathbb{P}(B)=\sum_{j=1}^{k}\mathbb{P}\left( B/E_{j}\right) \text{ }\mathbb{P}(E_{j}).  \label{propa01.cpi.FPT}
\end{equation}

\bigskip \noindent \textbf{Proof}. Let $B$ be an arbitrary event. By Formula (\ref{proba01.cpi.DecompCausis}) and by the additivity of the probability, we have

\begin{equation}
\mathbb{P}(B)=\sum_{j=1}^{k}\mathbb{P}(E_{j}\cap B).
\end{equation}

\noindent By applying the conditional  \index{conditional} probability as in Theorem 2 above, we have for each $j \in \{1,...,k\}$,

$$
\mathbb{P}(E_{j}\cap B)=\mathbb{P}(E_{j})\mathbb{P}(B/E_{j}).
$$

\bigskip \noindent By combining these two formulas, we arrive at

$$
\mathbb{P}(B)=\sum_{j=1}^{k}\mathbb{P}\left( B/E_{j}\right) \text{ }\mathbb{P}(E_{j}).
$$

\noindent QED.\\

\bigskip \noindent The total probability formula allows us to find the probability of a future event $B$, called \textit{effect}, by collecting the contributions of the causes  \index{causes} $E_i$, $i=1,...,k$.\\

\noindent The Bayes rule  \index{Bayes rule} intends to invert this process in the following sense : Given an event $B$ has occurred, what is the probable cause which made $B$ occur. The Bayes rule, in this case, computes the probability that each cause occurred prior to $B$.

\bigskip
\noindent \textbf{Bayes Theorem}.  \index{Bayes theorem} Suppose that \textit{prior probabilities}  \index{prior probabilities} are positive, that is $\mathbb{P}(E_{i})>0$ for each $1\leq i\leq k$. Then, for any $1\leq i\leq k$, for any event $B$, we have

\begin{equation}
\mathbb{P}\left( E_{i}/B\right)=\frac{\mathbb{P}(E_{i})\text{ }\mathbb{P}\left( B/E_{i}\right) }{%
\sum_{j=1}^{k}\mathbb{P}(E_{j})\text{ }\mathbb{P}\left( B/E_{j}\right)}.\label{proba01.cpi.Fbayes}
\end{equation}

\bigskip \noindent The formula computes the probability that the cause $E_i$ occurred given the effect $B$. We will come back to the important interpretations of this formula. Right now, let us give the proof.\\

\noindent \textbf{Proof}. Direct manipulations of the conditional  \index{conditional} probability lead to

\begin{eqnarray*}
\mathbb{P}\left( E_{i}/B\right) &=&\frac{\mathbb{P}(E_{i} \cap B)}{\mathbb{P}(B)}\\
&=&\frac{\mathbb{P}(E_{i})\mathbb{P}\left( B/E_{i}\right) }{\mathbb{\mathbb{P}(B)}}.
\end{eqnarray*}

\noindent We finish by replacing $\mathbb{P}(B)$ by its value using the total probability formula.\\

\noindent Now, let us mention some interpretations of this rule.\\

The Total Probability Formula shows how each cause contributes in forming the probability of future events called effects.\\

\noindent The Bayes rule  \index{Bayes rule} does the inverse way. From the effect, what are the probabilities that the causes  \index{causes} have occurred prior to the effect.\\

\noindent It is like we may invert the \textit{past} and the \textit{future}. But we must avoid to enter into philosophical problems regarding the past and the future. The context of the Bayes rule  \index{Bayes rule} is clear. All is about the past. The effect has occurred at a time \textbf{$t_{1}$} in the past. The first is the future of a second time \textbf{$t_{2}$} at which one of the causes  \index{causes} occurred. The application of the Bayes rule for the future leads to pure speculations.\\

\noindent Now, we need to highlight an interesting property of the Bayes rules  \index{Bayes rule} for two equally probable causes.  \index{causes} In this case, denote  $p=\mathbb{P}(E_{1})=\mathbb{P}(E_{2})$. We have

\begin{equation*}
\mathbb{P}\left( E_{1}/B\right) =\frac{p \ \mathbb{P}\left(B/E_{1}\right) }{\mathbb{P}(E_{1}) \ \mathbb{P}\left( B/E_{1}\right)
+\mathbb{P}(E_{2}) \ \mathbb{P}(B/E_{2})},
\end{equation*}

\Bin and

\begin{equation*}
\mathbb{P}\left( E_{2}/B\right) =\frac{p\text{ }\mathbb{P}\left(
B/E_{2}\right) }{\mathbb{P}(E_{1})\text{ }\mathbb{P}\left( B/E_{1}\right) +%
\mathbb{P}(E_{2})\text{ }\mathbb{P}(B/E_{2})},
\end{equation*}

\Bin and then

\begin{equation*}
\frac{\mathbb{P}\left( E_{1}/B\right) }{\mathbb{P}\left( E_{2}/B\right) }=
\frac{\mathbb{P}\left( B/E_{1}\right) }{\mathbb{P}\left( B/E_{2}\right) }.
\end{equation*}

\bigskip \noindent \textbf{Conclusion}. For two equi-probable  \index{equi-probable} causes,  \index{causes} the ratio of the conditional  \index{conditional} probabilities of the causes given the effect is the same as the ratio of the conditional probabilities of the effect given the causes.\\

\bigskip \noindent Both the Bayes Formula and the Total Probabilities Formula  \index{Total Probabilities Formula} are very useful in a huge number of real problems. Here are some examples.\\

\bigskip \noindent \textbf{Example}.\\

\noindent \textbf{Example 1}. (The umbrella problem).  \index{umbrella problem} An umbrella is in one the seven floors of a building with probability
$0\leq p\leq 1$. Precisely, it is in each floor with probability $p/7$. We searched it in the first six floors without success so that we are sure it is not in these first six floors. What is the probability that the umbrella is in the seventh floor?\\

\bigskip

\noindent \textbf{Solution}. Denote by $E_{i}$ the event : \textit{The umbrella is in the $i$-th floor}. The event
\begin{equation*}
E_{1}+...+E_{7}=E
\end{equation*}

\Bin
\noindent is the event : \textit{The umbrella is the building} and we have

\begin{equation*}
\mathbb{P}(E)=p.
\end{equation*}

\Bin
\noindent Denote by $F=E_{1}+ ....+E_{6}$ the event : \textit{The umbrella is in one of the six first floors}. We see that

$$
F^c=E_{1}^{c}\cap ....\cap E_{6}^{c}
$$

\Bin
\noindent is the event : \textit{The umbrella is not in the six first floors}.\\

\noindent The searched probability is the conditional  \index{conditional} probability :

\begin{equation*}
Q=\mathbb{P}\left( E_{7}/F^{c} \right) =\frac{\mathbb{P}(E_{7} \cap F^{c})}{\mathbb{P}(F^{c})}.
\end{equation*}

\Bin
\noindent We also have

$$
F^c=E_7 + E^c.
$$

\noindent Then
\begin{equation*}
E_{7} \subseteq F^c,
\end{equation*}

\Bin
\noindent which implies that

$$
E_{7} \cap F^{c} =E_{7}.
$$

\Bin
\noindent We arrive at

\begin{equation*}
Q=\mathbb{P}\left( E_{7}/F^{c} \right) =\frac{\mathbb{P}(E_{7})}{\mathbb{P}(F^{c})}.
\end{equation*}

\Bin
\noindent But
$$
\mathbb{P}(F^{c})=\mathbb{P}(E_7 + E^c)=(p/7)+(1-p)=(7-6p)/7.
$$

\Bin
\noindent Hence

\begin{equation}
Q=\mathbb{P}\left( E_{7}/F^{c} \right) =(p/7)/((7-6p)/7)=\frac{p}{7-6p}. \label{umbrella.prob}
\end{equation}

\bigskip \noindent We remark $Q=1$ for $p=1$. The interpretation is simple. If $p=1$, we are sure that the umbrella is in one of the seven floors. If it is not the six first floors, it is surely, that is with probability one, in the seventh.\\

\bigskip \noindent \textbf{Example 2}. (Disease test problem).  \index{disease test problem} In one farm, a medical test $T$ is set to detect infected animals by some Desease $D$.\\

\noindent We have the following facts :\\

\noindent (a) The probability that one animal infected by $D$ is $p=0.3$.\\

\noindent (b) For an infected animal, the probability that the test $T$ declares it positive is $r=0.9$.\\

\noindent (c) For a healthy animal, the probability that the test $T$ declares it negative is $s=0.8$.\\

\noindent  \textbf{Question}. An animal which is randomly picked has been tested and declared positive by the test $T$. What is the probability that it is really infected by $D$.\\

\Bin \textbf{Solution}. Let us introduce the following events:\\

\noindent  \textit{PR} : \textit{the animal positively responds to the test $T$}.\\

\noindent  \textit{NR} : \textit{the animal negatively responds to the test $T$}.\\

\noindent  \textit{D} : \textit{the animal is infected by $D$}.\\

\noindent  We are asked to find the number $\mathbb{P}(D/PR)$.\\

\noindent  We have $\Omega=D + D^{c}$. The two causes  \index{causes} are $D$ and $D^c$. We may use the Bayes rule  \index{Bayes rule} :

\begin{equation*}
\mathbb{P}\left( D/PR\right) =\frac{\mathbb{P}(D) \ \mathbb{P}\left(RP/D\right)}{\mathbb{P}(D) \ \mathbb{P}\left(PR/M\right)
+\mathbb{P}(D^{c}) \ \mathbb{P}\left(PR/D^{c}\right)}.
\end{equation*}

\bigskip \noindent  We are given above : $\mathbb{P}\left(D\right)=0.3$, $\mathbb{P}\left(PR/D\right)=0.9$, $\mathbb{P}\left(NR/D^c\right)=0.8.$\\

\noindent We infer that $\mathbb{P}\left(D^c\right)=0.7$ and

$$
\mathbb{P}\left( PR/D^{c}\right) =1-\mathbb{P}\left( NR/D^{c}\right)=1-0.8=0.2
$$

\noindent We conclude

\begin{equation*}
\mathbb{P}\left( M/RP\right) =\frac{(0.3\times 0.9)}{(0.3\times 0.9+0.7\times 0.2)}=0.7059
\end{equation*}

\bigskip

\noindent We are going to speak about the concept of independence,  \index{independence} that is closely related to what precedes.\\

\newpage

\bigskip

\noindent \textbf{Definition}. Let $A_{1}$, $A_{2}$, . . . , $A_{n}$ be events in a probability space  \index{probability space} $(\Omega , \mathcal{P}(\Omega), \mathbb{P})$. We have the following definitions.\\

\noindent (A)  The events $A_{1}$, $A_{2}$, . . . , $A_{n-1}$ and $A_{n}$ are pairwise independent  \index{pairwise independent} if and only if

\begin{equation*}
\mathbb{P}(A_{i}\cap A_{j})=\mathbb{P}(A_{i})\text{ }\mathbb{P}(A_{j}),\text{
for all }1\leq i\neq j\leq n.
\end{equation*}

\Bin
\noindent (B)  The events $A_{1}$, $A_{2}$, . . ., $A_{n-1}$ and $A_{n}$ are mutually independent  \index{mutually independent} if and only if for any subset $\left\{
i_{1,}i_{2,}...,i_{k}\right\} $ of  $\left\{ 1,2,...,n\right\}$, with $2\leq k\leq n$, we have

\begin{equation*}
\mathbb{P}(A_{i_{1}}\cap A_{i_{2}}\cap ...\cap A_{i_{k}})=\mathbb{P}(A_{i_{1}})\text{ }%
\mathbb{P}(A_{i_{2}}) ... \mathbb{P}(A_{i_{k}})\text{.}
\end{equation*}

\bigskip \noindent (C)  Finally, the events $A_{1}$, $A_{2}$, . . ., $A_{n-1}$ and $A_{n}$ satisfy the global factorization formula  \index{global factorization formula} if and only if

\begin{equation*}
\mathbb{P}(A_{1}\cap A_{2}\cap ...\cap A_{n})=\mathbb{P}(A_{1})\text{ }%
\mathbb{P}(A_{2}) ... \mathbb{P}(A_{n})\text{.}
\end{equation*}

\bigskip \noindent \textbf{Remarks}.\\

\noindent \textbf{(1)} For two events, the three definitions $(A)$, $(B)$ and $(C)$ coincide for $k=2$.\\

\noindent \textbf{(2)} For more that two events, independence  \index{independence} without any further indication, means  \index{mean} mutual independence.\\

\noindent \textbf{(3)} Formula (B) means  \index{mean} that the elements of $A_{1}$, $A_{2}$, . . ., $A_{n-1}$ and $A_{n}$ satisfy the factorization formula for any sub-collection of $A_{1}$, $A_{2}$, . . ., $A_{n-1}$ and $A_{n}$ of size $2$, $3$, ...., $n$.\\

\noindent \textbf{(P1)}  $A$ and $B$ are independent if and only if
$$
\mathbb{P}(A\cap B)=\mathbb{P}(A)\mathbb{P}(B).
$$

\bigskip \noindent \textbf{(P2)}  If $\mathbb{P}(A\cap B)=\mathbb{P}(A)P(B)$ and $\mathbb{P}(A)\neq 0$, then $\mathbb{P}\left( B/A\right) =\mathbb{P}(B)$.\\

\noindent In other words, if $B$ and $A$ are independent, the conditional  \index{conditional} probability $B$ given $A$, does not depend on $A$ : it remains equal to the unconditional probability of $B$.\\

\bigskip

Strictly speaking, proving the independence  \index{independence} requires checking the formula $(B)$. But in many real situations, the context itself
allows us to say that we intuitively have independence \index{independence} .\\

\noindent \textbf{Example.} We toss a die three times. We are sure that the outcomes from one tossing  \index{tossing} are independent of that of the two other tossing. This means  \index{mean} that the outcome of the second tossing is not influenced by the result of the first tossing nor does it influence the outcome of the third tossing. Let us consider the following events.\\

\noindent $A$ : The first tossing  \index{tossing} gives an even number.\\

\noindent $B$ : The number of the face occurring in the second tossing  \index{tossing} is different from $1$ and is a perfect square [that is, its square root is an integer number].\\

\noindent $C$ : The last occurring number is a multiple of 3.\\

\noindent The context of the current experience tells us that these events are independent and we have

\begin{equation*}
\mathbb{P}(A)=Card \{1,3,5\}/6=1/2,
\end{equation*}

\begin{equation*}
\mathbb{P}(B)=Card \{4\}/6=1/6,
\end{equation*}

\begin{equation*}
\mathbb{P}(C)=Card \{3,6\}/6=1/3,
\end{equation*}

\begin{equation*}
\mathbb{P}(A\cap B)=\mathbb{P}(A) \mathbb{P}(B)=(1/2)(1/6)=1/12,
\end{equation*}

\begin{equation*}
\mathbb{P}(A\cap C)=\mathbb{P}(A)  \mathbb{P}(C)=(1/2)(1/3)=1/6,
\end{equation*}

\begin{equation*}
\mathbb{P}(B\cap C)= \mathbb{P}(B)\ \mathbb{P}(C)=(1/6)(1/3)=1/18,
\end{equation*}

\Bin
\noindent and

\begin{equation*}
\mathbb{P}(A\cap B\cap C)= \mathbb{P}(B) \mathbb{P}(B) \mathbb{P}(C)=\frac{1}{2}\times \frac{1}{6}\times \frac{1}{3}=\frac{1}{36}.
\end{equation*}

\bigskip

\bigskip \noindent It is obvious that mutual independence  \index{independence} implies the pairwise disjoint  \index{joint} independence and the global factorization formula \index{global factorization formula} .\\

\noindent But neither of the pairwise independence  \index{independence} nor the global factorization formula  \index{global factorization formula} implies the mutual independence.\\

\noindent Here are counter-examples, from the book by \cite{stoyanov} \index{Stayonov} that is entirely devoted to counter-examples in Probability and Statistics \index{statistics} .\\

\noindent The following counter-example shows that the pairwise independence  \index{independence} does not imply the global factorization, and then does not imply the mutual independence.\\

\noindent A urn contains (4) four cards, respectively holding of the following numbers : $112$, $121$, $211$, and $222$.\\

\bigskip \noindent We want to pick one card at random. Consider the events.\\

\noindent $A_{i}$ : \textit{the number 1 (one) is at the $i$-th place in the number of the card}, $i = 1,2,3$.\\

\noindent  We have

\begin{equation*}
\mathbb{P}(A_{1})=\mathbb{P}(A_{2})=\mathbb{P}(A_{3})=\frac{2}{4}=\frac{1}{2}%
\text{,}
\end{equation*}

\begin{equation*}
\mathbb{P}(A_{1} \cap A_{2})=\mathbb{P}(A_{1}A_{3})=\mathbb{P}(A_{2}A_{3})=\frac{1%
}{4}\text{.}
\end{equation*}

\noindent Yet, we have
\begin{equation*}
\mathbb{P}(A_{1}\cap A_{2} \cap A_{3})=0\neq \mathbb{P}(A_{1})~P(A_{2})\mathbb{~P}%
(A_{3})=\frac{1}{8}\text{. }
\end{equation*}

\noindent We have pairwise independence  \index{independence} and not the mutual independance.\\

\bigskip

\noindent The following counter-example shows that the global factorization property does not imply the pairwise independence,  \index{independence} and then, the global factorization formula  \index{global factorization formula} does not imply the mutual independence.\\

\noindent {\bf Exemple}: We toss a die twice. The probability space  \index{probability space} is
$\Omega= \{1,2,...,6 \}^2$. Consider the events :\\

\noindent $A$ : \textit{The first tossing  \index{tossing} gives 1, 2, 3}.\\

\noindent $B$ : \textit{The second gives donne 4, 5, 6}.\\

\noindent and\\

\noindent $C$ : \textit{The sums of the two numbers is 9}.\\

\noindent We have :

\begin{equation*}
Card(A)=3\times 6=18\text{ \ and then  }\mathbb{P}(A)=\frac{18}{36}=\frac{1}{2},
\end{equation*}

\begin{equation*}
Card(B)=3\times 6=18\text{ and then }\mathbb{P}(B)=\frac{18}{36}=\frac{1}{2}
\end{equation*}

\noindent and
\begin{equation*}
Card(C)=4\ \ \text{and then  }\mathbb{P}(C)=\frac{4}{36}=\allowbreak \frac{1}{9}.
\end{equation*}

\bigskip \noindent We also have

$$
Card(A \cap  B) = 9,
$$

$$
Card (A \cap  C) = 1,
$$

$$
Card (B \cap  C) =3,
$$

\noindent  and

$$
A \cap  B \cap  C = B \cap  (A \cap  C) = \left\{ (3,6)\right\}.
$$

\noindent Then, we have

\begin{equation*}
\mathbb{P}(A\cap B\cap C)=\frac{1}{36}
\end{equation*}

\Bin \noindent and we have the global factorization

\begin{equation*}
\mathbb{P}(A)\mathbb{P}(B)\mathbb{P}(C)=\frac{1}{2}\times \frac{1}{2}\times
\frac{1}{9}=\frac{1}{36}=\mathbb{P}(A\cap B\cap C).
\end{equation*}

\bigskip \noindent But we do not have the pairwise independence  \index{independence} since

\begin{equation*}
\mathbb{P}(B\cap C)=\frac{3}{36}= \frac{1}{12} \ \ \neq \ \ \mathbb{P}(B)\mathbb{P}(C)=\frac{1}{2} \times \frac{1}{9}=\frac{1}{18}.
\end{equation*}

\bigskip \noindent \textbf{Conclusion}. The chapters 3 and 4 are enough to solve a huge number of problems in elementary Probability Theory  \index{probability theory} provided enough  mathematical tools of Analysis and Algebra are mastered. You will have the great and numerous opportunities to practice with the Exercises book related to this monograph.

 %proba01.cpi
\chapter{Random Variables} \label{proba01.rv}

\bigskip

\Bin

We begin with definitions, examples and notations.\\

\noindent \textbf{Definition 1}. Let $(\Omega, \mathcal{P}(\Omega), \mathbb{P})$ be a probability space.  \index{probability space} A random variable on $(\Omega, \mathcal{P}(\Omega), \mathbb{P})$ is an application $X$ from $\Omega$ to  $\mathbb{R}^{k}$.\\

\noindent If $k = 1$, we say that $X$ is a \textbf{real-valued random variable}, abbreviated in (\textit{rrv}).\\

\noindent If $k = 2$, $X$ is a \textbf{ random couple} or a \textbf{bi-dimensional random variable}.\\

\noindent If the general case, $X$ is a random variable of $k$ dimensions, or a \textbf{$k$-dimensional random variable}, or simply a \textbf{random vector}.\\

\noindent If $X$ takes a finite number of distinct values (points) in $\mathbb{R}^{k}$, $X$ is said to be a \textbf{random variable with finite number of values}.\\

\noindent If $X$ takes its values in a set $\mathcal{V}_X$ that can be written in an enumerable form :
$$
\mathcal{V}_X =\{ x_{1}, x_{2}, x_{3}, . . .\},
$$

\bigskip \noindent the random variable is said to be a \textbf{discrete random variable}. So, any random variable with finite number of values is a discrete random variable.\\

\noindent \textbf{Examples 1}. Let us toss twice a die whose faces are numbered from 1 to 6. Let $X$ be the addition of the two occuring numbers. Then $X$ is a real-valued random defined on
$$
\Omega =\left\{ 1,2,......,6\right\} ^{2}
$$

\noindent such that $X((i, j))=i + j$.\\

\noindent Let us consider the same experience and let $Y$ be the application on $\Omega$ defined by

$$
Y((i, j))=(i+2j, i/(1+j), (i+j)/(1+|i-j|).
$$

\Bin \noindent $Y$ is a random vector of dimension $k=3$. Both $X$ and $Y$ have a finite number of values.

\bigskip \noindent \textbf{Example 2.} Consider a passenger arriving at a bus station and beginning to wait for a bus to pick him up to somewhere. Let $\Omega$ be the set all possible passengers. Let $X$ be the time a passenger has to wait before the next bus arrives at the station.\\

\noindent This random variable is not discrete, since its values are positive real numbers $t>0$. This set of values is an interval of $\mathbb{R}$. It is not denumarable (See a course of calculus).\\

\noindent \textbf{Notations}. Let $X$ : $\Omega \mapsto \ \mathbb{R}^{k}$ be a random vector and $A \subseteq \mathbb{R}^{k}$, a subset of $\mathbb{R}^{k}$. Let us define the following events

\begin{equation}
(X\in A)=X^{-1}(A)=\left\{ \omega ~\in ~\Omega ,\text{ }X(\omega )~\in ~A\right\}. \label{proba01.rv.invA}
\end{equation}

\Bin
\noindent The set $A$ usually takes particular forms we introduce below :

\begin{equation*}
(X\leq x)=X^{-1}~\left( \left] -~\infty ,~x\right] ~\right) =\left\{ \omega
~\in ~\Omega ,\text{ }X(\omega )~\leq ~x\right\} ,
\end{equation*}

\begin{equation*}
(X~<~x)=X^{-1}\left( ~\left] -~\infty ,~x\right[ \right) ~=\left\{ \omega
~\in ~\Omega ,\text{ }X(\omega )~<~x\right\},
\end{equation*}

\begin{equation*}
(X\geq x)=X^{-1}~\left( \left[ x,~+\infty \right[ \right) =\left\{ \omega
~\in ~\Omega ,\text{ }X(\omega )~\geq ~x\right\},
\end{equation*}

\begin{equation*}
(X~>~x)=X^{-1}~\left( \left] ~x,~+\infty \right[ \right) =\left\{ \omega
~\in ~\Omega ,\text{ }X(\omega )~>~x\right\},
\end{equation*}

\begin{equation*}
(X=x)=X^{-1}(~\left\{ x\right\} ~)=\left\{ \omega ~\in ~\Omega ,\text{ }%
X(\omega )~=~x\right\}.
\end{equation*}

\Bin \noindent These notations are particular forms of Formula (\ref{proba01.rv.invA}). We suggest that you try to give more other specific forms by letting $A$ take particular sets.\\

\bigskip \noindent \textbf{Exercise}\\

\noindent (a) Let $X$ be a real random variable. Extend the list of the notation in the same form for

$$
A=[a,b], \ A=[a,b[, \ A=]a,b[, \ A=\{0\}, \ A=\{1,2\}, \ A=\{]1,2[ \cup \{5,10\}\}.
$$

\Bin (b) Let $(X,Y)$ be a random vectors. Define the sets $((X,Y)\in A)$ with\\
$$
A=\{(x,y), x\leq a, y\geq b\},
$$

$$
A=\{(x,y), y=x+1\},
$$

$$
A=\{(x,y), \frac{x}{1+x^2+y^2}\geq 1 \},
$$

$$
A=\{1\}\times [1,+\infty[,
$$

\noindent or

$$
A=\{1\}\times [1,10].
$$

\bigskip \noindent \textbf{Important Reminder}. We remind that the reciprocal image mapping (or inverse image mapping)  \index{inverse image mapping} $X^{-1}$ preserves all sets operations (See a course of general Algebra). In particular, we have

\begin{equation*}
X^{-1}(A)\cap X^{-1}(B)=X^{-1}(A\cap B),
\end{equation*}

\begin{equation*}
X^{-1}(A^{c})=X^{-1}(A)^{c},
\end{equation*}

\begin{equation*}
X^{-1}(A)\cup X^{-1}(B)=X^{-1}(A\cup B),
\end{equation*}

\begin{equation*}
X^{-1}\left( A\setminus B\right) =X^{-1}(A)\setminus X^{-1}(B).
\end{equation*}

\bigskip \noindent In Probability theory,  \index{probability theory} we have to compute the probability that events occur. But for random variables, we compute the probability that this random variable takes a particular value or fall in a specific region of $\mathbb{R}^{k}$. This is what we are going to do. But, in this first level of Probability Theory, we will mainly focus on the real case where $k=1$.

\noindent At this preliminary stage, we focus on discrete random variables.\\

\noindent \textbf{Definition 2}. Let $X$ be a discrete random variable defined on the probability space  \index{probability space} $\Omega$ onto $\mathcal{V} (X)=\left\{
x_{i},~i\in I\right\} \subseteq \mathbb{R}^{k}$, where $I$ is a subset of the set of non-negative integers $\mathbb{N}$. The probability law  \index{probability law} of $X$, is characterized by the numbers

\begin{equation*}
\mathbb{P}(X=x_{i}),~i\in I\text{.} \label{proba01.rv.dens1}
\end{equation*}

\noindent These numbers satisfy

\begin{equation}
\sum_{i\in I}~\mathbb{P}(X=x_{i})=1. \label{proba01.rv.dens2}
\end{equation}

\bigskip \noindent \textbf{Remarks}.\\

\noindent \textbf{(1)} We may rephrase the definition in these words : when you are asked to give the probability law  \index{probability law} of $X$, you have to provide the numbers in Formula (\ref{proba01.rv.dens1}) and you have to check that Formula (\ref{proba01.rv.dens2}) holds.\\

\noindent \textbf{(2)} The set $I$ is usually $\mathbb{N}$ or a finite set of it, of the form $I=\{1,2,...,m\}$ where $m$ is finite.\\

\bigskip \textbf{Direct Example}.\\

\noindent The probability law  \index{probability law} may be given by explicitly stating Formula (\ref{proba01.rv.dens1}). When the values set is finite and moderate, it is possible to use a table. For example, consider the real-valued random variable of the age of a student in the first year of high school. Suppose that $X$ takes the values 17, 18, 19, 20, 21 with probabilities 0.18, 0.22, 0.28, 0.22, 0.10.   We may represent its probability in the following table

\begin{table}[htbp]
\centering
\begin{tabular}{|l|c|c|c|c|c|c|}
\hline
Values $x_i$ 				& 17 							& 18 							& 19 							& 20 							& 21 							& Total\\
\hline
$\mathbb{P}(X=x_i)$ & $\frac{18}{100}$ 	& $\frac{22}{100}$ 	& $\frac{28}{100}$ 	& $\frac{22}{100}$ 	& $\frac{10}{100}$	& $\frac{100}{100}$\\
\hline
\end{tabular}
\end{table}

\noindent Formula (\ref{proba01.rv.dens2}), as suggests its name, enables the computation of the probability that $X$ fall in any specific area of $\mathbb{R}^k$. Indeed, we have the following result.\\

\begin{theorem} \label{proba01.rv.theo1} Let $X$ be a discrete random variable defined of $\Omega$ with values in $\mathcal{V}(X)=\left\{
x_{i},~i\in I\right\} \subseteq \mathbb{R}^{k}$. Then for any subset of $\mathbb{R}^{k}$, the probability that $X$ falls in $A$ is given by

\begin{equation*}
\mathbb{P}(X\in A)=\sum_{i\in I, x_{i}~\in ~A}\mathbb{P}(X=x_{i}).
\end{equation*}
\end{theorem}

\bigskip \noindent \textbf{Proof.} We have

\begin{equation*}
\Omega =\sum_{i~\in ~I}(X=x_{i}).
\end{equation*}

\noindent For any subset $A$ of  $\mathbb{R} ^k$, $(X \in  A)$ is a subset of $\Omega $, and then

\begin{equation}
(X\in A)=(X\in A)\cap \Omega =\sum_{i~\in ~I}(X\in A)\cap (X=x_{i}). \label{proba01.rv.decompEvent01}
\end{equation}

\noindent For each $i \in I$, we have :\\

\noindent \textbf{Either} $x_{i} \in  A$ and in this case,  $(X = x_{i}) \subseteq (X \in A)$,\\

\noindent \textbf{Or}, $x_{i}\notin  A$ and in this case, $(X \in  A) \cap  (X = x_{i}) = \emptyset $.\\

\noindent Then, we finally have

\begin{equation*}
(X\in A)=\sum_{x_{i}~\in ~A}(X=x_{i})\text{.}
\end{equation*}

\noindent We conclude the proof by applying the additivity of the probability measure  \index{probability measure} $\mathbb{P}$.\\

\bigskip \noindent We seize the opportunity to extend Formula (\ref{proba01.rv.decompEvent01}) into a general principle.\\

\noindent If $X$ is a discrete random variables taking the values $x_i$, $i \in I$, then for any event $B \subset \Omega$, we have the following decomposition formula

\begin{equation}
B =\sum_{i~\in ~I} B \cap (X=x_{i}). \label{proba01.rv.decompEvent02}
\end{equation}

\bigskip \noindent A considerable part of Probability Theory  \index{probability theory} consists of finding probability laws.  \index{probability law} So, we have to begin to learn the usual and unavoidable cases.\\

\Ni We are going to give a list of some probability laws  \index{probability law} including the very usual ones. For each case, we describe a random experience that generates such a random variable. It is very important to perfectly understand these experiences and to know the results by heart, if possible, for subsequent more advanced courses in Probability Theory \index{probability theory} .\\

\bigskip \noindent \textbf{(1) Degenerate or Constant Random variable}.\\

\noindent It is amazing to say that an element $c$ of $\mathbb{R}^{k}$ is a random variable. Indeed, the application $X : \Omega \mapsto \mathbb{R}^{k}$, such that

$$
\text{for any } \omega \in \Omega, \ \ X(\omega)=c,
$$

\Bin is the constant application that assigns to any $\omega$ the same value $c$. Even if $X$ is not \textit{variable} (in the common language), it is a random variable in the sense of an application.\\

\noindent We have in this case : $\mathcal{V}(X)=\{c\}$ and the probability law  \index{probability law} is given by

\begin{equation}
\mathbb{P}(X=c)=1. \label{proba01.rv.pl.degen}
\end{equation}

\bigskip \noindent \textbf{In a usual terminology, we say that $X$ is degenerated}.\\

\noindent In the remainder of this text, we focus on real-valued-random variables.\\

\bigskip \noindent \textbf{(2) Discrete Uniform Random variable \index{discrete uniform random variable} }.\\

\noindent A Discrete Uniform Random variable  \index{discrete uniform random variable} $X$ with parameter $n\geq 1$, denoted by $X \sim \mathcal{DU}(n)$, takes its values in a set of $n$ elements $\mathcal{V}(X)$ = $\left\{ x_{1},~x_{2},..,x_{n}\right\}$ such that

\begin{equation*}
\mathbb{P}(X=x_{i})=\frac{1}{n},~i=1,.....,n\text{.}
\end{equation*}
\bigskip

\noindent \textbf{Generation}. Consider a urn \textbf{U} in which we have $n$ balls that are numbered from 1 to $n$ and indistinguishable by hand. We want to draw at random on ball from the urn and consider the random number $X$ of the ball that will be drawn.\\

\noindent We identify the balls with their numbers and we have  $\mathcal{V}(X)=\left\{ 1,~2.......,n\right\}$.\\

\noindent We clearly see that all the $n$ balls will have the same probability of occurring, that is $1/n$. Then,

\begin{equation}
\mathbb{P}(X=i)=\frac{1}{n},~i=1,......,n. \label{proba01.rv.pw.du}
\end{equation}

\Bin
\noindent Clearly, the addition of the probabilities in (\ref{proba01.rv.pw.du}) gives the unity.\\

\bigskip

\bigskip \noindent \textbf{(3) Bernoulli Random variable  \index{Bernoulli random variable} $\mathcal{B}$(p), $0<\ p<\ 1$}.\\

\noindent A Bernoulli random variable  \index{Bernoulli random variable} $X$ with parameter $p$, $0< p< 1$, denoted by $X\sim \mathcal{B}(p)$, takes its values in  $\mathcal{V}(X)=\left\{ 0,~1\right\}$ and

\begin{equation}
\mathbb{P}(X=1)=p=1-\mathbb{P}(X=0). \label{proba01.rv.pw.bern}
\end{equation}

\Bin
\noindent * $X$ is called a success-failure random variable. The value one is associated with a success and the value zero with a failure.\\

\noindent \textbf{Generation}. A urn \textbf{U} contains $N$ black balls and $M$ red balls, all of them have the same form and are indistinguishable by hand touch. We draw at random one ball. Let $X$ be the random variable that takes the value $1$ (one) if the drawn ball is red and $0$ (zero) otherwise.\\

\noindent It is clear that we have : $X\sim \mathcal{B}\left( p\right)$ with $p=M/(N+M$).\\

\bigskip \noindent \textbf{(4) Binomial Random variable  \index{binomial random variable} $\mathcal{B}(n, p)$ with parameters $n\geq 1$ and $p\in ]0,1[$}.\\

\noindent A Binomial random variable  \index{binomial random variable} $X$ with parameters $n\geq 1$ and $p\in ]0,1[$, denoted by $X\sim \mathcal{B}(n,p)$, takes its values in $$
\mathcal{V}(X) = \left\{ 0,1,......,n\right\}
$$

\noindent and

\begin{equation}
\mathbb{P}(X=k)= \binom{n}{k} / p^{k} / (1-p)^{n-k}, \ \ 0\leq leq k \geq n. \label{proba01.rv.pl.bin}
\end{equation}

\bigskip

\noindent \textbf{Generation}. Consider a Bernoulli experience with probability $p$. For example, we may consider again a urn with $N+M$ similar balls, $N$ being black and $M$ red, with $p=M/(N+M)$. We are going to draw one ball at random, note down the color of the drawn ball and put it back in the urn. We repeat this experience $n$ times. Finally, we consider the number of times $X$ a red ball has been drawn.\\

\noindent In general, we consider a Bernoulli experience of probability of success $p$ and \textbf{independently} repeat this experience $n$ times. After $n$ repetitions, let $X$ be the number of successes.\\

\noindent It is clear that $X$ takes its values in

$$
\mathcal{V}(X)=\left\{ 0,1......,n\right\},
$$

\Bin \noindent and $X$ follows a Binomial law with parameters $n\geq 1$ and $p\in]0,1[$.\\

\noindent \textbf{Proof}. We have to compute  $\mathbb{P}(X = k)$ for $k=0,...,n$. We say that the event $(X = k)$ is composed by all permutations  \index{permutation} of $n$ elements with repetition of the following particular elementary event \index{elementary event} 

\begin{equation*}
\omega _{0}=\underset{k\text{ times }}{\underbrace{RR...R\text{ }}}\underset{%
n-k\text{ times}}{\text{ }\underbrace{NN...N\text{ }}},
\end{equation*}

\Bin
\noindent where $k$ elements are identical and $n-k$ others are identical and different from those of the first group. By independence,  \index{independence} all these permutations  \index{permutation} with repetition have the same individual probability

$$
p^{k}(1-p)^{n-k}.
$$

\Bin
\noindent And the number of permutations  \index{permutation} with repetition of $n$ elements with two distinguishable groups of indistinguishable elements of respective sizes of $k$ and $n-k$ is

\begin{equation*}
\frac{n!}{k!~(n-k)!}.
\end{equation*}

\Bin
\noindent Then, by the additivity of the probability measure,  \index{probability measure} we have

\begin{equation*}
\mathbb{P}(X=k)=\frac{n!}{k!~(n-k)!}  p^{k}(1-p)^{n-k},
\end{equation*}

\noindent which is Formula (\ref{proba01.rv.pl.bin}).

\bigskip \noindent \textbf{(5) Hypergeometric Random Variable  \index{geometric random variable}  \index{hypergeometric random variable} $\mathcal{H}(N, r,\theta)$}.\\

\noindent A Hypergeometric Random Variable  \index{geometric random variable}  \index{hypergeometric random variable} $X$ with parameter $N\geq 1$, $1\leq r \leq N$, and $\theta \in ]0,1[$ takes its values in

$$
\mathcal{V}(X)=\left\{ 0,1,.......,min(r,N\theta )\right\}
$$

\Bin \noindent and

\begin{equation}
\mathbb{P}(X=k)=\frac{\left(\begin{tabular}{c} $N\theta$ \\ $k$ \end{tabular} \right) \times \left(\begin{tabular}{c} $N(1-\theta)$ \\ $r-k$ \end{tabular} \right)}
{\left(\begin{tabular}{c} $N$ \\ $r$ \end{tabular} \right)}. \label{proba01.rv.pl.hyper}
\end{equation}

\bigskip

\noindent \textbf{Generation}. Consider a urn containing $N$ similar balls and exactly $M$ of these balls are black, with $0< M < N$. We want to draw without replacement $r\geq 1$ balls from the urn with $r\leq N$. Let $X$ be the number of black balls drawn out of the $r$ balls. Put

$$
\theta =M/N \in ]0,1[.
$$

\Bin
\noindent It is clear that $X$ is less than $r$ and less that $M$ since a drawn ball is not put back in the urn. So the values taken by $X$ is

$$
\mathcal{V}(X)=\left\{ 0,1,.......,min(r,N\theta )\right\}.
$$

\Bin
\noindent Now, when we do not take the order into account, drawing $r$ balls without replacement from the urn is equivalent to drawing at once the $r$ balls at the same time. So, our draw is simply a combination  \index{combination} and all the combinations of $r$ elements from $N$ elements have equal probabilities of occurring. Our problem becomes a counting one : In how many ways can we choose a combination of $r$ elements from the urn with exactly $k$ black balls? It is enough to choose $k$ balls among the black ones and $r-k$ among the non black. The total numbers of combination is

$$
\left(\begin{tabular}{c} $N$ \\ $r$ \end{tabular} \right)
$$

\Bin \noindent and the number of favorable cases  \index{number of favorable cases} to the event $(X=k)$ is

$$
\left(\begin{tabular}{c} $M$ \\ $k$ \end{tabular} \right) \times \left(\begin{tabular}{c} $N-M$ \\ $r-k$ \end{tabular} \right).
$$

\noindent So, we have

\begin{equation*}
\mathbb{P}(X=k)=\frac{\left(\begin{tabular}{c} $M$ \\ $k$ \end{tabular} \right) \times \left(\begin{tabular}{c} $N-M$ \\ $r-k$ \end{tabular} \right)}
{\left(\begin{tabular}{c} $N$ \\ $r$ \end{tabular} \right)}.
\end{equation*}

\bigskip \noindent We replace $M$ by $N\theta$ and get (\ref{proba01.rv.pl.hyper}).

\bigskip

\bigskip \noindent \textbf{(6) Geometric Random Variable  \index{geometric random variable} $\mathcal{G}(p)$}.\\

\noindent A Geometric random variable  \index{geometric random variable} $Y_1$ with parameter $p \in ]0,1[$, denoted by $Y_1 \sim \mathcal{G}(p)$, takes its values in

$$
\mathcal{V}(Y_1)=\{1,2, .....\}
$$

\noindent and

\begin{equation}
\mathbb{P}(Y_{1}=n)=p(1-p)^{n-1}, \ n \geq 1. \label{proba01.rv.pl.binneg1}
\end{equation}

\bigskip \noindent \textbf{Remark}. The Geometric random variable  \index{geometric random variable}  is the first example of random variable with an infinite number of values.\\

\bigskip \noindent \textbf{Generation}. Consider a Bernoulli experience with parameter $p \in ]0,1[$ which is independently repeated. We decide to repeat the experience until we get a success. Let $Y_1$ be the number of repetitions we need to have one (1) success.\\

\noindent Of course, we have to do at least one experience to hope to get a success. Hence,

$$
\mathcal{V}(X_1)=\{1,2, ...\}.
$$

\Bin \noindent Now, the event $(Y_1=n)$ requires that we have failures in first $(n-1)$ experiences and a success at the $n^{th}$ experience. We represent this event by

\begin{equation*}
(Y_{1}=k)=\underset{(n-1)~times}{\underbrace{FF \cdots F}~S}.
\end{equation*}

\noindent Then, by independence,  \index{independence} we have

\begin{equation*}
\mathbb{P}(Y_{1}=n)=p(1-p)^{n-1},~n=1,~2, \ \cdots
\end{equation*}

\bigskip \noindent Now, consider $X_1=Y_{1}-1$. We see that $X$ is the number of failures before the first success.\\

\noindent The values set of $X$ is $\mathbb{N}$ and

\begin{equation}
\mathbb{P}(X_1=n)=(1-p)^{n}p,~n=0,1,~2, \cdots, \label{proba01.rv.pl.nfailures}
\end{equation}

\noindent since

\begin{equation*}
\mathbb{P}(X_{1}=n)=\mathbb{P}(Y_{1}=n+1)=(1-p)^{n}p,~n=0,1,~2,\cdots
\end{equation*}

\bigskip \noindent \textbf{Nota-Bene}. Certainly, you are asking yourself why we used a subscript in $Y_1$. The reason is that we are working here for one success. In a coming subsection, we will work with $k>1$ successes and $Y_k$ will stand for the number of independent repetitions of the Bernoulli experience we need to ensure exactly $k$ successes.\\

\noindent Before we do this, let us learn a nice property of Geometric random variable \index{geometric random variable} s.\\

\noindent \textbf{Loss of memory property}.  \index{loss of memory property} Let $X+1$ follow a geometric law $\mathcal{G}(p)$. We are going to see that, for $m \geq 0$ and $n\geq 0$, the probability that $X=n+m$ given that $X$ has already exceeded $n$, depends only on $m$, and is still equal to the probability that $X=m$ :

\begin{equation*}
\mathbb{P}\left( X=n+m / X>n\right) =\mathbb{P}(X=>n)
\end{equation*}

\noindent Indeed, we have

\begin{eqnarray*}
\mathbb{P}(X\geq n)&=& \mathbb{P}(\bigcup_{r\geq n}(X=r))\\
&=& \sum_{r\geq n}p(1-p)^{r}\\
&=& p(1-p)^{n}\sum_{r\geq n}(1-p)^{r-n}\\
&=& p(1-p)^{n}\sum_{s\geq 0}(1-p)^{s} \text{ } (\text{ by change of variable } s=r-n)\\
&=& p(1-p)^{n}(1/p)=(1-p)^{n}.\\
\end{eqnarray*}

\noindent Then for $m \geq 0$ and $n\geq 0$, the event $(X=m+n)$ is included in $(X>n)$, and then $(X=m+n, X>n)=(X=m+n)$. By using the latter equality, we have

\begin{eqnarray*}
\mathbb{P}\left( X=n+m / X_{1}>n\right) &=&\frac{\mathbb{P}(X=m+n, X>n)}{\mathbb{P}(X>n)}\\
&=&\frac{\mathbb{P}(X=m+n)}{\mathbb{P}(X>n)}\\
&=&\frac{p(1-p)^{n+m}}{(1-p)^{n}}\\
&=&p(1-p)^{m}\\
&=&\mathbb{P}(X=m).
\end{eqnarray*}

\Bin  We say that $X$ has the property of forgetting the past or loosing the memory. This means  \index{mean} the following : if we continuously repeat the Bernoulli experience, then given any time $n$, the law of the number of failures before the first success after the time $n$, does not depend on $n$ and still has the probability law  \index{probability law} (\ref{proba01.rv.pl.nfailures}).\\

\bigskip \noindent \textbf{(7) Negative Binomial Random variable:  \index{negative binomial random variable}  \index{binomial random variable} $\mathcal{BN}(k, p)$}.\\

\noindent A Negative Binomial Random variable  \index{negative binomial random variable}  \index{binomial random variable} $Y_k$ with parameter $k\geq 1$ and $p \in ]0,1[$, denoted by $X \sim \mathcal{BN}(k, p)$, takes its values in

$$
\mathcal{V}(Y_k)=\{k,k+1, ...\}
$$

\noindent and

\begin{equation*}
\mathbb{P}(Y_{k}=n)= \binom{n-1}{k-1}\ p^{k}\ (1-p)^{n-k},~n=k,k+1,.......
\end{equation*}

\bigskip

\noindent \textbf{Generation}. Unlikely to the binomial case, we fix the number of successes. And we independently repeat a Bernoulli experience with probability $p \in ]0,1[$ until we reach $k$ successes. The number $Y_k$ of trials needed to reach $k$ successes is called a Negative Binomial Random variable  \index{negative binomial random variable}  \index{binomial random variable} or parameters $p$ and $k$.\\

\noindent (1) For Binomial random variable  \index{binomial random variable} $X \sim \mathcal{B}(n,p)$, we set the number of experiences to $n$ and search the probability that the number of success is equal to $k$ : $\mathbb{P}(X=k)$.\\

\noindent (1) For a Negative Binomial random variable  \index{negative binomial random variable}  \index{binomial random variable} $Y_k \sim \mathcal{NB}(n,p)$, we set the number of successes $k$ and search the probability that the number of experiences needed to have $k$ successes is equal to $n$ : $\mathbb{P}(Y_k=n)$.\\

\noindent Of course, we need at least $k$ experiences to hope to get $k$ successes. So, the values set of $Y_k$ is

$$
\mathcal{V}(Y_k)=\{k,k+1, ...\}.
$$

\Bin \noindent The event $(Y_{k} = n)$ occurs if and only if :\\

\noindent (1) The $n^{th}$ experience is successful.\\

\noindent (2) Among the $(n - 1)$ first experiences, we got exactly $(k - 1)$ successes.\\

\noindent By independence,  \index{independence} the probability of the event $(Y_{k}=n)$ is the product of the probabilities of the described events in (1) and (2), denoted $a \times b$. But $a=p$ and $b$ is the probability to get $k-1$ successes in $n-1$ independent Bernoulli experiences with parameter $p$. Then, we have

$$
b=\left(\begin{tabular}{c} $n-1$ \\ $k-1$ \end{tabular} \right) p^{k-1}~(1-p)^{(n-1)-(k-1)}.
$$

\Bin
\noindent By using, the product of probabilities $a\times b$, we arrive at

\begin{equation*}
\mathbb{P}(Y_{k}=n)= \left(\begin{tabular}{c} $n-1$ \\ $k-1$ \end{tabular} \right)
p^{k}~(1-p)^{n-k},~n\geq k\text{.}
\end{equation*}

\bigskip
\bigskip \noindent \textbf{(8) Poisson Random Variable  \index{Poisson random variable} $\mathcal{P}(\lambda)$}.\\

\noindent A Poisson Random Variable  \index{Poisson random variable} with parameter $\lambda>0$, denoted by $\mathcal{P}(\lambda)$ takes its values in $\mathbb{N}$ and

\begin{equation}
\mathbb{P}(X=k)=\frac{\lambda ^{k}}{k!} e^{-\lambda },~k=0,~1,~2,..., \label{proba01.rv.pl.pois}
\end{equation}

\Bin
\noindent where $\lambda$ is a positive number.\\

\noindent \textbf{Here, we are not able to generate the Poisson law by a simple experience}. Rather, we are going to show how to derive it by a limit procedure from binomial random variable \index{binomial random variable} s.\\

\noindent First, we have the show that Formula (\ref{proba01.rv.pl.pois}) satisfies (\ref{proba01.rv.dens2}). If this is true, we say that Formula (\ref{proba01.rv.pl.pois}) is the probability law  \index{probability law} of some random variable taking its values in $\mathbb{N}$.\\

\noindent We recall the  following expansion (from calculus courses)

\begin{equation}
e^x= \sum_{k=0}^{+\infty}  \frac{x^{k}}{k!}, x \in \mathbb{R}. \label{proba01.rv.expand.expo}\end{equation}

\bigskip \noindent By using (\ref{proba01.rv.expand.expo}), we get

\begin{eqnarray*}
\sum_{k=0}^{+\infty} \mathbb{P}(X=k)&=& \sum_{k=0}^{+\infty} \frac{\lambda ^{k}}{k!} e^{-\lambda}\\
&=& e^{-\lambda} \sum_{k=0}^{+\infty} \frac{\lambda ^{k}}{k!}\\
&=& e^{-\lambda}e^{+\lambda}=1.\\
\end{eqnarray*}

\noindent \textbf{A random variable whose probability law  \index{probability law} is given by (\ref{proba01.rv.pl.pois}}) is said to follow a Poisson law with parameter $\lambda$.\\

\noindent We are going to see how to derive it from a sequence binonial random variables.\\

\noindent Let $X_n$, a sequence of $\mathcal{B}(n,p_n)$ random variables such that, as $n \rightarrow +\infty$,

\begin{equation}
p_n \rightarrow 0 \text{ and } np_n \rightarrow \lambda, 0 <\lambda <+\infty. \label{proba01.rv.cond.bin.pois}
\end{equation}

\bigskip

\begin{theorem} \label{proba01.rv.theo2} Let

$$
f_{n}(k)= \left(\begin{tabular}{c} $n$ \\ $k$ \end{tabular} \right)p^{k}(1-p)^{n-k}
$$

\Bin
\noindent be the probability law  \index{probability law} of the random variable $X_n$ and the probability law of the Poisson random variable \index{Poisson random variable} 

\begin{equation*}
f(k)= \frac{e^{-\lambda }\lambda ^{k}}{k!}.
\end{equation*}

\Bin
\noindent Suppose that (\ref{proba01.rv.cond.bin.pois}) holds. Then, for any fixed $k\in \mathbb{N}$, we have, as  $n \rightarrow \infty$,

$$
f_{n}(k) \rightarrow f(k).
$$
\end{theorem}

\bigskip \noindent \textbf{Comments}. We say that the probability laws  \index{probability law} of binomial random variables  \index{binomial random variable} $X_n$ with parameters $n$ and $p_n$ converges to the probability law of a Poisson random variable  \index{Poisson random variable} $Z$ with parameter $\lambda$, $0<\lambda<+\infty$, as $n\rightarrow +\infty$, when $np_n \rightarrow \lambda$.\\

\noindent \textbf{Proof}. Since $np_{n} \rightarrow \lambda $, we may denote  $\varepsilon(n)=np_{n}-\lambda
\rightarrow 0$. Then $p_n=(\lambda +\varepsilon (n))/n$.\\

\noindent For $k$ fixed and $n>\ k$, we may write,

\begin{eqnarray*}
&&f_{n}(k,p_n)\\
&&=\frac{n!}{k!\times (n-k)!\times n^{k}}(\lambda +\varepsilon (n))^{k}\left(1-\frac{\lambda }{n}-\frac{\varepsilon (n)}{n}\right)^{n-k}.
\end{eqnarray*}

\noindent Hence,

\begin{eqnarray*}
f_{n}(k,~p)&=&\frac{\lambda ^{k}}{k!}\frac{(n-k+1)\times (n-k+2)\times \cdots \times (n-1)\times (n-0)}{n^{k}}\\
&\times& \left(1-\frac{\lambda +\varepsilon (n)}{n}\right)^{n}\left\{ \left(1+\varepsilon(n)/\lambda \right)^{k}\left(1-\frac{\lambda +\varepsilon (n)}{n}\right)^{-k}\right\}\\
&=&\frac{\lambda ^{k}}{k!}\times (1-\frac{k-1}{n})\times (1-\frac{k-2}{n})\times ...\times (1-\frac{1}{n})\times 1\\
&\times& \left(1-\frac{\lambda +\varepsilon (n)}{n}\right)^{n}\left\{ \left(1+\varepsilon(n)/\lambda \right)^{k}\left(1-\frac{\lambda +\varepsilon (n)}{n}\right)^{-k}\right\}\\
&\rightarrow& \frac{\lambda ^{k}}{k!}e^{-\lambda},
\end{eqnarray*}

\Bin
\noindent since

\begin{eqnarray*}
&&\left(1-\frac{k-1}{n}\right)\times \left(1-\frac{k-2}{n}\right)\times \cdots \times \left(1-\frac{1}{n}\right)\times 1\rightarrow 1,\\
&&\left\{ \left(1+\varepsilon (n)/\lambda \right)^{k}\left(1-\frac{\lambda +\varepsilon (n)}{n}\right)^{-k}\right\} \rightarrow 1\\
&& \text{and}\\
&&\left(1-\frac{\lambda +\varepsilon (n)}{n}\right)^{n}\rightarrow e^{-\lambda }.
\end{eqnarray*}

\bigskip

\bigskip

\noindent \textbf{Definition 3}. The discrete random variables $X_{1}$, $X_{2}$, . . ., $X_{n}$ are mutually independent  \index{mutually independent} if and only if
for any subset $\left\{ i_{1},i_{2},...,i_{k}\right\} $ of $\left\{1,2,......,n\right\} $ and for any $k$-tuple $\left(
s_{1},...,s_{k}\right)$, we have the factorizing formula

\begin{equation*}
\mathbb{P}(X_{i_{1}}=s_{1},...,X_{i_{k}}=s_{k})=\mathbb{P}(X_{i_{1}}=s_{1}) ... \mathbb{P}(X_{i_{k}}=s_{k}).
\end{equation*}

\bigskip \noindent The random variables $X_{1}$, $X_{2}$, . . . . . , $X_{n}$ are pairwise independent  \index{pairwise independent} if and only if for any
$1\leq i\neq j\leq n$, for any $s_i$ and $s_j$, we have

\begin{equation*}
\mathbb{P}(X_{i}=s_{i}~,~X_{j}=s_{j})=\mathbb{P}(X_{i}=s_{i})~\mathbb{P}%
(X_{j}=s_{j})\text{.}
\end{equation*}

\bigskip

\noindent \textbf{Example 4.} Consider a Bernoulli experience with parameter $p \in ]0,1[$ independently repeated. Set

\begin{equation*}
X_{i}=\left\{
\begin{array}{c}
\bigskip
1\text{ if the i}^{th}\text{ experience is successful} \\
0\text{ otherwise.}
\end{array}
\right.
\end{equation*}

\Bin Each  $X_{i}$ follows a Bernoulli law with parameter $p$, that is $X_{i}\sim \mathcal{B}(p)$.\\

\noindent The random variables $X_{i}$ are independent by construction. After $n$ repetitions, the number of successes $X$ is exactly

\begin{equation*}
X=\sum_{i=1}^{n}X_{i}.
\end{equation*}

\bigskip \noindent We have the following conclusion  : A Binomial random variable  \index{binomial random variable} $X \sim \mathcal{B}(n,p)$ is always the sum of $n$ independent random variables $X_1$, ..., $X_n$, following all a Bernoulli law $\mathcal{B}(p)$.\\

\bigskip \noindent \textbf{Example 6}.  Consider a Bernoulli experience with parameter $p \in ]0,1[$ independently repeated.\\

\noindent Let $T_{1}$ be the number of experiences needed to get one success.\\

\noindent Let $T_{2}$ be the number of experiences needed to get one success after the first success.\\

\noindent Generally, for $i>1$, let $T_{i}$ be the number of experiences needed to get one success after the $(i-1)^{th}$ success.\\

\noindent It is clear that

\begin{equation*}
T_{1}+...+T_{k}
\end{equation*}

\Bin
\noindent is the number of experiences needed to have $k$ successes, that is

\begin{equation*}
Y_k=T_{1}+...+T_{k}
\end{equation*}

\bigskip \noindent Because of the Memory loss property, the random variables $T_{i-1}$ are independent, so are the $T_i$. We get a  conclusion that is similar to binomial case.\\

\noindent  Conclusion : a Negative Binomial random variable  \index{negative binomial random variable}  \index{binomial random variable} $\mathcal{NB}(k,p)$ is a sum of $k$ independent random variables $T_1$, ...,$T_k$ following all a geometric law $\mathcal{G}(p)$.\\

\bigskip Because of the importance of such results, we are going to state them into lemmas. The first is the following.\\

\begin{lemma} \label{proba01.rv.pd.sumBern}
A Binomial random variable  \index{binomial random variable} $X \sim \mathcal{B}(n,p)$ is always the sum of $n$ independent random variables $X_1$, ..., $X_n$, following all a Bernoulli law $\mathcal{B}(p)$.
\end{lemma}

\Bin Here is the second lemma.\\

\begin{lemma} \label{proba01.rv.pd.sumGeo}
A Negative Binomial random variable  \index{negative binomial random variable}  \index{binomial random variable} $\mathcal{NB}(k,p)$ is a sum of $k$ independent random variables $T_1$, ...,$T_k$ following all a geometric law $\mathcal{G}(p)$.
\end{lemma}
 
 %proba01.rv
\chapter{Real-valued Random Variable Parameters} \label{proba01.param}

\Bin Before we begin, we strongly stress that the mathematical expectation  \index{mathematical expectation} is defined only for a real-valued random variable. If a random variable is not a real-valued one, we may have mathematical expectations of real-valued functions of it.

\Bin  Let us introduce the notion of mathematical expectation  \index{mathematical expectation} in an every day life example.\\

\noindent \textbf{Example 1.} Consider a class of $2$4 students. Suppose that $5$ of them have the age of $19$, $7$ the age of $20$
, $10$ the age of $23$, and $2$ the age of $17$. Let us denote by $m$ the average age of the class, that is :

\begin{equation*}
m=\frac{5\times 19~+~7\times 20~+~10\times 23~+~2\times 17}{24}=20.79.
\end{equation*}

\Bin \noindent Let $X$ be the random variable taking the distinct ages as values, that is $\mathcal{V}(X)=\{x_1,...,x_4\}$ with
$x_{1} = 19$, $x_{2} = 20$, $x_{3} = 23$, $x_{4} = 17$, with probability law  \index{probability law} :

\begin{equation*}
\mathbb{P}(X=x_{1})=\frac{5}{24},~\mathbb{P}(X=x_{2})=\frac{7}{24},~\mathbb{P(}X=x_{3})=\frac{10}{24},~\mathbb{P}(X=x_{4})=\frac{2}{24}.
\end{equation*}

\bigskip \noindent We may summarize this probability law  \index{probability law} in the following table :

\begin{equation*}
\begin{array}{|c|c|c|c|c|}
\hline
k & 19 & 20 & 23 & 17 \\
\hline
P(x=k) & \frac{5}{24} & \frac{7}{24} & \frac{10}{24} & \frac{2}{24}\\
\hline
\end{array}
\end{equation*}

\bigskip \noindent Set $\Omega$ as the class of these $24$ students. The following graph

\begin{equation*}
\begin{array}{cccc}
X: & \Omega & \rightarrow & \mathbb{R} \\
& \omega & \hookrightarrow & X(\omega)\text{= age of the student } \omega.
\end{array}
\end{equation*}

\bigskip \noindent defines a random variable. The probability that is used on $\Omega$ is the frequency  \index{frequency} of occurrence. For each $i\in \{1,2,3,4\}$, $\mathbb{P}(X=x_{i})$ is the frequency of occurrence of the age $x_{i}$ in the class.\\

\noindent Through these notations, we have the following expression of the mean  \index{mean} $m$ :

\begin{equation}
m=\sum_{i=1}^{4}x_{i}~\mathbb{P}(X=x_{i}). \label{math.expect.01}
\end{equation}

\bigskip
\noindent Formula (\ref{math.expect.01}) is the expression of the mathematical expectation  \index{mathematical expectation} of the random variable $X$. \textbf{Mean  \index{mean} values} in real life are particular cases of mathematical expectations in Probability Theory \index{probability theory} .

\Bin Let $X$ be a discrete real-valued random variable taking its values in $\mathcal{V}(X)=\left\{ x_{i}, \ i\in I\right\}$, where
$I$ is a subset of $\mathbb{N}$. If the quantity

\begin{equation}
\sum_{i\in I}x_{i}~\mathbb{P}(X=x_{i}) \label{math.expect.02}
\end{equation}

\Bin
\noindent exists in $\mathbb{R}$, we called it the mathematical expectation  \index{mathematical expectation} of $X$, denoted by $\mathbb{E}(X)$, and we write

\begin{equation}
\mathbb{E}(X)=\sum_{i\in I}x_{i}~\mathbb{P}(X=x_{i}). \label{math.expect}
\end{equation}

\Bin
\noindent If the quantity (\ref{math.expect.02}) is $-\infty$ or $+\infty$, we still name it the mathematical expectation \index{mathematical expectation} ,
with $\mathbb{E}(X)=-\infty$ or $\mathbb{E}(X)=+\infty$.\\

\bigskip \noindent We are now giving some examples.\\

\noindent \textbf{Example 2.}  Let $X$ be a degenerated random variable at $c$, that is
\begin{equation*}
\mathbb{P}(X=c)=1
\end{equation*}

\Bin
\noindent Then, Formula (\ref{math.expect}) reduces to :

\begin{equation*}
\mathbb{E}(X)=c\times \mathbb{P}(X=c)=c.
\end{equation*}

\bigskip \noindent Let us keep in mind that the mathematical expectation  \index{mathematical expectation} of a constant is itself.\\

\noindent \textbf{Example 3} Let $X$ be Bernoulli random variable,  \index{Bernoulli random variable} $X \sim \mathcal{B}(p)$, $p\in ]0,1[$. Formula (\ref{math.expect}) is

\begin{equation*}
\mathbb{E}(X)=1\times \mathbb{P}(X=1)+0\times \mathbb{P}(X=0)=p.
\end{equation*}

\bigskip \noindent \textbf{Example 4.} Let $X$ be a Binomial random variable,  \index{binomial random variable} that is  $X \sim \mathcal{B}(n, p)$, $n\geq 1$ and
$p\in ]0,1[$. Formula (\ref{math.expect}) gives

\begin{eqnarray*}
\mathbb{E}(X)&=&\sum_{k=0}^{n}k~\mathbb{P}(X=k)=\sum_{k=0}^{n}\frac{k~n!}{k!(n-k)!}p^{k}(1-p)^{n-k}\\
&=&\sum_{k=1}^{n}np~\frac{(n-1)!}{(k-1)!(n-k)!} \ p^{k-1}(1-p)^{n-k}\\
&=&np~\sum_{k=1}^{n}C_{n-1}^{k-1}~p^{k-1}(1-p)^{(n-1)-(k-1)}.
\end{eqnarray*}

\bigskip
\noindent Let us make the change of variable : $k^{\prime }=k-1$ and  $n^{\prime }=n-1$. By the Newton's formula  \index{Newton's formula} in Theorem \ref{theo2.4} in Chapter
\ref{proba01.combin},

\begin{eqnarray*}
\mathbb{E}(X)&=&np~\sum_{k^{\prime }=0}^{n^{\prime }} \left(\begin{tabular}{c} $n^{\prime}$ \\ $k^{\prime}$ \end{tabular} \right) ~p^{k^{\prime }}(1-p)^{n^{\prime }-k^{\prime }}\\
&=&np~\left(p+(1-p)\right) ^{n^{\prime }}=np.
\end{eqnarray*}

\bigskip \noindent \textbf{Exercise 1}. Find the mathematical expectation  \index{mathematical expectation} of the a hypergeometric random variable  \index{geometric random variable}  \index{hypergeometric random variable} $\mathcal{H}(N, r, \theta)$,
a Negative Binomial random variable  \index{negative binomial random variable}  \index{binomial random variable} $\mathcal{BN}(k, p)$ and a Poisson random variable$\mathcal{P}(\lambda)$,  \index{Poisson random variable}  with $N\geq 1$, $0<r\leq N$, $\theta\in ]0,1[$, $n\geq 1$, $p\in ]0,1[$, $\lambda>0$.\\

\bigskip

\bigskip \noindent \textbf{Answers}. We will find respectively : $r\theta$ ,$(1-\theta )\left( \frac{1-r}{N}%
\right) $ and $\lambda $.

\bigskip

\Bin Let $X$ and $Y$ be two discrete random variables with respective set of values $\left\{ x_{i}, \ i\in I\right\}$ and $\left\{ y_{j}, \ j\in J\right\}$. Let $g$ anf $h$ be two functions from $\mathbb{R}$ onto $\mathbb{R}$, $k $ a function from $\mathbb{R}^{2}$ onto  $\mathbb{R}$. Finally, $\lambda $ and $\mu
$ are real numbers.\\

\noindent We suppose that all the mathematical expectations  \index{mathematical expectation} used below exist. We have the following properties.\\

\noindent \textbf{(P1). Mathematical Expectation  \index{mathematical expectation} of a function of the random variable}.\\

\begin{equation}
\mathbb{E}(g(X))=\sum_{i\in I}~g(x_{i})~\mathbb{P}(X=x_{i}). \label{FIGD}
\end{equation}

\bigskip\noindent \textbf{(P2).  Mathematical expectation  \index{mathematical expectation} of a real function of a pair of real random variables}.

\begin{equation}
\mathbb{E}(k(X,Y))=\sum_{(i,~j)\in I\times J}~k(x_{i},~y_{j})~\mathbb{P}(X=x_{i},~Y=y_{j}). \label{FIGDDIM2}
\end{equation}

\bigskip\noindent \textbf{(P3). Linearity of the mathematical expectation  \index{mathematical expectation} operator}.\\

\bigskip\noindent The mathematical expectation  \index{mathematical expectation} operator : $X\rightarrow\mathbb{E}(X)$, is linear, that is

\begin{equation*}
\mathbb{E}(\lambda X+\mu Y)=\lambda \mathbb{E}(X)+\mu \mathbb{E}(Y).
\end{equation*}

\bigskip\noindent \textbf{(P4). Non-negativity of the mathematical expectation  \index{mathematical expectation} operator}.\\

\bigskip \noindent \textbf{(i)} The mathematical expectation  \index{mathematical expectation} operator is non-negative, that is

\begin{equation*}
X\geq 0\Rightarrow \mathbb{E}(X)\geq 0.
\end{equation*}

\Bin \noindent \textbf{(ii)} If $X$ is non-negative random variable, then

\begin{equation*}
(\mathbb{E}(X)=0)  \Leftrightarrow (X=0).
\end{equation*}

\Bin
\noindent \textbf{Rule}. A positive random variable has a positive mathematical expectation.  \index{mathematical expectation} A non-negative random variable is constant to zero if its mathematical expectation is zero.\\

\bigskip \noindent \textbf{(P5)}  The mathematical expectation  \index{mathematical expectation} operator is non-decreasing, that is

\begin{equation*}
X\leq Y\Rightarrow \mathbb{E}(X)\leq \mathbb{E}(Y).
\end{equation*}

\bigskip \noindent \textbf{(P6)}  Factorization Formula of the mathematical expectation  \index{mathematical expectation} of a product of real functions of two independent random variables : the random variables $X$ and $Y$ are independent if and only if for any real function $g$ of $X$ and $h$ of $Y$ such that the mathematical expectations of $X$ and of $Y$ exist, we have

\begin{equation*}
\mathbb{E}(g(X)\times h(Y))=\mathbb{E}(g(X))\times ~\mathbb{E}(h(Y))
\end{equation*}

\bigskip

\noindent \textbf{Remarks} : In Formulas (P1), (P2) and (P6), the random variables $X$ and $Y$ are not necessarily real-valued random variables. They may be of any kind : vectors, categorical objects, etc.\\

\bigskip \noindent \textbf{Proofs.}\\

\noindent {\bf Proof of (P1)}. Set $Z=g(X),~X\in \left\{ x_{i},~i\in I\right\}$. The values set of $Z$ are the distinct values of
$g(x_i)$, $i\in I$. Let us denote the values set by $\left\{ z_{k},~k\in K\right\}$, where necessarily, $Card(K) \leq Card(I)$, since

\begin{equation*}
\left\{ z_{k},~k\in K\right\} \subseteq \left\{ g(x_{i}),~i\in I\right\} .
\end{equation*}

\Bin
\noindent For each $k \in K$, we are going to define the class of all $x_i$ such that $g(x_i)=z_k$, that is
$$
I_{k}=\left\{ i\in I\text{ , }g(x_{i})=z_{k}\right\}.
$$

\Bin
\noindent It comes that for $k\in K$,

\begin{equation}
(Z=z_{k})=\underset{i\in I_{k}}{\cup }(X=x_{i}). \label{proba01.paramRealRv.class01}
\end{equation}

\noindent It is also obvious that the $I_{k},k\in K$, form a partition of $I$, that is

\begin{equation}
I= \sum_{k\in K} I_{k}. \label{proba01.paramRealRv.part01}
\end{equation}

\noindent Formula (\ref{proba01.paramRealRv.class01}) implies that

\begin{equation}
\mathbb{P}(Z=z_{k})=\sum_{i\in I_{k}}\mathbb{P}(X=x_{i}). \label{proba01.paramRealRv.class02}
\end{equation}

\noindent By definition, we have

\begin{equation}
\mathbb{E}(Z)=\sum_{k\in K}~z_{k}~\mathbb{P}(Z=z_{k}). \label{part6.efunction}
\end{equation}

\noindent Let us combine all these facts to conclude. By the partition (\ref{proba01.paramRealRv.part01}), we have

\begin{equation*}
\sum_{i\in I}~g(x_{i})~\mathbb{P}(X=x_{i})= \sum_{k\in K} \sum_{i\in I_k}~g(x_{i})~\mathbb{P}(X=x_{i}).
\end{equation*}

\noindent Now, we use the remark that $g(x_{i})=z_{k}$ on $I_k$, in the latter expression, to get

\begin{eqnarray*}
\sum_{k\in K} \sum_{i\in I_k}~g(x_{i})~\mathbb{P}(X=x_{i})&=& \sum_{k\in K} \sum_{i\in I_k}~z_k~\mathbb{P}(X=x_{i})\\
\sum_{k\in K} z_{k} \biggr( \sum_{i\in I_k}~~\mathbb{P}(X=x_{i})\biggr)\\
&=& \sum_{k\in K} z_{k} \mathbb{P}(Z=z_{k}),\\
\end{eqnarray*}

\noindent where we used (\ref{proba01.paramRealRv.class02}) in the last line. So we obtained (\ref{part6.efunction}). The proof is finished. $\square$\\

\bigskip
\noindent {\bf Proof of (P2)}. This proof of the property is very similar to that of (P1). Set $V=k(X,~Y),~X\in \left\{ x_{i},~i\in I\right\}$ and $Y\in \left\{ y_{j},~j\in J\right\}$. As in the proof of (P1), $V$ take the distinct values among the values $k(x_i,y_j)$, $(i,j)\in I\times J$, and we write

\begin{equation*}
\mathcal{V}_V =\left\{ v_{k},~k\in K\right\} \subseteq \left\{k(x_{i,}~y_{j}),~(i,~j)\in I\times J\right\} .
\end{equation*}

\bigskip \noindent Similarly, we partition $\mathcal{V}_V$ in the following way :

\begin{equation*}
I_{k}=\left\{ (i,~j)\in I\times J,k(x_{i},~y_{j})=v_{k}\right\} .
\end{equation*}

\Bin so that

\begin{equation}
I\times J= \sum_{k\in K} I_{k} \label{proba01.paramRealRv.part02}
\end{equation}

\Bin \noindent and

\begin{equation}
(V=v_{k})=\underset{(i,j) \in I_{k}}{\cup }(X=x_{i},Y=y_j). \label{proba01.paramRealRv.class03}
\end{equation}

\Bin \noindent Formula (\ref{proba01.paramRealRv.class02}) implies

\begin{equation}
\mathbb{P}(V=v_{k})=\sum_{(i,j) \in I_{k}}\mathbb{P}(X=x_{i},Y=y_j). \label{proba01.paramRealRv.class04}
\end{equation}

\Bin By definition, we have

\begin{equation}
\mathbb{E}(V)=\sum_{k\in K}~v_{k}~\mathbb{P}(V=v_{k}). \label{part6.efunction02}
\end{equation}

\Bin Next, by the partition (\ref{proba01.paramRealRv.part02}), we have that $k(x_{i},y_j)=v_{k}$ on $I_k$ and we are lead to

\begin{eqnarray*}
\sum_{(i,j)\in I\times J}~k(x_{i},y_j)~\mathbb{P}(X=x_{i},Y=y_j)&=& \sum_{k\in K} \sum_{(i,j)\in \in I_k}~ v_k~\mathbb{P}(X=x_{i},Y=y_j)\\
&=& \sum_{k\in K} v_k \ \biggr( \sum_{(i.j)\in I_k}\mathbb{P}(X=x_{i},Y=y_j) \biggr)\\
&=& \sum_{k\in K} v_k \ \mathbb{P}(V=v_k)
\end{eqnarray*}

\noindent where we used (\ref{proba01.paramRealRv.class04}) in the last line. We reached (\ref{part6.efunction02}). The proof is finished.$\square$\\

\bigskip
\noindent {\bf Proof of (P3)}. First, let us apply (P1) with $g(X)=\lambda X$, to get

\begin{equation*}
\mathbb{E}(\lambda X)=\sum_{i\in I}\lambda ~x_{i}~\mathbb{P}%
(X=x_{i})=\lambda ~\mathbb{E}(X).
\end{equation*}

\noindent Next, we apply (P2) with $k(x,~y)=x+y$, to get

\begin{eqnarray*}
\mathbb{E}(X+Y) &=&\sum_{(i,~j)\in I\times J}(x_{i}+y_{j})~\mathbb{P}(X=x_{i},~Y=y_{j}) \\
&=&\sum_{i\in I}x_{i}\left\{ \sum_{j\in J}\mathbb{P}(X=x_{i}, \ Y=y_{j})\right\} +\sum_{j\in J}y_{j}\left\{ \sum_{i\in I}\mathbb{P}(X=x_{i},~Y=y_{j})%
\right\} \text{.}
\end{eqnarray*}

\Ni But

\begin{equation*}
\Omega =\underset{j\in J}{\cup }(Y=y_{j})=\sum_{j\in J}(Y=y_{j}).
\end{equation*}

\Ni Thus,

\begin{equation*}
(X=x_{i})=(X=x_{i})\cap \Omega =\sum_{j\in J}(X=x_{i},~Y=y_{j}),
\end{equation*}

\Ni and then

\begin{equation*}
\mathbb{P}(X=x_{i})=\sum_{j\in J}\mathbb{P}(X=x_{i},~Y=y_{j}).
\end{equation*}

\noindent By the symmetry of the roles of $X$ and $Y$, we also have

\begin{equation*}
\mathbb{P}(Y=y_{j})=\sum_{i\in I}\mathbb{P}(X=x_{i},~Y=y_{j}).
\end{equation*}

\noindent Finally,

\begin{equation*}
\mathbb{E}(X+Y)=\sum_{i\in I}x_{i}~\mathbb{P}(X=x_{i})+\sum_{j\in J}y_{j}~\mathbb{P}(Y=y_{j})=\mathbb{E}(X)+\mathbb{E}(Y).
\end{equation*}

\Ni We just finished the proof that the mathematical expectation  \index{mathematical expectation} is linear.\\

\noindent {\bf Proof of (P4)}. \\

\noindent \textbf{Part (i)}.\\

\noindent Let  $X\in \left\{ x_{i},~i\in I\right\} $ and  $X\geq 0$. This means  \index{mean} that
the values of $X$, which are the $x_{i}$'s, are non-negative. Hence, we have

\begin{equation*}
\mathbb{E}(X)=\sum_{i\in I}x_{i}~\mathbb{P}(X=x_{i})\geq 0,
\end{equation*}

\Bin as a sum of the product of $x_{i}$ and $\mathbb{P}(X=x_{i})$ both non-negative.\\

\noindent \textbf{Part (ii)}.\\

\noindent Assume that $X$ is non-negative. This means  \index{mean} that the values  $x_{i}$ are non-negative. $\mathbb{E}(X)=0$ means that

\begin{equation*}
\mathbb{E}(X)=\sum_{i\in I}x_{i}~\mathbb{P}(X=x_{i})=0.
\end{equation*}

\Bin
\noindent Since the left-hand member of that equation is formed by a sum of non-negative terms, it is null if and only if each of the terms is null, that is for each  $i\in I$,

\begin{equation*}
x_{i}~\mathbb{P}(X=x_{i})=0.
\end{equation*}

\Bin
\noindent Hence, if $x_{i}$ is a value of $X$, we have $\mathbb{P}(X=x_{i})>0$, and then $x_{i}=0$. This means  \index{mean} that $X$ only takes the value $0$. Thus, we have $X=0$.

\bigskip \noindent {\bf Proof of (P5)}. Let $X\leq Y$. By linearity and, since the mathematical expectation  \index{mathematical expectation} is a non-negative operator,

\begin{equation*}
0\leq \mathbb{E}(Y-X)=\mathbb{E}(Y) - \mathbb{E}(X).
\end{equation*}

\noindent Hence $\mathbb{E}(Y) \geq \mathbb{E}(X)$.\\

\noindent {\bf Proof of  (P6)}. Suppose that $X$ and $Y$ are independent. Set
$k(x, y)= g(x)h(y)$. We may apply (P2) to get

\begin{equation*}
\mathbb{E}(g(X)h(Y))=\sum_{(i,j)}~g(x_{i})h(y_{j})~\mathbb{P}(X=x_{i},~Y=y_{j}).
\end{equation*}

\noindent But, by independence \index{independence} ,

\begin{equation*}
\mathbb{P}(X=x_{i},~Y=y_{j})=\mathbb{P}(X=x_{i})~\mathbb{P}(Y=y_{j}),
\end{equation*}

\Bin
\noindent for all $(i, j) \in I \times J$. Then, we have

\begin{eqnarray*}
\mathbb{E}(g(X)h(Y))&=&\sum_{j\in J}~\sum_{i\in I}g(x_{i})~\mathbb{P}(X=x_{i})~h(y_{j})~\mathbb{P}(Y=y_{j})\\
&=&\sum_{i\in I}g(x_{i})~\mathbb{P}(X=x_{i})\times \sum_{j\in J}h(y_{j})~\mathbb{P}(Y=y_{j})\\
&=&\mathbb{E}(g(X))~\mathbb{E}(h(Y)).
\end{eqnarray*}

\Bin
\noindent This proves the direct implication.\\

\noindent Now, suppose that

\begin{equation*}
\mathbb{E}(g(X)h(Y))=\mathbb{E}(g(X))\times ~\mathbb{E}(h(Y)),
\end{equation*}

\noindent for all functions $g$ and $h$ for which the mathematical expectations  \index{mathematical expectation} make sense. Consider two particular values in $I\times J$, $i_{0}\in I$ and $j_{0}\in J$. Set

\begin{equation*}
g(x)=\left\{
\begin{array}{c}
0~\text{if}~x\neq x_{i_{0}} \\
1~\text{if}~x=x_{i_{0}}
\end{array}
\right.
\end{equation*}

\Bin
\noindent and

\begin{equation*}
h(y)=\left\{
\begin{array}{c}
0\text{ if }y\neq y_{j_{0}} \\
1\text{ if }y=y_{j_{0}}
\end{array}
\right.
\end{equation*}

\Bin We have

\begin{eqnarray*}
\mathbb{E}(g(X)h(Y))&=&\sum_{(i,~j)\in I\times J}g(x_{i})h(y_{j})~\mathbb{P}(X=x_{i},~Y=y_{j})\\
&=&\mathbb{P}(X=x_{i_{0}},~Y=y_{j_{0}})\\
&=& \mathbb{E}g(X)\times \mathbb{E}h(Y)\\
&=&\sum_{i\in I}g(x_{i})~\mathbb{P}(X=x_{i})~\times \sum_{j\in J}h(y_{j})\mathbb{P}(Y=y_{j})\\
&=&\mathbb{P}(X=x_{i_{0}})~\mathbb{P}(Y=y_{j_{0}}).
\end{eqnarray*}

\Bin We proved that for any $(i_{0}, j_{0}) \in I \times J$, we have

\begin{equation*}
\mathbb{P}(X=x_{i_{0}},~Y=y_{j_{0}})=\mathbb{P}(X=x_{i_{0}})~\mathbb{P}(Y=y_{j_{0}}).
\end{equation*}

\Bin We then have the independence \index{independence} .

\bigskip

\Bin Let $X$ and $Y$ be a discrete real-valued random variables with respective values set $\left\{x_{i},~i\in I\right\} $ and $\left\{ y_{j},~j\in J\right\} $.\\

\noindent The most common parameters of such random variables are defined below.

\bigskip

\noindent {\bf $k^{th}$ non centered moments  \index{moments} or centered moments of order $k\geq 1$} :

\begin{equation*}
m_{k}=\mathbb{E}(X^{k})=\sum_{i\in I}x_{i}^{k}~\mathbb{P}(X=x_{i}).
\end{equation*}

\noindent {\bf Centered Moments  \index{moments} of order $k \geq 1$} :

\begin{equation*}
M_{k}=\mathbb{E}((X-\mathbb{E}(X))^{k})=\sum_{i\in I}(x_{i}-\mathbb{E}(X))^{k}~\mathbb{P}(X=x_{i}).
\end{equation*}

\noindent {\bf Variance  \index{variance} and Standard deviation \index{Standard deviation} }.\\

\noindent The variance  \index{variance} of $X$ is defined by :

\begin{equation*}
\sigma _{X}^{2}=Var(X)=\mathbb{E}((X-\mathbb{E}(X))^{2})=\sum_{i\in I}(x_{i}-\mathbb{E}(X))^{2}~\mathbb{P}(X=x_{i}).
\end{equation*}

\noindent The square root of the variance  \index{variance} of $X$, $\sigma _{X}=(Var(X))^{\frac{1}{2}}$, is the standard deviation  \index{Standard deviation} of $X$.\\

\noindent \textbf{Remark}. The variance  \index{variance} is the centered moment of order 2.\\

\bigskip \noindent {\bf Factorial moment  \index{factorial moment} of order two of $X$} :

\begin{equation*}
fm_{2}(X)=\mathbb{E}(X(X-1))=\sum_{i \in I} x_{i} \ (x_{i}-1) \ \mathbb{P}(X=x_{i}).
\end{equation*}

\Bin {\bf Covariance  \index{variance} between $X$ and $Y$} :

\begin{equation*}
Cov(X,Y)=\mathbb{E}((X-\mathbb{E}(X))(Y-\mathbb{E}(Y)))=\sum_{(i,~j) \in I\times J}(x_{i}-\mathbb{E}(X))(y_{j}-\mathbb{E}(Y))~\mathbb{P}(X=x_{i},~Y=y_{j}).
\end{equation*}

\bigskip \noindent {\bf Linear correlation coefficient between $X$ and $Y$}\\

\noindent If  $\sigma _{X}\neq 0$ and $\sigma _{Y}\neq 0$, we may define the number

\begin{equation*}
\sigma _{XY}=\frac{Cov(X,Y)}{\sigma _{X}~\sigma _{Y}}
\end{equation*}

\noindent as the linear correlation coefficient.\\

\noindent Now, we are going to review important properties of these parameters.\\

\noindent \textbf{(1) Other expression of the variance  \index{variance} and of the covariance}.\\

\noindent \textbf{(i)} The variance  \index{variance} and the covariance have the following two alternative expressions :

\begin{equation*}
Var(X)=\mathbb{E}(X^{2})-\mathbb{E}(X)^{2}
\end{equation*}

\Bin
\noindent and

\begin{equation*}
Cov(X,Y)=\mathbb{E}(XY)-\mathbb{E}(X)~\mathbb{E}(Y).
\end{equation*}

\Bin
\noindent \textbf{Rule}. The variance  \index{variance} is the difference of the non centered moment of order two \textit{and} the square of the mathematical expectation \index{mathematical expectation} .\\

\noindent \textbf{Rule}. The covariance  \index{variance} between $X$ and $Y$ is the difference of mathematical expectation  \index{mathematical expectation} of the product of $X$ and $Y$ \textit{and} the product of the mathematical expectations of $X$ and $Y$.\\

\noindent \textbf{(ii)} The variance  \index{variance} can be computed from the factorial moment  \index{factorial moment} of order 2 by

\begin{equation}
Var(X)=fm_{2}(X)+\mathbb{E}(X)-\mathbb{E}(X)^{2}. \label{varFactoMoment}
\end{equation}

\bigskip
\noindent \textbf{Remark}. In a number of cases, the computation of the second factorial moment  \index{factorial moment} is easier than that of the second moment, and Formula
(\ref{varFactoMoment}) becomes handy.\\

\bigskip
\noindent \textbf{(2). The standard deviation  \index{Standard deviation} is zero or the variance  \index{variance} is zero if and only if the random variable is constant and is equal to its mathematical expectation}  \index{mathematical expectation} :

\begin{equation*}
\sigma _{X}=0\text{ if and only if }X=\mathbb{E}(X).
\end{equation*}

\bigskip

\noindent \textbf{(3) Variance  \index{variance} of a linear combination  \index{combination} of random variables}.\\

\noindent \textbf{(i) Variance  \index{variance} of the sum of two random variables} :

\begin{equation*}
Var(X+Y)=Var(X)+Var(Y)+2Cov(X,Y).
\end{equation*}

\bigskip \noindent \textbf{(ii) Independence  \index{independence} and Covariance}.  \index{variance} If $X$ and  $Y$ are independent then

\begin{equation*}
Var(X+Y)=Var(X)+Var(Y)\text{ and }Cov(X,Y)=0\text{.}
\end{equation*}

\Bin
\noindent We conclude by :\\

\noindent (3ii-a) The variance  \index{variance} of a sum of independent random variables is the sum of their variance.\\

\noindent (3ii-b) If $X$ and $Y$ are independent, they are uncorrelated that is,

$$
Cov(X,Y)=0.
$$

\bigskip \noindent Two random variables $X$ and $Y$ are linearly uncorrelated if and only $Cov(X,Y)=0$. This is implied by the independence.  \index{independence} But lack of correlation does not imply independence. This may happen for special class of random variables, like Gaussian ones.\\

\bigskip \noindent \textbf{(iii) Variance  \index{variance} of linear combination  \index{combination} of random variables}.\\

\noindent Let $X_{i}$, $i\geq 1$, be a sequence of real-valued random variables. Let $a_{i}$, $i\geq 1$, be a sequence of real numbers. The variance  \index{variance} of the linear combination \index{combination} 
$$
\sum_{i=1}^{k}a_{i}X_{i}
$$

\noindent is given by :

\begin{equation*}
Var\left(\sum_{i=1}^{k} a_{i}X_{i}\right)=\sum_{i=1}^{k}a_{i}^{2}~Var(X_{i})+2\sum_{1%
\leq i<j\leq k}a_{i}a_{j}~Cov(X_{i},X_{j}).
\end{equation*}

\bigskip

\noindent \textbf{(4). Cauchy-Schwarz's inequality}.

\begin{equation*}
\left| Cov(X,Y)\right| \leq \sigma _{X}~\sigma _{Y}\text{.}
\end{equation*}

\bigskip

\noindent \textbf{(5). H\"{o}lder inequality}.\\

\noindent Let $p>1$ and $q>1$ such that  $1/p+1/q=1$. Suppose that

$$
\left\| X\right\| _{p}=\mathbb{E}(\left| X\right|^{p})^{1/p}
$$

\Bin
\noindent and

$$
\left\| Y\right\| _{q}=\mathbb{E}(\left| Y\right|^{q})^{1/q}
$$

\noindent exist. Then

\begin{equation*}
\left| \mathbb{E}(XY)\right| \leq \left\| X\right\|_{p}\times \left\|Y\right\|_{q}.
\end{equation*}

\bigskip \noindent

\noindent \textbf{(6)  Minkowski's Inequality}.\\

\noindent  For any $p\geq 1$,

\begin{equation*}
\left\| X+Y\right\| _{p}\leq \left\| X\right\| _{p}+\left\| Y\right\| _{p}.
\end{equation*}

\bigskip \noindent The number $\left\| X\right\| _{p}$ is called the $L_p$-norm of $X$.\\

\bigskip \noindent

\noindent \textbf{(7)  Jensen's Inequality}.\\

\noindent Let us consider a real-valued random variable taking the values $(x_i)_{i\in I}$ and let us denote

$$
\forall i \in I, \ p_i=\mathbb{P}(X=x_i).
$$

\Bin
\noindent  If $g : I \mapsto \mathbb{R}$ is a convex function  \index{convex function} on an interval including $\mathcal{V}_X$, and if

$$
\mathbb{E}(|X|)=\sum_{i\in I} p_{i} |x_i|<+\infty
$$

\noindent and

$$
\mathbb{E}(|g(X)|)=\sum_{i\in I} p_{i} |g(x_i)|<+\infty,
$$

\noindent then

$$
g(\mathbb{E}(X))\leq \mathbb{E}(g(X)).
$$

\bigskip \noindent \textbf{Proofs}. Let us provide the proofs of these properties one by one in the same order in which they are stated.\\

\noindent {\bf Proof of (1)}.\\

\noindent Part \textbf{(i)}. Let us apply the linearity of the mathematical expectation  \index{mathematical expectation} stated in Property (P3) above to the following development

$$
(X-\mathbb{E}(X))^{2}=X^{2}+2X\mathbb{E}(X)+\mathbb{E}(X)^{2}
$$

\noindent to get :

\begin{eqnarray*}
\mathbb{V}ar(X)=\mathbb{E}(X-\mathbb{E}(X))^{2}&=&\mathbb{E}(X^{2}-2~X~\mathbb{E}(X)+\mathbb{E}(X)^{2})\\
&=&\mathbb{E}(X^{2})-2\mathbb{E}(X)\mathbb{E}(X)+\mathbb{E}(X)^{2}\\
&=&\mathbb{E}(X^{2})-\mathbb{E}(X)^{2}.
\end{eqnarray*}

\Bin
\noindent We also have, by the same technique,

\begin{eqnarray}
Cov(X,Y)&=&\mathbb{E}\{(X-\mathbb{E}(X))(Y-\mathbb{E}(Y))\} \notag \\
&=&\mathbb{E}(XY-X\mathbb{E}(Y)-Y\mathbb{E}(X)+\mathbb{E}(X)~\mathbb{E}(Y)) \notag \\
&=&\mathbb{E}(XY)-\mathbb{E}(X)~\mathbb{E}(Y). \label{cov00}
\end{eqnarray}

\Bin
\noindent Let us notice that the following identity

\begin{equation*}
Var(X)=Cov(X,X)
\end{equation*}

\Bin \noindent and next remark that the second formula implies the first.\\

\noindent From the formula of the covariance,  \index{variance} we see by applying Property (P6), that the independence  \index{independence} of $X$ and $Y$ implies

\begin{equation*}
\mathbb{E}(XY)=\mathbb{E}(X)\times \mathbb{E}(Y)
\end{equation*}

\noindent and then

\begin{equation*}
\mathbb{C}ov(X,Y)=0
\end{equation*}

\Bin
\noindent Part \textbf{(ii)} It is enough to remark that

$$
fm_{2}(X)=\mathbb{E}(X(X-1))=\mathbb{E}(X^2)-\mathbb{E}(X),
$$

\Bin
\noindent and to use it the Formula of the variance  \index{variance} above.\\

\bigskip \noindent {\bf Proof of (2)}.  Let us begin to prove the direct implication. Assume $\sigma_{X}=0$. Then

\begin{equation*}
\sigma _{X}^{2}=\sum_{i}(x_{i}-\mathbb{E}(X))^{2}~\mathbb{P}(X=x_{i})=0.
\end{equation*}

\noindent Since this sum of non-negative terms is null, each of the terms is null, that is for each  $i\in I$,

\begin{equation*}
(x_{i}-E(X))^{2}~\mathbb{P}(X=x_{i})=0.
\end{equation*}

\Bin \noindent Hence, if $x_{i}$ is a value of $X$, we have $\mathbb{P}(X=x_{i}) > 0$, and then $x_{i}-E(X)=0$, that is $x_{i}=\mathbb{E}(X)$. This means  \index{mean} that $X$ only takes the value $\mathbb{E}(X)$.\\

\noindent Conversely, if $X$ takes one unique value $c$, we have $\mathbb{P}(X=c)=1$. Then for any $k\geq 1$,

$$
\mathbb{E}(X^{k})=c^{k}~\mathbb{P}(X=c)=c^{k}.
$$

\Bin \noindent Finally, the variance  \index{variance} is

\begin{equation*}
Var(X)=\mathbb{E}(X^{2})-\mathbb{E}(X)^{2}=c^{2}-c^{2}=0.
\end{equation*}

\bigskip \noindent {\bf Proof of (3)}. We begin by the general case of a linear combination  \index{combination} of random variables of the form
$Y=\sum_{i=1}^{k}a_{i}X_{i}$. By the linearity Property (P6) of the mathematical expectation,  \index{mathematical expectation} we have

\begin{equation*}
\mathbb{E}(Y)=\sum_{i=1}^{k}a_{i}\mathbb{E(}X_{i}).
\end{equation*}

\noindent It comes that
\begin{equation*}
Y-\mathbb{E}(Y)=\sum_{i=1}^{k}a_{i}(X_{i}-\mathbb{E(}X_{i})).
\end{equation*}

\noindent Now, we use the full development of products of two sums

\begin{equation*}
(Y-\mathbb{E}(Y))^{2}=\sum_{i=1}^{k}a_{i}^{2}(X_{i}-E(X_{i}))^{2}+2\sum_{1%
\leq i<j\leq k}a_{i}a_{j}~(X_{i}-\mathbb{E}(X_{i}))(X_{j}-\mathbb{E}(X_{j})).
\end{equation*}

\noindent We apply the mathematical expectation  \index{mathematical expectation} at both sides of the equality above to get

\begin{equation*}
Var(Y)=\sum_{i=1}^{k}a_{i}^{2}~Var(X_{i})+2\sum_{1\leq i<j\leq
k}a_{i}a_{j}~Cov(X_{i},X_{j}).
\end{equation*}

\noindent In particular, if $X_{1}$,...,$X_{k}$ are pairwise independent  \index{pairwise independent} (let alone mutually independent),  \index{mutually independent} we have

\begin{equation*}
Var(Y)=\sum_{i=1}^{k}a_{i}^{2}Var(X_{i}).
\end{equation*}

\bigskip \noindent {\bf Proof of (4)}. We are going to proof the Cauchy-Schwarz's inequality by studying the sign of the trinomial function
$VarX+\lambda ^{2}Var(Y)+2\lambda Cov(X,Y)$, where $\lambda \in \mathbb{R}$.\\

\noindent We begin by remarking that if  $Var(Y)=0$, then by Point (2) above, $Y$ is a constant, say $c=\mathbb{E}(Y)$. Since
$Y-\mathbb{E}(Y)=c-c$, this implies that

$$
\mathbb{C}ov(X,c)=\mathbb{E}(X-\mathbb{E}(X))(Y-\mathbb{E}(Y)))=\mathbb{E}(X-\mathbb{E}(X))(0)=0.
$$

\noindent In that case, the inequality

$$
|\mathbb{C}ov(X,Y)|\leq \sigma_{X} \sigma_{X},
$$

\Bin
\noindent is true, since the two members are both zero.\\

\noindent Now suppose that $Var(Y)>0$. By using Point (3) above, we have

\begin{equation*}
0 \leq Var(X+\lambda ~Y)=VarX+\lambda ^{2}Var(Y)+2\lambda \mathbb{C}ov(X,Y).
\end{equation*}

\noindent So the third order polynomial function $VarX+\lambda ^{2}Var(Y)+2\lambda Cov(X,Y)$ in $\lambda$ has the constant non-negative sign. This is possible if and only if the discriminant is non-positive, that is

\begin{equation*}
\Delta ^{\prime }=Cov(X,Y)^{2}-VarX\times VarY \leq 0.
\end{equation*}

\Bin \noindent This leads to\\

\begin{equation*}
\mathbb{C}ov(X,Y)^{2}\leq VarX\times VarY.
\end{equation*}

\Bin
\noindent By taking the square roots, we get

$$
|\mathbb{C}ov(X,Y)|\leq \sigma_{X} \sigma_{Y},
$$

\noindent which is the searched result.\\

\noindent {\bf Proof of (5)}. We are going to use the following  inequality on $\mathbb{R}$ :

\begin{equation}
\left| ab\right| \leq \frac{\left| a\right| ^{p}}{p}+\frac{\left| b\right|
^{q}}{q}. \label{ineqCauchyR}
\end{equation}

\noindent We may easily check if $p$ and $q$ are positive integers such that $1/p+1/q=1$, both $p$ and $q$ are greater than the unity and the formulae
\begin{equation*}
p=\frac{q}{q-1}\text{ and \ }q=\frac{p}{p-1}.
\end{equation*}

\noindent also hold. Assume first that we have : $\left\| X\right\| _{p}\neq 0$\ and $\left\| Y\right\| _{q}\neq
0$. Set

\begin{equation*}
a=\frac{X}{\left\| X\right\| _{p}}\text{ \ \ and \ \ \ }b=\frac{Y}{\left\|
Y\right\| _{q}}.
\end{equation*}

\Bin
\noindent By using Inequality (\ref{ineqCauchyR}), we have

\begin{equation*}
\frac{\left| X\times Y\right| }{\left\| X\right\| _{p}\times \left\|
Y\right\| _{q}}\leq \frac{\left| X\right| ^{p}}{p\times \left\| X\right\|
_{p}^{p}}+\frac{\left| Y\right| ^{q}}{q\times \left\| Y\right\| _{q}^{q}}.
\end{equation*}

\bigskip \noindent By taking the mathematical expectation,  \index{mathematical expectation} we obtain

\begin{equation*}
\frac{\mathbb{E}\left| X\times Y\right| }{\left\| X\right\| _{p}\times
\left\| Y\right\| _{q}}\leq \frac{\mathbb{E}\left| X\right| ^{p}}{p\times
\left\| X\right\| _{p}^{p}}+\frac{\mathbb{E}\left| Y\right| ^{q}}{q\times
\left\| Y\right\| _{q}^{q}}.
\end{equation*}

\Bin \noindent This yields

\begin{equation*}
\frac{\mathbb{E}\left| X\times Y\right| }{\left\| X\right\| _{p}\times
\left\| Y\right\| _{q}}\leq \frac{1}{p}+\frac{1}{q}=1.
\end{equation*}

\bigskip \noindent We proved the inequality if neither of $\left\| X\right\| _{p}$\ and $\left\| Y\right\| _{q}$ are null. If one of them is zero, say $\left\| X\right\| _{p}$, this mean  \index{mean} that

$$
\mathbb{E}(|X|^p)=0.
$$

\Bin
\noindent By Property Point (ii) of Property (p4) above, $|X|^p=0$ and hence $X=0$.\\

\noindent In this case, both members in the H\"{o}lder inequality are zero and then, it holds.\\

\noindent {\bf Proof of (6)}. We are going to use the H\"{o}lder inequality to establish the Minkowski's one.\\

\noindent First, if $\left\| X+Y\right\| _{p}=0$, there is nothing to prove. Now, suppose that $\left\| X+Y\right\| _{p} >0$.\\

\noindent We have

\begin{equation*}
\left| X+Y\right| ^{p}=\left| X+Y\right| ^{p-1}\left| X+Y\right| \leq \left|
X+Y\right| ^{p-1}\left| X\right| +\left| X+Y\right| ^{p-1}\left| Y\right|
\end{equation*}

\Bin
\noindent By taking the expectations, and by using the H\"{o}lder to each of the two terms the right-hand member of the inequality and by reminding that $p=p/(p-1)$, we get

\begin{eqnarray*}
\mathbb{E}\left| X+Y\right| ^{p} & \leq& \mathbb{E}\left| X+Y\right|^{p-1}\left| X\right| +\mathbb{E}\left| X+Y\right| ^{p-1}\left| Y\right|\\
&\leq& (\mathbb{E}\left| X\right| ^{p})^{1/p}\times \mathbb{E} \ (\left|X+Y\right| ^{q(p-1)})^{1/q}+(\mathbb{E}\left| Y\right| ^{p})^{1/p}\times (\mathbb{E}\left| X+Y\right| ^{q(p-1)})^{1/q}.
\end{eqnarray*}

\noindent Now, we divide the last formula by

\begin{equation*}
(\mathbb{E}\left| X+Y\right| ^{q(p-1)})^{1/q}=(\mathbb{E}\left| X+Y\right|
^{p})^{1/q},
\end{equation*}

\noindent which is not null, to get

\begin{equation*}
(\mathbb{E}\left| X+Y\right| ^{p})^{1-1/q}\leq \left\| X\right\|
_{p}+\left\| Y\right\| _{p}.
\end{equation*}

\bigskip \noindent This was the target.\\

\bigskip \noindent {\bf Proof of (7)}. Let $g : I \rightarrow \mathbb{R}$ be a convex function  \index{convex function} such that $I$ includes $\mathcal{V}_X$. We will see two particular cases and a general case.\\

\noindent Case 1 : $\mathcal{V}_X$ is finite. In that case, let us denote $\mathcal{V}_X=\{x_1,\cdots, x_k\}$, $k\geq 1$. If $g$ is convex, then by Formula
(\ref{proba01_convexK}) in Section \ref{proba01_appendix_convex} in the Appendix Chapter \ref{proba01_appendix}, Part (D-D1), we have

\begin{equation}
g(p_{1} x_{1}+...+p_{k} x_{k})\leq p_{1}g(x_{1})+ \cdots + p_{k}g(x_k).
\end{equation}

\Bin \noindent  By the property (P1) of the mathematical expectation  \index{mathematical expectation} above,  the left-hand member of this inequality is $g(\mathbb{E}(X))$ and the right-hand member is $\mathbb{E}(g(X))$.\\

\noindent Case 2 : $\mathcal{V}_X$ is countable infinite, $g$ is a bounded function or $I$ is a bounded and closed interval. In that case, let us denote $\mathcal{V}_X=\{x_i,i\geq 1\}$. If $g$ is convex, then by Formula (\ref{proba01_convexInf}) Section \ref{proba01_appendix_convex} in the Appendix Chapter \ref{proba01_appendix} , Part (D-D2), we have

\begin{equation}
g(p_{1} x_{1}+...+p_{k} x_{k})\leq p_{k}g(x_{1})+...+p_{k}g(x_k)
\end{equation}

\Bin
\noindent By the property (P1) of the mathematical expectation  \index{mathematical expectation} above,  the left-hand member of this inequality is $g(\mathbb{E}(X))$ and the right-hand member is $\mathbb{E}(g(X))$.\\

\noindent Case 3 : Let us suppose that $\mathcal{V}_X$ is countable infinite only. Denote $\mathcal{V}_X=\{x_i,i\geq 1\}$. Here, $g$ is defined on the whole real line $\mathbb{R}$. Let $J_n=[-n,n]$, $n\geq 1$. Put, for $n\geq 1$,

$$
I_{n}=\{ i, x_i\in J_n\},
$$

$$
\mathcal{V}_{n}=\{x_i, i\in I_n\}
$$

\Bin
\noindent and

$$
p(n)=\sum_{i \in I_n} p_i.
$$

\Bin \noindent Let $X_n$ be a random variable taking the values of $\mathcal{V}_{n}=\{x_i, i\in I_n\}$ with the probabilities

$$
p_{i,n}=p_{i}/p(n).
$$

\noindent By applying the Case 1, we have $g(\mathbb{E}(X_n))\leq \mathbb{E}g(X_n)$, which is

\begin{equation}
g\left(\sum_{i\in I_n} p_{i,n}x_i\right) \leq \sum_{i\in I_n} p_{i,n}g(x_i). \label{proba01_param_conv01}
\end{equation}

\Bin
\noindent Here we apply the dominated convergence  \index{dominated convergence} theorem in Section \ref{proba01_appendix_tcmtcd} in the Appendix Section \ref{proba01_appendix}. We have

$$
\left\vert \sum_{i\geq 1} p_{i}x_i 1_{(x_i \in J_n)}\right\vert \leq \sum_{i\geq 1} p_{i}|x_i| <+\infty
$$

\Bin
\noindent and each $p_{i}x_i 1_{(x_i \in J_n)}$ converges to $p_{i}x_i$ as $n\rightarrow +\infty$. Then, by the dominated convergence  \index{dominated convergence} theorem for series, we have

\begin{equation}
\sum_{i\geq 1} p_{i}x_i 1_{(x_i \in J_n)} \rightarrow \sum_{i\geq 1} p_{i}x_i. \label{proba01_param_conv02}
\end{equation}

\Bin The same argument based on the finiteness of $\sum_{i\geq 1} p_{i}|x_i|$ implies that

\begin{equation}
\sum_{i\geq 1} p_{i} g(x_i) 1_{(x_i \in J_n)} \rightarrow \sum_{i\geq 1} p_{i} g(x_i). \label{proba01_param_conv03}
\end{equation}

\Bin
\noindent By referring to Section \ref{proba01_appendix_convex} in the Appendix Chapter \ref{proba01_appendix}, Part B, we know that a convex function  \index{convex function} $g$ is continuous. We may also remark that $p(n)$ converges to the unity as $n\rightarrow +\infty$. By combining these facts with Formula \eqref{proba01_param_conv01}, we have

\begin{equation}
\sum_{i\geq 1} p_{i}x_i 1_{(x_i \in J_n)}+(1-p(n)) p_{i}x_i 1_{(x_i \in J_n)} \rightarrow \sum_{i\geq 1} p_{i}x_i. \label{proba01_param_conv04}
\end{equation}

\noindent and

$$
g\left(\sum_{i\geq 1} p_{i}x_i\right)=\lim_{n \rightarrow +\infty} g\left(\sum_{i\geq 1} p_{i}x_i 1_{(x_i \in J_n}+(1-p(n)) \sum_{i\geq 1} p_{i}x_i 1_{(x_i \in J_n)}\right)
$$

\Bin
\noindent which is

$$
g\left(\sum_{i\geq 1} p_{i}x_i\right)=\lim_{n \rightarrow +\infty} g\left(p(n) \sum_{i\geq 1} p_{i,n}x_i 1_{(x_i \in J_n}+(1-p(n)) \sum_{i\geq 1} p_{i}x_i 1_{(x_i \in J_n)}\right).
$$

\Bin
\noindent By convexity of $g$, we have

$$
g\left(\sum_{i\geq 1} p_{i}x_i\right)=\lim_{n \rightarrow +\infty} \left(p(n) g\left(\sum_{i\geq 1} p_{i,n}x_i 1_{(x_i \in J_n)}\right)+(1-p(n)) \
g\left( \sum_{i\geq 1} p_{i}x_i 1_{(x_i \in J_n)}\right)\right).
$$

\Bin
\noindent We may write it as

$$
g\left(\sum_{i\geq 1} p_{i}x_i\right)=\lim_{n \rightarrow +\infty} \biggr\{ p(n) g\left(\sum_{i\in I_n} p_{i,n}x_i \right)+(1-p(n))\biggr\}.
g\left( \sum_{i\in I_n} p_{i}x_i \right).
$$

\Bin
\noindent By using Formula \label{proba01_param_conv01} in the first term of the left-hand member of the latter inequality, we have

$$
g\left(\sum_{i\geq 1} p_{i}x_i\right)=\lim_{n \rightarrow +\infty} \left\{p(n) \sum_{i\in I_n} p_{i,n}g(x_i)+(1-p(n))
g\left(\sum_{i\in I_n} p_{i}x_i 1_{(x_i \in J_n)}\right)\right\}.
$$

\Bin
\noindent This gives

$$
g\left(\sum_{i\geq 1} p_{i}x_i\right)=\lim_{n \rightarrow +\infty} \sum_{i\in I_n} p_{i}g(x_i)+(1-p(n))
g\left(\sum_{i\in I_n} p_{i}x_i 1_{(x_i \in J_n)}\right).
$$

\noindent Finally, we have

$$
g\left(\sum_{i\geq 1} p_{i}x_i\right) \leq \sum_{i\geq 1} p_{i}g(x_i)+ 0 \times g\left(\sum_{i\geq 1} p_{i}x_i\right).
$$

\noindent which is exactly

$$
g(\mathbb{E}(X))\leq \mathbb{E}(g(X)).
$$

\newpage

\noindent We are going to find the most common parameters of usual real-valued random variables. It is recommended to know these results by heart.\\

\noindent In each example, we will apply Formula (\ref{FIGD}) and (\ref{FIGDDIM2}) to perform the computations.\\

\noindent  \textbf{(a) Bernoulli Random variable  \index{Bernoulli random variable} : $ X \sim \mathcal{B}(p)$}.\\

\noindent We have

\begin{equation*}
E(X)=1~\mathbb{P}(X=1)+0~\mathbb{P}(X=0)=p.
\end{equation*}

\begin{equation*}
\mathbb{E}(X^{2})=1^{2} \ \mathbb{P}(X=1)+0^2 \ \mathbb{P}(X=0)=p.
\end{equation*}

\Bin \noindent We conclude :

\begin{equation*}
Var(X)=p(1-p)=pq.
\end{equation*}

\bigskip
\noindent \textbf{(b) Binomial Random variable  \index{binomial random variable} : $Y_{n}\sim \mathcal{B}(n, p)$}.\\

\noindent We have, by Lemma \ref{proba01.rv.pd.sumBern} of Chapter \ref{proba01.rv},
\begin{equation*}
Y_{n}=\sum_{i=1}^{n}X_{i}
\end{equation*}

\Bin
\noindent where $X_{1},...,X_{n}$ are Bernoulli $\mathcal{B}$(p) random variables. We are going to use the known parameters of the Bernoulli random variables.  \index{Bernoulli random variable} By linearity property (P3), its comes that

\begin{equation*}
\mathbb{E}(Y_{n})=\sum_{i=1}^{n}\mathbb{E}(X_{i})=np
\end{equation*}

\noindent and by variance-covariance  \index{variance} properties in \textbf{Point (3)} above, the independence  \index{independence} of the $X_i$'s allows to conclude that

\begin{equation*}
Var(Y_{n})=\sum_{i=1}^{n}Var(X_{i})=npq.
\end{equation*}

\bigskip

\noindent  {\bf (c) Geometric Random Variable  \index{geometric random variable} : $X \sim \mathcal{G}(p)$, $p \in ]0,1[$}.\\

\noindent Remind that

\begin{equation*}
\mathbb{P}(X=k)=p(1-p)^{k-1},k=1,2....
\end{equation*}

\noindent Hence

\begin{eqnarray*}
\mathbb{E}(X)&=&\sum_{k=1}^{\infty }k~p(1-p)^{k-1}\\
&=&p\left( \sum_{k=1}^{\infty}kq^{k-1}\right)\\
&=&p\left( \sum_{k=0}^{\infty }q^{k}\right) ^{\prime} \text{ (use the primitive here) }\\
&=&p\left( \frac{1}{1-q}\right) ^{\prime }=\frac{p}{(1-q)^{2}}\\
&=&\frac{p}{p^{2}}=\frac{1}{p}.
\end{eqnarray*}

\Bin
\noindent Its factorial moment  \index{factorial moment} of second order is

\begin{eqnarray*}
\mathbb{E}(X(X-1))&=&p\sum_{k=1}^{\infty}k(k-1)q^{k-1}=pq\sum_{k=2}^{\infty}k(k-1)q^{k-2}\\
&=&pq\left( \sum_{k=0}^{\infty }q^{k}\right) ^{\prime \prime} \text{ (Use the primitive twice )}\\
&=&pq\left( \frac{1}{(1-q)^{2}}\right) ^{\prime}\\
&=&2pq(1-q)^{-3}=\frac{2pq}{p^{3}}=\frac{2q}{p^{2}}.
\end{eqnarray*}

\Bin
\noindent Thus, we get

\begin{equation*}
\mathbb{E}(X^{2} )=\mathbb{E}( X(X-1) )+\mathbb{E}(X)=\frac{2q}{p2}+\frac{1}{p},
\end{equation*}

\noindent and finally we arrive at

\begin{equation}
Var(X)=\frac{2q}{p^{2}} + \frac{1}{p} - \frac{1}{p^{2}} = \frac{q}{p^{2}}.
\end{equation}

\Bin
\noindent In summary :

\begin{equation*}
\mathbb{E}(X)=\frac{1}{p}
\end{equation*}

\noindent and

\begin{equation}
Var(X)=\frac{q}{p^{2}}.
\end{equation}

\bigskip
\noindent \textbf{(d) Negative Binomial Random Variable  \index{negative binomial random variable}  \index{binomial random variable} : $Y_{k} \sim \mathcal{NB}(k, p)$}.\\

\noindent By Lemma \ref{proba01.rv.pd.sumGeo} of Chapter \ref{proba01.rv}, $Y_{k}$ is the sum of $k$ independent geometric $\mathcal{G}(p)$ random variables. Like in the Binomial case, we apply the linearity property (P6) and the properties of the variance-covariance  \index{variance} in Point (3) above to get\\

\begin{equation*}
\mathbb{E}(X)=\frac{k}{p}
\end{equation*}

\noindent and

\begin{equation*}
Var(X)=\frac{kq}{p^{2}}.
\end{equation*}

\bigskip \noindent \textbf{(e) Hypergeometric Random Variable  \index{geometric random variable} \index{hypergeometric random variable} }

$$
Y_{r}\sim H(N,r,\theta), \ \ M=N\theta.\\$$

\noindent We remind that for $0\leq k\leq \min(r,M)$

\begin{equation*}
\mathbb{P}(Y_r=k)=\frac{\left(\begin{tabular}{c} $M$ \\ $k$ \end{tabular} \right) \times \left(\begin{tabular}{c} $N-M$ \\ $r-k$ \end{tabular} \right)}
{\left(\begin{tabular}{c} $N$ \\ $r$ \end{tabular} \right)}.
\end{equation*}

\Bin \noindent To make the computations simple, we assume that $r<M$. We have

\begin{eqnarray*}
\mathbb{E}(Y_{r})&=&\sum_{k=0}^{r}k\ \frac{\text{ \ \ \ \ }M!\text{ \ }(N-M)!%
\text{ \ \ \ \ }r!(N-r)!\text{ \ \ \ }}{k!(M-k)!(r-k)!(N-M-(r-k))!\text{ \ \
}N!}\\
&=&\sum_{k=1}^{r}\ \ \frac{rM\text{ \ }(M-1)!\text{ \ }\left[ (N-1)-(M-1)%
\right] !\text{ \ }}{N\text{ \ \ }(k-1)!\text{ \ }\left[ (M-1)-(k-1)\right] !%
\text{ \ }\left[ (r-1)-(k-1)\right] !\text{ \ }!(N-1)!}\\
&\times& \frac{\text{\ }(r-1)!\text{ \ }\left[ (N-1)-(r-1)\right] !}{\left\{ %
\left[ (N-1)-(M-1)\right] -\left[ (r-1)-(k-1)\right] \right\}! }
\end{eqnarray*}

\Bin
\noindent Let us make the following changes of variables  $k^{\prime }=k-1,$ $M^{\prime}=M-1,$ $r^{\prime }=r-1,$ $N^{\prime }=N-1$. We get

$$
\mathbb{E}(Y_{r})=\frac{rM}{N}\left[ \sum_{k^{\prime }=0}^{r^{\prime }}\frac{%
C_{M^{\prime }}^{k^{\prime }}C_{N^{\prime }-M^{\prime }}^{r^{\prime
}-k^{\prime }}}{C_{N^{\prime }}^{r^{\prime }}}\right] =\frac{rM}{N}=r\theta,
$$

\Bin
\noindent since the term in the bracket is the sum of probabilities of a Hypergeometric  random variable $\mathcal{H}(N-1,r-1,\theta ^{\prime})$, with $\theta ^{\prime }=(M-1)/(N-1)$, then it is equal to one.\\

\noindent The factorial moment  \index{factorial moment} of second order is

\begin{eqnarray*}
&&\mathbb{E}(Y_{r}(Y_{r}-1))=\\
&&\sum_{k=2}^{r}k(k-1)\ \frac{\text{ \ \ \ \ }M!%
\text{ \ \ \ \ }(N-M)!\text{ \ \ \ \ }r!(N-r)!\text{ \ \ \ }}{k!\text{ \ }%
(M-k)!\text{ \ }(N-M-(r-k))!\text{ \ }(r-k)!\text{ \ \ }N!}.
\end{eqnarray*}

\bigskip \noindent Then we make the changes of variables $k^{\prime \prime }=k-2,$ $M^{\prime \prime }=M-2,$ $N^{\prime \prime }=N-2,$ $r^{\prime \prime }=r-2$, to get

$$
\mathbb{E}(Y_{r}(Y_{r}-1))=\frac{Mr(M-1)(r-1)}{N(N-1)}\left[ \sum_{k^{\prime
\prime }=0}^{r^{\prime \prime }}\frac{C_{M^{\prime \prime }}^{k^{\prime
\prime }}C_{N^{\prime \prime }-M^{\prime \prime }}^{r^{\prime \prime
}-k^{\prime \prime }}}{C_{N^{\prime \prime }}^{r^{\prime \prime }}}\right].
$$

\Bin
\noindent Since the term in the bracket is the sum of probabilities of a Hypergeometric  random variable $\mathcal{H}(N-2,r-2,$ $\theta ^{\prime \prime }$),
with  $\theta ^{\prime \prime }=(M-2)/(N-2)$, and then is equal to one.\\

\noindent We get

$$
\mathbb{E}(Y_{r}(Y_{r}-1))=\frac{M(M-1)~r(r-1)}{N(N-1)},
$$

\noindent and finally,

\begin{eqnarray*}
\mathbb{V}ar(Y_{r})&=&\frac{M(M-1)\text{ }r(r-1)}{N(N-1)}+\frac{rM}{N}-\frac{r^{2}M^{2}}{N^{2}}\\
&=&r\frac{M}{N}\left(1+\frac{(M-1)(r-1)}{N-1}-\frac{rM}{N}\right)\\
&=&r\frac{M}{N}\left(\frac{N-M}{N}\right)\left(\frac{N-r}{N-1}\right).
\end{eqnarray*}

\noindent In summary, we have

$$
\mathbb{E}(Y_{r})=r\theta
$$

\noindent and

$$
\mathbb{V}ar (Y_{r})=r\theta (1-\theta)(1-f) (1-1/N)^{-1},
$$

\noindent where

$$
f=\frac{r}{N}
$$

\Bin is the sample ratio.

\bigskip \noindent  \textbf{(f) Poisson Random Variable  \index{Poisson random variable} : $X \sim \mathcal{P}(\lambda)$, $\lambda>0$}.\\

\noindent We remind that
\begin{equation*}
\mathbb{P}(X=k)=\frac{\lambda ^{k}}{k!}e^{-\lambda}
\end{equation*}

\noindent for $k=0,1,...$. Then

\begin{eqnarray*}
\mathbb{E}(X)&=&\sum_{k=0}^{\infty }k\times \frac{\lambda ^{k}}{k!}e^{-\lambda}\\
&=&\sum_{k=1}^{\infty }\frac{\lambda ^{k}}{(k-1)!}e^{-\lambda }\\
&=&\lambda \sum_{k=1}^{\infty }\frac{\lambda ^{k-1}}{(k-1)!}e^{-\lambda}.
\end{eqnarray*}

\bigskip \noindent By setting $k^{\prime }=k-1$, we have
\begin{eqnarray*}
\mathbb{E}(X)&=&\lambda e^{-\lambda }\sum_{k^{\prime }=0}^{\infty }\frac{\lambda ^{k^{\prime }}}{k^{\prime }!}\\
&=&\lambda e^{-\lambda }e^{\lambda}\\
&=&\lambda .
\end{eqnarray*}

\bigskip \noindent The factorial moment  \index{factorial moment} of second order is

\begin{eqnarray*}
\mathbb{E}(X(X-1))&=&\sum_{k=0}^{\infty }k(k-1)\times \frac{\lambda ^{k}}{k!}e^{-\lambda}\\
&=&\sum_{k=2}^{\infty }\frac{\lambda ^{k}}{(k-2)!}e^{-\lambda}\\
&=&\lambda ^{2}e^{-\lambda }\sum_{k=1}^{\infty }\frac{\lambda ^{k-2}}{(k-2)!}
\end{eqnarray*}

\bigskip \noindent By putting $k^{\prime \prime }=k-2$, we get

\begin{eqnarray*}
\mathbb{E}(X(X-1))&=&\lambda ^{2}e^{-\lambda }\sum_{k^{\prime \prime}=0}^{\infty }\frac{\lambda ^{k^{\prime }}}{k^{\prime}!}\\
&=&\lambda ^{2}.
\end{eqnarray*}

\bigskip \noindent Hence,

\begin{equation*}
\mathbb{E}(X^{2})=\mathbb{E}(X(X-1))+\mathbb{E}(X)=\lambda ^{2}+\lambda
\end{equation*}

\bigskip \noindent and then we have
\begin{equation*}
Var(X)=\lambda.
\end{equation*}

\bigskip \noindent \textbf{Remark}. For a Poisson random variable,  \index{Poisson random variable} the mathematical expectation  \index{mathematical expectation} and the variance  \index{variance} are both equal to the intensity $\lambda$.

\newpage

\noindent Let us denote the mean  \index{mean} of a random variable $X$ by $m$. The mean is meant to be a central parameter of $X$ and it is expected that the values of $X$ are near the mean $m$. This is not always true. Actually, the standard deviation  \index{Standard deviation} measures how the values are near the mean.\\

\noindent We already saw in Property (P3) that a null standard deviation  \index{Standard deviation} means  \index{mean} that there is no deviation from the mean and that $X$ is exactly equal to the mean.\\

\noindent Now, in a general case, Tchebychev's inequality  \index{Tchebychev's inequality} may be used to measure the deviation of a random variable $X$ from its mean.  \index{mean} Here is that inequality :

\bigskip

\begin{proposition} \label{prop.proba01.paramRealRv.01} (\textbf{Tchebychev Inequality}). if $\sigma_{X}^{2}=Var(X)>0$, then for any
$\lambda>0$, we have

$$
\mathbb{P}(\left| X-\mathbb{E}(X)\right| \text{ }>\text{ }\lambda \text{ }\sigma_X)\leq \frac{1}{\lambda ^{2}}.
$$

\end{proposition}

\bigskip \noindent \textbf{Commentary}. This inequality says that the probability that $X$ deviates from its mean  \index{mean} by at least $\lambda \sigma_{X}$ is less than one over $\lambda$.\\

\noindent In general terms, we say this : that $X$ deviates from its mean  \index{mean} by at least a multiple of $\sigma_{X}$ is  all the more unlikely that multiple is large. More precisely, the probability that $X$ deviates from its mean by at least a multiple of $\sigma_{X}$ is bounded by the inverse of the square of that multiple. This conveys the fact that the standard deviation  \index{Standard deviation} controls how the observations are near the mean.\\

\bigskip \noindent \textbf{Proof}. Let us begin to prove the Markov's inequality.  \index{Markov's inequality} Let $Y$ be a non-negative random variable
with values $\left\{ y_{j},\text{ }j\in J\right\}$. Then, for any $\lambda>0$, the following Markov Inequality

\begin{equation} \label{proba01.paramRealRv.MarkovIneq}
\mathbb{P}(Y\text{ }> \lambda \\mathbb{E}(Y))\text{ }<\text{ }\frac{1}{\lambda},
\end{equation}

\noindent holds. Before we give the proof of this inequality, we make a remark. By the formula in Theorem \ref{proba01.rv.theo1} on Chapter
\ref{proba01.rv}, we have for any number $y$,

\begin{equation}
\mathbb{P}(Y \leq y)=\sum_{y_{i}\leq y} \mathbb{P}(Y=y_{i}). \label{proba01.paramRealRv.A01}
\end{equation}

\Bin and

\begin{equation}
\mathbb{P}(Y > y)=\sum_{y_{i} > y} \mathbb{P}(Y=y_{i}). \label{proba01.paramRealRv.A01S}
\end{equation}

\Bin
\noindent Let us give now the proof of the inequality (\ref{proba01.paramRealRv.MarkovIneq}). Since, by the assumptions, $y_{j}\geq 0$ for all $j \in J$, we have

\begin{eqnarray*}
\mathbb{E}(Y)&=&\sum_{j\in J}y_{j}\mathbb{P}(Y=y_{j})\\
&=&\sum_{y_{j}\leq \lambda \mathbb{E}(Y)} y_{j}\mathbb{P}(Y=y_{i})+\sum_{y_{j}>\lambda \mathbb{E}(Y)}
y_{j}\mathbb{P}(Y=y_{j}).
\end{eqnarray*}

\Bin
\noindent Then by getting rid of the first term and use Formula (\ref{proba01.paramRealRv.A01}) to get

\begin{eqnarray*}
\mathbb{E}(Y) &\geq& \sum_{y_{j}>\lambda \mathbb{E}(Y)}y_{j}\mathbb{P}(Y=y_{j})\\
&\geq& (\lambda \mathbb{E}(Y)) \times \sum_{y_{j}>\lambda \mathbb{E}(Y)}\mathbb{P}(Y=y_{j})\\
&=&(\lambda \mathbb{E}(Y)) \times \mathbb{P}(Y>\lambda \mathbb{E}(Y)).
\end{eqnarray*}

\Bin
\noindent By comparing the two extreme members of this latter formula, we get the Markov inequality in Formula
(\ref{proba01.paramRealRv.MarkovIneq}).\\

\noindent Now, let us apply the Markov's inequality  \index{Markov's inequality} to the random variable

\begin{equation*}
Y=(X-\mathbb{E}(X))^{2}.
\end{equation*}

\noindent By definition, we have  $\mathbb{E}(Y)=Var(X)=\sigma ^{2}$. The application of Markov's inequality  \index{Markov's inequality} to $Y$ leads to

\begin{eqnarray*}
\mathbb{P}(\left| X-\mathbb{E}(X)\right| \text{ }>\text{ }\lambda\sigma)&=&\mathbb{P}(\left| X-\mathbb{E}(X)\right| ^{2}\text{ }>\lambda ^{2}\sigma ^{2})\\
&\leq& \frac{1}{\lambda ^{2}}.
\end{eqnarray*}

\Bin This puts an end to the proof.\\

\bigskip

\noindent \textbf{Application of the Tchebychev's Inequality  \index{Tchebychev's inequality} : Confidence intervals  \index{confidence interval} of $X$ around the mean \index{mean} }.\\

\noindent The Chebychev's inequality allows to describe the deviation of the
random variable X from mean  \index{mean} $m_{X}$ in the following probability covering

\begin{equation}
\mathbb{P}\left( m_{X}-\Delta \leq X\leq m_X+\Delta \right) \geq 1-\alpha .
\label{proba01.paramRealRv.ici.01}
\end{equation}

\Bin
\noindent This relation says that we are $100(1-\alpha )\%$ \textbf{confident} that $X$ lies in the interval \ $I(m_X)=[m_X-\Delta, \ m_X+\Delta]$. If Formula (\ref{proba01.paramRealRv.ici.01}) holds, $I(m_{X})=[m_{X}-\Delta ,m_{X}+\Delta ]$ is called a $(1-\alpha )$-confidence interval  \index{confidence interval} of $X$.\\

\noindent By using the Chebychev's inequality, $I(m_{X},\alpha )=[m_{X}-\sigma _{X}/%
\sqrt{\alpha },m_{X}+\sigma _{X}/\sqrt{\alpha }]$ is a $(1-\alpha )$%
-confidence interval  \index{confidence interval} of $X,$ that is

\begin{equation}
\mathbb{P}\left( m_{X}-\text{ }\frac{\sigma _{X}}{\sqrt{\alpha }}\leq X\leq
m_{X}+\text{ }\frac{\sigma _{X}}{\sqrt{\alpha }}\right) \geq 1-\alpha \text{.%
}  \label{proba01.paramRealRv.ici.02}
\end{equation}

\noindent Let us make two remarks.\\

\noindent \textbf{(1)} The confidence intervals  \index{confidence interval} are generally used for small values of $\alpha ,
$ and most commonly for $\alpha =5\%.$ And we have a $95\%$ confidence
interval of $X$ in the form (since $1/\sqrt{\alpha }=4.47)$

\begin{equation}
\mathbb{P}\left( m_{X}-4.47\text{ }\sigma _{X}\leq X\leq m_{X}+4.47\sigma
_{X}\right) \geq 95\%\text{.}
\end{equation}

\bigskip \noindent $95\%$-confidence intervals  \index{confidence interval} are used for comparing different groups with
respect to the same character $X$.\\

\bigskip \noindent \textbf{(2)} Denote by
\begin{equation*}
RVC_{X}=\left\vert \frac{\sigma }{m}\right\vert .
\end{equation*}

\bigskip \noindent Formula \eqref{proba01.paramRealRv.ici.02} gives for $0<\alpha <1,$

\begin{equation*}
\mathbb{P}\left( 1 -  \ \frac{RVC_{X}}{\sqrt{\alpha }}\leq \left\vert
\frac{X}{m}\right\vert \leq 1+\text{ }\frac{RVC_{X}}{\sqrt{\alpha }}\right)
\geq 1-\alpha.
\end{equation*}

\bigskip \noindent The coefficient $RVC_{X}$ is called the Relative Variation Coefficient of $X$ which indicates how small the variable $X$ deviates from its mean \index{mean} .
Generally, it is expected that $RVC_{X}$ is less $30\%$ for a homogeneous variable.\\

\noindent \noindent (3) The confidence intervals  \index{confidence interval} obtained by the Tchebychev's inequality  \index{Tchebychev's inequality} are not  so precise as we may need them. Actually, they are mainly used for theoretical
purposes. In Applied Statistics  \index{statistics} analysis, more sophisticated confidence
intervals are used. For example, normal confidence interval  \index{confidence interval} are
systematically exploited. But for learning purposes,  Tchebychev confidence
intervals are good tools.\\

\bigskip \noindent Let us now prove Formula (\ref{proba01.paramRealRv.ici.02}).\\

\noindent \textbf{Proof of Formula \ref{proba01.paramRealRv.ici.02}}. Fix $0<\alpha <1$ and $\varepsilon >0.$
Let us apply Tchebychev's inequality  \index{Tchebychev's inequality} to get

\begin{equation*}
\mathbb{P}(\left\vert X-\mathbb{E}(X)\right\vert \text{ }>\text{ }%
\varepsilon )=\mathbb{P}(\left\vert X-E(X)\right\vert \text{ }>\text{ }((%
\frac{\varepsilon }{\sigma _{X}})\text{ }\sigma ))\leq \frac{\sigma _{X}^{2}%
}{\varepsilon ^{2}}\text{.}
\end{equation*}

\noindent Let
\begin{equation*}
\alpha =\frac{\sigma _{X}^{2}}{\varepsilon ^{2}}.
\end{equation*}

\noindent We get
\begin{equation*}
\varepsilon =\frac{\sigma _{X}}{\sqrt{\alpha }}
\end{equation*}

\noindent and then
\begin{equation*}
\mathbb{P}\left( \left\vert X-m_{X}\right\vert \text{ }>\text{ }\frac{\sigma
_{X}}{\sqrt{\alpha }}\right) \leq \alpha \text{,}
\end{equation*}

\noindent which equivalent to

\begin{equation*}
\mathbb{P}\left( \left\vert X-m_{X}\right\vert \text{ }\leq \text{ }\frac{%
\sigma _{X}}{\sqrt{\alpha }}\right) \geq 1-\alpha ,
\end{equation*}

\noindent and this is Formula (\ref{proba01.paramRealRv.ici.01}).\\

 %proba01.param
\chapter{Random pairs}  \index{random pair} \label{proba01.randomcouples}

\noindent
This chapter is devoted to an introduction to the study of probability law \index{probability law} s
of random pairs  \index{random pair} $(X,Y)$, \textit{i.e.} two-dimensional random vectors, their usual parameters and related concepts. Its
follows the lines of Chapters \ref{proba01.rv} and \ref{proba01.param} which focused on real-valued random variables.
As in the aforementioned chapter, here again, we focus on discrete random
pairs.

\bigskip

\noindent A discrete random pair  \index{random pair} $(X,Y)$ takes its values in a set of the form%

\begin{equation*}
\mathcal{SV}_{(X,Y)}=\{(a_{i},b_{j}),(i,j)\in K\},
\end{equation*}

\Bin \noindent where $K$ is an enumerable set with%

\begin{equation}
\forall (i,j)\in K, \ \mathbb{P}\biggr((X,Y)=(a_{i},b_{j})\biggr)>0.
\label{probab01.rcouples.strictDom}
\end{equation}

\bigskip \noindent We say that $\mathcal{SV}_{(X,Y)}$ is a \textbf{strict} support or a domain or a values
set of the random pair  \index{random pair} $(X,Y)$, if and only if, all its elements are taken by $(X,Y)$ with non-zero probabilities, as in (\ref{probab01.rcouples.strictDom}).\\

\noindent We want the reader to remark for once that adding supplementary points $(x,y)$ which are not taken by $(X,Y)$ [meaning  \index{mean} that $\mathbb{P}((X,Y)=(x,y))=0$] to $\mathcal{SV}_{(X,Y)}$, does not change anything regarding computations using the probability law  \index{probability law} of $(X,Y)$. If we have such points in a values set, we call this latter an \textbf{extended values set}.\\

\bigskip \noindent We consider the first projections  \index{projection} of the pair values defined as follows:

\begin{equation*}
(a_{i},b_{j})\hookrightarrow a_{i}.
\end{equation*}

\bigskip
\noindent By forming a set of these projections,  \index{projection} we obtain a set
\begin{equation*}
\mathcal{V}_{X}=\{x_{h},h\in I\}.
\end{equation*}

\noindent For any $i\in I$, $x_{i}$ is equal to one of the projections  \index{projection} $a_{i}$ for which
there exists at least one $b_{j}$ such that $(a_{i},b_{j})\in \mathcal{SV}_{(X,Y)}.$ It
is clear that $\mathcal{V}_{X}$ is the strict values set of $X$.\\

\bigskip \noindent We have to pay attention to the possibility of having some values $x_{i}$ that are repeated while taking all the projections  \index{projection} $(a_{i},b_{j}) \hookrightarrow a_{i}$. Let us give an example.

\bigskip \noindent Let

\begin{equation}
\mathcal{SV}_{(X,Y)}=\{(1,2),(1,3),(2,2),(2,4)\}  \label{proba01.rcouples.02}
\end{equation}

\bigskip \noindent with $\mathbb{P}((X,Y)=(1,2))=0.2$, $\mathbb{P}((X,Y)=(1,3))=0.3$, $\mathbb{P}((X,Y)=(2,2))=0.1$ and $\mathbb{P}((X,Y)=(2,4))=0.4$.\\

\noindent We have
\begin{equation*}
\mathcal{V}_{X}=\{1,2\}.
\end{equation*}

\Bin  \noindent Here the projections  \index{projection} $(a_{i},b_{j})\hookrightarrow a_{i}$ give the value $1$
two times and the value $2$ two times.\\

\noindent
We may and do proceed similarly for the second projection \index{projection} s

\begin{equation*}
(a_{i},b_{j})\hookrightarrow b_{j},
\end{equation*}

\Bin \noindent to define
\begin{equation*}
\mathcal{V}_{Y}=\{y_{j},j \in J\},
\end{equation*}

\noindent as the strict values set of $Y$.\\

\bigskip \noindent We express a strong warning to not think that the strict values set of the
pair is the Cartesian product of the two strict values sets of $X$ and $Y$, that is

\begin{equation*}
\mathcal{SV}_{(X,Y)}=\mathcal{V}_{X}\times \mathcal{V}_{Y}.
\end{equation*}

\Bin \noindent In the example of (\ref{proba01.rcouples.02}), we may check that

\begin{eqnarray*}
\mathcal{V}_{X}\times \mathcal{V}_{Y} &=&\{1,2\}\times \{2,3,4\} \\
&=&\{(1,2),(1,3),(1,4),(2,2),(2,3),(2,4)\} \\
&\neq &\{(1,2),(1,3),(2,2),(2,4)\} \\
&=&\mathcal{SV}_{(X,Y)}.
\end{eqnarray*}

\Bin
\noindent Now, even we have this fact, we may use

\begin{equation*}
\mathcal{V}_{(X,Y)}=\mathcal{V}_{X}\times \mathcal{V}_{Y},
\end{equation*}

\Bin
\noindent as an extended values set and remind ourselves that some of the elements of $%
\mathcal{V}_{(X,Y)}$ may be assigned null probabilities by $\mathbb{P}_{(X,Y)}$.\\

\noindent In our example (\ref{proba01.rcouples.02}), using the extended values sets leads to the probability law  \index{probability law} represented in Table \ref{proba01.rcouples.03}.

\begin{table}[htbp]
\centering
\caption{Simple example of a probability law  \index{probability law} in dimension 2}
\begin{tabular}{|l|l|l|l|l|}
\hline
$X \ \setminus \ Y$ & $2$ & $3$ & $4$ & $X$ \\ \hline
$1$ & $0.2$ & $0.3$ & $0$ & $0.5$ \\ \hline
$2$ & $0$ & $0.1$ & $0.4$ & $0.5$ \\ \hline
$Y$ & $0.2$ & $0.4$ & $0.4$ & $1$ \\ \hline
\end{tabular}
\label{proba01.rcouples.03}
\end{table}

\bigskip \noindent In this table,\\

\noindent (1) We put the probability $\mathbb{P}((X,Y)=(x,y))$ at the intersection of
the value $x$ of $X$ (in columns) and the value $y$ of $Y$ (in lines).\\

\noindent (2) In the column $X$, we put, at each line corresponding to a value $x$ of $%
X$, the sum of all the probability of that line. This column represents the
probability law  \index{probability law} of $X$, as we will see it soon.\\

\noindent (3) In the line $Y$, we put, at each column corresponding to a value $y$ of $Y$, the sum of all the probability of that
column. This line represents the probability law  \index{probability law} of $Y$.\\

\bigskip \noindent \textbf{Conclusion}. In our study of random pairs  \index{random pair} $(X,Y)$, we may always use an
extended values set of the form

\begin{equation*}
\mathcal{V}_{(X,Y)}=\mathcal{V}_{X}\times \mathcal{V}_{Y},
\end{equation*}

\bigskip \noindent where $\mathcal{V}_{X}$ is the strict values set of $X$ and $\mathcal{V}_{Y}$ the
strict values set of $Y$. If $\mathcal{V}_{X}$  and $\mathcal{V}_{Y}$ are finite with
relatively small sizes, we may appeal to tables like Table \ref{proba01.rcouples.03} to represent the probability law  \index{probability law} of $(X,Y)$.

\noindent Let $(X,Y)$ be a discrete random pair  \index{random pair} such that $\mathcal{V}_{X}=\{x_{i}, \ i\in I\}$ is the strict values set of $X$ and $\mathcal{V}%
_{Y}=\{y_{j},j\in J\}$ the strict values set of $Y$, with $I\subset \mathbb{N%
}$ and $J\subset \mathbb{N}.$ The probability law  \index{probability law} of $(X,Y)$ is given on the
extended values set

\begin{equation*}
\mathcal{V}_{(X,Y)}=\mathcal{V}_{X}\times \mathcal{V}_{Y},
\end{equation*}

\noindent by

\begin{equation}
\mathbb{P}(X=x_{i},\text{ }Y=y_{j}),(i,j)\in I\times J. \label{proba01.rcouples.jpl01}
\end{equation}

\bigskip \noindent where the event $(X=x_{i},\text{ }Y=y_{j})$ is a notation of the intersection

$$
(X=x_{i}) \cap (Y=y_{j}).
$$

\Bin \noindent This probability law  \index{probability law} may be summarized in Table \ref{proba01.rcouples.t04}.

\begin{table}[htbp]
\centering
\begin{tabular}{|l|l|l|l|l|l|}
\hline
X $\diagup$  Y & $y_{1}$ & $\cdots$ & $y_{j}$ & $\cdots$ & X \\ \hline
$x_{1}$ & $\mathbb{P}(X=x_1,\text{ }Y=y_1)$ & $\cdots$ & $\mathbb{P%
}(X=x_1,\text{ }Y=y_j)$ & $\cdots$ & $\mathbb{P}(X=x_1)$ \\ \hline
$x_{2}$ & $\mathbb{P}(X=x_2,\text{ }Y=y_1)$ & $\cdots$ & $\mathbb{P%
}(X=x_2,\text{ }Y=y_j)$ & $\cdots$ & $\mathbb{P}(X=x_2)$ \\ \hline
$\cdots$ & $\cdots$ & $\cdots$ & $\cdots$ & $\cdots$ & $\cdots$\\ \hline
$x_{i}$ & $\mathbb{P}(X=x_i,\text{ }Y=y_1)$ & $\cdots$ & $\mathbb{P}%
(X=x_i,\text{ }Y=y_j)$ & $\cdots$ & $\mathbb{P}(X=x_i)$ \\ \hline
$\cdots$ & $\cdots$ & $\cdots$ & $\cdots$ & $\cdots$ & $\cdots$\\ \hline
Y & $\mathbb{P}(Y=y_1)$ & $\cdots$ & $\mathbb{P}(Y=y_j)$ &
$\cdots$ & $100\%$ \\ \hline
\end{tabular}
\caption{Probability Law  \index{probability law} table of a random pair  \index{random pair} $(X,Y)$}
\label{proba01.rcouples.t04}
\end{table}

\bigskip \noindent This table is formed as follows :\\

\noindent (1) We put the probability $\mathbb{P}((X,Y)=(x_i,y_j))$ at the
intersection of the value $x_i$ of $X$ (in line) and the value $y_j$ of $Y$
(in column).\\

\noindent (2) In the column $X,$ we put, at each line the corresponding to a value $x_i$
of $X$, the sum of all the probability of that line. As we will see it later, this
column represents the probability law  \index{probability law} of $X$.\\

\noindent (3) In the line $Y$, we put, at each column corresponding to a value $y$ of $Y$, the sum of all the probability of that
column. This line represents the probability law  \index{probability law} of $Y$.\\

\bigskip  \noindent Let us introduce the following terminology.\\

\bigskip \noindent \textbf{Joint  \index{joint} Probability Law \index{probability law} }.\\

\noindent The joint  \index{joint} probability law  \index{probability law} of $(X,Y)$ is simply the probability law of the pair, which is  given in (\ref{proba01.rcouples.jpl01}).\\

\bigskip \noindent \textbf{Marginal  \index{marginal} probability law \index{probability law} s}.\\

\noindent The probability laws  \index{probability law} of $X$ and $Y$ are called marginal  \index{marginal} probability laws.\\

\noindent Here is how we get the marginal  \index{marginal} probability law  \index{probability law} of $X$ and $Y$ from the joint  \index{joint} probability law. Let us begin by $X$.\\

\noindent  We already knew, from Formula (\ref{proba01.rv.decompEvent02}), that any event $B$ can be decomposed into

\begin{equation*}
B= \sum_{j~\in ~J} B \cap (Y=y_{j}).
\end{equation*}

\Bin
\noindent So, for a fixed $i \in I$, we have

\begin{equation*}
(X=x_i)= \sum_{j~\in ~J} (X=x_i) \cap (Y=y_{j})=\sum_{j~\in ~J} (X=x_i, \ \ Y=y_{j}).
\end{equation*}

\Bin
\noindent By applying the additivity of the probability measure,  \index{probability measure} it comes that

\begin{equation}
\mathbb{P}(X=x_i)= \sum_{j~\in ~J} \mathbb{P}(X=x_i, \ \ Y=y_{j}), \ \ i\in I. \label{proba01.rcouples.mplX}
\end{equation}

\Bin
\noindent By doing the same of $Y$, we get

\begin{equation}
\mathbb{P}(Y=y_j)= \sum_{i~\in ~I} \mathbb{P}(X=x_i, \ \ Y=y_{j}),\ \ j\in J. \label{proba01.rcouples.mplY}
\end{equation}

\noindent We may conclude as follows.\\

\noindent (1) The marginal  \index{marginal} probability law  \index{probability law} of $X$, given in (\ref{proba01.rcouples.mplX}), is obtained by summing the joint  \index{joint} probability law over the values of $Y$.\\

\noindent In Table \ref{proba01.rcouples.t04}, the marginal  \index{marginal} probability law  \index{probability law} of $X$ is represented by the last column, in which each line corresponding to a value $x_i$ of $X$ is the sum of all the probability of that line.\\

\noindent (2) The marginal  \index{marginal} probability law  \index{probability law} of $Y$, given in (\ref{proba01.rcouples.mplY}), is obtained by summing the joint  \index{joint} probability law over the values of $X$.\\

\noindent In Table \ref{proba01.rcouples.t04}, the marginal  \index{marginal} probability law  \index{probability law} of $Y$ is represented by the last line, in which each column corresponding to a value $y_j$ of $Y$ is the sum of all the probability of that column.\\

\noindent (3) The box at the intersection of the last line and the last column contains the unit value, that is, at the same time, the sum of all the joint  \index{joint} probabilities, all marginal  \index{marginal} probabilities of $X$ and all marginal probabilities of $Y$.\\

\noindent Usually, it is convenient to use the following notations :\\

\begin{equation}
p_{i,j}=\mathbb{P}(X=x_{i},\text{ }Y=y_{j}), \ (i,j)\in I\times J. \label{proba01.rcouples.jpl02}
\end{equation}

\begin{equation}
p_{i,\bullet}=\mathbb{P}(X=x_{i}), \ \ i\in I. \label{proba01.rcouples.mplX1}
\end{equation}

\Bin
\noindent and

\begin{equation}
p_{\bullet,j}=\mathbb{P}(Y=y_{j}), \ \ j \in J. \label{proba01.rcouples.mplY1}
\end{equation}

\bigskip
\noindent We have the formulas

\begin{eqnarray}
p_{i,\bullet}&=&\mathbb{P}(X=x_i, Y \in \Omega) \label{proba01.rcouples.mplX2} \\
&=& \sum_{j~\in ~J} p_{i,j}, \ \ i\in I. \notag
\end{eqnarray}

\noindent and

\begin{eqnarray}
p_{\bullet,j}&=&\mathbb{P}(X \in \Omega, Y=y_j) \label{proba01.rcouples.mplY2} \\
&=& \sum_{i~\in ~I} p_{i,j}, \ \ j\in J. \notag
\end{eqnarray}

\bigskip

\noindent We come back to the concept of conditional  \index{conditional} probability and independence  \index{independence} that was introduced in Chapter \ref{proba01.cpi}. This section is developed for real-valued random variables. But we should keep in mind that the theory is still valid for $E$-valued random variables, where $E=\mathbb{R}^k$. When a result hold only for real-valued random variables, we will clearly specify it.\\

\bigskip \noindent \textbf{ I - Independence \index{independence} }.\\

\bigskip \noindent \textbf{(a) Definition}.\\

\noindent We remind that, by definition (see Property (P6) in Chapter \ref{proba01.param}), $X$ and $Y$ are independent if and only if one of the following assertions holds.\\

\noindent \textbf{(CI1)} For all $i\in J$ and $j\in J$,

$$
p_{ij}=p_{i,\bullet} \times p_{\bullet,j}.
$$

\Bin
\noindent \textbf{(CI2)} For any subsets $A$ and $B$ of $E=\mathbb{R}$,

$$
\mathbb{P}(X \in A, Y\in B)= \mathbb{P}(X \in A) \times \mathbb{P}(Y\in B).
$$

\bigskip
\noindent \textbf{(CI3)} For any non negative functions $f,g : \Omega \mapsto E$,

$$
\mathbb{E}(h(X)g(Y))= \mathbb{E}(h(X)) \times \mathbb{E}(g(X))
$$

\bigskip
\noindent Actually, Property (P6) in Chapter \ref{proba01.param}, only establishes the equivalence between (CI1) and (CI3). But clearly, (CI3) implies (CI2) by using the following functions
$$
h=1_{A} \text{ and } g=1_{B}.
$$

\bigskip
\noindent Also (CI2) implies (CI1) by using

$$
A=\{x_i\} \text{ and } B=\{y_j\}.
$$

\bigskip
\noindent We get the equivalence between the three assertions by the circular argument : \\

$$
(CI1) \Rightarrow (CI3) \Rightarrow (CI2) \Rightarrow (CI1).
$$

\bigskip \bigskip \noindent \textbf{(b) New tools for assessing the independence \index{independence} }.\\

\noindent \textbf{For a real-valued random variable} $X$, we may define the two following functions depending on the law of $X$.\\

\noindent \textbf{(MGF1)} The first moment generating function  \index{moment generating function} of $X$ :

$$
\Phi_{X}(s)=\mathbb{E}(\exp(sX))=\sum_{i \in I} \mathbb{P}(X=x_i) \exp(sx_i), \ \ s\in \mathbb{R}.
$$

\Bin
\noindent \textbf{(MGF2)} The second moment generating function  \index{moment generating function} of $X$ :

$$
\Psi_{X}(s)=\mathbb{E}(s^{X})=\sum_{i \in I} \mathbb{P}(X=x_i) s^{x_i}, s\in ]0,1].
$$

\bigskip \noindent It is clear that we have the following relation between these two functions:

\begin{equation}
\Psi_{X}(s)=\Phi_{X}(\log s) \ \  and \ \ \Psi_{X}(1)=\Phi_{X}(0). \label{proba01.rcouples.MGF12}
\end{equation}

\bigskip \noindent So, we do not need to study them separately. It is enough to study one of them and to transfer the obtained properties to the other. Since, the first is more popular, we will study it. Consequently, the first form is simply called the moment generating function  \index{moment generating function} (m.g.f). When we use the form given in (MGF2), we will call it by its full name as the \textit{second m.g.f}.\\

\noindent \textbf{Characterization of probability laws}.  \index{probability law} We admit that the \textit{m.g.f} of $X$, when it is defined, is characteristic of the probability law of $X$, meaning  \index{mean} that two real-valued random variables having the same \textit{m.g.f} have the same law. The same characterization is valid for the second \textit{m.g.t}, because of the relation between the two \textit{m.g.f} forms through Formula (\ref{proba01.rcouples.MGF12}). The proof of this result, which is beyond the level of this book, will be presented in the monograph of \textit{Mathematical Foundations of Probability Theory}  \index{probability theory} of this series.\\

\noindent Each form of \textit{m.g.f} has its own merits and properties. We have to use them in smart ways depending on the context.\\

\noindent The \textit{m.g.f} for real-valued random variables may be extended to a pair of random variables. If $X$ and $Y$ are both real applications, we may define the bi-dimensional \textit{m.g.f} of the pair $(X,Y)$ by

\begin{equation}
\Phi_{(X,Y)}(s,t)=\mathbb{E}( \exp(sX+tY))=\sum_{(i,j) \in I \times J} \mathbb{P}(X=x_i,Y=y_i) \exp(sx_i+ty_j), \tag{MGFC}
\end{equation}

\noindent for $(s,t)\in \mathbb{R}^2$.\\

\noindent The \textit{m.g.f} has a nice and simple affine transformation formula. Indeed, if $a$ and $b$ are real numbers, we have for any $s \in \mathbb{R}$

$$
\exp((aX+b)s))=e^{bs} \exp((as)X)
$$

\Bin
\noindent and then, we have

\begin{eqnarray*}
\Phi_{Z}(s)&=&\mathbb{E}(\exp((aX+b)s)\\
&=& \mathbb{E}(e^{bs} \exp((as)X)\\
&=& e^{bs}  \mathbb{E}(\exp((as)X)\\
&=& e^{bs}  \Phi_{X}(as).
\end{eqnarray*}

\bigskip \noindent We obtained this remarkable formula : for any real numbers $a$ and $b$, for any $s \in \mathbb{R}$, we have

\begin{equation}
\Phi_{aX+b}(s)= e^{bs}  \Phi_{X}(as). \label{proba01.rcouples.MGFTransform}
\end{equation}

\bigskip \noindent Formula (\ref{proba01.rcouples.MGFTransform}) is a useful tool in finding news probability laws  \index{probability law} for known ones.\\

\bigskip \noindent We are going to see two other interesting properties of these functions.\\

\bigskip \noindent \textbf{A - Moment generating Functions  \index{moment generating function} and Independence \index{independence} }.\\

\noindent If $X$ and $Y$ are real-valued independent random variables, then we have : \\

\bigskip
\noindent For any $s \in \mathbb{R}$,

\begin{equation}
\Phi_{X+Y}(s)=\Psi_{X}(s) \times \Phi_{Y}(s). \label{proba01.rcouples.MGFFACT1}
\end{equation}

\Bin
\noindent For any $(s,t)\in \mathbb{R}^2$,

\begin{equation}
\Phi_{(X,Y)}(s,t)=\Psi_{X}(s) \times \Phi_{Y}(t). \label{proba01.rcouples.MGFFACT2}
\end{equation}

\bigskip \noindent Before we give the proofs, let us enrich our list of independent conditions.\\

\noindent \textbf{(CI4)} \textbf{Two real-valued random variables $X$ and $Y$ are independent if and only Formula (\ref{proba01.rcouples.MGFFACT2}) holds for any $(s,t)\in \mathbb{R}^2$}.\\

\noindent This assertion also is beyond the current level. We admit it. We will prove it in the book of \textit{Mathematical Foundations of probability Theory}  \index{probability theory} in this series.\\

\noindent \textbf{Warning}. Formula (\ref{proba01.rcouples.MGFFACT2}) characterizes the independence  \index{independence} between two random variables, but not (\ref{proba01.rcouples.MGFFACT1}). Let us show it with this counter-example in \cite{stoyanov} \index{Stayonov}, page 62, 2nd Edition.\\

\noindent Consider the random pair  \index{random pair} $(X,Y)$ taking with domain $\mathcal{V}_{(X,Y)}=\{1,2,3\}^2$, and whose probability law  \index{probability law} is given by

\begin{table}[htbp]
\centering
\begin{tabular}{|l|l|l|l|l|}
\hline
$X / Y$ & $1$ & $2$ & $3$ &  $X$ \\ \hline
$1$ & $\frac{2}{18}$ & $\frac{1}{18}$ & $\frac{3}{18}$ & $\frac{6}{18}$ \\ \hline
$2$ & $\frac{3}{18}$ & $\frac{2}{18}$ & $\frac{1}{18}$ & $\frac{2}{18}$ \\ \hline
$3$ & $\frac{1}{18}$ & $\frac{3}{18}$ & $\frac{2}{18}$ & $\frac{6}{18}$ \\ \hline
$Y$ & $\frac{6}{18}$ & $\frac{6}{18}$ & $\frac{6}{18}$ & $\frac{18}{18}$ \\ \hline
\end{tabular}
\vskip 0.5cm
\caption{Counter-example : (\ref{proba01.rcouples.MGFFACT1}) does not imply independence \index{independence} }
\label{proba01.rcouples.06}

\end{table}

\noindent We may check that this table gives a probability laws.  \index{probability law} Besides, we have :\\

\noindent (1) $X$ and $Y$ have a uniform distribution on $\{1,2,3\}$ with

$$
\mathbb{P}(X=1)=\mathbb{P}(X=2)=\mathbb{P}(X=3)=\frac{6}{18}=\frac{1}{3}
$$

\noindent and

$$
\mathbb{P}(Y=1)=\mathbb{P}(Y=2)=\mathbb{P}(Y=3)=\frac{6}{18}=\frac{1}{3}.
$$

\Bin
\noindent (2) $X$ and $Y$ have the common \textit{m.g.f}

$$
\Phi_{X}(s)=\Phi_{Y}(s)=\frac{1}{3}(e^{s}+e^{2s}+e^{3s}),\ \ s\in \mathbb{R}.
$$

\Bin
\noindent (3) We have, for $s\in \mathbb{R}$,

\begin{eqnarray*}
\Phi_X(s) \times \Phi_Y(s) &=&\frac{1}{9}(e^{2s}+2e^{3s}+3e^{4s}+2e^{5s}+e^{6s})\\
&=&\frac{1}{18}(2e^{2s}+4e^{3s}+6e^{4s}+4e^{5s}+2e^{6s})
\end{eqnarray*}

\noindent (4) The random variable $Z=X+Y$ takes its  values in $\{2,3,4,5\}$. The events $(X=k)$, $2\leq k \leq 6$ may be expressed  with respect to the values of $(X,Y)$ as follows

\begin{eqnarray*}
(Z=2)&=&(X=1,Y=1)\\
(Z=3)&=&(X=1,Y=2)+(X=2,Y=1)\\
(Z=4)&=&(X=1,Y=3)+(X=2,Y=2)+(X=3,Y=1)\\
(Z=5)&=&(X=2,Y=3)+(X=3,Y=2)\\
(Z=6)&=&(X=3,Y=3)
\end{eqnarray*}

\Bin
\noindent By combining this and the joint  \index{joint} probability law  \index{probability law} of $(X,Y)$ given in Table \ref{proba01.rcouples.06}, we have

\begin{eqnarray*}
\mathbb{P}(Z=2)&=&\frac{2}{18}\\
\mathbb{P}(Z=3)&=&\frac{4}{18}\\
\mathbb{P}(Z=4)&=&\frac{6}{18}\\
\mathbb{P}(Z=5)&=&\frac{4}{18}\\
\mathbb{P}(Z=6)&=&\frac{2}{18}\\
\end{eqnarray*}

\noindent (5) The \textit{m.g.f} of $Z=X+Y$ is

\begin{equation}
\Phi_{X+Y}(s)=\frac{1}{18}(2e^{2s}+4e^{3s}+6e^{4s}+4e^{5s}+2e^{6s})
\end{equation}

\noindent (6) By combining Points (3) and (5), we see that Formula (\ref{proba01.rcouples.MGFFACT1}) holds.\\

\noindent (7) Yet, we do not have independence,  \index{independence} since, for example

$$
\mathbb{P}(X=2,Y=1)=\frac{3}{18} \ \ \neq \ \ \frac{2}{18}=\mathbb{P}(X=2)\times \mathbb{P}(Y=1).
$$

\noindent $\square$\\

\bigskip \noindent \textbf{Proofs of Formulas (\ref{proba01.rcouples.MGFFACT1}) and (\ref{proba01.rcouples.MGFFACT2})}.\\

\noindent We simply use Property (CI3) below to get for any $s \in \mathbb{R}$,

\begin{eqnarray*}
\Phi_{X+Y}(s)&=&\mathbb{E}(\exp(sX+sY)) \\
&=&\mathbb{E}(\exp(sX) \times \exp(sY)) \ \ \text{(Property of the exponential)}\\
&=&\mathbb{E}(\exp(sX)) \times \mathbb{E}(\exp(sY)) \ \ \text{(By Property (CI3))}\\
&=& \Phi_{X}(s) \times \Phi_{Y}(s)
\end{eqnarray*}

\bigskip \noindent and for any $(s,t)\in \mathbb{R}^2$,

\begin{eqnarray*}
\Phi_{(X,Y)}(s,t)&=&\mathbb{E}(\exp(sX+tY)) \\
&=&\mathbb{E}(\exp(sX) \times \exp(tY)) \ \ \text{(Property of the exponential)}\\
&=&\mathbb{E}(\exp(sX)) \times \mathbb{E}(\exp(tY)) \ \ \text{(By Property (CI3))}\\
&=& \Psi_{X}(s) \times \Psi_{Y}(t).
\end{eqnarray*}

\bigskip \noindent \textbf{B - Moment generating Functions  \index{moment generating function} and Moments \index{moments} }.\\

\noindent Here, we do not deal with the existence of the mathematical expectation  \index{mathematical expectation} nor its finiteness. Besides, we admit that we may exchange the $\mathbb{E}$ symbol and the differentiation symbol $d()/ds$ in the form

\begin{equation}
\frac{d}{ds}(\mathbb{E}(f(s,X))=\mathbb{E}(\frac{d}{ds}f(s,X)), \label{proba01.rcouples.SWAPID}
\end{equation}

\Bin
\noindent where $f(s,X)$ is a function of real numbers $t$ for a fixed value $X(\omega)$. The validity of this operation is completely settled in the course of Measure and Integration of this series.\\

\noindent Let us successively differentiate $\Phi_X(s)$ with respect to $s$. Let us begin by a simple differentiation. We have :

\begin{eqnarray*}
\left( \Phi_{X}(s) \right)^{\prime}&=& \left(\mathbb{E}(\exp(sX)) \right)^{\prime}\\
&=& \mathbb{E}\left( \left(\exp(sX\right)^{\prime} \right)\\
&=& \mathbb{E}\left( X \exp(sX) \right).
\end{eqnarray*}

\Bin
\noindent Now, by iterating the differentiation, the same rule applies again and we have

\begin{eqnarray*}
\left( \Phi_{X}(s) \right)^{\prime}&=& \mathbb{E}\left( X \exp(sX) \right)\\
\left( \Phi_{X}(s) \right)^{\prime \prime}&=& \mathbb{E}\left( X^2 \exp(sX) \right)\\
\left( \Phi_{X}(s) \right)^{(3)}&=& \mathbb{E}\left( X^3 \exp(sX) \right)\\
\left( \Phi_{X}(s) \right)^{(4)}&=& \mathbb{E}\left( X^{(4)}\exp(sX) \right)\\
\end{eqnarray*}

\noindent We see that a simple induction leads to

\begin{eqnarray*}
\left(\Phi_{X}(s)\right)^{(k)}&=& \mathbb{E}\left( X^k \exp(sX) \right),
\end{eqnarray*}

\noindent where $\left(\Phi_{X}(s)\right)^{(k)}$ stands for the $k^{th}$ derivative of $\Phi_{X}(s)$. Applying this formula to zero gives

\begin{eqnarray*}
\left( \Phi_{X}(0) \right)^{(k)}&=& \mathbb{E}\left( X^k\right).\\
\end{eqnarray*}

\noindent We conclude : The moments  \index{moments} of a real-valued random variable, if they exist, may be obtained from the derivative of the \textit{m.g.f} of $X$ applied to zero :

\begin{equation}
m_{k}(X)=\mathbb{E}\left( X^k\right)=\left( \Phi_{X}(0) \right)^{(k)}, \ \  k \geq 1. \label{proba01.rcouples.MGFMmentsA}
\end{equation}

\bigskip
\noindent \textbf{Remark}. Now, we understand why the function is called the generation moment functions of $X$, since it leads to all the finite moments  \index{moments} of $X$.\\

\bigskip \noindent \textbf{C - The Second Moment generating Function  \index{moment generating function} and Moments \index{moments} }.\\

\noindent First define for an integer $h \geq 1$, following functions from $\mathbb{R}$ to $\mathbb{R}$ :

$$
(x)_{h}=x \times (x-1) \times .... \times (x-h+1), x\in \mathbb{R}.
$$

\bigskip \noindent For example, we have : $(x)_{1}=x$, $(x)_{2}=x(x-1)$, $(x)_{3}=x(x-1)(x-2)$. We already encountered these functions, since
$\mathbb{E}((X)_{2})$ is the factorial moment  \index{factorial moment} of second order. As well, we define

$$
fm_{h}(X)=\mathbb{E}((X)_{h}), \ \ h \geq 1,
$$

\bigskip \noindent as the factorial moment  \index{factorial moment} of order $h$. We are going to see how to find these factorial moments  \index{moments} from the second \textit{m.g.f}.\\

\noindent \noindent Using Formula (\ref{proba01.rcouples.SWAPID}) and the iterated differentiation of $\Psi_{X}$, together, lead to

\begin{eqnarray*}
\left( \Psi_{X}(s) \right)^{\prime}&=& \mathbb{E}\left( X s^{X-1}) \right)\\
\left( \Psi_{X}(s) \right)^{\prime \prime}&=& \mathbb{E}\left( X(X-1) s^{X-2} \right)\\
\left( \Psi_{X}(s) \right)^{(3)}&=& \mathbb{E}\left( X(X-1)(X-2) s^{X-3} \right)\\
\left( \Psi_{X}(s) \right)^{(4)}&=& \mathbb{E}\left( (X)_{(4)} s^{X-4} \right).\\
\end{eqnarray*}

\noindent We see that a simple induction leads to

\begin{eqnarray*}
\left(\Psi_{X}(s)\right)^{(k)}&=& \mathbb{E}\left( (X)_{(k)} s^{X-k} \right),\\
\end{eqnarray*}

\noindent where $\left(\Psi_{X}(s)\right)^{(k)}$ stands for the $k^{th}$ derivative of $\Psi_{X}(s)$. Applying this formula to the unity in $\mathbb{R}$  gives

\begin{eqnarray*}
\left( \Psi_{X}(1) \right)^{(k)}&=& \mathbb{E}\left( (X)_{(k)} \right).\\
\end{eqnarray*}

\bigskip \noindent We conclude as follows. The factorial moments  \index{factorial moment}  \index{moments} of a real-valued random variable, if they exist, may be obtained from the derivatives of \textit{m.g.f} of $X$ applied to the unity :

\begin{equation}
fm_{k}(X)=\mathbb{E}\left(  (X)_{(k)} \right)=\left( \Psi_{X}(1) \right)^{(k)}, \ \  k \geq 1. \label{proba01.rcouples.MGFMmentsB}
\end{equation}

\bigskip
\noindent \textbf{Remark}. This function also yields the moments  \index{moments} of $X$, indirectly through the factorial moment \index{factorial moment} s.\\

\bigskip \noindent \textbf{D - Probability laws  \index{probability law} convolution \index{convolution} }.\\

\noindent \textbf{ (1) Definition}.\\

\noindent A similar factorization formula (\ref{proba01.rcouples.MGFFACT1}) for independent real random variables $X$ and $Y$ is easy to get when we use the second \textit{m.g.f}, in virtue of the relation \ref{proba01.rcouples.MGF12}. Indeed, if $X$ and $Y$ are independent, we have for any
$s \in ]0,1]$,

$$
\Psi_{X+Y}(s)=\Phi_{X+Y}(\log s)=\Phi_{X}(\log s) \Phi_{Y}(\log s)=\Psi_{X}(s) \Psi_{Y}(s).
$$

\bigskip \noindent We may write this as follows. If $X$ and $Y$ are independent real random variables, then for any $s \in ]0,1]$,

\begin{equation}
\Psi_{X+Y}(s)=\Psi_{X}(s) \Psi_{Y}(s). \label{proba01.rcouples.MGFFACT1X}
\end{equation}

\bigskip
\noindent However, many authors use another approach to prove the latter result. That approach may be useful  in a number of studies. For example, it an instrumental tool for the study of Markov chains.\\

\noindent We are going to introduce it here.\\

\noindent \textbf{Convolutions  \index{convolution} of the laws of $X$ and $Y$}. The convolution of two probability laws  \index{probability law} on $\mathbb{R}$ is the law of the addition of two real random variables following these probability laws.\\

\noindent If $X$ and $Y$ are independent real random variables of respective probability laws  \index{probability law} $\mathbb{P}_X$ and $\mathbb{P}_Y$, the probability law of $X+Y$, $\mathbb{P}_{X+Y}$, is called the convolution  \index{convolution} product of the probability laws $\mathbb{P}_X$ and $\mathbb{P}_Y$, denoted by

$$
\mathbb{P}_{X} * \mathbb{P}_Y.
$$

\bigskip \noindent In the remainder of this part (D) of this section, we suppose for once that $X$ and $Y$ are independent real random variables.\\

\noindent Let us find the law of $Z=X+Y$. The domain of $Z$, defined by $\mathcal{V}_{Z}=\{z_k, \ k\in K\}$, is formed by the distinct values of the numbers $x_i+y_j, \ (i,j) \in I\times J$. Based on the following decomposition of $\Omega$,

$$
\Omega=\sum_{i\in I, j\in J} (X=x_i,Y=y_j),
$$

\noindent we have for $k \in K$,

$$
(X+Y=z_k)=\sum_{i\in I, j\in J} (X+Y=z_k)\cap (X=x_i,Y=y_j).
$$

\noindent We may begin by summing over $i \in I$, and then, we have $y_j=z_k-x_i$ on the event $(X+Y=x_k, X=x_i, Y=y_j)$. Thus we have

$$
(X+Y=k)=\sum_{i\in I} (X=x_i, Y=z_k-x_i).
$$

\noindent Then, the independence  \index{independence} between $X$ and $Y$ implies

\begin{equation}
\mathbb{P}(X+Y=z_k)=\sum_{i\in I} \mathbb{P}(X=x_i) \mathbb{P}(Y=z_k-x_i). \label{proba01.rcouples.conv01}
\end{equation}

\noindent \textbf{Warning}. When we use this formula, we have to restrict the summation over $i\in I$ to the values of $i$ such that the event $(Y=z_k-x_i)$ is not empty.\\

\noindent \textbf{ (2) Example}.\\

\noindent Consider the random variables $X$ and $Y$ that both follow a Poisson law with respective parameters $\lambda>0$ and $\mu>0$. We are going to use Formula
(\ref{proba01.rcouples.conv01}). Here for a fixed $k \in \mathbb{N}=\mathcal{V}_{X+Y}$, the summation over $i \in \mathbb{N}=\mathcal{V}_{X}$ will be restricted to $0\leq i \leq k$, since the events $(Y=k-i)$ are impossible for $i>k$. Then we have :

\begin{eqnarray}
\mathbb{P}(X+Y=k)&=&\sum_{i=0}^{k} \mathbb{P}(X=i) \mathbb{P}(Y=k-i)\\
&=& \sum_{i=0}^{k} \frac{\lambda^{i}e^{-\lambda}}{i!}\times \frac{\mu^{k-i}e^{-\mu}}{(k-i)!}\\
&=& \sum_{i=0}^{k} \frac{\lambda^{i}e^{-\lambda}}{i!}\times \frac{\mu^{k-i}e^{-\mu}}{(k-i)!}\\
&=& \frac{e^{-(\lambda+\mu)}}{k!} \biggr( \sum_{i=0}^{k} \frac{k!}{i!(k-i)!} \lambda^{i}\mu^{k-i} \biggr).
\end{eqnarray}

\bigskip
\noindent Now, we apply the Newton's Formula  \index{Newton's formula} (\ref{proba01.ForBin}) in Chapter \ref{proba01.combin} to get

\begin{eqnarray}
\mathbb{P}(X+Y=k)=\frac{e^{-(\lambda+\mu)}(\lambda+\mu)^k}{k!}.
\end{eqnarray}

\bigskip
\noindent \textbf{Conclusion}. By the characterization of probability laws  \index{probability law} by their first or second \textit{m.g.f}'s, whom we admitted in the beginning of this part I, we infer from the latter result, that the sum of two independent Poisson random variables  \index{Poisson random variable} with parameters $\lambda>0$ and $\mu>0$ is a Poisson random variable with parameter $\lambda+\mu$. We may represent that rule by

$$
\mathcal{P}(\lambda)*\mathcal{P}(\mu)=\mathcal{P}(\lambda+\mu).
$$

\bigskip \noindent \textbf{ (2) Direct Factorization Formula}.\\

\noindent In the particular case where the random variables $X$ and $Y$ take the non-negative integers values, we may use Analysis results to directly establish Formula (\ref{proba01.rcouples.MGFFACT1X}).\\

\noindent Suppose that $X$ and $Y$ are independent and have the non-negative integers as values. Let $\mathbb{P}(X=k)=a_k$, $\mathbb{P}(Y=k)=b_k$ and $\mathbb{P}(X+Y=k)=c_k$, for $k\geq 0$. Formula (\ref{proba01.rcouples.conv01}) becomes

\begin{equation}
c_n=\sum_{k=}^{n} a_{k} b_{n-k}, \ n\geq 0. \label{proba01.rcouples.conv02}
\end{equation}

\noindent \textbf{Definition}. The convolution  \index{convolution} of the two sequences of real numbers $(a_n)_{n\geq 0}$ and $(b_n)_{n\geq 0}$ is the sequence of real numbers $(c_n)_{n\geq 0}$ defined in (\ref{proba01.rcouples.conv02}).\\

\noindent If the sequence $(a_n)_{n\geq 0}$ is absolutely convergent, that is,

$$
\sum_{n\geq 0} |a_n| < +\infty,
$$

\noindent we may define the function

$$
a(s)=\sum_{n\geq 0} a_{n} s^{n} < +\infty, \ \ 0<s\leq 1.
$$

\bigskip \noindent We have the following property.\\

\begin{lemma} \label{proba01.rcouples.conv03} Let $(a_n)_{n\geq 0}$ and $(b_n)_{n\geq 0}$ be two absolutely convergent sequences of real numbers, and $(c_n)_{n\geq 0}$ be their convolution,  \index{convolution} defined in (\ref{proba01.rcouples.conv02}). Then we have

$$
c(s)=a(s) b(s), 0<s\leq 1.
$$
\end{lemma}

\bigskip \noindent \textbf{Proof}. We have

\begin{equation*}
c(s)=\sum_{n\geq 0}c_{n}s^{n}=\sum_{n\geq 0}\left(
\sum_{k=0}^{n}a_{k}b_{n-k}\right) s^{n}=\sum_{n=0}^{\infty
}\sum_{k=0}^{n}\left( a_{k}s^{k}\right) \left( b_{n-k}s^{n-k}\right)
\end{equation*}
\begin{equation*}
=\sum_{n=0}^{\infty }\sum_{k=0}^{\infty }\left( a_{k}s^{k}\right) \left(
b_{n-k}s^{n-k}\right) 1_{(k\leq n)}
\end{equation*}

\Bin
\noindent Let us apply Fubini's  \index{Fubini's} property for convergent sequences of real numbers by exchanging the two summation symbols, to get

\begin{eqnarray*}
c(s)&=&\sum_{k=0}^{\infty }\sum_{n=0}^{\infty }\left( a_{k}s^{k}\right) \left(b_{n-k}s^{n-k}\right) 1_{(k\leq n)}\\
c(s)&=&\sum_{k=0}^{\infty }\sum_{n=k}^{\infty }\left( a_{k}s^{k}\right) \left(b_{n-k}s^{n-k}\right)\\
&=& \sum_{k=0}^{\infty }\left( a_{k}s^{k}\right) \left\{\sum_{n=k}^{\infty }\left( b_{n-k}s^{n-k}\right) \right\}.
\end{eqnarray*}

\bigskip \noindent The quantity in the brackets in the last line is equal to $b(s)$. To see this, it suffices to make the change of variables $\ell=n-k$ and $\ell$ runs from $0$ to $+\infty$. QED.\\

\bigskip \noindent By applying this lemma to the discrete probability laws  \index{probability law} with non-negative integer values, we obtain the factorization formula (\ref{proba01.rcouples.MGFFACT1X}).\\

\bigskip \noindent We concluding by saying that : it is much simpler to work with the first \textit{m.g.f}. But in some situations, like when handling the Markov chains, the second \textit{m.g.f} may be very useful.\\

\bigskip \noindent \textbf{F - Simple examples of \textit{m.g.f}}.\\

\noindent Here, we use the usual examples of discrete probability law  \index{probability law} that were reviewed in Section \ref{proba01.rv.review} of Chapter \ref{proba01.rv} and the properties of the mathematical expectation  \index{mathematical expectation} stated in Chapter \ref{proba01.param}

\bigskip \noindent (1) Constant $c$ with $\mathbb{P}(X=c)=1$ :\\

\begin{eqnarray*}
\Phi_{X}(s)&=&\mathbb{E}\left( \exp(sX) \right)=1  \times e^{cs}\\
&=& e^{cs}, s\in \mathbb{R}.
\end{eqnarray*}

\bigskip  \noindent (2) Bernoulli law $\mathcal{B}(p)$, $p\in]0,1[$.\\

\begin{eqnarray*}
\Phi_{X}(s)&=&\mathbb{E}\left( \exp(sX) \right)=q \times e^{0s}+p e^{1s}\\
&=& q+pe^s, s\in \mathbb{R}.
\end{eqnarray*}

\bigskip \noindent (3) Geometric law $\mathcal{G}(p)$, $p\in]0,1[$.\\

\begin{eqnarray*}
\Phi_{X}(s)&=&\mathbb{E}\left( \exp(sX) \right)\\
&=& \sum_{n\geq 1} q^{n-1}p e^{ns}\\
&=& (pe^s) \sum_{n\geq 1} q^{n-1} e^{(n-1)s}\\
&=& (pe^s) \sum_{k\geq 0} q^{k} e^{ks}, \text{ (by change of variable: } k=n-1) \\
&=& (pe^s) \sum_{k \geq 0} (qe^{s})^k \\
&=& \frac{pe^{s}}{1-qe^{s}}, \text{ for } 0<qe^{s}<1, \text{ i.e. } s<-\log q.
\end{eqnarray*}

\bigskip \noindent (4) $X$ is the number of failures before the first success $X \sim \mathcal{NF}(p)$, $p\in]0,1[$.\\

\noindent If $X$ is the number of failures before the first success, then  $Y=X+1$ follows a geometric law $\mathcal{G}(p)$. Thus, we may apply Formula (\ref{proba01.rcouples.MGFTransform}) and the \textit{m.g.f} of $Y$ given by Point (3) above to have

$$
\Phi{X}(s)=\Phi_{Y-1}(s)=e^{-s} \times \Phi_{Y}(s)=e^{-s} \times \frac{pe^s}{1-qe^{s}}=\frac{p}{1-qe^{s}}.
$$

\bigskip \noindent (5) Discrete Uniform Law on {1,2,...,n} $\mathcal{DU}(n)$, $n\geq 1$.\\

\noindent We have

\begin{eqnarray*}
\Phi_{X}(s)&=&\mathbb{E}\left( \exp(sX) \right)\\
&=& \sum_{1 \leq i \leq n} \frac{1}{n} e^{is}\\
&=& \frac{e^{s}(1+e^{s}+e^{2s}+...+e^{(n-1)s}}{n} \\
&=& \frac{(e^{(n+1)s}-e^{s}}{n(e^{s}-1)}. \\
\end{eqnarray*}

\bigskip \noindent (6) Binomial Law $\mathcal{B}(n,p)$, $p \in ]0,1[$, $n\geq 1$.\\

\noindent Since a $\mathcal{B}(n,p)$  random variable has the same law as the sum of $n$ independent Bernoulli $\mathcal{B}(p)$ random variables, the factorization formula (\ref{proba01.rcouples.MGFFACT1}) and the expression the \textit{m.g.f} of a $\mathcal{B}(p)$ given in Point (2) above lead to

$$
\Phi_{X}(s)=(q+pe^s)^n, s\in \mathbb{R}.
$$

\Bin
\bigskip \noindent (7) Negative Binomial Law $\mathcal{NB}(k,p)$, $p \in ]0,1[$, $k\geq 1$.\\

\noindent Since a $\mathcal{NB}(k,p)$ random variable has the same law as the sum of $k$ independent Geometric $\mathcal{G}(p)$ random variables, the factorization formula (\ref{proba01.rcouples.MGFFACT1}) and the expression the \textit{m.g.f} of a $\mathcal{G}(p)$ given in Point (3) above lead to

$$
\Phi_{X}(s)=\left(\frac{pe^{s}}{1-qe^{s}}\right)^k, \ , \text{ for } 0<qe^{s}<1, \text{ i.e. } s<-\log q.
$$

\bigskip \noindent (8) Poisson Law $\mathcal{P}(\lambda)$, $\lambda >0$. We have

\begin{eqnarray*}
\Phi_{X}(s)&=&\mathbb{E}\left( \exp(sX) \right)\\
&=& \sum_{k\geq 0} \frac{\lambda^{k}e^{\lambda}}{k!} e^{ks}\\
&=& e^{-\lambda} \sum_{k\geq 0} \frac{(\lambda e^{s})^{k}}{k!} \\
&=& e^{-\lambda} \exp(\lambda e^{s}) \text{ (Exponential Expansion) }\\
&=& \exp(\lambda(e^{s}-1)).\\
\end{eqnarray*}

\newpage

\Bin \textbf{F - Table of some usual discrete probability law \index{probability law} s}.\\

\noindent \textbf{Abbreviations}.\\

\begin{table}[htbp]
\centering
\begin{tabular}{|l|l|l|l|}
\hline
Name & Symbol & Parameters & Domain\\
\hline
Const. & $c$ & $c\in \mathbb{R}$ & $\mathcal{V}={c}$ \\
\hline
Bernoulli & $\mathcal{B}(p)$ & $0<p<1$ &$\mathcal{V}=\{0,1\}$\\
\hline
Discrete Uniform 	 &	$\mathcal{DU}(n)$ & $n\geq 1$ &  $\mathcal{V}=\{1,...,n\}$ \\
\hline
Geometric &	$\mathcal{G}(p)$ & $0<p<1$ & $\mathcal{V}={1,2,...}$ \\
\hline
Binomial &	$\mathcal{B}(n,p)$&  $n\geq 1$, $0<p<1$ & $\mathcal{V}={0,1,...,n}$\\
\hline
Negative Binomial &	$\mathcal{NB}(k,p)$& $k\geq 1$, $0<p<1$& $\mathcal{V}=\{k,k+1,...1\}$\\
\hline
Poisson  &	$\mathcal{P}(\lambda)$ & $\lambda>0$ & $\mathcal{V}=\mathbb{N}=\{0,1,...\}$\\
\hline
\end{tabular}
\vskip 0.5cm
\caption{Abbreviations and names of of some usual probability Laws  \index{probability law} on $\mathbb{R}$}
\label{propb01.rcouples.Tab01}
\end{table}

\newpage
\noindent \textbf{Probability Law \index{probability law} s}

\begin{table}[htbp]
\centering
\begin{tabular}{|l|l| l| l| l|}
\hline
Symbol & Probability law.  \index{probability law} & (ME) & Variance  \index{variance} & \textit{m.g.t} in (s)\\
\hline
$c$ & $1$, $k=1$ & $c$ & 0 & $\exp(cs)$ \\
\hline
$\mathcal{B}(p)$  & $p$, $1-p$, $k=0,1$        & $p$ 						& $p(1-p)$ & $q+pe^{s}$\\
\hline
$\mathcal{DU}(n)$& $\frac{1}{n},1\leq i\leq n$  & $\frac{n+1}{2}$ &  &$\frac{(e^{(n+1)s}-e^{s}}{n(e^{s}-1)}$ \\
\hline
$\mathcal{G}(p)$ & $p(1-p)^{n-1},n\geq 1$  & $\frac{p}{q}$ & $\frac{p}{q^2}$ & $\frac{pe^s}{1-e^s}$ \\
\hline
$\mathcal{B}(n,p)$& $\left(\begin{tabular}{c} $n$ \\ $k$ \end{tabular} \right) p^k (1-p)^{n-k},1\leq k\leq n$  & $np$ & $np(1-p)$ & $(q+pe^{s})^n$\\
\hline
$\mathcal{NB}(n, \ k)$& $\left(\begin{tabular}{c} $n-1$ \\ $k-1$ \end{tabular} \right) p^k(1-p)^{n-k},n\geq k$  & $\frac{kp}{q}$ & $\frac{kp}{q^2}$ &	$(\frac{pe^s}{1-e^s})^k$\\
\hline
$\mathcal{P}(\lambda)$ & $\frac{\lambda \exp(-\lambda)}{k!},k\geq 0$  & $\lambda$ & $\lambda$ & $\exp(\lambda(e^s-1))$\\
\hline
\end{tabular}
\vskip 0.5cm
\caption{Probability Laws  \index{probability law} of some usual real-valued discrete random variables. (ME) stands for Mathematical Expectation \index{mathematical expectation} .}
\label{propb01.rcouples.tab01}
\end{table}

\bigskip \noindent \textbf{I - Conditional  \index{conditional} Probability Law \index{probability law} s}.\\

\noindent The conditional  \index{conditional} probability law  \index{probability law} of $X$ given $Y=y_j$ is given by the conditional probability formula

\begin{equation}
\mathbb{P}(X=x_i/Y=y_j)=\frac{\mathbb{P}(X=x_i, Y=y_j)}{\mathbb{P}(Y=y_j)}. \label{proba01.rcouples.cp01}
\end{equation}

\Bin \noindent If we use an extended domain of $Y$, this probability is $0$ if $\mathbb{P}(Y=y_j)=0$. We adopt the following. For $i\in I,j\in J$, set

$$
p_{i}^{(j)}= \mathbb{P}(X=x_i/Y=y_j).
$$

\Bin
\noindent According to the notation already introduced, we have for a fixed $j \in J$

\begin{equation}
p_{i}^{(j)}=\frac{p_{i,j}}{p_{\bullet,j}}, \ i\in I. \label{proba01.rcouples.cp01a}
\end{equation}

\noindent We may see, by Formula (\ref{proba01.rcouples.mplY2}), that for each $j\in J$, the numbers $p_{i}^{(j)}, \ i\in I$, stands for a discrete probability measure.  \index{probability measure} It is called the probability law  \index{probability law} of $X$ given $Y=y_j$.\\

\noindent We may summarize the conditional  \index{conditional} laws in Table \ref{proba01.rcouples.t05}.

\begin{table}[htbp]
\centering
\begin{tabular}{|l|l|l|l|l|}
\hline
X 	& $y_{1}$ 			& $\cdots$ & $y_{j}$ 			& $\cdots$ 					\\
\hline
$x_{1}$ & $p_{1}^{(1)}$ 	& $\cdots$ & $p_{1}^{(j)}$ & $\cdots$ 					\\
\hline
$x_{2}$ & $p_{2}^{(1)}$ 	& $\cdots$ & $p_{2}^{(j)}$ & $\cdots$  				\\
\hline
$\vdots$& $\vdots$ 			& $\vdots$ & $\vdots$ 		& $\vdots$ 				\\
\hline
$x_{i}$ & $p_{i}^{(1)}$ 	& $\cdots$ & $p_{i}^{(j)}$ & $\cdots$  			\\
\hline
$\vdots$ & $\vdots$ 		& $\vdots$ & $\vdots$ 		& $\vdots$  \\
\hline
Total    & $100\%$  	 & $\cdots$  &  $100\%$ 		& $\cdots$ \\
\hline
\end{tabular}
\vskip 0.5cm
\caption{Conditional  \index{conditional} Probability Laws  \index{probability law} of $X$}
\label{proba01.rcouples.t05}
\end{table}

\noindent Each column, with the exception of the first, represents a conditional  \index{conditional} probability law  \index{probability law} of $X$. This table is read column by column.\\

\noindent If we have considered the conditional  \index{conditional} probability laws  \index{probability law} of $Y$ given some value $x_i$ of $X$, we would have a table to be read by lines, each line representing a conditional probability law of $Y$. We only study the conditional probability laws of $X$. The results and formulas we obtained are easy to transfer to the conditional laws of $Y$.\\

\noindent The most remarkable conclusion of this point is the following :\\

\noindent \textbf{The random variables (not necessarily real-valued) $X$ and $Y$ are independent if and only if any conditional  \index{conditional} probability law  \index{probability law} of $X$ given a value of $Y$ is exactly the marginal  \index{marginal} probability law of $X$ (called its unconditional probability law).}\\

\noindent This is easy and comes from the following lines. $X$ and $Y$ are independent if and only if for any $j\in J$, we have

$$
p_{i,j}=p_{i,\bullet} \times p_{\bullet,j}, \ i\in I.
$$

\Bin  This assertion is equivalent to : for any $j\in J$, we have

$$
p_{i}^{(j)}=\frac{p_{i,j}}{p_{\bullet,j}}=\frac{p_{i,\bullet} \times p_{\bullet,j}}{p_{\bullet,j}}=p_{i,\bullet}, \ i\in I.
$$

\Bin
\noindent Then, $X$ and $Y$ are independent if and only if for any $j\in J$,

$$
p_{i}^{(j)}=p_{i,\bullet}, \ i\in I.
$$

\Bin
\noindent We denote by $X^{(j)}$ a random variable taking the same values as $X$ whose probability law  \index{probability law} is the conditional  \index{conditional} probability law of $X$ given $Y=y_j$.\\

\noindent \textbf{II - Computing the mathematical expectation  \index{mathematical expectation} of a function of $X$}.\\

\noindent Let $h(X)$ be a real-valued function of $X$. We define the mathematical expectation  \index{mathematical expectation} of $X$ given $Y=y_j$, denoted by
$\mathbb{E}(h(X)/Y=y_{j})$ as the mathematical expectation  \index{mathematical expectation} of $h(X^{(j)})$. It is given by

$$
\mathbb{E}(h(X^{(j)}))=\sum_{i \in I} h(x_i) p_{i}^{(j)}.
$$

\noindent This gives the two expressions:

\begin{equation}
\mathbb{E}(h(X)/Y=y_{j})=\sum_{i \in I} h(x_i) p_{i}^{(j)} \label{proba01.rcouples.cp03}
\end{equation}

\noindent or

\begin{equation}
\mathbb{E}(h(X)/Y=y_{j})=\sum_{i \in I} h(x_i) \mathbb{P}(X=x_i/Y=y_j). \label{proba01.rcouples.cp03a}
\end{equation}

\noindent It is clear that, by construction, the mathematical expectation  \index{mathematical expectation} given $Y=y_j$ is linear. If $Z$ is another discrete random variable and $\ell(Z)$ is a real-valued function of $Z$, if $a$ and $b$ are two real numbers, we have, for each $j\in J$,

$$
\mathbb{E}(ah(X)+b\ell(Z)/Y=y_{j})=a\mathbb{E}(h(X)/Y=y_{j})+ b \mathbb{E}(\ell(Z)/Y=y_{j}).
$$

\bigskip \noindent We have the interesting and useful formula for computing the mathematical expectation  \index{mathematical expectation} of $h(X)$.\\

\begin{equation}
\mathbb{E}(h(X))=\sum_{j \in J} \mathbb{P}(Y=y_{j}) \mathbb{E}(h(X)/Y=y_{j}). \label{proba01.rcouples.cp10a}
\end{equation}

\noindent Let us give a proof of this formula. We have,

\begin{eqnarray*}
&& \sum_{j \in J} \mathbb{P}(Y=y_{j}) \mathbb{E}(h(X)/Y=y_{j})\\
&=& \sum_{j \in J} \mathbb{P}(Y=y_{j}) \sum_{i \in I} h(x_i) p_{i}^{(j)} \text{( by Forumla (\ref{proba01.rcouples.cp03}))}\\
&=&  \sum_{j \in J} \sum_{i \in I} h(x_i)  \mathbb{P}(Y=y_{j}) p_{i}^{(j)}\\
&=&  \sum_{j \in J} \sum_{i \in I} h(x_i)  \mathbb{P}(X=x_i, \ Y=y_{j})  \text{( by definition in Formula (\ref{proba01.rcouples.cp01}) )}\\
&=&  \sum_{i \in I} h(x_i) \sum_{j \in J}   \mathbb{P}(X=x_i,  Y=y_{j})  \text{( by Fubini's  \index{Fubini's} rule)}\\
&=&  \sum_{i \in I} h(x_i) \mathbb{P}(X=x_i)\\
&=&  \mathbb{E}(h(X)).
\end{eqnarray*}

\bigskip \noindent \textbf{A simple Exercise}. Let $X_1$, $X_2$, ..., be a sequence of real-valued discrete random variables with common mean  \index{mean} $\mu$ and let $Y$ a non-negative integer real-valued discrete random variable whose mean is $\lambda$. Define

$$
S=X_1+X_2+ \cdots + X_{Y}=\sum_{i=1}^{Y} X_{i}.
$$

\noindent Show that $\mathbb{E}(S)=\mu\lambda$.\\

\noindent \textbf{Solution}. Let us apply Formula (\ref{proba01.rcouples.cp10a}). First, we have for any $j \in J$

\begin{eqnarray*}
\mathbb{E}(S/Y=y_j)&=&\mathbb{E}(X_1+X_2+ \cdots + X_{Y}/Y=y_j)\\
&=&\mathbb{E}(X_1+X_2+ \cdots + X_{Y}/Y=y_j)\\
&=&\mathbb{E}(X_1+X_2+ \cdots + X_{y_j}) \ \ L3\\
&=& \mu y_{j}
\end{eqnarray*}

\Bin where, in Line L3, we used the fact that the number of terms is nonrandom at that step. Now, the application of Formula (\ref{proba01.rcouples.cp10a}) gives

\begin{eqnarray*}
\mathbb{E}(S)&=&\sum_{j \in J} \mathbb{P}(Y=y_{j}) \mathbb{E}(Z/Y=y_{j})\\
&=& \mu \sum_{j \in J} \mathbb{P}(Y=y_{j}) y_{j}\\
&=& \mu \mathbb{E}(Y)=\mu \lambda.
\end{eqnarray*}

\bigskip \noindent \textbf{Conditional  \index{conditional} Expectation Random Variable}.\\

\noindent Now, Let us consider $\mathbb{E}(h(X)/Y=y_{j})$ as a function of $y_j$ denoted $g(y_j)$. So, $g$ is a real-valued function which is defined on a domain including the values set of $Y$, such that for any $j \in J$

$$
\mathbb{E}(h(X)/Y=y_{j})=g(y_j).
$$

\Bin This function if also called the regression function of $X$ given $Y$.\\

\bigskip \noindent \textbf{Definition}. The random variable $g(Y)$ is the conditional  \index{conditional} expectation of $h(X)$ given $Y$ and is denoted by
$\mathbb{E}(h(X)/Y)$.\\

\noindent \textbf{Warning}. The conditional  \index{conditional} expectation of a real-valued function of $h(X)$ given $Y$ is a random valued.\\

\noindent \textbf{Properties}.\\

\noindent \textbf{(1)} \textit{How to use the conditional  \index{conditional} expectation to compute the mathematical expectation \index{mathematical expectation} ?}\\

\noindent We have an extrapolation of Formula (\ref{proba01.rcouples.cp10a}) in the form :

$$
\mathbb{E}(h(X))= \mathbb{E}\biggr(\mathbb{E}(h(X)/Y) \biggr).
$$

\bigskip
\noindent Let us consider the exercise just given above. We had for each $j \in J$

$$
\mathbb{E}(S/Y=y_j)=g(y_{j}),
$$

\bigskip
\noindent with $g(y)=\mu y$. Then we have $\mathbb{E}(S/Y)=\mu Y$. Next, we get

$$
\mathbb{E}\biggr(\mathbb{E}(h(S)/Y) \biggr)=\mathbb{E}(\mu Y)=\mu \mathbb{E}(Y)=\mu \lambda=\mathbb{E}(S).
$$

\bigskip \noindent \textbf{(2)} \textit{Case where $X$ and $Y$ are independent}.\\

\noindent We have

$$
\mathbb{E}(h(X)/Y)=\mathbb{E}(h(X)).
$$

\bigskip \noindent \textbf{Proof}. Suppose that $X$ and $Y$ are independent. So for any $(i,j) \in I \times J$, we have

$$
\mathbb{P}(X=x_i/Y=y_j)=\mathbb{P}(X=x_i).
$$

\Bin \noindent Then, by Formula (\ref{proba01.rcouples.cp03a})

\begin{eqnarray*}
\mathbb{E}(h(X)/Y=y_{j})&=&\sum_{i \in I} h(x_i) \mathbb{P}(X=x_i/Y=y_j)\\
&=&\sum_{i \in I} h(x_i) \mathbb{P}(X=x_i)\\
&=& \mathbb{E}(h(X)).\\
\end{eqnarray*}

\noindent QED.\\

\bigskip
\noindent \textbf{(3)} \textit{Conditional  \index{conditional} Expectation of a multiple of a function of $Y$}.\\

\noindent Let $\ell$ be a real-value function whose domain contains the values of $Y$. Then, we have

$$
\mathbb{E}(\ell(Y) h(X)/Y)=\ell(Y) \mathbb{E}(h(X)/Y).
$$

\bigskip \noindent \textbf{Proof}. We have

\begin{eqnarray*}
\mathbb{E}(\{\ell(Y)h(X)\}/\{Y=y_{j}\})&=&\mathbb{E}(\ell(y_j)h(X)/Y=y_{j})\\
&=&\ell(y_j) \mathbb{E}(h(X)/Y=y_{j}.\\
\end{eqnarray*}

\noindent So, if $\mathbb{E}(h(X)/Y=y_{j})=g(y_j)$, the latter formula says that

$$
\mathbb{E}(\ell(Y)h(X)/Y)=\ell(Y)g(Y)=\ell(Y) \mathbb{E}(h(X)/Y).
$$

\bigskip \noindent which was the target. QED.\\

\bigskip

\noindent This gives the rule : when computing a conditional  \index{conditional} expectation given $Y$, may get the factors that are functions of $Y$ out of the conditional expectation.\\

\noindent In particular, if $h$ is the function constantly equal to one, we have

$$
\mathbb{E}(\ell(Y)/Y)=\ell(Y).
$$

\bigskip \noindent The conditional  \index{conditional} expectation of a function of $Y$, given $Y$, is itself.

\noindent
 
 %proba01.rcouples
\chapter{Continuous Random Variables}  \index{continuous Random Variable} \label{proba01.crv}

\Bin Until now, we exclusively dealt with discrete random variables. Introducing probability theory  \index{probability theory} in discrete space is a good pedagogical approach. Beyond this, we will learn in advanced courses that the general Theory of Measure and Integration and Probability Theory are based on the method of discretization. General formulas depending on measurable mappings and/or probability measure  \index{probability measure} are extensions of the same formulas established for discrete applications and/or discrete probability measures. This means  \index{mean} that the formulas stated in this textbook, beyond their usefulness in real problems, actually constitute the foundation of the Probability Theory.\\

\noindent In this chapter we will explain the notion of continuous random variables  \index{continuous Random Variable} as a consequence of the study of the cumulative distribution function  \index{cumulative Distribution Function} (\textit{cdf}) properties. We will provide a list of a limited number of examples as an introduction to a general chapter on probability law \index{probability law} s.\\

\noindent Especially, we will see how our special guest, the \textbf{normal} or \textbf{Gaussian} probability law,  \index{probability law} has been derived from the historical works of \textit{de Moivre}  \index{de Moivre} (1732) and Laplace  \index{Laplace} (1801) using elementary real calculus courses.\\

\Bin Let us begin by the following definition.\\

\noindent \textbf{Definition}. Let $X$ be a \textit{rrv}. The following function defined from $\mathbb{R}$ to $[0,1]$ by

\begin{equation*}
F_{X}(x)=\mathbb{P}(X\leq x)\text{ for }x\in \mathbb{R},
\end{equation*}

\Bin
\noindent is called the \textbf{cumulative distribution function}  \index{cumulative Distribution Function} \textit{cdf} of the random variable $X$.

\bigskip \noindent Before we proceed any further, let us consider two very simple cases of \textit{cdf}'s.\\

\noindent \textbf{Example 1. Cumulative distribution function  \index{cumulative Distribution Function} of a constant \textit{rrv}}. Let $X$ be the constant \textit{rrv} at $a$, that is

\begin{equation*}
\mathbb{P}(X=a)=1.
\end{equation*}

\Bin
\noindent It is not difficult to see that we have the following facts : $(X\leq x)=\emptyset$ for $x<a$ and $(X\leq x)=\Omega$ for $x\geq a$. Based on these facts, we see that \textit{cdf} of $X$ is given by \\

$F_{X}(x)=\left\{
\begin{tabular}{lll}
$1$ & if & $x\geq a$ \\
$0$ & if & $x<a$%
\end{tabular}%
\right. $.

\bigskip
\bigskip \noindent \textbf{Example 2. Cumulative Distribution Function  \index{cumulative Distribution Function} of a Bernoulli $\mathcal{B}(p)$ \textit{rrv}}.\\

\noindent Remark that $(X\leq x)=\emptyset$ for $x<0$, $(X\leq x)=(X=0)$ for $0\leq x <1$ and $(X\leq x)=(X=1)$ for $x\geq 1$. From these remarks, we derive that:\\

$F_{X}(x)=\left\{
\begin{tabular}{lll}
$1$ &  & $x\geq 1$ \\
$1-p$ &  & $0\leq x<1$ \\
$0$ &  & $x<0$%
\end{tabular}%
.\right. $ \ \ \ \

\bigskip \noindent We may infer from the two previous example the general method of computing the \textit{cdf} of a discrete \textsl{rrv}.

\noindent Suppose that $X$ is a discrete \textit{rrv} and takes its values in $\mathcal{V}_X=\left\{x_{i},\text{ \ \ }i\in I\right\}$. We suppose that the $x_{i}$'s are listed in an increasing order :

$$
\left\{ x_{1}\leq x_{2}\leq ...\leq x_{j}...\right\},
$$

\Bin \noindent with convention that $x_{0}=-\infty$. Let us denote

\begin{equation*}
p_{k}=\mathbb{P}(X=x_{k}), \ \  k\geq 1.
\end{equation*}

\noindent As well, we may define the cumulative probabilities

\begin{equation*}
p_{1}^{\ast}=p_{1},
\end{equation*}

\noindent and

\begin{equation*}
p_{k}^{\ast}=p_{1}+...+p_{k-1}+p_{k}.
\end{equation*}

\noindent We say that $p_{k}^{\ast}$ is the cumulative probabilities up to $k\geq 1$.\\

\noindent The cumulative distribution function  \index{cumulative Distribution Function} of $X$ is characterized by

\begin{equation}
\mathbb{P}(X\leq x)=\left\{
\begin{tabular}{lll}
$p_{k}^{\ast }$ & if & $x_{k}\leq x\text{ }<\text{ }x_{k+1},\text{ for some }%
k\text{ }\geq 1$ \\
$0$ & if & otherwise%
\end{tabular}%
.\right.   \label{proba01.crv.fr1}
\end{equation}

\Bin
\noindent In short, we may write

\begin{equation}
\mathbb{P}(X\leq x)=p_{j}^{\ast }\text{ if }x_{j}\leq x\text{ }<\text{ }x_{j+1}, \label{proba01.crv.fr2}
\end{equation}

\Bin \noindent  while keeping in mind that $F_{X}(x)$ is \textit{zero} if $x$ is strictly less than all the elements of $\mathcal{V}_X$(case where $\mathcal{V}_X$ is bounded below) and \textit{one} if $x$ is greater or equal to all the elements of $\mathcal{V}_X$ (case where $\mathcal{V}_X$ is bounded above).\\

\noindent Be careful about the fact that a discrete \textit{rrv} may have an unbounded values both below and above. For example, let us consider a random variable $X$ with values $\mathcal{V}_X=\{ j \in \mathbb{Z} \}$ such that

$$
\mathbb{P}(X=j)=\frac{1}{(2e-1)|j|!}.
$$

\bigskip
\noindent It is clear that $\mathcal{V}_X$ is unbounded. Proving that the sum of the probabilities given above is \textit{one} is left to the reader an exercise. He or She is suggested to use Formula (\ref{proba01.rv.expand.expo}) of Chapter \ref{proba01.rv} in his or her solution.

\bigskip

\noindent \textbf{Proof of (\ref{proba01.crv.fr1})}. We already know that

\begin{equation*}
\ (X\leq x)=\bigcup_{x_{i}\leq x}(X=x_{i}).
\end{equation*}

\Bin
\noindent Next, if $x_{k} \leq x < x_{k+1}$, the values $x_j$ less or equal to $x$ are exactly : $x_1$, $x_2$, ..., $x_k$. Combining these two remarks leads to

\begin{equation*}
(X\leq x)=\bigcup_{i=1}^{i=k} (X=x_{i}).
\end{equation*}

\noindent By Theorem \ref{proba01.rv.theo1} in Chapter \ref{proba01.rv}, we have

\begin{equation*}
\mathbb{P}(X\leq x)=\sum_{i=1}^{k}\mathbb{P}(X=x_{i})=p_{k}^{\ast},
\end{equation*}

\noindent for $x_{k} \leq x < x_{k+1}$. Further, if $x <x_1$, the event $(X \leq x)$ is impossible and $\mathbb{P}(X\leq x)=0$. This proves Formula \ref{proba01.crv.fr1}.\\

\bigskip

\noindent The \textit{cdf} of discrete random variables are piecewise constant functions, meaning  \index{mean} that they take constant values on segments of $\mathbb{R}$. Here are some simple illustrations for finite values sets.\\

\noindent \textbf{Example 3}. Let $X$ be a \textit{rrv} whose values set and probability law  \index{probability law} are defined in the following table.

\bigskip

\begin{center}
\begin{tabular}{|l|l|l|l|l|}
\hline
$X$ & $0$ & $1$ & $2$ & $3$ \\ \hline
$\mathbb{P}(X=k)$ & 0.125 & 0.375 & 0.375 & 0.125 \\ \hline
$p_{k}^{\ast }$ & 0.125 & 0.5 & 0.875 & 1 \\ \hline
\end{tabular}
\end{center}

\bigskip \noindent The graph of the associated \textit{cdf} defined by\\

\begin{center}
\begin{tabular}{|l|l|l|l|l|l|}
\hline
$x$ & $x<0$ & $0\leq x<1$ & $1\leq x<2$ & $2\leq x<3$ & $3x \geq 3$ \\ \hline
$\mathbb{P}(X\leq x)=p_{k}^{\ast }$ & 0 & 0.125 & 0.5 & 0.875 & 1 \\ \hline
\end{tabular}
\end{center}

\Bin \noindent can be found in Figure \ref{proba01.crv.fig1}

\begin{figure}[htbp]
\centering
\includegraphics[width=0.90\textwidth]{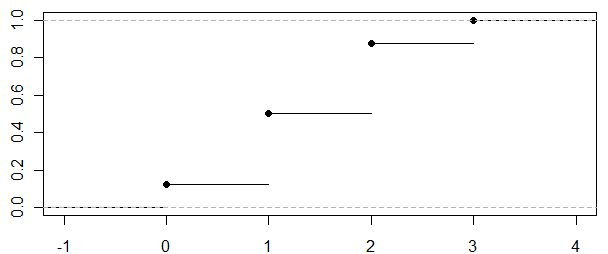}
\caption{DF}
\label{proba01.crv.fig1}
\end{figure}

\bigskip \noindent In the next section, we are going to discover properties of the \textit{cdf} in the discrete approach. Next, we will take these properties as the definition of a \textit{cdf} of an arbitrary \textit{rrv}. By this way, it will be easy to introduce continuous \textit{rrv}.

\noindent The properties of \textit{cdf}'s are given in the following proposition.\\

\begin{proposition} \label{proba01.crv.prop01} Let $F_{X}$ be the \textit{cdf} of a discrete random variable $X$. Then, $F_{X}$ fulfills the following properties :\\

\noindent \textbf{(1)}  $0\leq F_{X}\leq 1.$\\

\noindent \textbf{(2)} $\lim_{x\rightarrow -\infty}F_{X}(x)=0$ and $\lim_{x\rightarrow +\infty}F_{X}(x)=1$.\\

\noindent \textbf{(3)}  $F_{X}$ is right-continuous.\\

\noindent \textbf{(4)}  $F_{X}$ is $1$-non-decreasing, that is for $a\leq b$, $(a,b)\in \mathbb{R}^2$,

\begin{equation}
\Delta_{a,b}F_X=F_X(b)-F_X(a) \geq 0. \label{propa01.crv.DefDelta01}
\end{equation}

\bigskip
\noindent \textbf{Other terminology}. We will also say that $F_X$ assigns non-negative lengths to intervals $[a,b]$, $a\leq b$. We will come back to this terminology.
\end{proposition}

\bigskip \noindent \textbf{Proof.}\\

\noindent \textbf{A general remark}. In the proof below, we will have to deal with limits of $F_X$. But $F_X$ will be proved to be non-decreasing in Point (4). By Calculus Courses, general limits of non-decreasing functions, whenever they exist, are the same as monotone limits. So, when dealing with limits of $F_X$, we restrict ourselves to monotone limits.\\

\noindent In the sequel, a broadly increasing (resp. decreasing) sequence to something \textit{should be interpreted as} a non strictly increasing (resp. decreasing) to something.\\

\noindent Consequently, we begin by the proof of Point (4).\\

\noindent {\bf Proof of (4)}. Let $a$ and $b$ be real numbers such that  $a \leq b$, then  $(X\leq a) \subset (X\leq b)$. Since probability measures  \index{probability measure} are non-decreasing, we get

$$
F_{X}(a)=\mathbb{P}(X\leq a)\leq \mathbb{P}(X\leq b)=F_{X}(b),
$$

\Bin \noindent and then

\begin{equation*}
\Delta_{a,b}F_X=F_X(b)-F_X(a) \geq 0,
\end{equation*}

\noindent \bigskip which is the targeted Formula \ref{propa01.crv.DefDelta01}.

\bigskip \noindent {\bf Proof of (1)}. The proof of \textbf{(1)} is immediate since for any $x\in \mathbb{R}$, $F_X(x)$ is probability value.\\

\bigskip \noindent {\bf Proof of  (2)}.\\

\noindent \textbf{First}, consider a sequence $t(n), n\geq 1$, that broadly decreases to  $-\infty$ as $n\uparrow +\infty$, that is $t(n)\downarrow -\infty$. Let $A_{t(n)}=(X\leq t(n))$. Then, the sequence of events $(A_{t(n)})$ broadly decreases to $\emptyset$, that is

$$
\bigcap_{n\geq 1} A_{t(n)}=\emptyset.
$$

\bigskip \noindent By the continuity of Probability measures  \index{probability measure} (See Properties (D) and (E) in Chapter \ref{proba01.pm}), we have

$$
0=\mathbb{P}(\emptyset)=\lim_{n\uparrow +\infty} \mathbb{P}(A_{t(n)}),
$$

\noindent that is, $F_{X}(t(n)) \downarrow 0$ for any sequence $t(n)\downarrow -\infty$ as $n\uparrow +\infty$. By the previous remark on monotone limits, this means  \index{mean} that

$$
\lim_{x\rightarrow -\infty}F_{X}(x)=0.
$$

\Bin
\noindent Similarly, consider a broadly increasing sequence $t(n), n\geq 1$, to $+\infty$ as $n\uparrow +\infty$, that is $t(n)\uparrow +\infty$. Let $A_{t(n)}=(X\leq t(n))$. Then, the sequence of events $(A_{t(n)})$ broadly increases to $\Omega$, that is

$$
\bigcup_{n\geq 1} A_{t(n)}=\Omega.
$$

\Bin
\noindent By the continuity of Probability measures  \index{probability measure} (See Properties (D) and (E) in Chapter \ref{proba01.pm}), we have

$$
1=\mathbb{P}(\Omega)=\lim_{n\uparrow +\infty} \mathbb{P}(A_{t(n)}),
$$

\noindent that is, $F_{X}(t(n)) \uparrow 1$ for any sequence $t(n)\uparrow +\infty$ as $n\uparrow +\infty$. By the previous remark on monotone limits, this means  \index{mean} that

$$
\lim_{x\rightarrow +\infty}F_{X}(x)=1.
$$

\Bin
\noindent The proof of Point (2) is over.\\

\bigskip \noindent {\bf Proof of 3}. Let us show that $F_{X}$ is right-continuous at any point $x \in \mathbb{R}$. Let $h(n),n\geq 1$ a sequence of real and non-negative numbers decreasing to $0$ in a broad sense, that is $h(n)\downarrow 0$ as $n\uparrow +\infty$. We clearly have (use a drawing if not clear) that

\begin{equation*}
\left(-\infty,x+h(n)\right] \downarrow \left(-\infty,x\right] \ \ as \ \  n\uparrow +\infty.
\end{equation*}

\Bin
\noindent By using the fact that an inverse image mapping  \index{inverse image mapping} $X^{-1}$ preserves sets operations, we get that, as $n\uparrow +\infty$,

$$
A_{x+h(n)}=(X \leq x+h(n))=X^{-1}(\left(-\infty,x+h(n)\right]) \downarrow X^{-1}(\left(-\infty,x\right])=(X\leq x)=A_{x},
$$

\Bin
\noindent By the continuity of Probability measures  \index{probability measure} (See Properties (D) and (E) in Chapter \ref{proba01.pm}), we have

$$
\mathbb{P}(A_{x+h(n)}) \downarrow \mathbb{P}(A_{x}),
$$

\Bin
\noindent that is $F_X(x+h(n)) \downarrow F_X(x)$, for any sequence $h(n)\downarrow 0$ as $n\uparrow +\infty$. From Calculus  Courses, this means  \index{mean} that

$$
F_X(x+h) \downarrow F_X(x) \ \ as \ \ 0<h \rightarrow 0
$$

\Bin
\noindent Hence $F_X$ is right-continuous at $x$.

\bigskip
\noindent The previous properties allow us to state a general definition of the notion of cumulative distribution function  \index{cumulative Distribution Function} and to unveil continuous
real-valued random variables.

\noindent Let us begin by the general definition of a  cumulative distribution function  \index{cumulative Distribution Function} on $\mathbb{R}$.\\

\noindent \textbf{Definition}. A real-valued function $F$ defined on $\mathbb{R}$ is a cumulative distribution function  \index{cumulative Distribution Function} if and only if \\

\noindent \textbf{(1)}  $0\leq F\leq 1.$\\

\noindent \textbf{(2)}  $\lim_{x\rightarrow -\infty}F(x)=0$ and  $\lim_{x\rightarrow +\infty}F(x)=1$\\

\noindent \textbf{(3)}  $F$ is right-continuous.\\

\noindent \textbf{(4)}  $F$ is $1$-non-decreasing.\\

\bigskip \noindent \textbf{Random Variable associated with a \textit{cdf}}.\\

\noindent Before we begin, let us make some remarks on the domain of a probability measure.  \index{probability measure} So far, we used discrete probability measures defined of the power set $\mathcal{P}(\Omega)$ of the sample set $\Omega$.\\

\noindent But it is possible to have a probability measure  \index{probability measure} that is defined only on a subclass $\mathcal{A}$ of $\mathcal{P}(\Omega)$. In this case, we have to ensure the following points for the needs of the computations.\\

\noindent (1) The subclass $\mathcal{A}$ of $\mathcal{P}(\Omega)$ should contain the impossible event $\emptyset$ and the sure event $\Omega$.\\

\noindent (2) Since we use the sets operations (union, intersection, difference of sets, complements of sets, etc.), The subclass $\mathcal{A}$ of $\mathcal{P}(\Omega)$ should be stable under sets operations, but only when they are applied at most a countable number of times.\\

\noindent (3) For a real-valued random variable $X$,  we need to define its \textit{cdf} $F_X(x)$, $x\in \mathbb{R}$, we have to ensure that the events $(X\leq x)$ lie in $\mathcal{A}$ for all $x\in \mathbb{R}$.\\

\noindent The first two points hold for discrete probability measures  \index{probability measure} when $\mathcal{A}=\mathcal{P}(\Omega)$ and the third point holds for a discrete real-valued random variable when $\mathcal{A}=\mathcal{P}(\Omega)$.\\

\noindent In general, a subclass $\mathcal{A}$ of $\mathcal{P}(\Omega)$ fulfilling the first two points is called a \textit{fields of subsets}
or a \textit{$\sigma$-algebra of subsets} of $\Omega$. A real-value application $X$ defined on  $\Omega$ such that events $(X\leq x)$ lie in $\mathcal{A}$ for all $x\in \mathbb{R}$, is called a \textit{measurable} application, or a random variable, with respect to $\mathcal{A}$.\\

\noindent \textbf{The consequence is that, you may see the symbol $\mathcal{A}$ and read about measurable application in the rest of this chapter and subsequent chapters. But, we do not pay attention to them. We only focus on the probability computations, assuming everything is well defined and works well, at least at this level. We will have time and space for such things in the mathematical course on Probability Theory \index{probability theory} }.\\

\noindent Now, it is time to introduce our definitions.\\

\bigskip \noindent \textbf{Definition}. Given a \textit{cdf} $F$, we say that a random variable $X$ defined on some probability space  \index{probability space} $(\Omega, \mathcal{A}, \mathbb{P})$ possesses
(or, is associated to) the \textit{cdf} $F$, if and only if

\begin{equation}
\mathbb{P}(X \leq x)=F(x), \ \ for \ all \ \ x\in \mathbb{R}, \label{proba01.rcv.Kolmogorov01}
\end{equation}

\Bin
\noindent that is $F_X=F$. Each particular \textit{cdf} is associated to a name of probability law  \index{probability law} that applies both to the \textit{cdf} and to the \textit{rrv}.\\

\bigskip \noindent \textbf{Immediate examples}\\

\noindent Let us give the following examples.\\

\noindent \textbf{Exponential law with parameter $\lambda>0$, denoted $\mathcal{E}(\lambda)$}.\\

\noindent A \textit{rrv} $X$ follows an Exponential law with parameter $\lambda >0$, denoted $X \sim \mathcal{E}(\lambda)$ if and only if its \textit{cdf} is given by

\begin{equation*}
F(x)=\left\{
\begin{array}{c}
1-e^{-\lambda x}\text{ \ \ \ }x\geq 0 \\
0\text{ \ \ \ \ \ \ \ \ \ \ \ \ \ }x<0.
\end{array}
\right.
\end{equation*}

\bigskip \noindent \textbf{Uniform law on $(a,b)$, denoted $\mathcal{U}(a,b)$, $a<b$}.\\

\noindent A \textit{rrv} $X$ follows a uniform law on $[a,b]$, $a<b$, denoted $X \sim \mathcal{U}(a,b)$, if and only if its \textit{cdf} is given by\\

$F(x)=\left\{
\begin{tabular}{lll}
$1$ & $if$ & $x>b$ \\
$(x-a)/(b-a)$ & $if$ & $a\leq x\leq b$ \\
$0$ & $if$ & $x<a.$
\end{tabular}%
.\right. $ \ \ \ \

\bigskip \noindent \textbf{These two functions} are clearly \textit{cdf}'s. Besides, they are continuous \textit{cdf}'s.\\

\noindent Let us begin by definitions and examples.\\

\noindent Let $F_{X}$ be the \textit{cdf} of a \textit{rrv} $X$. Suppose that $F_X$ satisfies :

\noindent \textbf{(a)}
\begin{equation*}
\text{ \ \ \ \ }f_{X}=F_{X}^{\prime }=\frac{dF_{X}}{dx}\text{ exists on }\mathbb{R}
\end{equation*}

\noindent and \\

\noindent \textbf{(b)}
\begin{equation*}
f_{X}=F_{X}^{\prime }\text{ is integrable on intervals }[s,t], \ s<t.
\end{equation*}

\bigskip \noindent Then for any $b<x\in \mathbb{R}$

\begin{equation*}
F_{X}(x)-F_{X}(b)=\int_{b}^{x}F_{X}^{\prime}(t)dt, \ \ x \in  \mathbb{R}.
\end{equation*}

\Bin
\noindent When $b\downarrow -\infty$, we get for  $x\in \mathbb{R}$,

\begin{equation}
F_{X}(x)=\int_{-\infty }^{x}f_{X}(t)\text{ }dt.  \label{proba01.crv.densite01}
\end{equation}

\noindent We remark that $f_{X}$ is non-negative as the derivative of a non-decreasing function. Next, by applying Formula
(\ref{proba01.crv.densite01}) to $+\infty$, we get that

$$
1=\int_{-\infty }^{+\infty}f_{X}(t)dt.
$$

\Bin
\noindent In summary, under our assumptions, there exists a function $f_X$ satisfying

\begin{equation}
f_X(x) \geq 0, \text{for all } x\in \mathbb{R}, \label{proba01.crv.densite02}
\end{equation}

\Bin
\noindent and

\begin{equation}
\int_{-\infty }^{+\infty}f_{X}(t)dt=1.  \label{proba01.crv.densite03}
\end{equation}

\bigskip \noindent We say that $X$ is \textbf{an absolutely continuous} random variable and that $f_X$ is its \textbf{probability density function}  \index{probability density function} (\textit{pdf}). This leads us to the following definition.\\

\noindent \textbf{Definition}. A function $f:\mathbb{R}\rightarrow \mathbb{R}$ is a probability density function  \index{probability density function} (\textit{pdf}) if and only if  Formulas (\ref{proba01.crv.densite02}) and (\ref{proba01.crv.densite03}) hold.

\bigskip \noindent \textbf{Definition}. A \textit{rrv} is said to have an absolutely continuous probability law  \index{probability law} if and only if there exists a \textit{pdf} denoted by $f_X$ such that (\ref{proba01.crv.densite01}) holds.\\

\noindent The name \textit{absolutely continuity} is applied to the probability law  \index{probability law} and to the random variable itself and to the \textit{cdf}.\\

\bigskip \noindent \textbf{Definition}: The \textbf{domain} or the \textbf{support} of an absolutely continuous random variable  \index{continuous Random Variable} is defined by

\begin{equation}
\mathcal{V}_X=\{x \in \mathbb{R}, f_X(x) \neq 0\}. \label{proba01.crv.densite04}
\end{equation}

\Bin \textbf{Remark for advanced readers}. Actually, we need to take the closure of that set. But, at this earlier stage, we simplify things and fortunately, such simplifications do not impact the type of problems we treat here.

\noindent We saw that knowing the values set of a discrete \textit{rvv} is very important for the computations. So is it for absolutely continuous \textit{rvv}.\\

\noindent \textbf{Examples}. Let us return to the examples of the Exponential law and the uniform law. By deriving the \textbf{cdf}, we see that :\\

\noindent An exponential \textit{rvv} with parameter $\lambda >0$, $X \sim \mathcal{E}(\lambda)$, is absolutely continuous with \textit{pdf}

\begin{equation*}
f(x)=\lambda exp(-\lambda x) 1_{(x \geq  0)},
\end{equation*}

\Bin
\noindent with support

$$
\mathcal{V}_X=\mathbb{R}_{+}=\{x \in \mathbb{R}, \  x\geq 0\}.
$$

\bigskip \noindent A uniform (or uniformly distributed) random variable on $[a,b]$, $a<b$, $X \sim \mathcal{U}(a,b)$, is absolutely continuous with \textit{pdf}

\begin{equation*}
f(x)=\frac{1}{b-a} 1_{(a\leq x\leq b)},
\end{equation*}

\Bin
\noindent and with domain

$$
\mathcal{V}_X=[a,b].
$$

\bigskip

\noindent In this subsection, we are going to present two important classes of distribution functions : the \textbf{Gamma} and the \textbf{Gaussian} classes. But we begin by showing a general way of finding \textit{pdf}'s.\\

\noindent \textbf{I - A general way of finding \textit{pdf}'s.}\\

\noindent A general method of finding a \textit{pdf} is to consider a non-negative function $g\geq 0$ which is integrable on $\mathbb{R}$ such that the integral on $%
\mathbb{R}$ is strictly positive, that is

\begin{equation*}
C=\int_{\mathbb{R}}g(x)dx>0.
\end{equation*}

\Bin
\noindent We get a \textit{pdf} $f\geq 0$ of the form%

\begin{equation*}
f(x)=C^{-1} \ g(x), \ x\in \mathbb{R}.
\end{equation*}

\bigskip
\noindent Here are some examples.\\

\bigskip \noindent \textbf{II - Gamma Law  \index{Gamma Law} $\gamma (a,b)$ with parameters $a>0$ and $b>0$}.\\

\noindent \textbf{ (a) Definition}. Let
\begin{equation*}
\Gamma (a)=\int_{0}^{+\infty }x^{a-1}e^{-x}dx.
\end{equation*}

\Bin
\noindent The function

\begin{equation*}
f_{\gamma (a,b)}(x)=\frac{b^{a}}{\Gamma (a)}x^{a-1}e^{-bx}1_{(x\geq 0)}
\end{equation*}

\bigskip \noindent is a \textit{pdf}. A non-negative random variable $X$ admitting $f_{\gamma (a,b)}$ as
a pdf is said to follow a Gamma law  \index{Gamma Law} $\gamma (a,b)$\ with parameter $a>0$ and $b>0,$ denoted as $X\sim \gamma (a,b)$.\\

\noindent \textbf{Important remark}. An exponential $\mathcal{E}(\lambda)$ random variable, with $\lambda >0$
is also a $\gamma (1,\lambda)$ \textit{rrv}.\\

\begin{figure}[htbp]
\centering
\includegraphics[width=0.90\textwidth]{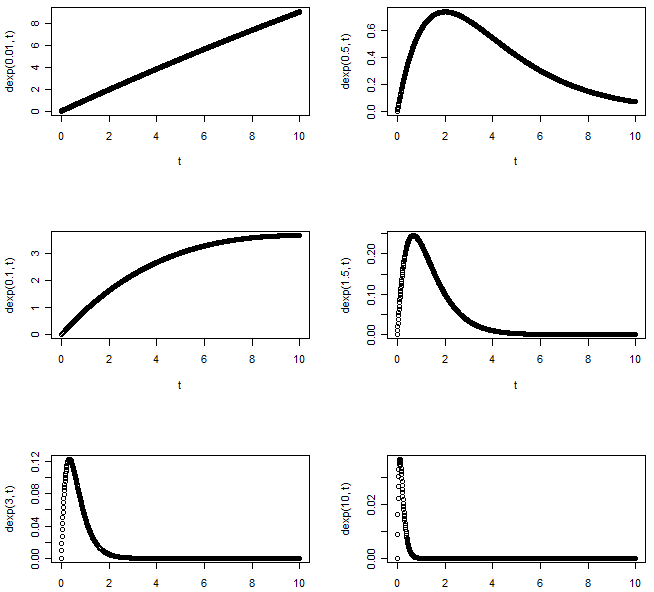}
\caption{Exponential Probability Density Function \index{probability density function} s}
\label{proba01.crv.fig2}
\end{figure}

\begin{figure}
\centering
\includegraphics[width=0.90\textwidth]{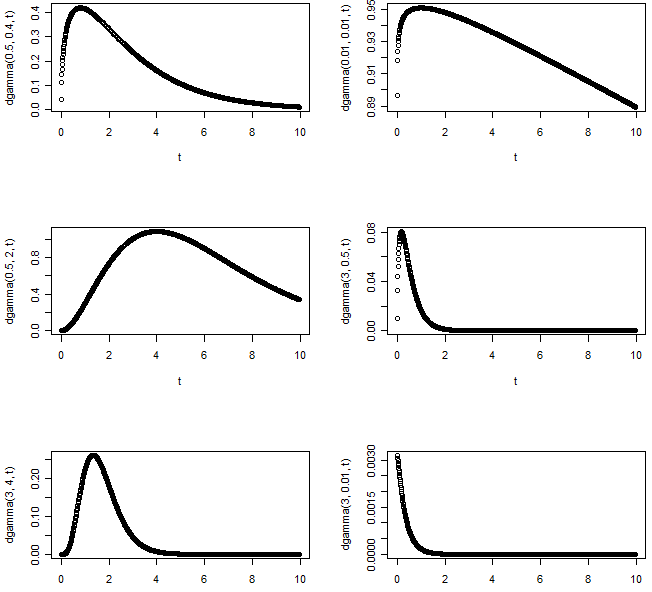}
\caption{$\gamma$-Probability Density Function \index{probability density function} s}
\label{proba01.crv.fig3}
\end{figure}

\noindent We have in Figure \ref{proba01.crv.fig2} graphical representations of probability density functions  \index{probability density function} depending of the parameter $\lambda$. In Figure
\ref{proba01.crv.fig3}, probability density functions  \index{probability density function} of $\gamma$ laws are illustrated.

\bigskip \noindent \textbf{Justification}. From Calculus courses, we may justify that the following function

\begin{equation*}
x^{a-1}e^{-bx}1_{(x>0)}
\end{equation*}

\Bin
\noindent is integrable on its domain $\mathbb{R}_{+}$ and we denote

\begin{equation}
\Gamma (a,b)=\int_{0}^{+\infty }x^{a-1}e^{-bx}dx>0.  \label{gammaAB}
\end{equation}

\Bin
\noindent We obtain a \textit{pdf}

\begin{equation*}
\frac{1}{\Gamma (a,b)}x^{a-1}e^{-bx}1_{(x>0)}.
\end{equation*}

\Bin
\noindent If $b=1$ in (\ref{gammaAB}), we denote

\begin{equation*}
\Gamma (a)=\Gamma (a,1)=\int_{0}^{+\infty }x^{a-1}e^{-x}dx.
\end{equation*}

\Bin
\noindent The function $\Gamma (\circ )$ satisfies the relation%

\begin{equation}
\forall (a>1),\text{ \ }\Gamma (a)=(a-1)\Gamma (a-1).  \label{gammaA01}
\end{equation}

\Bin \noindent This comes from the following partial integration

\begin{eqnarray*}
\Gamma (a) &=&\int_{0}^{+\infty }x^{a-1}d(-e^{-x}) \\
&=&\left[ -x^{a-1}e^{-x}\right] _{0}^{+\infty }+(a-1)\int_{0}^{+\infty
}x^{a-2}e^{-x}dx \\
&=&0+(a-1)\Gamma (a-1).
\end{eqnarray*}

\noindent We also remark that

\begin{equation}
\Gamma (1)=1,  \label{gammaA02}
\end{equation}

\Bin \noindent since

\begin{equation*}
\Gamma (1)=\int_{0}^{+\infty }e^{-x}dx=\int_{0}^{+\infty }d(-e^{-x})=\left[
-e^{-x}\right] _{0}^{+\infty }=1.
\end{equation*}

\Bin \noindent So, by applying (\ref{gammaA01}) to integers $a=n\in \mathbb{N}$, $n \geq 1$, and by taking (\ref{gammaA02}) into account, we get by induction that

\begin{equation*}
\forall (n\in \mathbb{N}\setminus \{0\}),\text{ }\Gamma (n)=(n-1)!.
\end{equation*}

\Bin \noindent Finally, we have
\begin{equation}
\forall (a>0,b>0),\text{ }\Gamma (a,b)=\frac{\Gamma (a)}{b^{a}}. \label{gammaAAB}
\end{equation}

\Bin \noindent This comes from the change of variables $u=bx$ in (\ref{gammaAB}) which gives%

\begin{eqnarray*}
\Gamma (a,b) &=&\int_{0}^{+\infty }\left( \frac{u}{b}\right)
^{a-1}e^{-u}\left( \frac{du}{b}\right)  \\
&=&\frac{1}{b^{a}}\int_{0}^{+\infty }u^{a-1}e^{-u}du \\
&=&\frac{\Gamma (a)}{b^{a}}.
\end{eqnarray*}

\noindent The proof is concluded by putting together (\ref{gammaAB}), and (\ref{gammaA01}), and (\ref{gammaA02}), and (\ref{gammaAAB}).\\

\bigskip \noindent \textbf{III - Gaussian law  \index{Gaussian law} $\mathcal{N}(m,\sigma ^{2})$ with parameters $m\in \mathbb{R}$ and $\sigma ^{2}>0$}
.\\

\noindent \textbf{(a) Definition}. A random variable $X$ admitting the \textit{pdf}

\begin{equation}
f(x)=\frac{1}{\sigma \sqrt{2\pi }}\exp \left( -\frac{(x-m)^{2}}{2\sigma ^{2}}%
\right), \ \ x\in \mathbb{R}, \label{proba01.crv.pdfGaussGen}
\end{equation}

\bigskip \noindent is said to follow a \textit{normal }or a \textit{Gaussian}  law of
parameters $m\in \mathbb{R}$ and $\sigma ^{2}>0,$ denoted by

$$
X\sim \mathcal{N}(m,\sigma ^{2})
$$

\begin{figure}[htbp]
\centering
\includegraphics[width=0.90\textwidth]{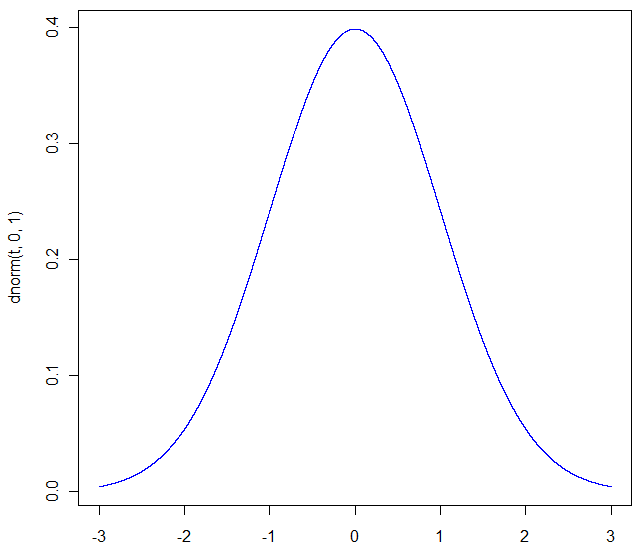}
\caption{Probability Density Function  \index{probability density function} of a Standard Gaussian Random Variable}
\label{proba01.cvr.fig4}
\end{figure}

\bigskip \noindent The probability density function  \index{probability density function} of the standard Gaussian Normal variable is illustrated in Figure \ref{proba01.cvr.fig4}.\\

\bigskip \noindent \textbf{(b) Justification}. The justification comes from the computation \bigskip of

\begin{equation*}
I=\int_{-\infty }^{+\infty }\exp \left( -\frac{(x-m)^{2}}{2\sigma ^{2}}%
\right) dx.
\end{equation*}

\noindent By the change of variable

\begin{equation*}
u=\frac{x-m}{\sigma \sqrt{2}},
\end{equation*}

\Bin \noindent we get

\begin{equation*}
I=\sigma \sqrt{2}\int_{-\infty }^{+\infty }e^{-u^{2}}du.
\end{equation*}

\Bin \noindent Advanced readers already know the value of $I$. Other readers, especially those in the first year of University, have to know that a forwarding Calculus course on multiple integration will show

\begin{equation*}
\int_{-\infty }^{+\infty }e^{-u^{2}}du=\sqrt{\pi},
\end{equation*}

\Bin
\noindent so that we arrive at

\begin{equation*}
I=\sigma \sqrt{2\pi }.
\end{equation*}

\noindent This justifies that
\begin{equation*}
\frac{1}{\sigma \sqrt{2\pi }}\exp \left( -\frac{(x-m)^{2}}{2\sigma ^{2}}%
\right) ,x\in \mathbb{R},
\end{equation*}

\Bin
\noindent is a \textit{pdf}.\\

\noindent But in this chapter, we will provide elementary methods based on the earlier works of \textit{de} Moivre (1732) and Laplace  \index{Laplace} (1801) to directly proof that

\begin{equation*}
J=\int_{-\infty}^{+\infty} \frac{1}{\sqrt{2\pi }}\exp \left( -t^{2}/2\right) \ dt=1,
\end{equation*}

\Bin \noindent which proves that (\ref{proba01.crv.pdfGaussGen}) is a \textit{pdf} for $m=0$ and $\sigma=1$. For proving the general case of arbitrary $m$ and $\sigma^2 >0$, we can use the change of variable $u=(x-m)/\sigma$ to see that $I=\sigma \sqrt{1\pi}$.\\

\bigskip \noindent \textbf{Commentary}. This Gaussian law  \index{Gaussian law} is one of the most important probability laws  \index{probability law} both in Probability Theory  \index{probability theory} and in Statistics.  \index{statistics} It will be a key tool in the next courses on Probability Theory and Mathematical Statistics.\\

\noindent By now, we want to devoted a full section on some historical facts and exhibit the first central limit theorem through the \textit{de} Moivre and Laplace  \index{Laplace} theorems.\\

\bigskip \noindent \textbf{Warning}. The reader is not obliged to read this section. The results of this section are available in modern books with elegant proofs. But, those who want acquire a deep expertise in this field are invited to read such a section because of its historical
aspects. They are invited to discover the great works of former mathematicians.\\

\noindent A very important example of absolutely continuous \textit{cdf} is the Gaussian \textit{rrv}. There are many ways to present it. Here, in this first course, we want to show the reader the very historical way to derive it from the Binomial law, through the \textit{de} Moivre and Laplace  \index{Laplace} theorems (1732 - 1801).\\

\noindent  Here, the standard Gaussian law  \index{Gaussian law} with $m=0$ and $\sigma=1$ is derived as an approximation of centered and normalized Binomial random variable \index{binomial random variable} s.\\

\noindent Let $(X_{n})_{n\geq 1}$ be a sequence of Bernoulli random variables  \index{Bernoulli random variable} of respective parameters $n$ and $p$, $\ $that is, for each $n\geq 1$, $X_{n}\sim \mathcal{B}(n,p)$. We recall the following parameters of such \textit{rrv}'s :

\begin{equation*}
\mathbb{E}(X_{n})=np\text{ and }Var(X_{n})=npq.
\end{equation*}

\Bin
\noindent Put

\begin{equation*}
Z_{n}=\frac{X_{n}\text{ }-\text{ }np}{\sqrt{npq}}.
\end{equation*}

\Bin
\noindent We have $\mathbb{E}(Z_{n})=0$ and  $Var(Z_{n})=1$. Here are historical results.

\bigskip

\begin{theorem} \noindent (\textbf{\textit{de} Moivre, 1732, see  \cite{loeve} \index{Lo\`eve}, page 23}) \label{proba01.crv.LGM}.\\

\noindent The equivalence

\begin{equation*}
\mathbb{P}(X_{n}=j)\sim \frac{e^{-\frac{x^{2}}{2}}}{\sqrt{2^{\pi }npq}},
\end{equation*}

\noindent holds uniformly in

\begin{equation*}
x=\frac{(j\text{ }-\text{ }np)}{\sqrt{2^{\pi }npq}}\in \left[ a,\text{ }b%
\right] ,a<b.
\end{equation*}

\noindent as $n\rightarrow \infty .$
\end{theorem}

\bigskip \noindent Next, we have:\\

\bigskip

\begin{theorem} (Laplace  \index{Laplace} 1801, see \cite{loeve} \index{Lo\`eve}, page 23) \label{proba01.crv.LAP}.\\

\noindent (1) For any $b<a$, we have
\begin{equation*}
\mathbb{P}(b\leq Z_{n}\leq a)\rightarrow \int_{b}^{a}\frac{1}{\sqrt{2\pi }}%
\text{ }e^{-\frac{x^{2}}{2}}dx.
\end{equation*}

\noindent (2) For any $x\in \mathbb{R}$, we have
\begin{equation*}
\mathbb{P}(Z_{n}\leq x)\rightarrow \int_{-\infty }^{x}\frac{1}{\sqrt{2\pi }}%
\text{ }e^{-\frac{x^{2}}{2}}dx.
\end{equation*}

\noindent (3) The function
\begin{equation*}
\frac{1}{\sqrt{2\pi }}\text{ }e^{-\frac{x^{2}}{2}}dx,\ x\in \mathbb{R},
\end{equation*}

\Bin
\noindent is a probability distribution function.
\end{theorem}

\bigskip \noindent \textbf{Proofs of the results}.\\

\noindent \textbf{(A) Proof of Theorem \ref{proba01.crv.LGM}}.\\

\noindent At the beginning, let us recall some needed tools.\\

\noindent \textbf{(i)} $f_{n}(x)\rightarrow f(x)$ uniformly in $x\in A$, as $%
n\rightarrow \infty $, means  \index{mean} that

\begin{equation*}
\lim_{n\rightarrow \infty }\sup_{x\in A}\text{ }\left| f_{n}(x)-f(x)\right|
=0.
\end{equation*}

\Bin
\noindent \textbf{(ii)} $f_{n}(x)\sim f(x)$, as $n\rightarrow \infty $, for
a fixed $x,$ mean \index{mean} s

\begin{equation*}
\lim_{n\rightarrow \infty }\text{ }\frac{f_{n}(x)}{f(x)}=1.
\end{equation*}

\Bin
\noindent \textbf{(iii)} $f_{n}(x)\sim f(x)$ uniformly in $x\in A$, as $%
n\rightarrow \infty $, means  \index{mean} that

\begin{equation*}
\lim_{n\rightarrow \infty }\sup_{x\in A}\text{ }\left| \left( \frac{f_{n}(x)%
}{f(x)}\right) -1\right| =0.
\end{equation*}

\bigskip \noindent \textbf{(iv)} $a_{n}=0(b_{n})$ as  $n\rightarrow \infty $, mean \index{mean} s
that

\begin{equation*}
\exists \text{ \ }M\text{ }>\text{ }0,\text{ }\exists \text{ \ }N,\text{ \ }%
\forall \text{ }n\geq N,\text{ }\left| \frac{a_{n}}{b_{n}}\right| \text{ }<%
\text{ }M.
\end{equation*}

\Bin
\noindent \textbf{(v)} $a_{n}(x)=0(b_{n}(x))$ uniformly in  $x\in A,$ as $%
n\rightarrow \infty ,$ means  \index{mean} that

\begin{equation*}
\exists \text{ \ }M\text{ }>\text{ }0,\text{ }\exists \text{ \ }N,\text{ }%
\forall \text{ }n\geq N,\text{ }\forall \text{ }x\in A,\text{ }\left| \frac{%
a_{n}(x)}{b_{n}(x)}\right| \text{ }<\text{ }M.
\end{equation*}

\Bin
\noindent \textbf{(vi)} We will need the Sterling Formula we studied in Chapter \ref{proba01.combin}:

\begin{equation*}
n!=\sqrt{2 \pi n}\text{ }n^{n}\text{ }e^{-n}\text{ }e^{\theta _{n}},
\end{equation*}

\Bin
\noindent where for any $\eta>0$, we have for $n$ large enough,

$$
\left| \theta _{n}\right| \leq \frac{1+\eta}{(12n)}.
$$

\Bin
\noindent We begin the proof by denoting  $k=n-j,$ $q=1-p,$ and by writing

\begin{equation}
S_{n}(j)=\mathbb{P}(X_{n}=j)=\left( \frac{1}{\sqrt{2\pi }}\right) \sqrt{(%
\frac{n}{jk})}\left( \frac{np}{j}\right) ^{j}\text{ }\left( \frac{nq}{k}%
\right) ^{k}\text{ }e^{\theta _{n}}e^{-\theta _{j}}e^{-\theta _{k}}.
\label{gauss1}
\end{equation}

\Bin
\noindent Hence

\begin{equation}
x=\frac{(j-np)}{\sqrt{2\pi npq}}\in \left[ a,\text{ }b\right]
\Leftrightarrow np-a\sqrt{npq}\leq j\leq np-b\sqrt{npq}  \label{gauss2}
\end{equation}

\begin{equation}
\Leftrightarrow nq-b\sqrt{npq}\leq k\leq nq-a\sqrt{npq}.  \label{gauss3}
\end{equation}

\bigskip \noindent Hence  $j$ and  $k$ tend to $+\infty $ uniformly in  $x=\frac{(j-np)%
}{\sqrt{2\pi npq}}\in \left[ a,\text{ }b\right] .$ Then  $e^{-\theta
_{j}}\sim 1$ and $e^{-\theta _{k}}\sim 1$, uniformly in  $x=\frac{(j-np)}{%
\sqrt{2\pi npq}}\in \left[ a,\text{ }b\right] $ as $n\rightarrow \infty .$
Let us use the following second order expansion of $\log (1+u):$

\begin{equation}
\log (1+u)=u-\frac{1}{2}u^{2}+O(u^{3}),\text{ as }u\rightarrow 0.
\label{gauss4}
\end{equation}

\bigskip \noindent Then, there exist a constant $C>0$  and a number $u_{0}$
$>$ $0$ such that

\begin{equation*}
\forall (0\leq u\leq u_{0}),\text{ }\left\vert O(u^{3})\right\vert \leq C\
u^{3}.
\end{equation*}

\Bin \noindent We have,

\begin{equation*}
\frac{j}{(np)}=\left( 1+x\sqrt{\frac{q}{(np)}}\right)
\end{equation*}

\begin{equation*}
\left\vert x\sqrt{\frac{q}{(np)}}\right\vert \leq \max (\left\vert
a\right\vert ,\left\vert b\right\vert )\sqrt{\frac{q}{(np)}}\rightarrow 0,
\end{equation*}

\Bin
\noindent uniformly in  $\ x\in \left[ a,\text{ }b\right] .$ Combining this with (\ref%
{gauss4}) leads to

\begin{equation*}
\left( \frac{np}{j}\right) ^{j}=\exp \left( -\text{ }j\text{ }\log (\frac{j}{%
np})\right) .\text{ \ \ \ }
\end{equation*}

\Bin
\noindent But we also have

\begin{equation*}
\log \left( \frac{j}{np}\right) =\log \left( 1+x\sqrt{\frac{q}{(np)}}\right)
\end{equation*}
\begin{equation}
=x\sqrt{\frac{q}{(np)}}-x^{2}\text{ }\frac{q}{\left[ 2(np)\right] }+O(n^{-%
\frac{1}{2}}).  \label{gauss6}
\end{equation}

\bigskip \noindent We remark that the big $O,$ which does not depend on $n$
nor on  $x,$ is uniformly bounded. By using  $j=np+x\sqrt{\frac{q}{(np)}}$
in (\ref{gauss6}), and by operating the product, we get

\begin{equation*}
j\text{ }\log \left( \frac{j}{np}\right) =x\sqrt{npq}+x^{2}\frac{q}{2}+0(n^{-%
\frac{1}{2}})
\end{equation*}

\Bin
\noindent and then,

\begin{equation}
\left( \frac{np}{j}\right) ^{j}=\exp \left( \text{ }-\text{ }\sqrt{npq}-x^{2}%
\frac{q}{2}+0(n^{-\frac{1}{2}})\right) .  \label{gauss7}
\end{equation}

\bigskip \noindent We just proved that

\begin{equation}
\left( \frac{np}{k}\right) ^{j}=\exp \left( \text{ }+\sqrt{npq}-x^{2}\frac{p%
}{2}+0(n^{-\frac{1}{2}})\right) .  \label{gauss8}
\end{equation}

\Bin
\noindent The product gives
\begin{equation}
\left( \frac{np}{j}\right) ^{j}\left( \frac{np}{k}\right) ^{j}=\exp \left(
\text{ }-\frac{1}{2}x^{2}\right) \exp \left( \text{ }0(n^{-\frac{1}{2}%
})\right) ,  \label{gauss78}
\end{equation}

\Bin
\noindent where the big $0(n^{-\frac{1}{2}})$ uniformly holds in
\begin{equation*}
x=\frac{(j\text{ }-\text{ }np)}{\sqrt{2^{\pi }npq}}\in \left[ a,\text{ }b%
\right] ,a<b.
\end{equation*}

\bigskip \noindent The proof is finished by combining (\ref{gauss78}) and
consequences of  (\ref{gauss2}) and  (\ref{gauss3}).

\bigskip \noindent \textbf{Proof of Theorem \ref{proba01.crv.LAP}}.\\

\noindent Let  $j_{0}=np+a\sqrt{npq}.$ Let  $j_{h},$ $h=1,...,m_{n-1},$ be the
intergers between  $j_{0}$ and $np+b\sqrt{npq}=j_{m_{n}}$ and denote
\begin{equation*}
x_{h}=\frac{(j_{h}-np)}{\sqrt{npq}}.
\end{equation*}

\noindent It is immediate that

\begin{equation}
\max_{0\leq h\leq m_{n-1}}\text{ }j_{h+1}-j_{h}\leq \frac{1}{\sqrt{npq}}.%
\text{ }  \label{gauss9}
\end{equation}

\bigskip \noindent By the Laplace-Moivre-Gauss  \index{Laplace} Theorem \ref{proba01.crv.LGM}, we
have

\begin{equation}
\mathbb{P}(X_{n}=j_{h})=\frac{e\frac{-x_{h}^{z}}{2}}{\sqrt{2\pi npq}}%
(1+\varepsilon _{n}(h)),\text{ }h=1,......,m_{n-1}  \label{gauss10}
\end{equation}

\Bin
\noindent with

\begin{equation}
\underset{1\leq h\leq m_{n-1}}{\max }\left\vert \varepsilon
_{n}(h)\right\vert =\varepsilon \rightarrow 0  \label{gauss10a}
\end{equation}

\Bin
\noindent as  $n\rightarrow +\infty.$ But,

\begin{equation*}
\mathbb{P}(a<Z_{n}<b)=\sum_{1\leq h\leq m_{n-1}}\mathbb{P}(X_{n}=j_{h}).\
\end{equation*}
\begin{equation}
=\sum_{1\leq h\leq m_{n-1}}\frac{e\frac{-x_{h}^{z}}{2}}{\sqrt{2\pi npq}}%
+\sum_{1\leq h\leq m_{n-1}}\frac{e\frac{-x_{h}^{z}}{2}}{\sqrt{2\pi npq}}%
\text{ }\varepsilon _{n}(h).  \label{gauss11}
\end{equation}
\begin{equation*}
=\sum_{1\leq h\leq m_{n-1}}\frac{e\frac{-x_{h}^{z}}{2}}{\sqrt{2\pi npq}}%
(1+\varepsilon _{n}(h))
\end{equation*}

\noindent Also, we have

\begin{equation*}
\underset{n\rightarrow \infty }{\lim }\ \sum_{1\leq h\leq m_{n-1}}\frac{e%
\frac{-x_{h}^{z}}{2}}{\sqrt{2\pi npq}}=\sum_{0\leq h\leq m_{n}}\frac{e%
\frac{-x_{h}^{z}}{2}}{\sqrt{2\pi npq}}
\end{equation*}

\bigskip \noindent since each of the two terms that were added converges to 0. The latter expression is a Riemann sum on $\left[ a,\text{ }b\right] $ based on a
uniform subdivision and associated to the function  $x\mapsto e^{-\frac{x^{2}%
}{2}}.$ This function is continuous on $[a,b],$ and therefore, is integrable
on $[a,b].$ Therefore the Riemann sums  \index{Riemann sums} converge to the integral

\begin{equation*}
\frac{1}{\sqrt{2\pi }}\int_{a}^{b}e^{-\frac{x^{2}}{2}}dx.
\end{equation*}

\Bin
\noindent This means  \index{mean} that the first term of (\ref{gauss11}) satisfies

\begin{equation}
\sum_{1\leq h\leq m_{n-1}}\frac{e\frac{-x_{h}^{z}}{2}}{\sqrt{2^{\pi }npq}}%
\rightarrow \left( \frac{1}{\sqrt{2^{\pi }}}\right) \int_{a}^{b}e^{-\frac{%
x^{2}}{2}}dx.  \label{gauss12}
\end{equation}

\Bin
\noindent The second term (\ref{gauss11}) satisfies

\begin{equation}
0\leq \sum_{1\leq h\leq m_{n-1}}\frac{e^{-\frac{x_{h}^{2}}{2}}}{\sqrt{2^{\pi
}npq}}\varepsilon _{n}(h)\leq \varepsilon _{n}\sum_{1\leq h\leq m_{n-1}}%
\frac{e^{-\frac{x_{h}^{2}}{2}}}{\sqrt{2^{\pi }npq}}.  \label{gauss13}
\end{equation}

\noindent Because, the sequence
\begin{equation*}
\varepsilon _{n}=\sup \{\left\vert \varepsilon _{n}(h)\right\vert ,0\leq
h\leq m_{n}\}
\end{equation*}

\noindent goes to zero based on (\ref{gauss10a}) above, and
of (\ref{gauss12}), it comes that second term of (\ref{gauss11}) tends to
zero. This establishes Point (1) of the theorem. Each of the signs $<$ and $>
$ may be replaced $\leq $ or $\geq .$ The effect in the proof would result
in adding of at most two  indices \ in the sequence $i_{1},...,i_{m_{n}-1}$
( $j_{0}$ and/or $j_{m_{n}})$ and the corresponding terms of this or these
two indices in (\ref{gauss11})\ tend to zero and the result remains valid.\\

\Bin
\noindent Let us move to the two other points. We have
\begin{equation*}
\underset{t\rightarrow \infty }{\lim }t^{2}e^{-\frac{t^{2}}{2}}=0.
\end{equation*}

\Bin
\noindent Thus, there exists $A>0$ such that
\begin{equation}
\left\vert t\right\vert \geq A\Rightarrow 0\leq e^{-\frac{t^{2}}{2}}<\frac{1%
}{t^{2}}.  \label{gauss001}
\end{equation}

\Bin
\noindent Then, for $a>A$ and $b<-A$

\begin{equation*}
\int_{b}^{a}e^{-\frac{t^{2}}{2}}dt\leq \int_{-A}^{A}e^{-\frac{t^{2}}{2}%
}dt+\int_{b}^{-A}e^{-\frac{t^{2}}{2}}dt+\int_{A}^{a}e^{-\frac{t^{2}}{2}}dt
\end{equation*}

\begin{eqnarray*}
\  &\leq &\int_{-A}^{A}e^{-\frac{t^{2}}{2}}dt\ +\int_{-\infty }^{-A}e^{-%
\frac{t^{2}}{2}}dt+\int_{A}^{+\infty }e^{-\frac{t^{2}}{2}}dt \\
&\leq &\int_{-A}^{A}e^{-\frac{t^{2}}{2}}dt\ +2\int_{A}^{+\infty }e^{-\frac{%
t^{2}}{2}}dt \\
&\leq &\int_{-A}^{A}e^{-\frac{t^{2}}{2}}dt\ +2\int_{A}^{+\infty }\frac{1}{%
t^{2}}dt\text{ \ \ (By \ref{gauss001})} \\
&=&\int_{-A}^{A}e^{-\frac{t^{2}}{2}}dt\ +\frac{2}{A}.
\end{eqnarray*}

\Bin \noindent The inequality

\begin{equation}
\int_{b}^{a}e^{-\frac{t^{2}}{2}}dt\leq \int_{-A}^{A}e^{-\frac{t^{2}}{2}}dt\ +%
\frac{2}{A},  \label{gauss002}
\end{equation}

\noindent remains true if $a\leq A$ or $b\geq -A$. Thus, Formula (\ref{gauss002}) holds for any
real numbers $a$ and $b$ with $a<b$. Since

\begin{equation*}
\int_{b}^{a}e^{-\frac{t^{2}}{2}}dt
\end{equation*}

\noindent is increasing as $a\uparrow \infty $ and $b\downarrow -\infty ,$ it comes
from Calculus courses on limits that its limit as $a\rightarrow \infty $ and $%
b\rightarrow -\infty $, exists if and only if it is bounded. This is the
case with \ref{gauss002}. Then

\begin{equation*}
\lim_{a\rightarrow \infty ,b\rightarrow -\infty }\int_{b}^{a}e^{-\frac{t^{2}%
}{2}}dt=\int_{-\infty }^{\infty }e^{-\frac{t^{2}}{2}}dt\in \mathbb{R}\
\end{equation*}

\noindent We denote

\begin{equation*}
G(x)=\int_{-\infty }^{x}e^{-\frac{t^{2}}{2}}dt\in \mathbb{R}
\end{equation*}

\Bin \noindent and

\begin{equation*}
F_{n}(x)=\mathbb{P}(Z_{n}\leq x),x\in \mathbb{R}.
\end{equation*}

\bigskip
\noindent It is clear that $G(x)\rightarrow 0$ as $x\rightarrow -\infty $. \
Thus, for any  $\varepsilon >0$, there exists $b_{1}<0$ such that

\begin{equation*}
\forall \text{ }b\leq b_{1},\text{ }G(b)\leq \frac{\varepsilon }{3}
\end{equation*}

\Bin \noindent For any $b<0$, we apply the Markov's inequality,  \index{Markov's inequality} to have

\begin{equation*}
\mathbb{P}(Z_{n}\leq b)\leq \mathbb{P}\left( \left\vert Z_{n}\right\vert
\geq \left\vert b\right\vert \right) \leq \frac{\mathbb{E}\left\vert
Z_{n}\right\vert ^{2}}{\left\vert b\right\vert ^{2}}=\frac{1}{\left\vert
b\right\vert ^{2}},
\end{equation*}

\Bin
\noindent by the fact that $\mathbb{E}Z_{n}=0,$%

\begin{equation*}
\mathbb{E}\left\vert Z_{n}\right\vert ^{2}=\mathbb{E}Z_{n}^{2}=var(Z_{n})=1.
\end{equation*}

\noindent Then, for any  $\varepsilon >0$, there exist $b_{2}<b_{1}$ such that

\begin{equation*}
\forall \text{ }b\leq b_{2},\text{ }F_{n}(b)=\mathbb{P}(Z_{n}\leq b)\leq
\frac{\varepsilon }{3}
\end{equation*}

\Bin \noindent By $Point$ (1), for any $a\in \mathbb{R},$ and for $b\leq \min (b_{2},a),$
there exists $n_{0}$ such that for $n\geq n_{0},$

\begin{equation*}
\left\vert (F_{n}(a)-F_{n}(b))-(G(a)-G(b))\right\vert \leq \frac{\varepsilon
}{3}\text{.}
\end{equation*}

\Bin
\noindent We combine the previous facts to get that for any $\varepsilon >0,$ for any
$a\in \mathbb{R}$, there exists $n_{0}$ such that for $n\geq n_{0},
$

\begin{eqnarray*}
\left\vert F_{n}(a)-G(a)\right\vert  &\leq &\left\vert
(F_{n}(a)-F_{n}(b_{0}))-(G(a)-G(b_{0}))\right\vert +\left\vert
F_{n}(b_{0})-G(b_{0})\right\vert  \\
&\leq &\frac{\varepsilon }{3}+\left\vert F_{n}(b_{0})\right\vert +\left\vert
G(b_{0})\right\vert \leq \frac{\varepsilon }{3}+\frac{\varepsilon }{3}+\frac{%
\varepsilon }{3}=\varepsilon \text{.}
\end{eqnarray*}

\Bin
\noindent We conclude that for any  $a\in \mathbb{R},$

\begin{equation*}
\mathbb{P}(Z_{n}\leq a)\rightarrow G(a)=\frac{1}{\sqrt{2^{\pi }}}%
\int_{-\infty }^{a}e^{-\frac{t^{2}}{2}}dt\text{. }
\end{equation*}

\Bin
\noindent This gives Point (2) of the theorem.\\

\bigskip \noindent The last thing to do to close the proof is to show Point (3) by establishing
that $G(\infty )=1$. By the Markov's inequality  \index{Markov's inequality} and the remarks done
above, we have for any  $a>0$,

\begin{equation*}
\mathbb{P}(Z_{n}>a)\leq \mathbb{P}(\left\vert Z_{n}\right\vert >a)<\frac{1}{%
a^{2}}
\end{equation*}

\Bin
\noindent Then for any $\varepsilon >0,$ there exists $a_{0}>0$ such that

\begin{equation*}
\mathbb{P}(Z_{n}>a)<\varepsilon \text{.}
\end{equation*}

\Bin
\noindent This is equivalent to saying that there exists $a_{0}>0$ such that for any $a\geq
a_{0}$ and for any  $n\geq 1,$

\begin{equation*}
1-\varepsilon \leq F_{n}(a)\leq 1\text{. }
\end{equation*}

\Bin
\noindent By letting $n\rightarrow \infty $ and by using Point (2), we get for any $%
a\geq a_{0},$

\begin{equation*}
1-\varepsilon \leq G(a)\leq 1\text{.}
\end{equation*}

\Bin
\noindent By the monotonicity of $G,$ we have that  $G(\infty )$ is the limit of $G(a)$
as $a\rightarrow \infty $. So, by letting $a\rightarrow \infty $ first and
next $\varepsilon \downarrow 0$ in the last formula, we arrive at

\begin{equation*}
1=G(+\infty )=\frac{1}{\sqrt{2^{\pi }}}\int_{-\infty }^{+\infty }e^{-\frac{%
t^{2}}{2}}dt.
\end{equation*}

\Bin \noindent This is Point (3).

\newpage

\Bin We need to make computations in solving problems based on probability theory.  \index{probability theory} Especially, when we use the cumulative distribution functions  \index{cumulative Distribution Function} $F_X$ of some real-valued random variable, we need some times to find the numerical value of $F(x)$ for a specific value of $x$. Sometimes, we need to know the quantile function (to be defined in the next lines). Some decades ago, we were using probability tables, that you can find in earlier versions of many books on probability theory, even in some modern books.\\

\noindent Fortunately, we now have free and powerful software packages for probability and statistics  \index{statistics} computations. One of them is the software \textbf{R}. You may find at \textit{www.r-project.org} the latest version of \textbf{R}. Such a software offer almost everything we need here.\\

\noindent Let us begin by the following problem.\\

\noindent \textbf{Example}. \label{lazyStudent} The final exam for a course, say of Geography, is composed by $n=20$ Multiple Choice Questions (MCQ). For each equation, $k$ answers are proposed and only  one them is true. The students are asked to answer each question by choosing one the $k$ proposed answers. At the end of the test, the grade of each student is the total number of correct answers. A student passes the exam if and only if he has at least the grade $x=10$ over $n$.\\

\noindent One of the students, say Jean, decides to answer at random. Let $X$ be his possible grade before the exam. Here $X$ follows a binomial law with parameters
$n=20$ and $p=1/4$. So the probability he passes is

$$
p_{pass}=\mathbb{P}(X\geq x)=\sum_{x\leq j \leq n}=\mathbb{P}(X=j).
$$

\noindent We will give the solution as applications of the material we are going to present. We begin by introducing the quantile function.\\

\noindent \textbf{A - Quantile function}.\\

\noindent Let $F : \mathbb{R} \rightarrow [0,1]$ be a cumulative distribution function.  \index{cumulative Distribution Function} Since $F$ is already non-decreasing, it is invertible if it is continuous and increasing. In that case, we may use the inverse function $F^{-1}$, defined by

$$
\forall s\in [0,1], \ \forall x\in \mathbb{R}, \ F(x)=s \Leftrightarrow x=F^{-1}(s).
$$

\Bin
\noindent In the general case, we may also define the generalized inverse

$$
\forall s\in [0,1], \ F^{-1}(s)=\inf \{x\in \mathbb{R}, \ F(x) \geq s \}.
$$

\Bin
\noindent we have the following interesting properties

\begin{eqnarray*}
&&\forall s\in [0,1],  \forall x\in \mathbb{R}, \ (F^{-1}(s)\leq x) \Leftrightarrow (s\leq F(x))\\
&&\forall s\in [0,1],  \forall x\in \mathbb{R}, \ (F^{-1}(s)> x) \Leftrightarrow (s > F(x))\\
&&\forall s\in [0,1],  \  F(F^{-1}(s)-0) \leq s \leq F(F^{-1}(s)),
\end{eqnarray*}

\noindent where

$$
F(F^{-1}(s)+0)=\lim_{h\downarrow 0} F(F^{-1}(s)-h).
$$

\Bin
\noindent The generalized inverse is called the quantile function of $F$. Remember that the quantile function is the inverse function of $F$ if $F$ is invertible.\\

\noindent \textbf{B - Installation of R}.\\

\begin{figure}[htbp]
\centering
\includegraphics[width=0.99\textwidth]{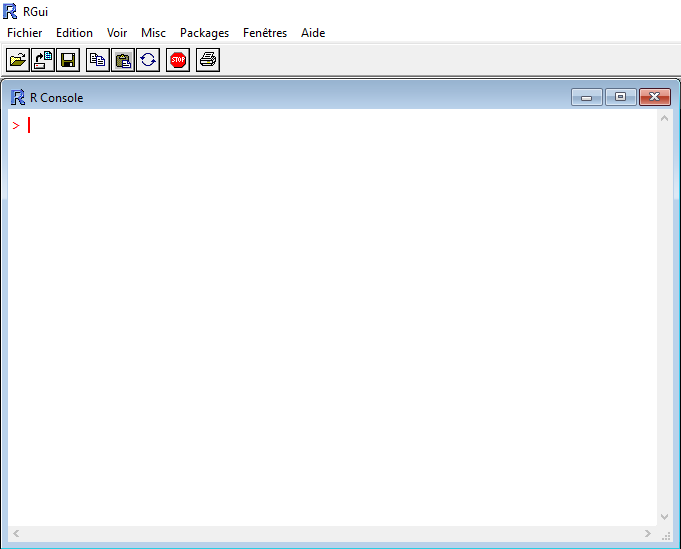}
\label{figWindow}
\end{figure}

\Bin Once you have installed \textbf{R} on your computer, and you launch, you will see the window given in Figure \ref{figWindow}. You will write your command and press the ENTER button of your keyboard. The result of the graph is displayed automatically.\\

\Bin \noindent \textbf{C - Numerical probabilities in R \index{Numerical probabilities in R} }.\\

In \textit{R}, each probability law  \index{probability law} has its name. Here are the most common ones in Table \ref{tabProbabR}.

\begin{table}
\centering
\begin{tabular}{lll}
\hline
Distribution & R name & additional arguments\\
\hline
beta &beta &shape1, shape2, ncp\\
\hline
binomial &binom &size, prob\\
\hline
Cauchy &cauchy &location, scale\\
\hline
chi-squared &chisq &df, ncp\\
\hline
exponential &exp &rate\\
\hline
F& f &df1, df2, ncp\\
\hline
gamma &gamma &shape, scale\\
\hline
geometric& geom &prob\\
\hline
hypergeometric &hyper& m, n, k\\
\hline
log-normal &lnorm &meanlog,  \index{mean} sdlog\\
\hline
logistic& logis &location, scale\\
\hline
\hline
negative binomial &nbinom& size, prob\\
\hline
normal &norm &mean,  \index{mean} sd\\
\hline
Poisson &pois &lambda\\
\hline
Student’s t& t &df, ncp\\
\hline
uniform& unif &min, max\\
\hline
Weibull &weibull& shape, scale\\
\hline
Wilcoxon &wilcox& m, n\\			
\end{tabular}
\label{tabProbabR}
\caption{Names of Probability laws  \index{probability law} in \textbf{R}}
\end{table}

\Bin
\noindent \textbf{How to use each probability law \index{probability law} ?}\\

\noindent In column 2, under \textbf{R name} of Table \ref{tabProbabR}, are listed the names of some random variables as they re used the Software \textbf{R}. In the last column, the names ate parameters are given.\\

\noindent Those names are used as follows :\\

\noindent (a) Add the letter \textbf{p} before the name to have the cumulative distribution function,  \index{cumulative Distribution Function} that is called as follows:

$$
pname(x,param1,param2,param3).
$$

\Bin \noindent (b) Add the letter \textbf{d} before the name to have the probability density function  \index{probability density function} (discrete or absolutely continuous), that is called as follows :

$$
dname(x,param1,param2,param3).
$$

\Bin \noindent (c) Add the letter \textbf{q} before the name to have the quantile function, that is called as follows

$$
qname(x,param1,param2,param3).
$$

\Bin \noindent \textbf{Example}. Let us use the Normal Random variable of parameters : mean  \index{mean} ($m=0$) and variance  \index{variance} ($sd=1$).\\

\noindent For $X \sim \mathcal{N}(0,1)$, we computed
$$
p=\mathbb{P}(X \leq 1.96)
$$

\begin{figure}
\centering
\includegraphics[width=0.99\textwidth]{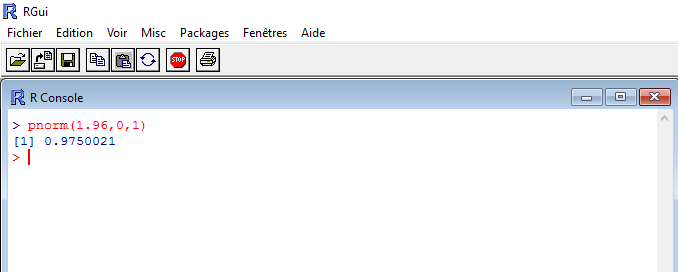}
\label{figNormal196}
\end{figure}

\noindent We may read in Figure \ref{figNormal196} the value : $p=0.975$.

\bigskip \noindent We may do more and get more values in Table \ref{tabNorm01}.

\begin{table}[htbp]
\centering
\begin{tabular}{llll}
\hline \hline
Probability 								& argument 	& R command 				& values\\
$\mathbb{P}(X \leq -2)$ 		& -2	 			& pnorm(-2,0,1)			&	 0.02275013		\\
$\mathbb{P}(X \leq -1.95)$ 	& -1.96 		& pnorm(-1.96,0,1)	&  0.02499790		\\
$\mathbb{P}(X \leq -0.8)$ 	& -0.8			& pnorm(-0.8,0,1)		&	 0.02275013		\\
$\mathbb{P}(X \leq 0)$ 			&	0				 	& pnorm(0,0,1)			&	 0.5					\\
$\mathbb{P}(X \leq 0.8)$ 		&	0.8			 	& pnorm(0.2,0,1)		&	 0.02275013		\\
$\mathbb{P}(X \leq -1.96)$ 	& 1.96			& pnorm(1.96,0,1)		&	 0.975				\\
$\mathbb{P}(X \leq 2)$ 			&	2			 		& pnorm(2,0,1)			&	 0.02275013	\\
\hline \hline
\end{tabular}
\label{tabNorm01}
\end{table}

\noindent \textbf{Exercise}. Use the quantile function applied to the probabilities above to find the arguments. For example, set $u=0.02275013$, and use the \textbf{R} command

$$
qnorm(u,0,1)
$$

\noindent to find again $x=-2$.\\

\noindent Follow that example to learn how to use all the probability law.  \index{probability law} Now, we return back to our binomial example (page : \pageref{lazyStudent}).\\

\noindent \textbf{D - Solution of the Lazy student problem using R}.\\

\noindent We already know that the probability that he passes by answering at random is

$$
p_{pass}=\mathbb{P}(X\geq x)=\sum_{x\leq j \leq n}=\mathbb{P}(X=j).
$$

\noindent The particular values are : $n=20$, $p=0.25$, $x=10$. From Table \ref{tabProbabR}, we see that the \textbf{R} name for the binomial
probability is \textit{binom} and its parameters are $n=size$  and $prob=p$. By using \textbf{R}, this probability is

\begin{eqnarray*}
p_{pass}&=&\mathbb{P}(X \geq 10)\\
&=&1- \mathbb{P}(X \leq 10)\\
&=&1-pbinom(9,20,0.25)= 0.01386442.
\end{eqnarray*}

\Bin \noindent  So, he has only 13 chances over one thousand to pass.\\

 %proba01.crv
\chapter{Appendix : Elements of Calculus} \label{proba01_appendix}

\noindent Calculus is fundamental to Probability Theory  \index{probability theory} and Statistics.  \index{statistics} Especially, the notions on limits in $\mathbb{R}$ are extremely important. The current section allows the reader to revise these notions and to complete his knowledge on this subject through exercises whose solutions are given in detail \index{tail} s.\\

\noindent \textbf{Definition}: $\ell \in \overline{\mathbb{R}}$ is an accumulation point  \index{accumulation point} of a sequence
$(x_{n})_{n\geq 0}$ of real numbers finite or infinite, in $\overline{\mathbb{R}}$, if and only if there exists a sub sequence $(x_{n(k)})_{k\geq 0}$ of
$(x_{n})_{n\geq 0}$ such that $%
x_{n(k)}$ converges to $\ell $, as $k\rightarrow +\infty $.\newline

\noindent \textbf{Exercise 1: } Set $y_{n}=\inf_{p\geq n}x_{p}$ and $%
z_{n}=\sup_{p\geq n}x_{p} $ for all $n\geq 0$. Show that :\newline

\noindent \textbf{(1)} $\forall n\geq 0,y_{n}\leq x_{n}\leq z_{n}$\newline

\noindent \textbf{(2)} Justify the existence of the limit of $y_{n}$ called limit inferior  \index{limit inferior} of the sequence $(x_{n})_{n\geq 0}$, denoted by $%
\liminf x_{n}$ or $\underline{\lim }$ $x_{n},$ and that it is equal to the following%

\begin{equation*}
\underline{\lim }\text{ }x_{n}=\lim \inf x_{n}=\sup_{n\geq 0}\inf_{p\geq n}x_{p}.
\end{equation*}

\Bin \noindent \textbf{(3)} Justify the existence of the limit of $z_{n}$ called limit superior  \index{limit superior} of the sequence $(x_{n})_{n\geq 0}$ denoted by $%
\lim \sup x_{n}$ or $\overline{\lim }$ $x_{n},$ and that it is equal%

\begin{equation*}
\overline{\lim }\text{ }x_{n}=\lim \sup x_{n}=\inf_{n\geq 0} \ \sup_{p\geq n}x_{p}.
\end{equation*}

\bigskip
\noindent \textbf{(4)} Establish that

\begin{equation*}
-\liminf x_{n}=\limsup (-x_{n})\noindent \text{ \ \ and \ }-\limsup
x_{n}=\liminf (-x_{n}).
\end{equation*}

\Bin  \noindent \textbf{(5)} Show that the limit superior  \index{limit superior} is sub-additive and the limit inferior  \index{limit inferior} is super-additive, i.e. :  for two sequences
$(s_{n})_{n\geq 0}$ and $(t_{n})_{n\geq 0}$
\begin{equation*}
\limsup (s_{n}+t_{n})\leq \limsup s_{n}+\limsup t_{n}
\end{equation*}

\Bin and

\begin{equation*}
\lim \inf (s_{n}+t_{n})\geq \lim \inf s_{n}+\lim \inf t_{n}
\end{equation*}

\Bin \noindent  \textbf{(6)} Deduce from (1) that if%

\begin{equation*}
\lim \inf x_{n}=\lim \sup x_{n},
\end{equation*}

\Bin then $(x_{n})_{n\geq 0}$ has a limit and

\begin{equation*}
\lim x_{n}=\lim \inf x_{n}=\lim \sup x_{n}
\end{equation*}

\bigskip
\noindent \textbf{Exercise 2.} Accumulation points  \index{accumulation point} of $(x_{n})_{n\geq 0}$.\newline

\noindent \textbf{(a)} Show that  $\ell _{1}$=$\lim \inf x_{n}$ and $\ell_{2}=\lim \sup x_{n}$ are accumulation points  \index{accumulation point} of $(x_{n})_{n\geq 0}$. Show one case and deduce the second using point (3) of Exercise 1.\\

\noindent \textbf{(b)} Show that $\ell _{1}$ is the smallest accumulation point  \index{accumulation point} of $(x_{n})_{n\geq 0}$ and $\ell _{2}$ is the biggest.
(Similarly, show one case and deduce the second using point (3) of exercise 1).\newline

\noindent \textbf{(c)} Deduce from (a) that if $(x_{n})_{n\geq 0}$ has
a limit $\ell$, then it is equal to the unique accumulation point  \index{accumulation point} and so,%

\begin{equation*}
\ell =\overline{\lim }\text{ }x_{n}=\lim \sup x_{n}=\inf_{n\geq 0}\sup_{p\geq n}x_{p}.
\end{equation*}

\Bin
\noindent \textbf{(d)} Combine his result with point \textbf{(6)} of Exercise 1 to show that a sequence $(x_{n})_{n\geq 0}$ of $\overline{\mathbb{R}}
$ has a limit $\ell $ in $\overline{\mathbb{R}}$ if and only if\ $\lim \inf
x_{n}=\lim \sup x_{n}$ and then

\begin{equation*}
\ell =\lim x_{n}=\lim \inf x_{n}=\lim \sup x_{n}.
\end{equation*}

\Bin \noindent \textbf{Exercise 3}. Let $(x_{n})_{n\geq 0}$ be a non-decreasing sequence
of $\overline{\mathbb{R}}$. Study its limit superior  \index{limit superior} and its limit inferior  \index{limit inferior} and deduce that%
\begin{equation*}
\lim x_{n}=\sup_{n\geq 0}x_{n}.
\end{equation*}

\noindent Deduce that for a non-increasing sequence $(x_{n})_{n\geq 0}$
of $\overline{\mathbb{R}},$%
\begin{equation*}
\lim x_{n}=\inf_{n\geq 0}x_{n}.
\end{equation*}

\bigskip

\noindent \textbf{Exercise 4.} (Convergence criteria \index{convergence criteria} )\newline

\noindent \textbf{Criterion 1.} Let $(x_{n})_{n\geq 0}$ be a sequence of $\overline{%
\mathbb{R}}$ and a real number $\ell \in \overline{\mathbb{R}}$ such that: Every sub-sequence of $(x_{n})_{n\geq 0}$
also has a sub-sequence ( that is a sub-sub-sequence of $(x_{n})_{n\geq 0}$ ) that converges to $\ell .$
Then, the limit of $(x_{n})_{n\geq 0}$ exists and is equal $\ell .$\newline

\noindent \textbf{Criterion 2. } Upcrossings and downcrossings.  \index{Upcrossings and downcrossings} \newline

\noindent Let $(x_{n})_{n\geq 0}$ be a sequence in $\overline{\mathbb{R}}$ and two real numbers $a$ and $b$ such that $a<b.$
We define

\begin{equation*}
\nu _{1}=\left\{
\begin{array}{cc}
\inf  & \{n\geq 0,x_{n}<a\} \\
+\infty  & \text{if (}\forall n\geq 0,x_{n}\geq a\text{)}%
\end{array}%
\right. .
\end{equation*}

\Bin If $\nu _{1}$ is finite, let

\begin{equation*}
\nu _{2}=\left\{
\begin{array}{cc}
\inf  & \{n>\nu _{1},x_{n}>b\} \\
+\infty  & \text{if } (\forall (n>\nu _{1}), \ x_{n}\leq b)
\end{array}%
\right. .
\end{equation*}

\Bin \noindent As long as the $\nu _{j}'s$ are finite, we can define for $\nu
_{2k-2}(k\geq 2)$

\begin{equation*}
\nu _{2k-1}=\left\{
\begin{array}{cc}
\inf  & \{n>\nu _{2k-2},x_{n}<a\} \\
+\infty  & \text{ if } ( (\forall n>\nu _{2k-2}), \ (x_{n}\geq a))
\end{array}
\right.
\end{equation*}

\Bin and for $\nu _{2k-1}$ finite,

\begin{equation*}
\nu _{2k}=\left\{
\begin{array}{cc}
\inf  & \{n>\nu _{2k-1},x_{n}>b\} \\
+\infty  & \text{if }((n>\nu _{2k-1}), \ (x_{n}\leq b))
\end{array}
\right. .
\end{equation*}

\noindent We stop once one $\nu _{j}$ is $+\infty$. If $\nu_{2j}$ is finite, then
\begin{equation*}
x_{\nu _{2j}}-x_{\nu _{2j-1}}>b-a.
\end{equation*}

\Bin
\noindent We then say : by that moving from $x_{\nu _{2j-1}}$ to $x_{\nu
_{2j}},$ we have accomplished a crossing (toward the up) of the segment $[a,b]$
called \textit{up-crossings}. Similarly, if one $\nu _{2j+1}$
is finite, then the segment $[x_{\nu _{2j}},x_{\nu _{2j+1}}]$ is a crossing downward (downcrossing) of the segment $[a,b].$ Let%

\begin{equation*}
D(a,b)=\text{ number of upcrossings of the sequence of the segment }[a,b]\text{.}
\end{equation*}

\bigskip
\noindent \textbf{(a)} What is the value of $D(a,b)$ if \ $\nu _{2k}$ is finite and $\nu
_{2k+1}$ infinite.\newline

\noindent \textbf{(b)} What is the value of $D(a,b)$ if \ $\nu _{2k+1}$ is finite and $\nu
_{2k+2}$ infinite.\newline

\noindent \textbf{(c)} What is the value of $D(a,b)$ if \ all the $\nu _{j}$'s are finite.
\newline

\noindent \textbf{(d)} Show that $(x_{n})_{n\geq 0}$ has a limit if and only if for all $a<b,$ $D(a,b)<\infty.$\newline

\noindent \textbf{(e)} Show that $(x_{n})_{n\geq 0}$ has a limit if and only if for all $a<b,$ $(a,b)\in \mathbb{Q}^{2},D(a,b)<\infty .$\newline

\bigskip

\noindent \textbf{Exercise 5}. (Cauchy Criterion). Let $(x_{n})_{n\geq 0} \subset \mathbb{R}$ be a sequence of (\textbf{real numbers}).\newline

\noindent \textbf{(a)} Show that if $(x_{n})_{n\geq 0}$ is Cauchy,
then it has a unique accumulation point  \index{accumulation point} $\ell \in
\mathbb{R}$ which is its limit.\newline

\noindent \textbf{(b)} Show that if a sequence $(x_{n})_{n\geq 0}\subset
\mathbb{R}$ \ converges to $\ell \in \mathbb{R},$ then, it is Cauchy.%
\newline

\noindent \textbf{(c)} Deduce the Cauchy criterion for sequences of real numbers.

\newpage

\begin{center}
\textbf{SOLUTIONS}
\end{center}

\noindent \textbf{Exercise 1}.\newline

\noindent \textbf{Question (1) :}. It is obvious that :%
\begin{equation*}
\underset{p\geq n}{\inf }x_{p}\leq x_{n}\leq \underset{p\geq n}{\sup }x_{p},
\end{equation*}

\noindent since $x_{n}$ is an element of $\left\{
x_{n},x_{n+1},...\right\} $ on which we take the supremum or the infinimum.%
\newline

\noindent \textbf{Question (2) :}. Let $y_{n}=\underset{p\geq n}{\inf }x_{p}=\inf A_{n}$,  where $A_{n}=\left\{
x_{n},x_{n+1},...\right\} $ is a non-increasing sequence of sets : $\forall n\geq 0$,

\begin{equation*}
A_{n+1}\subset A_{n}.
\end{equation*}

\Bin \noindent So the infinimum on $A_{n}$ increases. If $y_{n}$ increases in $%
\overline{\mathbb{R}},$ its limit is its upper bound, finite or infinite. So%
\begin{equation*}
y_{n}\nearrow \underline{\lim }\text{ }x_{n},
\end{equation*}%

\Bin is a finite or infinite number.\newline

\noindent \textbf{Question (3) :}. We also show that $z_{n}=\sup A_{n}$ decreases and $z_{n}\downarrow \overline{\lim }$ $x_{n}$.\newline

\noindent \textbf{Question (4) \label{qst4}:}. We recall that

\begin{equation*}
-\sup \left\{ x,x\in A\right\} =\inf \left\{ -x,x\in A\right\}.
\end{equation*}

\Bin \noindent Which we write as

\begin{equation*}
-\sup A=\inf (-A).
\end{equation*}

\Bin \noindent Thus,

\begin{equation*}
-z_{n}=-\sup A_{n}=\inf -A_{n} = \inf \left\{-x_{p},p\geq n\right\}.
\end{equation*}

\Bin
\noindent The right hand term tends to $-\overline{\lim}\ x_{n}$ and the left hand to $\underline{\lim} \ (-x_{n})$ and so

\begin{equation*}
-\overline{\lim}\ x_{n}=\underline{\lim }\ (-x_{n}).
\end{equation*}

\bigskip \noindent Similarly, we show:

\begin{equation*}
-\underline{\lim } \ (x_{n})=\overline{\lim} \ (-x_{n}).
\end{equation*}

\noindent

\noindent \textbf{Question (5)}. These properties come from the formulas, where $A\subseteq \mathbb{R}$, $B\subseteq \mathbb{R}$:

\begin{equation*}
\sup \left\{ x+y, \ x \in A, \ y \in B\right\} \leq \sup A+\sup B.
\end{equation*}

\noindent In fact :
\begin{equation*}
\forall x\in \mathbb{R},x\leq \sup A
\end{equation*}

\noindent and
\begin{equation*}
\forall y\in \mathbb{R},y\leq \sup B.
\end{equation*}

\noindent Thus
\begin{equation*}
x+y\leq \sup A+\sup B,
\end{equation*}

\noindent where

\begin{equation*}
\underset{x\in A,y\in B}{\sup} \ (x+y) \leq \sup A+\sup B.
\end{equation*}

\Bin  Similarly,

\begin{equation*}
\inf (A+B\geq \inf A+\inf B.
\end{equation*}

\noindent In fact :

\begin{equation*}
\forall (x,y)\in A\times B,x\geq \inf A\text{ and }y\geq \inf B.
\end{equation*}

\noindent Hence

\begin{equation*}
x+y\geq \inf A+\inf B.
\end{equation*}

\noindent Thus
\begin{equation*}
\underset{x\in A,y\in B}{\inf }(x+y)\geq \inf A+\inf B
\end{equation*}

\Bin \noindent \textbf{Application}.\newline

\begin{equation*}
\underset{p\geq n}{\sup } \ (x_{p}+y_{p})\leq \underset{p\geq n}{\sup } \ x_{p}+\underset{p\geq n}{\sup } \ y_{p}.
\end{equation*}

\Bin \noindent All these sequences are non-increasing. Taking infimum, we obtain the limits superior :

\begin{equation*}
\overline{\lim }\text{ }(x_{n}+y_{n})\leq \overline{\lim }\text{ }x_{n}+%
\overline{\lim }\text{ }x_{n}.
\end{equation*}

\bigskip

\noindent \textbf{Question (6)}: Set

\begin{equation*}
\underline{\lim } \ x_{n}=\overline{\lim } \ x_{n},
\end{equation*}

\Bin
\noindent Since

\begin{equation*}
\forall n\geq 1,\text{ }y_{n}\leq x_{n}\leq z_{n},
\end{equation*}

\begin{equation*}
y_{n}\rightarrow \underline{\lim} \ x_{n}
\end{equation*}

\noindent and

\begin{equation*}
z_{n}\rightarrow \overline{\lim } \ x_{n},
\end{equation*}

\Bin \noindent we apply Sandwich Theorem to conclude that the limit of $x_{n}$ exists and :

\begin{equation*}
\lim \text{ }x_{n}=\underline{\lim }\text{ }x_{n}=\overline{\lim }\text{ }%
x_{n}.
\end{equation*}

\bigskip
\noindent \textbf{Exercise 2}.\newline

\noindent \textbf{Question (a).}\\

\noindent Thanks to question (4) of exercise 1, it suffices to show this property for one of the limits. Consider the limit superior  \index{limit superior} and the three cases:\\

\noindent \textbf{The case of a finite limit superior}  \index{limit superior} :

\begin{equation*}
\underline{\lim \text{ }}x_{n}=\ell, \ \ell \ finite.
\end{equation*}

\noindent By definition,
\begin{equation*}
z_{n}=\underset{p\geq n}{\sup} \ x_{p}\downarrow \ell .
\end{equation*}

\noindent So:

\begin{equation*}
\forall \varepsilon >0,\exists (N(\varepsilon )\geq 1),\forall p\geq
N(\varepsilon ),\ell -\varepsilon <x_{p}\leq \ell +\varepsilon .
\end{equation*}

\Bin
\noindent Take less than that:

\begin{equation*}
\forall \varepsilon >0,\exists n_{\varepsilon }\geq 1:\ell -\varepsilon
<x_{n_{\varepsilon }}\leq \ell +\varepsilon.
\end{equation*}

\Bin
\noindent We shall construct a subsequence converging to $\ell$.\\

\noindent Let $\varepsilon =1$:

\begin{equation*}
\exists N_{1}:\ell -1<z_{N_{1}}=\underset{p\geq N_1}{\sup }x_{p}\leq \ell +1.
\end{equation*}

\noindent But if

\begin{equation}
z_{N_{1}}=\underset{p\geq N-1}{\sup } \ x_{p}>\ell -1, \label{cc}
\end{equation}

\Bin
\noindent there surely exists an $n_{1}\geq N_{1}$ such that%
\begin{equation*}
x_{n_{1}}>\ell -1.
\end{equation*}

\Bin
\noindent If not, we would have

\begin{equation*}
( \forall p\geq N_{1},x_{p}\leq \ell -1\ ) \Longrightarrow \sup \left\{
x_{p},p\geq N_{1}\right\} =z_{N_{1}}\leq \ell -1,
\end{equation*}

\Bin which is contradictory with (\ref{cc}). So, there exists $n_{1}\geq N_{1}$ such that

\begin{equation*}
\ell -1<x_{n_{1}}\leq \underset{p\geq N_{1}}{\sup }x_{p}\leq \ell +1.
\end{equation*}

\Bin
\noindent i.e.

\begin{equation*}
\ell -1<x_{n_{1}}\leq \ell +1.
\end{equation*}

\Bin
\noindent We move to step $\varepsilon =\frac{1}{2}$ and we consider the sequence%
$(z_{n})_{n\geq n_{1}}$ whose limit remains $\ell$. So, there exists $N_{2}>n_{1}:$%

\begin{equation*}
\ell -\frac{1}{2}<z_{N_{2}}\leq \ell +\frac{1}{2}.
\end{equation*}

\Bin
\noindent We deduce like previously that $n_{2}\geq N_{2}$ such that%
\begin{equation*}
\ell -\frac{1}{2}<x_{n_{2}}\leq \ell +\frac{1}{2}
\end{equation*}

\noindent with $n_{2}\geq N_{1}>n_{1}$.\\

\noindent Next, we set $\varepsilon =1/3,$ there will exist $N_{3}>n_{2}$ such that%

\begin{equation*}
\ell -\frac{1}{3}<z_{N_{3}}\leq \ell -\frac{1}{3}
\end{equation*}

\Bin \noindent and we could find an $n_{3}\geq N_{3}$ such that

\begin{equation*}
\ell -\frac{1}{3}<x_{n_{3}}\leq \ell +\frac{1}{3}.
\end{equation*}

\Bin
\noindent Step by step, we deduce the existence of $%
x_{n_{1}},x_{n_{2}},x_{n_{3}},...,x_{n_{k}},...$ with $n_{1}<n_{2}<n_{3}%
\,<...<n_{k}<n_{k+1}<...$ such that

$$
\forall k\geq 1, \ell -\frac{1}{k}<x_{n_{k}}\leq \ell -\frac{1}{k},
$$

\noindent i.e.

\begin{equation*}
\left\vert \ell -x_{n_{k}}\right\vert \leq \frac{1}{k}.
\end{equation*}

\Bin
\noindent Which will imply:

\begin{equation*}
x_{n_{k}}\rightarrow \ell
\end{equation*}

\Bin \noindent Conclusion : $(x_{n_{k}})_{k\geq 1}$ is very well a subsequence since $n_{k}<n_{k+1}$ for all $k \geq 1$
and it converges to $\ell$, which is then an accumulation point \index{accumulation point} .\\

\noindent \textbf{Case of the limit superior  \index{limit superior} equal $+\infty$} :

$$
\overline{\lim} \text{ } x_{n}=+\infty.
$$

\Bin \noindent Since $z_{n}\uparrow +\infty ,$ we have : $\forall k\geq 1,\exists
N_{k}\geq 1,$

\begin{equation*}
z_{N_{k}}\geq k+1.
\end{equation*}

\Bin \noindent For $k=1$, let $z_{N_{1}}=\underset{p\geq N_{1}}{\inf }%
x_{p}\geq 1+1=2.$ So there exists

\begin{equation*}
n_{1}\geq N_{1}
\end{equation*}

\Bin such that:

\begin{equation*}
x_{n_{1}}\geq 1.
\end{equation*}

\Bin \noindent For $k=2:$ consider the sequence $(z_{n})_{n\geq n_{1}+1}.$
We find in the same manner

\begin{equation*}
n_2 \geq n_{1}+1
\end{equation*}%

\Bin \noindent and

\begin{equation*}
x_{n_{2}}\geq 2.
\end{equation*}

\Bin \noindent Step by step, we find for all $k\geq 3$, an $n_{k}\geq n_{k-1}+1$ such that

\begin{equation*}
x_{n_{k}}\geq k,
\end{equation*}

\Bin \noindent which leads to $x_{n_{k}}\rightarrow +\infty $ as $k\rightarrow +\infty $.\\

\noindent \textbf{Case of the limit superior  \index{limit superior} equal $-\infty$} :

$$
\overline{\lim }x_{n}=-\infty.
$$

\Bin
\noindent This implies : $\forall k\geq 1,\exists N_{k}\geq 1,$ such that%

\begin{equation*}
z_{n_{k}}\leq -k.
\end{equation*}

\Bin \noindent For $k=1,\exists n_{1}$ such that%

\begin{equation*}
z_{n_{1}}\leq -1.
\end{equation*}

\Bin But

\begin{equation*}
x_{n_{1}}\leq z_{n_{1}}\leq -1
\end{equation*}

\noindent Let $k=2$. Consider $\left( z_{n}\right) _{n\geq
n_{1}+1}\downarrow -\infty .$ There will exist $n_{2}\geq n_{1}+1:$%

\begin{equation*}
x_{n_{2}}\leq z_{n_{2}}\leq -2
\end{equation*}

\Bin
\noindent Step by step, we find $n_{k+1}<n_{k+1}$ in such a way that $x_{n_{k}}<-k$ for all $k$ bigger that $1$. So%
\begin{equation*}
x_{n_{k}}\rightarrow -\infty, \ as k\rightarrow +\infty,
\end{equation*}

\bigskip
\noindent \textbf{Question (b).}\\

\noindent Let $\ell$ be an accumulation point  \index{accumulation point} of $(x_n)_{n \geq 1}$, the limit of one of its subsequences $(x_{n_{k}})_{k \geq 1}$. We have

$$
y_{n_{k}}=\inf_{p\geq n_k} \ x_p \leq x_{n_{k}} \leq  \sup_{p\geq n_k} \ x_p=z_{n_{k}}
$$

\Bin
\noindent The left hand side term is a subsequence of $(y_n)$ tending to the limit inferior  \index{limit inferior} and the right hand side is a
subsequence of $(z_n)$ tending to the limit superior.  \index{limit superior} So we will have:

$$
\underline{\lim} \ x_{n} \leq \ell \leq \overline{\lim } \ x_{n},
$$

\Bin
\noindent which shows that $\underline{\lim} \ x_{n}$ is the smallest accumulation point  \index{accumulation point} and $\overline{\lim } \ x_{n}$ is the largest.\\

\noindent \textbf{Question (c).} If the sequence $(x_n)_{n \geq 1}$ has a limit $\ell$, it is the limit of all its subsequences,
so subsequences tending to the limits superior and inferior. Which answers question (b).\\

\noindent \textbf{Question (d).} We answer this question by combining point (d) of this exercise and point (\textbf{6}) of the exercise \textbf{1}.\\

\noindent \textbf{Exercise 3}. Let $(x_{n})_{n\geq 0}$ be a non-decreasing sequence, we have:%

\begin{equation*}
z_{n}=\underset{p\geq n}{\sup} \ x_{p}=\underset{p\geq 0}{\sup} \ x_{p},\forall
n\geq 0.
\end{equation*}

\Bin \noindent Why? Because by increasingness,

\begin{equation*}
\left\{ x_{p},p\geq 0\right\} =\left\{ x_{p},0\leq p\leq n-1\right\} \cup
\left\{ x_{p},p\geq n\right\}
\end{equation*}

\bigskip
\noindent Since all the elements of $\left\{ x_{p},0\leq p\leq
n-1\right\} $ are smaller than that of $\left\{ x_{p},p\geq n\right\} ,$
the supremum is achieved on $\left\{ x_{p},p\geq n\right\} $ and so

\begin{equation*}
\ell =\underset{p\geq 0}{\sup } \ x_{p}=\underset{p\geq n}{\sup } \ x_{p}=z_{n}.
\end{equation*}

\Bin Thus

\begin{equation*}
z_{n}=\ell \rightarrow \ell .
\end{equation*}

\noindent We also have $y_n=\inf \left\{ x_{p},p\geq n\right\}=x_n$ which is a non-decreasing sequence and so converges to
$\ell =\underset{p\geq 0}{\sup } \ x_{p}$. \\

\bigskip

\noindent \textbf{Exercise 4}.\\

\noindent Let $\ell \in \overline{\mathbb{R}}$ having the indicated property. Let $\ell ^{\prime }$ be a given accumulation point \index{accumulation point} .%

\begin{equation*}
\left( x_{n_{k}}\right)_{k\geq 1} \subseteq \left( x_{n}\right) _{n\geq 0}%
\text{ such that }x_{n_{k}}\rightarrow \ell ^{\prime}.
\end{equation*}

\noindent By hypothesis this subsequence $\left( x_{n_{k}}\right) $
has in turn a subsubsequence $\left( x_{n_{\left( k(p)\right) }}\right)_{p\geq 1} $ such that $x_{n_{\left( k(p)\right) }}\rightarrow
\ell $ as $p\rightarrow +\infty $.\newline

\noindent But as a subsequence of $\left( x_{n_{\left( k\right)}}\right)$,

\begin{equation*}
x_{n_{\left( k(\ell )\right) }}\rightarrow \ell ^{\prime }.
\end{equation*}

\Bin Thus

\begin{equation*}
\ell =\ell ^{\prime}.
\end{equation*}

\Bin \noindent Applying that to the limit superior  \index{limit superior} and limit inferior,  \index{limit inferior} we have:%

\begin{equation*}
\overline{\lim} \ x_{n}=\underline{\lim}\ x_{n}=\ell.
\end{equation*}

\noindent And so $\lim x_{n}$ exists and equals $\ell$.\\

\noindent \textbf{Exercise 5}.\\

\noindent \textbf{Question (a)}. If $\nu _{2k}$ finite and $\nu _{2k+1}$ infinite, it then has exactly $k$ up-crossings :
$[x_{\nu_{2j-1}},x_{\nu _{2j}}]$, $j=1,...,k$ : $D(a,b)=k$.\\

\noindent \textbf{Question (b)}. If $\nu _{2k+1}$ finite and $\nu _{2k+2}$ infinite, it then has exactly $k$ up-crossings:
$[x_{\nu_{2j-1}},x_{\nu_{2j}}]$, $j=1,...,k$ : $D(a,b)=k$.\\

\noindent \textbf{Question (c)}. If all the $\nu_{j}$'s are finite, then, there are an infinite number of up-crossings :
$[x_{\nu_{2j-1}},x_{\nu_{2j}}]$, $j\geq k$ : $D(a,b)=+\infty$.\\

\noindent \textbf{Question (d)}. Suppose that there exist $a < b$ rationals such that $D(a,b)=+\infty$.
Then all the $\nu _{j}$'s are finite. The subsequence $x_{\nu_{2j-1}}$ is strictly below $a$.
So its limit inferior  \index{limit inferior} is below $a$. This limit inferior is an accumulation point  \index{accumulation point} of the sequence $(x_n)_{n\geq 1}$,
so is more than $\underline{\lim}\ x_{n}$, which is below $a$.\\

\noindent Similarly, the subsequence $x_{\nu_{2j}}$ is strictly below $b$. So the limit superior  \index{limit superior} is above $b$.
This limit superior  \index{limit superior} is an accumulation point  \index{accumulation point} of the sequence $(x_n)_{n\geq 1}$, so it is below $\overline{\lim}\ x_{n}$,
which is directly above $b$. Which leads to:

$$
\underline{\lim}\ x_{n} \leq a < b \leq \overline{\lim}\ x_{n}.
$$

\Bin
\noindent That implies that the limit of $(x_n)$ does not exist. In contrary, we just proved that the limit of $(x_n)$ exists,
meanwhile  \index{mean} for all the real numbers $a$ and $b$ such that $a<b$, $D(a,b)$ is finite.\\

\noindent Now, suppose that the limit of $(x_n)$ does not exist. Then,

$$
\underline{\lim}\ x_{n} < \overline{\lim}\ x_{n}.
$$

\noindent We can then find two rationals $a$ and $b$ such that $a<b$ and a number $\epsilon>0$ such that

$$
\underline{\lim}\ x_{n} < a-\epsilon < a < b < b+\epsilon <  \overline{\lim}\ x_{n}.
$$

\Bin
\noindent If $\underline{\lim}\ x_{n} < a-\epsilon$, we can return to question \textbf{(a)} of exercise \textbf{2} and construct a subsequence of $(x_n)$
which tends to $\underline{\lim}\ x_{n}$ while remaining below $a-\epsilon$. Similarly, if $b+\epsilon < \overline{\lim}\ x_{n}$,
we can create a subsequence of $(x_n)$ which tends to $\overline{\lim}\ x_{n}$ while staying above $b+\epsilon$.
It is evident with these two sequences that we could define with these two sequences all $\nu_{j}$ finite and so $D(a,b)=+\infty$.\\

\noindent We have just shown by contradiction that if all the $D(a,b)$ are finite for all rationals $a$ and $b$ such that $a<b$,
then, the limit of $(x_n)_ {n\geq 0}$ exists.\\

\noindent \textbf{Exercise 5}. Cauchy criterion in $\mathbb{R}$.\\

\noindent Suppose that the sequence $(x_n)_ {n\geq 0}$ is Cauchy, $i.e.$,

$$
\lim_{(p,q)\rightarrow (+\infty,+\infty)} \ (x_p-x_q)=0.
$$

\Bin \noindent Then let $x_{n_{k,1}}$ and $x_{n_{k,2}}$ be two subsequences converging respectively to $\ell_1=\underline{\lim}\ x_{n}$ and $\ell_2=\overline{\lim}\ x_{n}$. So

$$
\lim_{(p,q)\rightarrow (+\infty,+\infty)} \ (x_{n_{p,1}}-x_{n_{q,2}})=0.
$$

\Bin
\noindent By first letting $p\rightarrow +\infty$, we have

$$
\lim_{q\rightarrow +\infty} \ \ell_1-x_{n_{q,2}}=0,
$$

\noindent which shows that $\ell_1$ is finite, else $\ell_1-x_{n_{q,2}}$ would remain infinite and would not tend to $0$.
By interchanging the roles of $p$ and $q$, we also have that $\ell_2$ is finite.\\

\noindent Finally, by letting $q\rightarrow +\infty$, in the last equation, we obtain
$$
\ell_1=\underline{\lim}\ x_{n}=\overline{\lim}\ x_{n}=\ell_2,
$$

\noindent which proves the existence of the finite limit of the sequence $(x_n)_ {n\geq 0}$.\\

\noindent Now suppose that the finite limit $\ell$ of $(x_n)_ {n\geq 0}$ exists. Then

$$
\lim_{(p,q)\rightarrow (+\infty,+\infty)} \ (x_p-x_q)=\ell-\ell=0,
$$

\Bin
\noindent  which shows that the sequence is Cauchy.\\

\bigskip
\noindent \textbf{Bravo! Grab this knowledge. Limits in $\overline{\mathbb{R}}$ have no more secrets for you!}.

\newpage

Convex functions  \index{convex function} play an important role in real analysis, and in probability theory  \index{probability theory} in particular.\\

\noindent A convex function  \index{convex function} is defined as follows.\\

\noindent \textbf{A - Definition}. A real-valued function $ g:I\longrightarrow \mathbb{R}$ defined on an interval on $\mathbb{R}$ is convex if and only if, for any $0<\alpha<1$, for any $(x,y)\in I^2$, we have

\begin{equation}
g(\alpha x+(1-\alpha)y)\leq \alpha g(x)+(1-\alpha)g(y). \label{proba01_convex_def}
\end{equation}

\Bin \noindent The first important thing to know is that a convex function  \index{convex function} is continuous.\\

\bigskip \noindent \textbf{B - A convex function  \index{convex function} is continuous}.\\

\noindent Actually, we have more.

\begin{proposition} If $g:I\longrightarrow \mathbb{R}$ is convex on the interval $I$, then g admits a right-derivative and a left-derivative at each point of $I$. In particular, $g$ is continuous on $I$.
\end{proposition}

\bigskip \noindent \textbf{Proof}. In this proof, we take $I=\mathbb{R}$, which is the most general case. Suppose that $g$ is convex. We begin to prove this formula :

\begin{equation*}
\forall s < t < u, \frac{g(t)-g(s)}{t-s} \leq \frac{g(u)-g(s)}{u-s} \leq
\frac{g(u)-g(t)}{u-t}. \ \ \  \text{(FCC1)}
\end{equation*}

\Bin
\noindent We check that, for $s < t < u$, we have

\begin{equation*}
t=\frac{u-t}{u-s} s + \frac{t-s}{u-s} u = \lambda s + (1-\lambda)t,
\end{equation*}

\Bin
\noindent where $\lambda=(u-t)/(u-s) \in ]0,1[$ and $1-\lambda=(t-s)/(u-s)$. By convexity, we have

\begin{equation*}
g(t)\leq \frac{u-t}{u-s} g(s) + \frac{t-s}{u-s} g(u). \ \  \text{(FC1a)}.
\end{equation*}

\bigskip
\noindent Let us multiply all members of (FC1) by $(u-s)$ to get

\begin{equation*}
(u-s) g(t)\leq (u-t) g(s) + (t-s) g(u). \ \ \  \text{(FC1b)}.
\end{equation*}

\bigskip \noindent First, we split $(u-t)$ into $(u-t)=(u-s)-(t-s)$ in the first term in
the left-hand member to have

\begin{eqnarray*}
&&(u-s) g(t)\leq (u-s) g(s) - (t-s) g(s) +(t-s) g(u)\\
&\Rightarrow& (u-s) (g(t)-g(s)) \leq (t-s) (g(u)-g(s)).
\end{eqnarray*}

\bigskip
\noindent This leads to

\begin{equation*}
\frac{g(t)-g(s)}{t-s} \leq \frac{g(u)-g(s)}{u-s}. \ \ \  \text{(FC2a)}.
\end{equation*}

\Bin
\noindent Next, we split $(t-s)$ into $(t-s)=(u-s)-(u-t)$ in the second term in the left-hand member to have

\begin{eqnarray*}
&&(u-s) g(t)\leq (u-t) g(s) + (u-s) g(u) - (u-t) g(u)\\
&\Rightarrow& (u-s) (g(t)-g(u)) \leq (u-t) (g(s)-g(u)).
\end{eqnarray*}

\bigskip
\noindent Let us multiply the last inequality by $-1$ to get

\begin{equation*}
(u-s) (g(u)-g(t)) \geq (u-t) (g(u)-g(s))
\end{equation*}

\Bin
\noindent and then

\begin{equation*}
\frac{g(u)-g(s)}{u-s} \leq \frac{g(u)-g(t)}{u-t} \text{ (FC2b) }
\end{equation*}

\bigskip
\noindent Formulas \textit{(FC2a)} and \textit{(FC2b)} together prove \textit{(FCC1)}.\newline

\noindent Let us write \textit{(FCC1)} in the following form

\begin{equation*}
\forall t_1 < t_2 < t_3, \frac{g(t_2)-g(t_1)}{t_2 - t_1} \leq \frac{g(t_3)-g(t_1)}{t_3-t_1} \leq \frac{g(t_3)-g(t_2)}{t_3-t_2}. \ \ \  \text{(FCC2)}
\end{equation*}

\bigskip
\noindent We also may apply \textit{(FCC1)} to get :\newline

\begin{equation*}
\forall r < s < t, \frac{g(s)-g(r)}{s-r} \leq \frac{g(t)-g(r)}{t-r} \leq \frac{g(t)-g(s)}{t-s}. \ \ \  \text{(FCC3)}
\end{equation*}

\bigskip
\noindent We are on the point to conclude. Fix $r < s$. From \textit{(FFC1)}, we may
see that the function

\begin{equation*}
G(v)=\frac{g(v)-g(t)}{v - t}, v >t,
\end{equation*}

\bigskip
\noindent is increasing since in \textit{(FCC1)}, $G(t) \leq G(u)$ for $t<u$. By (FCC2), $G(v)$ is bounded below by $G(r)$ (that is fixed with $r$ and $t$).\\

\noindent We conclude $G(v)$ decreases to a real number as $v \downarrow t$ and this limit is the right-derivative of $g$ at $t$ :

\begin{equation*}
\lim_{v \downarrow t} \frac{g(v)-g(t)}{v - t}=g^{\prime}_r(t).  \ \  \text{(RC)}
\end{equation*}

\Bin
\noindent Since $g$ has a right-hand derivative at each point, it is right-continuous at each point since, by (RC),

\begin{equation*}
g(t+h)-g(t)= h (G(t+h)=h(g^{\prime}_r(t)) + o(1)) \ \ as \ \ h\rightarrow 0.
\end{equation*}

\bigskip \noindent \textbf{Extension}. To prove that the left-hand derivative exists, we use a very similar way. We conclude that $g$ is left-continuous and
finally, $g$ is continuous.$\square$\newline

\bigskip \noindent \textbf{C - A weakened convexity condition}.\\

\begin{lemma} \label{proba01_convex_def2} Suppose that $g:I\longrightarrow \mathbb{R}$ is continuous. Then the mapping $g$ is convex whenever we have for any $(x,y)\in I^2$,
$$
g\left(\frac{x+y}{2}\right) \leq \frac{g(x)+g(y)}{2}.
$$
\end{lemma}

\Bin
\noindent So, for a continuous function, we may only check the condition of convexity formula (\ref{proba01_convex_def}) for $\alpha=1/2$.\\

\bigskip \noindent \textbf{Proof of Lemma \ref{proba01_convex_def2}}. Let us suppose that the assumptions of the lemma hold. We are going to exploit the density of dyadic numbers

$$
D_{[0,1]}=\left\{\frac{k}{2^n}, \ n\geq 0, \ 0\leq k\leq 2^n\right\}
$$

\Bin in $[0,1]$ and proceed in two steps.\\

\Ni \textbf{Step 1}. We have to prove that for any dyadic number $r \in [0,1]$, for $(x,y) \in I^2$, we have

$$
g(rx + (1-r)y)\leq r g(x) + (1-r) g(y). \ \ \ (CV)
$$

\Bin We fix $(x,y)$ in the rest of the proof. Let us put, for $n\geq 0$,

$$
\mathcal{P}_n=\biggr(\forall k \in \{0,\cdots,2^n\}, \ g\left(\frac{k}{n}x + \left(1-\frac{k}{n}\right)y\right)\leq \frac{k}{n} g(x) + \left(1-\frac{k}{n}\right) g(y) \ \biggr)
$$

\Bin Let us prove that $\mathcal{P}_n$ holds for $n\geq 0$. For $n=0$, we have $2^n=1$ and hence $k=0$ or $k=1$. We respectively get

$$
g(0 \ x+1 \ y)=g(y)=0 \ g(x) + 1 \ g(y) \ and \ g(1 \ x+0 \ y)=g(x)=1 \  g(x) + 0 \ g(y).
$$

\Bin For $n=1$, $2^n=2$ and hence $k \in \{0,1,2\}$. The cases $k \in \{0,1\}$ are handled previously. For $k=1$, $\mathcal{P}_2$ reduces to the hypothesis

$$
g\left(\frac{x+y}{2}\right) \leq \frac{1}{2} g(x) + \frac{1}{2} g(y).
$$

\Bin Now let us prove $\{\mathcal{P}_0, \mathcal{P}_1, \cdots, \mathcal{P}_n\} \Rightarrow \mathcal{P}_{n+1}$ for $n\geq 2$. Let us suppose that that the
$\mathbb{P}_j$'s hold true for $0\leq j\leq n$. Let $k \in \{0,\cdots,n+1\}$, we have

\begin{eqnarray*}
g\left( \frac{k}{2^{n+1}} x + \left(1 - \frac{k}{2^{n+1}}\right) y\right)
&=&g\left(  \frac{1}{2} \frac{k}{2^{n}} x + \frac{1}{2} \left(2 - \frac{k}{2^{n}}\right) \right)\\
&=&g\left( \frac{1}{2} \frac{k}{2^{n}} x +  \frac{1}{2}\left(1 - \frac{k}{2^{n}}\right) y + \frac{1}{2} y\right)\\
&=&g\left( \frac{1}{2} \biggr\{\frac{k}{2^{n}} x +  \left(1 - \frac{k}{2^{n}}\right) y\biggr\} + \frac{1}{2} y\right)\\
&=:&g\left( \frac{1}{2} z + \frac{1}{2} y\right),\\
\end{eqnarray*}

\Bin with

$$
z=\biggr\{\frac{k}{2^{n}} x + \left(1 - \frac{k}{2^{n}}\right) y\biggr\}.
$$

\Bin We use $\mathcal{P}_2$ to get

\begin{eqnarray*}
g\left( \frac{k}{2^{n+1}} x + \left(1 - \frac{k}{2^{n+1}}\right) y\right)
&\leq& \frac{1}{2} g(z) + \frac{1}{2} g(y)\\
&=& \frac{1}{k} g\left(\frac{k}{2^{n}} x + \frac{1}{2} \left(1 - \frac{k}{2^{n}}\right)y\right) + \frac{1}{2} g(y).
\end{eqnarray*}

\Bin Now, we use $\mathcal{P}_n$ to conclude

\begin{eqnarray*}
g\left( \frac{2}{2^{n+1}} x + \left(1 - \frac{k}{2^{n+1}}\right) y\right)&=& \frac{1}{2} \frac{k}{2^{n}} g(x) + \frac{1}{2} \left(1 - \frac{k}{2^{n}}\right) g(y) + \frac{1}{2} g(y)\\
&=&\frac{2}{2^{n+1}} g(x) +  \left(1 - \frac{k}{2^{n+1}}\right) g(y).
\end{eqnarray*}

\Bin We finished the proof by induction. So, Formula \textit{CV} of that step is established.\\

\Bin \textbf{Step 2}. Let $r \in [0,1]$. By the density of $D_{[0,1]}$ in $[0,1]$, there exists $(r_p)_{p\geq 1} \subset D_{[0,1]}$ such that $r_p \rightarrow to r$.
We have for all $p\geq 1$

$$
g(r_p x + (1-r_p) y)\leq r_p g(x) + (1-r_p) g(y).
$$

\Bin By letting $p\rightarrow +\infty$ and by the continuity of $g$, we get

$$
g(r x + (1-r) y)\leq r g(x) + (1-r) g(y).
$$

\Bin The proof is over. $\square$\\

\bigskip \noindent Finally, we are going to give generalization of Formula (\ref{proba01_convex_def}) to more than two terms.

\bigskip \noindent \textbf{D - General Formula of Convexity}.\\

\noindent \textbf{D1 - A convexity formula with an arbitrary finite number of points}.\\

\noindent Let $g$ be convex and let $k\geq 2$. Then for any $\alpha _{i}$, $0<\alpha _{i}<1$, $1\leq i\leq k$ such that

\begin{equation*}
\alpha _{1}+...+\alpha _{k}=1
\end{equation*}

\Bin \noindent and for any $(x_{1},...,x_{k})\in I^{k},$ we have

\begin{equation}
g(\alpha _{1}x_{1}+...+\alpha _{k}x_{k})\leq \alpha _{1}g(x_{1})+...+\alpha
_{k}g(x_{k}).  \label{proba01_convexK}
\end{equation}

\bigskip \noindent \textbf{Proof}. Let us use an induction argument. The formula is true for $k=2$ by
definition. Let us assume it is true for up to $k$. And we have to prove it
is true for $k+1.$ Set $\alpha _{i}$, $0<\alpha _{i}<1$, $1\leq i\leq
k+1$ such that
\begin{equation*}
\alpha _{1}+...+\alpha _{k}+\alpha _{k+1}=1
\end{equation*}

\bigskip
\noindent and let $(x_{1},...,x_{k},x_{k})\in I^{k+1}.$ By denoting

\begin{equation*}
\beta _{k+1}=\alpha _{1}+...+\alpha _{k}
\end{equation*}%
\begin{equation*}
\alpha _{i}^{\ast }=\alpha _{i}/\beta _{k+1},1\leq i\leq k,
\end{equation*}

\Bin
\noindent and
\begin{equation*}
y_{k+1}=\alpha _{1}^{\ast }x_{1}+...+\alpha _{k}^{\ast }x_{k},
\end{equation*}

\Bin
\noindent we have
\begin{eqnarray*}
g(\alpha _{1}x_{1}+...+\alpha _{k}x_{k}+\alpha _{k+1}x_{k+1}) &=&g\left((\alpha _{1}+...+\alpha _{k})(\alpha _{1}^{\ast }x_{1}+...+\alpha _{k}^{\ast
}x_{k})+\alpha _{k+1}x_{k+1})\right)  \\
&=&g\left( \beta _{k+1}y_{k+1}+\alpha _{k+1}x_{k+1}\right).
\end{eqnarray*}

\Bin
\noindent Since $\beta _{k+1}>0,\alpha _{k+1}>0$ and $\beta _{k+1}+\alpha _{k+1}=1,$
we have by convexity

\begin{equation*}
g(\alpha _{1}x_{1}+...+\alpha _{k}x_{k}+\alpha _{k+1}x_{k+1})\leq \beta
_{k+1}g(y_{k+1})+\alpha _{k+1}g(x_{k+1}).
\end{equation*}

\bigskip
\noindent But, we also have

\begin{equation*}
\alpha _{1}^{\ast }+...+\alpha _{k}^{\ast }=1,\text{ }\alpha _{i}^{\ast
}>0,1\leq i\leq k.
\end{equation*}

\bigskip
\noindent Thus, by the induction hypothesis, we have

\begin{equation*}
g(y_{k+1})=g(\alpha _{1}^{\ast }x_{1}+...+\alpha _{k}^{\ast }x_{k})\leq
\alpha _{1}^{\ast }g(x_{1})+...+\alpha _{k}^{\ast }g(x_{k}).
\end{equation*}

\bigskip
\noindent By combining these facts, we get

\begin{eqnarray*}
g(\alpha _{1}x_{1}+...+\alpha _{k}x_{k}+\alpha _{k+1}x_{k+1}) &\leq &\alpha
_{1}^{\ast }\beta _{k+1}g(x_{1})+...+\alpha _{k}^{\ast }\beta
_{k+1}g(x_{k})+\alpha _{k+1}g(x_{k+1}) \\
&=&\alpha _{1}g(x_{1})+...+\alpha _{k}g(x_{k})+\alpha _{k+1}g(x_{k+1}),
\end{eqnarray*}

\bigskip
\noindent since $\alpha _{i}^{\ast }\beta _{k+1}=\alpha _{i},1\leq i\leq k.$ Thus, we have proved

\begin{equation*}
g(\alpha _{1}x_{1}+...+\alpha _{k}x_{k}+\alpha _{k+1}x_{k+1})\leq \alpha
_{1}g(x_{1})+...+\alpha _{k}g(x_{k})+\alpha _{k+1}g(x_{k+1}).
\end{equation*}

\bigskip
\noindent We conclude that Formula (\ref{proba01_convexK}) holds for any $k\geq 2.$\\

\bigskip
\noindent \textbf{D2 - A convexity formula with a infinite and countable arbitrary finite number of points}.\\

\noindent Let $g$ be bounded on $I$ or $I$ be a closed bounded interval. Then for an infinite and countable number of coefficients $(\alpha _{i})_{i\geq 0}
$ with

\begin{equation*}
\sum_{i\geq 0}\alpha _{i}=1,\text{ \ and  }\forall (i\geq 0),\text{ }\alpha
_{i}>0
\end{equation*}

\noindent and for any family of $(x_{i})_{i\geq 0}$ of points in $I$, we have

\begin{equation}
g\left( \sum_{i\geq 0}\alpha_{i}x_{i}\right) \leq \sum_{i\geq 0} \alpha_{i}g(x_{i}). \label{proba01_convexInf}
\end{equation}

\bigskip \noindent \textbf{Proof}. Assume we have the hypotheses and the notations of the assertion to be proved. Either $g$ is bounded or $I$ in a bounded interval. In Part $A$ of this section, we saw that a convex function  \index{convex function} is continuous. If $I$ is a bounded closed interval, we get that $g$ is bounded on $I$. So, by the hypotheses of the assertion, $g$ is bounded. Let $M$ be a bound of $|g|$.\\

\noindent Now, fix $k\geq 1$ and denote
\begin{equation*}
\gamma _{k}=\sum_{i=0}^{k}\alpha _{i}\text{ and }\alpha _{i}^{\ast }=\alpha
_{i}/\gamma _{k}\text{ for }i=0,...,k
\end{equation*}

\noindent and

\begin{equation*}
\beta _{k}=\sum_{i\geq k+1}\alpha _{i}\text{ and }\alpha _{i}^{\ast \ast
}=\alpha _{i}/\beta _{k}\text{ for }i\geq k+1.
\end{equation*}

\noindent We have
\begin{eqnarray*}
g\left( \sum_{i\geq 0}\alpha _{i}x_{i}\right)  &=&g\left(
\sum_{i=0}^{k}\alpha _{i}x_{i}+\sum_{i\geq k+1}\alpha _{i}x_{i}\right)  \\
&=&g\left( \gamma _{k}\sum_{i=0}^{k}\alpha _{i}^{\ast }x_{i}+\beta
_{k}\sum_{i\geq k+1}\alpha _{i}^{\ast \ast }x_{i}\right)  \\
&\leq &\gamma _{k}g\left( \sum_{i=0}^{k}\alpha _{i}^{\ast }x_{i}\right)
+\beta _{k}g\left( \sum_{i\geq k+1}\alpha _{i}^{\ast \ast }x_{i}\right) ,
\end{eqnarray*}

\noindent by convexity. By using Formula (\ref{proba01_convexK}), and by using the
fact that $\gamma _{k}\alpha _{i}^{\ast }=\alpha _{i}$ for $i=1,...,k$, we
arrive at

\begin{equation}
g\left( \sum_{i\geq 0}\alpha _{i}x_{i}\right) \leq \sum_{i=0}^{k}\alpha
_{i}g(x_{i})+\beta _{k}g\left( \sum_{i\geq k+1}\alpha _{i}^{\ast \ast
}x_{i}\right) .  \label{proba01_convex4}
\end{equation}

\Bin \noindent Now, we use the bound $M$ of $g$ to have

\begin{equation*}
\left\vert \beta _{k}g\left( \sum_{i\geq k+1}\alpha _{i}^{\ast \ast
}x_{i}\right) \right\vert \leq M\beta _{k}\rightarrow 0\text{ as }
k\rightarrow \infty .
\end{equation*}

\Bin
\noindent Then, by letting $k\rightarrow \infty $ in (\ref{proba01_convex4}), we get

\begin{equation}
g\left( \sum_{i\geq 0}\alpha _{i}x_{i}\right) \leq \sum_{i=0}^{+infty}\alpha
_{i}g(x_{i}).
\end{equation}

\noindent $\blacksquare$\\

\bigskip

\noindent Later, a course on Measure Theory will allow to have powerful convergence theorems. But many people use probability theory  \index{probability theory} and/or statistical mathematics without taking a course such a course.\\

\noindent To these readers, we may show how to stay at the level of elementary calculus and to have workable versions of powerful tools when dealing with discrete probability measure \index{probability measure} s.\\

\noindent We will have two results of calculus that will help in that sense. In this section, we deal with real numbers series of the form

$$
\sum_{n\in \mathbb{Z}} a_n.
$$

\noindent The sequence $(a_n)_{n \in \mathbb{Z}}$ may be considered as a function $f : \mathbb{Z} \mapsto \overline{\mathbb{R}}$ such that, any $n \in \mathbb{Z}$, we have for $f(n)=a_n$, and the series becomes

$$
\sum_{n\in \mathbb{Z}} f(n).
$$

\Bin
\noindent We may speak of a sequence $(a_n)_{n \in \mathbb{Z}}$ or simply of a function $f$ on $\mathbb{Z}$.\\

\bigskip \noindent \textbf{A - Monotone Convergence  \index{monotone convergence} Theorem of series}.\\

\begin{theorem} Let $f : \mathbb{Z} \mapsto \overline{\mathbb{R}}_{+}$ be non-negative function and let $f_p : \mathbb{Z} \mapsto \overline{\mathbb{R}}_{+}$,  $p\geq 1$, be a sequence of non-negative functions increasing to $f$ in the following sense :

\begin{equation}
\forall (n\in \mathbb{Z}), \ \ 0\leq f_{p}(n)\uparrow f(n)\text{ as }p\uparrow\infty   \label{comptage01}
\end{equation}

\bigskip
\noindent Then, we have

\begin{equation*}
\sum_{n\in \mathbb{Z}}f_{p}(n)\uparrow \sum_{n\in \mathbb{Z}}f(n)
\end{equation*}
\end{theorem}

\bigskip \noindent \textbf{Proof}. \noindent Let us study two cases.\\

\noindent \textbf{Case 1}. There exists $n_{0} \in \mathbb{Z}$ such that $f(n_{0})=+\infty$. Then

\begin{equation*}
\sum_{n\in \mathbb{Z}}f(n)=+\infty
\end{equation*}

\Bin \noindent Since $f_{p}(n_{0})\uparrow f(n_{0})=+\infty$, then for any $M>0$, there exits $P_0\geq 1$, such that

\begin{equation*}
p>P_0 \Rightarrow f_{p}(n_{0})>M
\end{equation*}

\Bin \noindent Thus, we have

\begin{equation*}
\forall (M>0),\exists (P_0 \geq 1),\text{ }p>P_0 \Rightarrow \sum_{n\in \mathbb{Z}%
}f_{p}(n)>M.
\end{equation*}

\Bin \noindent By letting $p\rightarrow \infty$ in the latter formula, we have :

\begin{equation*}
\forall M>0, \ \sum_{n\in \mathbb{Z}}f_{p}(n) \geq M.
\end{equation*}

\Bin By letting $M \rightarrow +\infty$, we get

\begin{equation*}
\sum_{n\in \mathbb{Z}}f(n)\uparrow +\infty =\sum_{n\in \mathbb{Z}}f(n)
\end{equation*}

\Bin
\noindent The proof is complete of Case 1.\\

\bigskip \noindent \textbf{Case 2}. Suppose that $f(n)<\infty$ for any $n\in \mathbb{Z}$. By \eqref{comptage01}, it is clear that

\begin{equation*}
\sum_{n\in \mathbb{Z}}f_{p}(n)\leq \sum_{n\in \mathbb{Z}}f(n).
\end{equation*}

\Bin
\noindent The left-hand member is nondecreasing in $p$. So its limit is a monotone limit and it always exists in $\overline{\mathbb{R}}$ and we have

\begin{equation*}
\lim_{p\rightarrow +\infty}\sum_{n\in \mathbb{Z}}f_{p}(n)\leq \sum_{n\in \mathbb{Z}}f(n).
\end{equation*}

\Bin \noindent Now, fix an integer $N>1$. For any $\varepsilon >0$, there exists $P_{N}$ such that

\begin{equation*}
p>P_{N}\Rightarrow (\forall (-N\leq n\leq N),\text{ }f(n)-\varepsilon
/(2N+1)\leq f_{p}(n)\leq f(n)+\varepsilon /(2N+1)
\end{equation*}

\Bin \noindent and thus,

\begin{equation*}
p>P_{N}\Rightarrow \sum_{-N\leq n\leq N}f_{p}(n)\geq \sum_{-N\leq n\leq N}f(n)-\varepsilon.
\end{equation*}

\Bin \noindent Thus, we have

\begin{equation*}
p>P_{N}\Rightarrow \sum_{n\in \mathbb{Z}}f_{p}(n)\geq \sum_{-N\leq n\leq N}f_{p}(n)\geq
\sum_{-N\leq n\leq N}f(n)-\varepsilon,
\end{equation*}

\Bin  \noindent meaning  \index{mean} that for $p>P_{N}$,

\begin{equation*}
\sum_{-N\leq n\leq N}f(n)-\varepsilon \leq \sum_{-N\leq n\leq N}f_{p}(n).
\end{equation*}

\Bin \noindent  and then, for $p>P_{N}$,

\begin{equation*}
\sum_{-N\leq n\leq N}f(n)-\varepsilon \leq \sum_{n\in \mathbb{Z}}f_{p}(n).
\end{equation*}

\Bin
\noindent By letting  $p\uparrow \infty$ first and next, $N\uparrow \infty$, we get for any $\varepsilon >0,$%

\begin{equation*}
\sum_{n\in \mathbb{Z}}f(n)-\varepsilon \leq \lim_{p\rightarrow +\infty}\sum_{n\in \mathbb{Z}}f_{p}(n).
\end{equation*}

\Bin \noindent Finally, by letting $\varepsilon \downarrow 0$, we get

\begin{equation*}
\sum_{n\in \mathbb{Z}}f(n)\leq \lim_{p\rightarrow +\infty}\sum_{n\in \mathbb{Z}}f_{p}(n).
\end{equation*}

\Bin \noindent We conclude that

\begin{equation*}
\sum_{n\in \mathbb{Z}}f(n)=lim_{p\uparrow+\infty} \sum_{n\in \mathbb{Z}}f_{p}(n)
\end{equation*}

\Bin \noindent Case 2 is now closed and the proof of the theorem is complete.\\

\bigskip \noindent \textbf{B - Fatou-Lebesgues Theorem or Dominated Convergence  \index{dominated convergence} Theorem for series}.\\

\begin{theorem} Let $f_p : \mathbb{Z} \mapsto \overline{\mathbb{R}}$,  $p\geq 1$, be a sequence of functions.\\

\noindent (1) Suppose there exists a function $g : \mathbb{Z} \mapsto \mathbb{R}$ such that :\\

\noindent (a) $\sum_{n \in \mathbb{Z}} |g(n)| <+\infty$,\\

\noindent and\\

\noindent (b) For any $n \in \mathbb{Z}$, $g(n)\leq f_p(n)$.\\

\noindent Then we have

\begin{equation}
\sum_{n\in \mathbb{Z}} \liminf_{p\rightarrow +\infty} f_p(n) \leq \liminf_{p\rightarrow +\infty} \sum_{n\in \mathbb{Z}} f_p(n).
\end{equation}

\bigskip
\noindent (2) Suppose there exists a function $g : \mathbb{Z} \mapsto \mathbb{R}$ such that :\\

\noindent (a) $\sum_{n \in \mathbb{Z}} |g(n)| <+\infty$,\\

\noindent and\\

\noindent (b) For any $n \in \mathbb{Z}$, $f_p{n} \leq g(n)$.\\

\noindent Then we have

\begin{equation}
\limsup_{p\rightarrow +\infty} \sum_{n\in \mathbb{Z}} f_p(n) \leq \sum_{n\in \mathbb{Z}} \limsup_{p\rightarrow +\infty} f_p(n).
\end{equation}

\bigskip
\noindent (3) Suppose that there exists a function $f : \mathbb{Z} \mapsto \overline{\mathbb{R}}$ such that the sequence of functions
$(f_p)_{p\geq 1}$ point-wisely converges to $f$, that is :

\begin{equation}
\forall (n\in \mathbb{Z}), \ \  f_{p}(n)\rightarrow f(n) \  as \ p\rightarrow \infty.   \label{comptage03}
\end{equation}

\bigskip
\noindent And suppose that there exists a function $g : \mathbb{Z} \mapsto \mathbb{R}$ such that :

\noindent (a) $\sum_{n \in \mathbb{Z}}|g(n)| <+\infty$,\\

\noindent and\\

\noindent (b) For any $n \in \mathbb{Z}$, $|f_p(n)| \leq |g(n)|$.\\

\noindent Then we have
\begin{equation}
\lim_{p\rightarrow +\infty} \sum_{n\in \mathbb{Z}} f_p(n)= \sum_{n\in \mathbb{Z}} f(n) \in \mathbb{R}.
\end{equation}
\end{theorem}

\bigskip \noindent \textbf{Proof}. We proceed with three parts.\\

\noindent Part (1). Under the hypotheses of this part, we have that $g(n)$ is finite for any $n\in \mathbb{Z}$ and then $(f_p-g)_{p\geq 1}$ is a sequence nonnegative function defined on $\mathbb{Z}$ with values in $\overline{\mathbb{R}}_{+}$. The sequence of functions

$$
h_p=\inf_{k\geq p} (f_k-g)
$$

\Bin
\noindent defined by

$$
h_p(n)=\inf_{k\geq p} (f_k(n)-g(n)), \ n\in \mathbb{Z}
$$

\noindent is a sequence non-negative functions defined on $\mathbb{Z}$ with values in $\overline{\mathbb{R}}_{+}$.\\

\noindent By reminding the definition of the inferior limit (See the first section of that Appendix Chapter), we see that for any
$n \in \mathbb{Z}$,

$$
h_p(n) \uparrow \liminf_{p\rightarrow +\infty} (f_p(n)-g(n)).
$$

\bigskip
\noindent Recall that, for any fixed $n\in \mathbb{Z}$,

$$
\liminf_{p\rightarrow +\infty} (f_p(n)-g(n))= \left( \liminf_{_p \rightarrow +\infty}  f_p(n) \right)  - g(n).
$$

\bigskip
\noindent On one side, we may apply the Monotone Convergence  \index{monotone convergence} Theorem to get, as $p \uparrow +\infty$,

\begin{equation}
\sum_{n\in \mathbb{Z}} h_p(n) \uparrow \sum_{n\in \mathbb{Z}} \liminf_{p\rightarrow +\infty} f_p(n) - \sum_{n\in \mathbb{Z}} g(n).
\label{comptage05}
\end{equation}

\bigskip
\noindent We also have, for any $n\in \mathbb{Z}$, for any $k\geq p$,

$$
h_p(n) \leq (f_k(n)-g(n))
$$

\Bin
\noindent and for any $k\geq p$

$$
\sum_{n\in \mathbb{Z}} h_p(n) \leq \sum_{n\in \mathbb{Z}} (f_k(n)-g(n))
$$

\begin{equation}
\sum_{n\in \mathbb{Z}} h_p(n) \leq \inf_{k\geq p} \sum_{n\in \mathbb{Z}} (f_k(n)-g(n)). \label{comptage06}
\end{equation}

\Bin
\noindent Remark that

$$
\inf_{k\geq p} \sum_{n\in \mathbb{Z}} (f_k(n)-g(n))= \left\{\inf_{k\geq p} \sum_{n\in \mathbb{Z}} f_k(n) \right\} - \sum_{n\in \mathbb{Z}} g(n)
$$

\noindent to get

\begin{equation}
\sum_{n\in \mathbb{Z}} h_p(n) \leq \left\{\inf_{k\geq p} \sum_{n\in \mathbb{Z}} f_k(n)\right\} - \sum_{n\in \mathbb{Z}} g(n). \label{comptage07}
\end{equation}

\bigskip
\noindent By letting $p\rightarrow +\infty$, by using (\ref{comptage05}) and by reminding the definition of the inferior limit, we get

\begin{equation}
\sum_{n\in \mathbb{Z}} \liminf_{p\rightarrow +\infty} f_p(n) - \sum_{n\in \mathbb{Z}} g(n) \leq \left(\liminf_{p\rightarrow+\infty} \sum_{n\in \mathbb{Z}} f_k(n)\}\right) - \left( \sum_{n\in \mathbb{Z}} g(n)\right). \label{comptage08}
\end{equation}

\Bin
\noindent Since $\sum_{n\in \mathbb{Z}} g(n)$ is finite, we may add it to both members to have

\begin{equation}
\sum_{n\in \mathbb{Z}} \liminf_{p\rightarrow +\infty} f_p(n) \leq \liminf_{p\rightarrow} \sum_{n\in \mathbb{Z}} f_p(n)
\end{equation}.

\noindent Part (2). By taking the opposite functions $-f_p$, we find ourselves in the Part 1. In applying this part, the inferior limits are converted into superior limits when the minus sign is taken out of the limits, and next eliminated.\\

\noindent Part (3). In this case, the conditions of the two first parts hold. We have the two conclusions we may combine in

\begin{eqnarray*}
&& \sum_{n\in \mathbb{Z}} \liminf_{p\rightarrow +\infty} f_p(n)\\
&\leq& \liminf_{p\rightarrow+\infty} \sum_{n\in \mathbb{Z}} f_p(n)\\
&\leq& \limsup_{p\rightarrow+\infty} \sum_{n\in \mathbb{Z}} f_p(n)\\
&\leq&  \sum_{n\in \mathbb{Z}} \limsup_{p\rightarrow+\infty} f_p(n).
\end{eqnarray*}

\Bin
\noindent Since the limit of the sequence of functions $f_p$ is $f$, that is for any $n\in \mathbb{Z}$,

$$
\liminf_{p\rightarrow+\infty} f_p(n)=\limsup_{p\rightarrow+\infty} f_p(n)=f(n),
$$

\noindent and since

$$
\left\vert \sum_{n\in \mathbb{Z}} f_p(n) \right\vert \leq \mathbb{Z} |f_p(n)| <+\infty,
$$

\Bin \noindent  we get

$$
\lim_{p\rightarrow+\infty} \sum_{n\in \mathbb{Z}} f_p(n)=\sum_{n\in \mathbb{Z}} f(n)
$$

and

$$
\left\vert \sum_{n\in \mathbb{Z}} f(n) \right\vert<+\infty.
$$

\newpage

\noindent We are going to provide a proof of this formula. One can find several proofs (See \cite{feller1} \index{Feller}, page 52, for example). Here, we  give the proof in \cite{valiron} \index{Valiron}, pp. 167, that is based on Wallis integrals.  \index{Wallis integral} We think that a student in first year of University will be interested by an application the Riemann Integration course as below.\\

\noindent \textbf{A - Wallis Formula}.\\

\noindent We have the Wallis Formula

\begin{equation}
\pi =\lim_{n\rightarrow +\infty }\left\{ \frac{2^{4n}\times (n!)^{4}}{%
n\left\{ (2n)!\right\} ^{2}}\right\} .  \label{proba01_wallis_02}
\end{equation}

\noindent \textbf{Proof}. Define for $n\geq 1,$
\begin{equation}
I_{n}=\int_{0}^{\pi /2}\sin ^{n}xdx.  \label{proba01_wallis_01}
\end{equation}

\noindent On one hand, let us remark that for $0<x$ $<\pi /2,$ we have by the
elementary properties of the sine function that
\begin{equation*}
0<\sin x<1<(1/\sin x).
\end{equation*}

\noindent Then by multiplying all members of this double inequality by $\sin ^{2n}x$
for $0<x$ $<\pi /2$ and for $n\geq 1,$ we have%

\begin{equation*}
\sin ^{2n+1}x<\sin ^{2n}x<\sin ^{2n-1}x.
\end{equation*}

\noindent By integrating the functions in the inequality over $[0,\pi /2],$ and by
using (\ref{proba01_wallis_01}), we obtain for $n\geq 1.$%
\begin{equation}
I_{2n+1}<I_{2n}<I_{2n-1}. \label{proba01_wallis_DI}
\end{equation}

\noindent On another hand, we may integrate by part and get for $n\geq 2,$

\begin{eqnarray*}
I_{n} &=&\int_{0}^{\pi /2}\sin ^{2}x\sin ^{n-2}xdx \\
&=&\int_{0}^{\pi /2}(1-\cos ^{2}x)\sin ^{n-2}x\text{ }dx \\
&=&I_{n-2}-\int_{0}^{\pi /2}\cos x\text{ }(\cos x\sin ^{n-2}x)\text{ }dx \\
&=&I_{n-2}-\frac{1}{n-1}\int_{0}^{\pi /2}\cos x\text{ }d(\sin ^{n-1}x) \\
&=&I_{n-2}-\left[ \frac{\cos x\text{ }\sin ^{n-1}x}{n-1}\right] _{0}^{\pi
/2}-\frac{1}{n-1}\int_{0}^{\pi /2}\sin ^{n}x\text{ }dx \\
&=&I_{n-2}-\frac{1}{n-1}\int_{0}^{\pi /2}\sin ^{n}x\text{ }dx \\
&=&I_{n-2}-\frac{1}{n-1}I_{n}.
\end{eqnarray*}

\bigskip
\noindent Then for $n\geq 2,$

\begin{equation*}
I_{n}=\frac{n-1}{n}I_{n-2}.
\end{equation*}

\bigskip
\noindent  We apply this to an even number $2n$, $n\geq 1,$ to get by induction%
\begin{eqnarray*}
I_{2n} &=&\frac{2n-1}{2n}I_{2n-2}=\frac{2n-1}{2n}\times \frac{2n-3}{2n-2}%
I_{2n-4} \\
&=&\frac{2n-1}{2n}\frac{2n-3}{2n-2}\times ...\times \frac{2n-2p+1}{2n-2p+2}%
I_{2n-2p}.
\end{eqnarray*}

\bigskip
\noindent For $p=n\geq 1,$ we get

\begin{equation*}
I_{2n}=I_{0}\frac{1\times 3\times 5\times ...\times (2n-1)}{2\times 4\times
...\times 2n}.
\end{equation*}

\bigskip
\noindent
For an odd number $2n+1$, $n\geq 1,$ we have
\begin{eqnarray*}
I_{2n+1} &=&\frac{2n}{2n+1}I_{2n-1}=\frac{2n}{2n+1}\times \frac{2n-2}{2n-1}%
I_{2n-3} \\
&=&\frac{2n}{2n+1}\times \frac{2n-2}{2n-1}I_{2n-3}\times ...\times \frac{%
2n-2p}{2n-2p+1}I_{2n-2p+1}.
\end{eqnarray*}

\bigskip
\noindent For $p=n-1\geq 0,$ we have
\begin{equation*}
I_{2n+1}=I_{1}\frac{2\times 4\times ...\times 2n}{3\times 5\times ...\times
(2n+1)}.
\end{equation*}

\bigskip
\noindent We easily check that
\begin{equation*}
I_{0}=\pi /2\text{ and }I_{1}=1.
\end{equation*}

\bigskip
\noindent Thus the inequality $I_{2n+1}<I_{2n}$ of (\ref{proba01_wallis_DI}), for $n\geq 1$, yields

\begin{equation*}
\frac{\left\{ 2\times 4\times ...\times 2n\right\} ^{2}}{\left\{ 3\times
5\times ...\times (2n-1)\right\} ^{2}(2n+1)}<\frac{\pi }{2}
\end{equation*}

\Bin
\noindent and $I_{2n}<I_{2n-1}$ in (\ref{proba01_wallis_DI}), for $n\geq 1$, leads to
\begin{equation*}
\frac{\pi }{2}<\frac{\left\{ 2\times 4\times ...\times 2n\right\} ^{2}}{%
\left\{ 3\times 5\times ...\times (2n-1)\right\} ^{2}(2n)}.
\end{equation*}

\noindent We get the double inequality

\begin{equation*}
\frac{\left\{ 2\times 4\times ...\times 2n\right\} ^{2}}{\left\{ 3\times
5\times ...\times (2n-1)\right\} ^{2}(2n+1)}<\frac{\pi }{2}<\frac{\left\{
2\times 4\times ...\times 2n\right\} ^{2}}{\left\{ 3\times 5\times ...\times
(2n-1)\right\} ^{2}(2n)},
\end{equation*}

\bigskip \noindent which leads to,

\begin{equation*}
\frac{2n}{2n+1}<\left\{ \pi /2\right\} \left\{ \frac{\left\{ 2\times 4\times
...\times 2n\right\} ^{2}}{\left\{ 3\times 5\times ...\times (2n-1)\right\}
^{2}(2n)}\right\} ^{-1} \leq 1,\text{ }n\geq 1.
\end{equation*}

\Bin
\noindent In particular, this implies that

\begin{equation}
\pi /2=\lim_{n\rightarrow +\infty }\left\{ \frac{\left\{ 2\times 4\times
...\times 2n\right\} ^{2}}{\left\{ 3\times 5\times ...\times (2n-1)\right\}
^{2}(2n)}\right\} .  \label{proba01_wallis_03}
\end{equation}

\Bin
\noindent Put

\begin{equation*}
A_{n}=\frac{\left\{ 2\times 4\times ...\times 2n\right\} ^{2}}{\left\{
3\times 5\times ...\times (2n-1)\right\} ^{2}(2n)}.
\end{equation*}

\Bin
\noindent Let us transform its the numerator to get

\begin{eqnarray*}
A_{n} &=&\frac{\left\{ 2\times 4\times ...\times 2n\right\} ^{2}}{\left\{
3\times 5\times ...\times (2n-1)\right\} ^{2}(2n)} \\
&=&\frac{\left\{ (2\times 1)\times (2\times 2)\times ...\times (2\times
n)\right\} ^{2}}{\left\{ 3\times 5\times ...\times (2n-1)\right\} ^{2}(2n)}
\\
&=&\frac{\left\{ 2^{n}\times n!\right\} ^{2}}{\left\{ 3\times 5\times
...\times (2n-1)\right\} ^{2}(2n)}.
\end{eqnarray*}

\bigskip \noindent Next, we transform the denominator to get

\begin{eqnarray*}
A_{n} &=&\left\{ 2^{n}\times n!\right\} ^{2}\frac{\left\{ 2\times 4\times
...\times 2(n-1)\right\} ^{2}}{\left\{ (2n-1)!\right\} ^{2}(2n)} \\
&=&\left\{ 2^{n}\times n!\right\} ^{2}\frac{\left\{ 2\times 4\times
...\times 2(n-1)\right\} ^{2}(2n)^{2}}{\left\{ (2n-1)!\right\} ^{2}(2n)^{2}}%
\times \frac{1}{2n} \\
&=&\left\{ 2^{n}\times n!\right\} ^{2}\frac{\left\{ 2\times 4\times
...\times 2(n-1)(2n)\right\} ^{2}}{\left\{ (2n)!\right\} ^{2}}\times \frac{1%
}{2n} \\
&=&\frac{1}{2}\times \frac{2^{4n}\times (n!)^{4}}{n\left\{ (2n)!\right\} ^{2}%
}.
\end{eqnarray*}

\Bin \noindent Combining this with (\ref{proba01_wallis_03}) leads to the Wallis formula (\ref{proba01_wallis_02}).\\

\noindent The Wallis formula is now proved. Now, let us move to the Stirling Formula.\\

\bigskip \noindent \textbf{Stirling's Formula \index{Stirling's formula} }.\\

\noindent We have for $n\geq 1$
\begin{equation*}
\log (n!)=\sum_{p=1}^{p=n}\log p
\end{equation*}

\Bin
\noindent and

\begin{eqnarray*}
\int_{0}^{n}\log x\text{ }dx &=&\int_{0}^{n}d(x\log x-x)=\left[ x\log x-x%
\right] _{0}^{n} \\
&=&n\log n-n=n\log (n/e).
\end{eqnarray*}

\Bin
\noindent Let consider the sequence

\begin{equation*}
s_{n}=\log (n!)-n(\log n-1).
\end{equation*}

\Bin
\noindent By using an expansion of the logarithm function in the neighborhood of $1$ wa have for $n>1,$
\begin{eqnarray}
u_{n} &=&s_{n}-s_{n-1}=1+(n-1)\log \frac{n-1}{n} \notag \\
&=&1+(n-1)\log (1-\frac{1}{n})  \notag \\
&=&1+(n-1)\left( -\frac{1}{n}-\frac{1}{2n^{2}}-\frac{(1+\varepsilon _{n}(1))}{3n^{3}} \right)  \label{proba01_wallis_log_01} \\
&=&\frac{1}{2n}+\frac{1}{6n^{2}}-\frac{\varepsilon _{n}(1)}{3n^{2}}+\frac{(1+\varepsilon _{n}(1))}{3n^{3}},  \notag
\end{eqnarray}

\Bin
\noindent where

\begin{eqnarray*}
0 \leq \varepsilon _{n}(1)&=&3 \sum_{p=4} \frac{1}{pn^{p-3}}\\
&\leq & \frac{3}{4n} \sum_{p=4} \frac{4}{pn^{p-4}}\\
&\leq & \frac{3}{4n} \sum_{p=4} \frac{1}{n^{p-4}}\\
&=& \frac{3}{4n} \rightarrow 0 \text{ as } n\rightarrow +\infty.
\end{eqnarray*}

\Bin
\noindent Set

\begin{eqnarray}
S_{n} &=&s_{n}-\frac{1}{2}\log n  \notag \\
&=&\log (-n!)-n\log (n/e)-\frac{1}{2}\log n.  \label{proba01_wallis_log_02}
\end{eqnarray}

\Bin
\noindent We have
\begin{eqnarray*}
v_{n} &=&S_{n}-S_{n-1} \\
&=&s_{n}-s_{n-1}+\frac{1}{2}\log \frac{n-1}{n}.
\end{eqnarray*}

\bigskip
\noindent By using again the same expansion for%
\begin{equation*}
\log \frac{n-1}{n}=-\log \frac{n}{n-1},
\end{equation*}

\noindent and by combining with (\ref{proba01_wallis_log_01}), we get
\begin{equation*}
v_{n}=-\frac{1}{12n^{2}}-\frac{\varepsilon _{n}(1)}{3n^{2}}+\frac{%
(1+\varepsilon _{n}(1))}{6n^{3}}.
\end{equation*}

\noindent We see that the series $S=\sum_{n}v_{n}$ is finite and we have

\begin{equation}
R_{n}=S-S_{n}=\sum_{p=n+1}^{\infty }\left\{ -\frac{1}{12p^{2}}-\frac{%
\varepsilon _{p}(1)}{3p^{2}}+\frac{(1+\varepsilon _{p}(1))}{6p^{3}}%
)\right\} . \label{proba01_wallis_rn}
\end{equation}

\bigskip \noindent Is is clear that $R_{n}$ may be bounded by the sum of three remainders of convergent series, for $n$, large enough. So
$R_{n}$ goes to zero  as $n\rightarrow +\infty$. We will come back to a finer analysis to $R_n$.\\

\noindent From (\ref{proba01_wallis_log_02}), we have for $n\geq 1,$%
\begin{eqnarray*}
n! &=&n^{1/2}(n/e)^{n}\exp (-(S-S_{n})) \\
&=&n^{1/2}(n/e)^{n}e^{S}\exp (-R_{n}).
\end{eqnarray*}

\bigskip
\noindent Let us use the Wallis Formula

\begin{eqnarray*}
\pi &=&\lim_{n\rightarrow +\infty }\frac{2^{4n}\times (n!)^{4}}{n\left\{
(2n)!\right\} ^{2}} \\
&=&\lim_{n\rightarrow +\infty }\frac{2^{4n}\times n^{2}(n/e)^{4n}e^{S}\exp
(-4R_{n})}{n\left\{ (2n)^{1/2}(2n/e)^{2n}e^{4S}\exp (-R_{2n})\right\} ^{2}}
\\
&=&\lim_{n\rightarrow +\infty }\frac{2^{4n}\times n^{2}(n/e)^{4n}e^{4S}\exp
(-4R_{n})}{n(2n)(2n/e)^{4n}e^{2S}\exp (-2R_{2n})} \\
&=&\lim_{n\rightarrow +\infty }\frac{1}{2}\frac{2^{4n}\times
n^{2}(n/e)^{4n}e^{4S}}{2^{4n}n^{2}(n/e)^{4n}e^{2S}}\exp (-4R_{n}-2R_{n}) \\
&=&\frac{1}{2}e^{2S},
\end{eqnarray*}

\bigskip \noindent since $\exp (-4R_{n}-2R_{n})\rightarrow 0$ as $n\rightarrow +\infty$.\\

\bigskip
\noindent Thus, we arrive at

\begin{equation*}
e^{S}=\sqrt{2\pi }.
\end{equation*}

\noindent We get the formula

\begin{equation*}
n!=\sqrt{n2\pi }(n/e)^{n}\exp (-R_{n}).
\end{equation*}

\bigskip \noindent Now let us make a finer analysis to $R_n$. Let us make use of the classical comparison between series of monotone sequences and integrals.\\

\noindent Let $b\geq 1$. By comparing the area under the curve of $f(x)=x^{-b}$ going from $j$ to $k$
and that of the rectangles based on the intervals $[h,h+1],$ $h=1,..,k-1,$ we obtain

\begin{equation*}
\sum_{h=j+1}^{k}h^{-b}\leq \int_{j}^{k-1}x^{-b} \ dx \leq \sum_{h=j}^{k-1}h^{-b}.
\end{equation*}

\noindent This implies that

\begin{equation}
\int_{j}^{k}x^{-b}dx+(k)^{-b}\leq \sum_{h=j}^{k}h^{-b}\leq
\int_{j}^{k}x^{-b}dx+j^{-b}.  \label{fgintegn}
\end{equation}

\bigskip \noindent By applying this and letting $k$ go to $+\infty$ gives

$$
\frac{1}{(n+1)} \leq \sum_{p\geq n+1}^{k} p^{-2} \leq \frac{1}{n+1}+\frac{1}{(n+1)}
$$

\Bin and

$$
\frac{1}{2(n+1)^2} \leq \sum_{p\geq n+1}^{k} p^{-3} \leq \frac{1}{2(n+1)^2}+\frac{1}{(n+1)^3}.
$$

\Bin
\noindent Let $r>0$ be an arbitrary real number. For $n$ large enough, we will have $0\geq \varepsilon_{p}(1)\leq r$ for $p\geq (n+1)$. Then by combining the previous inequalities, we have for $n$ large enough,

\begin{eqnarray*}
-R_n&\geq &\frac{1}{12(n+1)}-\frac{(1+r)}{12(n+1)^2}-\frac{(1+r)}{6(n+1)^3}\\
&\geq &\frac{1}{(n+1)}\left( \frac{1}{12}-\frac{(1+r)}{12(n+1)}-\frac{(1+r)}{6(n+1)^2} \right).\\
\end{eqnarray*}

\Bin
\noindent Then $-R_n\geq 0$ for $n$ large enough. Next

\begin{eqnarray*}
-R_n &\leq & \frac{1}{12(n+1)}+\frac{1}{12(n+1)^2}+\frac{r}{3(n+1)}+\frac{r}{3(n+1)^2}-\frac{(1+r)}{12(n+1)^3}\\
&=&  \frac{1}{(n+1)} \left(\frac{1+4r}{12}+\frac{1}{12(n+1)}+\frac{r}{3(n+1)}-\frac{(1+r)}{12(n+1)^3}\right).
\end{eqnarray*}

\noindent And clearly, for $n$ large enough, we have

$$
-R_n \leq \frac{1+5r}{12(n+1)}.
$$

\Bin
\noindent Since $r>0$ is arbitrary, we have for any $\eta>$, for $n$ large enough,

$$
|R_n| \leq \frac{1+\eta}{12n}.
$$

\bigskip

\noindent  This finishes the proof of the Stirling Formula.

\include{99_proba01_cours_biblio_ang}

\newpage
\printindex

%\includ7e{pe_chap2_ang}
%\include{pe_chap3_ang}
%\include{pe_chap4_ang}
%\include{pe_chap5_ang}
%\include{pe_chap6_ang}
%\include{pe_chap7_ang} %couple
%\include{pe_chap8_ang}
%\include{pe_chap9_ang} % Resum\'e
%\include{pe_chap10_ang} % Software R and probability
%\include{pe_biblio}
\end{document}